\documentclass[12pt,centertags,oneside]{amsart}

\usepackage{tikz}

 \usetikzlibrary{decorations.pathreplacing,calligraphy}
\usetikzlibrary{arrows}
\usetikzlibrary{positioning,decorations.text}
\usetikzlibrary{matrix}

\usepackage{xcolor}
\usetikzlibrary{arrows.meta}

\usepackage{lmodern}
\usepackage[T1]{fontenc}

\usepackage[mathscr]{eucal}
\usepackage{asymptote}
\usepackage{caption}
\usetikzlibrary{arrows}
\usepackage{import}
\usepackage{graphicx,graphics,array,booktabs}

\usepackage{amsmath,amstext,amsthm,amscd,typearea,hyperref}
\usepackage{amssymb}
\usepackage{a4wide}
\usepackage[mathscr]{eucal}
\usepackage{mathrsfs}
\usepackage{typearea}
\usepackage{charter}
\usepackage{pdfsync}
\usepackage{multirow}
\usepackage{array}
\usepackage{amsmath, amssymb, latexsym, amscd, amsthm,amsfonts,amstext}
\usepackage[mathscr]{eucal}

\usepackage{amsmath,amstext,amsthm,amscd,typearea,stmaryrd,mathabx}
\usepackage{amssymb}
\usepackage{a4wide}
\usepackage[mathscr]{eucal}
\usepackage{mathrsfs}
\usepackage{typearea}
\usepackage{charter}
\usepackage{pdfsync}
\usepackage[a4paper,width=16.2cm,top=3cm,bottom=3cm,lines=50]{geometry}

\numberwithin{equation}{section}

\newtheorem{theorem}{Theorem}[section]
\newtheorem{definition}[theorem]{Definition}
\newtheorem{proposition}[theorem]{Proposition}
\newtheorem{corollary}[theorem]{Corollary}
\newtheorem{lemma}[theorem]{Lemma}
\newtheorem{remark}[theorem]{Remark}

\newtheorem{notation}[theorem]{Notation}

\newcommand{\cali}[1]{\mathscr{#1}}
\newcommand{\Tan}{\mathop{\mathrm{Tan}}\nolimits}

\newcommand{\comp}{{\mathop{\mathrm{comp}}\nolimits}}
\usepackage{graphicx} 

\newcommand{\AB}{{\rm AB}}

\newcommand{\Tube}{{\rm Tube}}

\newcommand{\supp}{{\rm supp}}

\newcommand{\dist}{{\rm dist}}

\newcommand{\ver}{{\rm ver}}
\newcommand{\hor}{{\rm hor}}

\renewcommand{\top}{{top}}
\newcommand{\ddc}{{dd^c}}

\newcommand{\ddcv}{{dd^c_\ver}}
\newcommand{\dc}{{d^c}}
\newcommand{\dbar}{{\overline\partial}}
\newcommand{\ddbar}{{\partial\overline\partial}}

\newcommand{\SH}{{\rm SH}}
\newcommand{\PH}{{\rm PH}}
\newcommand{\CL}{{\rm CL}}
\newcommand{\GL}{{\rm GL}}

\newcommand{\tot}{{\rm tot}}
\newcommand{\ind}{{\bf 1}}

\newcommand{\bfx}{{\bf x}}
\newcommand{\id}{{\rm id}}

\newcommand{\Vol}{{\rm Vol}}
\newcommand{\bfc}{{\rm \mathbf{c}}}

\newcommand{\bfr}{{\rm \mathbf{r}}}

\newcommand{\bfi}{{\rm \mathbf{i}}}

\newcommand{\bfj}{{\rm \mathbf{j}}}
\newcommand{\bfU}{{\rm \mathbf{U}}}
\newcommand{\bfW}{{\rm \mathbf{W}}}
\newcommand{\bfe}{{\rm \mathbf{e}}}
\newcommand{\upm}{{\mathrm{ \overline{m}}}}
\newcommand{\lowm}{{\mathrm{ \underline{m}}}}

\newcommand{\Cc}{\cali{C}}
\newcommand{\Dc}{\cali{D}}

\newcommand{\Fc}{\cali{F}}

\newcommand{\Hc}{\cali{H}}
\newcommand{\Kc}{\cali{K}}
\newcommand{\Lc}{\cali{L}}
\renewcommand{\Mc}{\cali{M}}
\newcommand{\Oc}{\cali{O}}

\newcommand{\Uc}{\cali{U}}
\newcommand{\Vc}{\cali{V}}

\newcommand{\Ic}{\cali{I}}
\newcommand{\Jc}{\cali{J}}
\newcommand{\Nc}{\cali{N}}
\newcommand{\FS}{{\rm FS}}
\newcommand{\DS}{{\rm DS}}

\newcommand{\C}{\mathbb{C}}
\newcommand{\D}{\mathbb{D}}
\newcommand{\E}{\mathbb{E}}
\renewcommand{\H}{\mathbb{H}}
\newcommand{\N}{\mathbb{N}}

\newcommand{\R}{\mathbb{R}}
\newcommand{\T}{\mathbb{T}}
\newcommand{\B}{\mathbb{B}}

\newcommand{\U}{\mathbb{U}}

\renewcommand{\P}{\mathbb{P}}
\newcommand{\X}{\mathbb{X}}

\def\mD{\mathcal{D}}

\def\mU{\mathcal{U}}

\def\mR{\mathcal{R}}

\title[]{The generalized Lelong numbers and intersection theory}

\author{Vi{\^e}t-Anh Nguy{\^e}n}
\address{Universit\'e de Lille, 
Laboratoire de math\'ematiques Paul Painlev\'e, 
CNRS U.M.R. 8524,  
59655 Villeneuve d'Ascq Cedex, 
France. }

\address{and Vietnam Institute for Advanced Study in Mathematics (VIASM),  157 Chua Lang Street, Hanoi, Vietnam.
}

\email{Viet-Anh.Nguyen@univ-lille.fr, {\tt   https://pro.univ-lille.fr/viet-anh-nguyen/}}

\date{January 04, 2025}

\begin{document}


\begin{abstract}

   Let $X$ be a  complex  manifold of dimension $k,$ and  $(V,\omega)$ be  a K\"ahler submanifold of dimension $l$ in $X,$
   and $B\Subset V$ be a domain with  $\Cc^2$-smooth boundary. 
   Let $T$ be a positive   plurisubharmonic  current  on $X$  
   such that $T$ satisfies a reasonable  approximation  condition on $X$ and near $\partial B.$ 
   In our previous  work  \cite{Nguyen21}   we
introduce    the   concept of  the generalized  Lelong numbers  $\nu_j(T,B)\in\R$   of $T$  along $B$  for $0\leq j\leq l.$  When $l=0,$  $V=B=\{ \text{a single point  $x\in X$}  \},$   $\nu_0(T,B)$ is  none  other  than  the   classical  Lelong number $\nu(T,x)$ of $T$ at $x.$

This  article has  five purposes.

Firstly,  we    formulate the notion of 
the generalized Lelong number of $T$ associated to every closed smooth $(j,j)$-form $\omega^{(j)}$ on $V.$ 
This  concept extends  the previous notion of the generalized Lelong numbers.   
We also establish their basic properties.

Secondly, we   define the horizontal dimension $\hbar$  of  such a current $T$ along $B.$  
Next,  we characterize   $\hbar$ in terms of  the generalized  Lelong numbers.
We also establish a Siu's  upper-semicontinuity type  theorem  for the  generalized Lelong numbers.
The horizontal dimension was first introduced  by Dinh-Sibony \cite{DinhSibony18}
in the   context  where the current  $T$ is   positive closed and   $\supp(T)\cap V$ is compact in $V.$

In their above-mentioned  context,
Dinh and  Sibony  introduced  the  cohomology classes $\bfc^\DS_j(T,B)\in H^{2l-2j}(B)$'s  which  may be regarded
as  their analogues of the  classical Lelong numbers. Our third objective is  to    generalize
 their notion  to  the  broader context  where $T$ is (merely) positive pluriharmonic.
Moreover,  we also  establish a formula  relating   Dinh-Sibony  classes  and the generalized 
Lelong numbers.  As  a  consequence,  the generalized  Lelong numbers  are totally intrinsic in this context.

Fourthly, we    obtain an effective sufficient condition for defining   the intersection $T_1\curlywedge \ldots\curlywedge  T_m$
in  the sense of Dinh-Sibony's theory  of tangent currents  \cite{DinhSibony18} of  $m\ (\geq 2)$ positive closed  currents $T_i$ of bidegree $(p_i,p_i)$  for $1\leq i\leq m$ on a compact K\"ahler  manifold $X.$

Finally,   we establish an effective sufficient condition  for   $m\ (\geq 2)$  sequences of   positive closed  currents $(T_{i,n})_{n\in\N}$  of bidegree $  (p_i,p_i)$ for $1\leq i\leq m$ on a compact K\"ahler  manifold $X$ to satisfy the  following    equality:
$$
\lim_{n\to\infty} \big( T_{1,n}\curlywedge \ldots\curlywedge  T_{m,n}\big)=(\lim_{n\to\infty}  T_{1,n})\curlywedge \ldots\curlywedge  (\lim_{n\to\infty} T_{m,n}).
$$

\end{abstract}

\maketitle

\medskip\medskip


 \medskip



\section{Introduction} \label{S:Intro}


 Let $X$ be a   complex   manifold of dimension $k.$ Let   $d,$ $\dc$ denote the  real  differential operators on $X$  defined by
$d:=\partial+\overline\partial,$  $\dc:= {1\over 2\pi i}(\partial -\overline\partial) $  with $i:=\sqrt{-1}$ so that 
$\ddc={i\over \pi} \partial\overline\partial.$  We start with the  following basic definitions.

\begin{definition}\label{D:positive-currents}\rm  A  $(p,p)$-current  $T$ defined on $X$  is  said to be
 positive  if $\langle T,\Phi\rangle \geq 0$  for every smooth  test form  $\Phi$  of the type
$$
\Phi=f\cdot (i\phi_1\wedge \bar \phi_1)\wedge  \ldots\wedge  (i\phi_q\wedge \bar\phi_q), 
$$
 where $q:=k-p$ and $\phi_1,\ldots, \phi_q$ are smooth $(1,0)$-forms on $X,$ $f$ is  a non-negative smooth function compactly supported on $X.$  
 \end{definition}
 
 \begin{definition}\label{D:3-classes-of-currents} \rm A   $(p,p)$-current  $T$ defined on $X$  is  said to be
 {\it closed}  (resp. {\it  pluriharmonic\footnote{Some  authors
use the  terminology  {\it harmonic} instead of {\it pluriharmonic}.}}), (resp. {\it plurisubharmonic})
if $dT=0$  (resp. $\ddc T=0$), (resp. $\ddc T$  is a positive current).  
\end{definition}
 
 If $V$ is a  subvariety of $X$ of pure codimension $p,$  then  the current of integration $[V]$ 
 associated to $V$ is a positive closed $(p,p)$-current. Moreover, if $T$ is  positive  closed, then
 it is clearly positive pluriharmonic. If  $T$ is  positive  pluriharmonic, then
 it is clearly positive plurisubharmonic.

 \subsection{Main purpose of this work}
 Let  $(V,\omega)$ be  a Hermitian submanifold of dimension $l$ in $X$ such that the Hermitian form is K\"ahler. But we point out that  the K\"ahlerian assumption  of $\omega$  could  be somehow  relaxed. 
   Let  $B\Subset V$ be a domain with  piecewise $\Cc^2$-smooth boundary. 
   Let $T$ be a positive   plurisubharmonic  current  on $X$  
   such that $T$ satisfies a reasonable  approximation  condition near $B$  and near $\partial B$ in  $X.$
   In our previous  work  \cite{Nguyen21}   we
introduce    the   concept of  the generalized  Lelong numbers  $\nu_j(T,B)\in\R$   of $T$  along $B$  for $0\leq j\leq l.$  
When $l=0,$  $V=B =\{ \text{a single point  $x\in X$}\},$   $\nu_0(T,B)$ becomes the   classical  Lelong number $\nu(T,x)$ of $T$ at $x.$   

Let us  first  talk about the particular case of  positive  closed  currents.
In this case the  Lelong numbers   play a fundamental role and have numerous  applications 
in  Complex  Analysis, Complex  Geometry, Algebraic Geometry and Complex  Dynamics, see, for example,  Griffiths-Harris \cite{GriffithsHarris},   Demailly  \cite{Demailly12} and  Dinh-Sibony  \cite{DinhSibony10b}, Lazarsfeld  \cite{Lazarsfeld-I,Lazarsfeld-II} as well as   the references therein.
 The readers could  find some systematic  developments on Lelong numbers  for  positive closed  currents in Siu \cite{Siu}  and Demailly
\cite{Demailly93,Demailly} etc. 

The  theory of intersection  for positive  closed $(1,1)$-currents  is by now   well-developed and well-understood.
The key point here is that positive closed $(1, 1)$-currents can be locally written as the $\ddc$ of plurisubharmonic (psh) functions
which are unique up to pluriharmonic functions. Therefore, the study of positive closed $(1, 1)$-currents can be systematically reduced to the study of psh
functions, see  the pioneering works of  Bedford-Taylor \cite{BedfordTaylor}, Forn{\ae}ss-Sibony \cite{FornaessSibony95},  Demailly \cite{Demailly92,Demailly93,Demailly12} 
etc. This theory  has a lasting  impact and finds many application in the contemporary study of degenerate  Monge-Amp\`ere  equations  which is  an important part of the  interface between Complex Analysis and Differential  Geometry etc.

However, the case of arbitrary bi-degree currents is still far from being well understood because  their local potentials may differ by singular currents.
 In \cite[p. 146]{Demailly92}, Demailly raised the problem of developing a theory of intersection for positive closed currents of higher bidegree. In  \cite{DinhSibony09,DinhSibony10,DinhSibony18} Dinh and  Sibony  gave  two remarkable answers to   this  problem via their two new theories, namely,  the super-potential theory and the theory of tangent  currents. More   concretely,  given  $m\ (\geq 2)$ positive closed  currents $T_i$ of bidegree $(p_i,p_i)$  for $1\leq i\leq m$ on a compact K\"ahler  manifold $X,$ under  theoretically reasonable assumptions,  Dinh-Sibony  can define the intersection $T_1\wedge \ldots\wedge T_m$ in the sense of the super-potential theory  (see  Definition \ref{D:DS-wedge-prod-super} below),  and  the intersection   $T_1\curlywedge \ldots\curlywedge  T_m$ in the sense of the   theory of tangent  currents (see  Definition \ref{D:DS-wedge-prod} below).    Dinh-Nguyen-Vu  \cite{DinhNguyenVu}   study the  relation between these  two  theories.
Ahn, Bayraktar, De Th\'elin, Dinh, Huynh, Kaufmann, Truong,Vigny, Vu  and  several other authors    apply these two theories to  many interesting  problems,  see \cite{AhnVu,Bayraktar,
 DeThelinVigny,DinhNguyenTruong15, DinhNguyenTruong17,HuynhKaufmannVu,HuynhVu,Kaufmann,KaufmannVu,Vigny,Vu16,Vu21a,Vu21b,Vu22}  etc.
In continuation  of Dinh-Sibony's pioneering works, 
we address, in this  article,  the following basic questions which  seem to play a very important role in  applications: 

\smallskip

\noindent{\bf Problem 1.} {\it  Find an   effective sufficient condition so that  $m\ (\geq 2)$ positive closed  currents $T_i$ of bidegree $(p_i,p_i)$  for $1\leq i\leq m$ on a compact K\"ahler  manifold $X$   are  wedgeable  in the sense  of the   theory of tangent  currents, that is, the current  $T_1\curlywedge \ldots\curlywedge  T_m$ exists.}

\smallskip

\noindent{\bf Problem 2.} {\it 
 For  $m\geq 2$ and    for $1\leq i\leq m,$ let     $(T_{i,n})_{n\in\N}$   be   a  sequence of   positive closed  currents   of bidegree $  (p_i,p_i)$  on a compact K\"ahler  manifold $X$  so  that $\lim_{n\to\infty}  T_{i,n}$ exists  in the sense of currents.  Find an   effective sufficient condition imposed on these $m$ sequences  so that the following equality holds:
$$
\lim_{n\to\infty} \big( T_{1,n}\curlywedge \ldots\curlywedge  T_{m,n}\big) =(\lim_{n\to\infty}  T_{1,n})\curlywedge \ldots\curlywedge  (\lim_{n\to\infty} T_{m,n}).
$$}

\smallskip
We refer the  readers to  Coman-Marinescu \cite{ComanMarinescu} for a connection between these two  problems and the equidistribution of  the zero-divisors of
random sequences of holomorphic sections in high tensor powers of a holomorphic line
bundle.

Next, we speak about  the Lelong numbers of positive  non-closed  currents.
The theory of singular holomorphic  foliations has been undergoing    several developments  during  the last 25 years.
One of its main development is the advent of the  ergodic  theory which emphasizes the global  study of singular holomorphic foliations
using     the analysis of positive currents and the  tools from  geometric complex    analysis.  These  techniques  are  developed by     Berndtsson-Sibony \cite{BerndtssonSibony}  and Forn{\ae}ss-Sibony \cite{FornaessSibony05,FornaessSibony08,FornaessSibony10}.
One major goal of the latter theory is to  
study the  global behaviour of  a  generic  leaf of a foliation $\Fc.$  Here, the  term {\it ``generic''} is measured  with respect to   some new objects called {\it  directed positive harmonic  currents} (see for example \cite[Definition 7]{FornaessSibony08} etc. for a definition of this  notion). Although these currents are positive  pluriharmonic,  the  problem is that  they 
are in general not closed. On the other hand, 
the      analysis  of these  currents 
are often  difficult when one approaches  the singularities of  $\Fc.$ The classcial Lelong numbers are only suitable   to the case where these singularities are  isolated. When the set of singularities is of positive  dimension, the  analysis  becomes much harder.  We expect  that 
the toolbox  that we develop  in \cite{Nguyen21} and  in the present work may  provide an effective way to handle  this  situation.

 As for   positive  plurisubharmonic currents,  the reader could  consult Alessandrini-Bassanelli  \cite{AlessandriniBassanelli96}.

   \smallskip
   
   {\bf  Notation.} Throughout  the article,  we denote by
   \begin{itemize}
   
 \item   $\D$  the unit disc in $\C;$
 \item  $\C^*$ the punctured complex plane $\C\setminus \{0\};$  
 \item  $\R^+:=[0,\infty)$ and  $\R^+_*:=(0,\infty);$
 \item  $\partial B$  the boundary  of an open  set $B$ in a manifold $Y.$

   \end{itemize}

   If $X$ is an oriented manifold, denote by $H^* (X, \C)$ the de Rham cohomology
group of $X$ and $H^*_\comp (X, \C)$ the de Rham cohomology group deﬁned by forms
or currents with compact support in $X.$ If $V$ is a submanifold of $X,$ denote
by $H^*_V (X, \C)$ the de Rham cohomology group deﬁned in the same way using
only forms or currents on $X$ whose supports intersect $V$ in a compact set.

If $T$ is either  a closed current on $X$  or a $\ddc$-closed current on a compact K\"ahler manifold $X,$ denote by $\{T \}$
its class in $H^* (X, \C).$ When
$T$ is supposed to have compact support, then $\{T \}$ denotes the class of $T$
in $H_\comp^* (X, \C).$ If we only assume that $\supp(T ) \cap V$ is compact, then $\{T \}$ denotes the class of $ T$ in $H_V^* (X, \C).$ The current of integration on an oriented
submanifold or a  complex  variety $Y$ is denoted by $[Y ].$ Its class is denoted by $\{Y \}.$

  For a differentiable map $\pi: X\to Y$ between  manifolds,  $\pi^*$ (resp. $\pi_*$) denotes the pull-back (resp. the push-forward) operator
  acting  on forms and currents  defined on  $Y$ (resp. on $X$). These operators induce  natural maps on cohomological levels:
  $\pi^*:\  H^*(Y,\C)\to  H^*(X,\C)$ and $\pi_*:\  H^*(X,\C)\to  H^*(Y,\C).$

 In the next  three subsections  we will recall known results. 
\subsection{Classical Lelong numbers}

Let $T$ be  a positive plurisubharmonic  $(p,p)$-current  defined on $X$  and $x\in X$   a point.
We first recall the  notion  of  Lelong number  $\nu(T,x)$ of $T$ at  $x.$
This  notion  was first introduced  by   Lelong in \cite{Lelong}  for the class of positive closed currents.
It was later  formulated  by Skoda in \cite{Skoda} for the wider class of positive  plurisubharmonic currents.  

Choose a local holomorphic coordinate system $z$ near $x$ such that $x = 0$ in
these coordinates. The Lelong number $\nu(T,x)$ of $T$ at $x$ is the limit of the normalized
mass of $\|T \|$  on the ball $\B(0, r)$ of center $0$ and radius $r$ when $r$ tends to $0.$ More
precisely, we have
\begin{equation}\label{e:Lelong-number-point}
\nu(T , x) := \lim\limits_{r\to 0}\nu(T,x,r),\quad\text{where}\quad \nu(T,x,r):=
{\sigma_T(\B(0,r))\over 
(2\pi)^{k-p} r^{ 2k-2p}}. 
\end{equation}
Here, $\sigma_T:={1\over (k-p)!}\, T\wedge ({i\over 2}\ddbar \|z\|^2)^{k-p}$ is  the trace measure of $T.$
Note that $(2\pi)^{k-p} r^{2k-2p}$  is the mass on $\B(0, r)$ of the $(p , p  )$-current of
integration on a linear subspace of dimension $k - p$ through $0.$  When $T$ is a positive closed current, Lelong establishes  in \cite{Lelong} (see also \cite{Lelong68})
that 
the  {\it  average mean} $\nu(T,x,r)$  is a  non-negative-valued  increasing function in the  radius $r.$
So  the limit \eqref{e:Lelong-number-point} always exists. 
Skoda \cite{Skoda} proves  the  same result  for  positive  plurisubharmonic currents.   Thie \cite{Thie} shows
that when $T$ is given by an
analytic set this number is the multiplicity of this set at $x.$ Siu proves  that when $T$ is  a positive  closed  current, the
limit \eqref{e:Lelong-number-point} does not depend on the choice of coordinates.

There is  another equivalent  logarithmic definition of the Lelong number of a positive closed  current that we want to discuss in this   work. More specifically,
consider
\begin{equation}\label{e:Lelong-number-point-loga}
 \kappa^\bullet (T,x,r):= \int_{\B(0,r)\setminus \{0\}} T(z)\wedge (\ddc\log{(\|z\|^2}))^{k-p}\quad\text{and}\quad
                        \kappa (T,x,r):= \int_{\B(0,r)} T(z)\wedge (\ddc\log{(\|z\|^2}))^{k-p}         \,\cdot
\end{equation}
The {\it  logarithmic means}  $\kappa^\bullet (T,x,r)$ and   $\kappa(T,x,r)$ are  non-negative-valued  increasing functions in the  radius $r.$  
Observe  that  in the expression of  $\kappa(T,x,r)$  in \eqref{e:Lelong-number-point-loga}, the  wedge-product of currents  is only well-defined  outside the origin $0$ because
the  second factor $ (\ddc\log{(\|z\|^2}))^{k-p}$ is  only smooth there. 
Here, we consider  two simple interpretations of  \eqref{e:Lelong-number-point-loga}  which correspond to
regularizing either  the first  or the second factor of the  wedge-product of currents in the  expression of
$ \kappa(T,x,r).$
 The  first interpretation concerns the  notion of approximation of currents.  By a standard regularization (e.g.  a convolution),  we see that there is a sequence of
 positive smooth closed  $(p,p)$-form  on   $\B(0,r+\epsilon)$ for some $\epsilon>0$  such that $T_n$  converges weakly  to $T.$ 
 The first  interpretation   
 of the integral on the RHS  of  \eqref{e:Lelong-number-point-loga}  is formulated as follows:
\begin{equation}\label{e:Lelong-number-point-bisbis(1)}
 \int_{\B(0,r)} T(z)\wedge (\ddc\log{(\|z\|^2}))^{k-p}:=\lim\limits_{n\to \infty}   \int_{\B(0,r)} T_n(z)\wedge (\ddc\log{(\|z\|^2)})^{k-p} \,\cdot
\end{equation}
provided that the limit exists. In fact, this  is indeed the case.
The  second  interpretation  consists in  regularizing  the integral kernel $(\ddc\log{(\|z\|^2}))^{k-p}$ in a    canonical way:
\begin{equation}\label{e:Lelong-number-point-bisbis(2)}
 \int_{\B(0,r)} T(z)\wedge (\ddc\log{(\|z\|^2)})^{k-p}:=\lim\limits_{\epsilon\to 0+}   \int_{\B(0,r)} T(z)\wedge (\ddc\log{(\|z\|^2+\epsilon^2)})^{k-p} \,\cdot
\end{equation}
provided that the limit exists. In fact, this is always  the case.

We record here the basic identities of the   logarithmic means:
\begin{equation}\label{e:Lelong-number-point-bis}
\lim\limits_{r\to 0} \kappa^\bullet(T,x,r)=0\quad\text{and }\quad \nu(T , x) := \lim\limits_{r\to 0} \kappa(T,x,r).
\end{equation}

\subsection{Tangent  currents     for $l > 0$  and Dinh-Sibony theory}

Next,  we deal with  the  situation where    the single point $x$ in the  previous  subsection  is  replaced  by   a submanifold $V\subset X$ of positive dimension $l$  ($1\leq l< k$).
Only recently,   Dinh and  Sibony \cite{DinhSibony18} have  developed  a  satisfactory  theory of tangent  currents and  density currents   for  positive  closed currents  in  this
context.

  Let $\E$ be the normal
vector bundle to $V$ in $X$ and $\pi:\ \E\to V$ be the  canonical projection.
Consider a point $x\in V.$ If $\Tan_x(X)$  and $\Tan_x(V)$ denote, respectively, the tangent spaces of $X$ and of $V$ at $x,$
the  fiber $\E_x$ of $\E$ over $x$ is canonically  identified  with the quotient space $\Tan_x(X)/\Tan_x(V).$

For $\lambda \in \C^\ast ,$ let $A_\lambda :\ \E \to  \E$ be the multiplication by $\lambda$ in fibers of $\E,$
that is, 
\begin{equation}\label{e:A_lambda} A_\lambda(y):=\lambda y\qquad\text{for}\qquad y\in \E.
\end{equation}
A  current $T$ on $\E$ is  said to be {\it $V$-conic} if $T$ is  invariant under  the action of  $A_\lambda,$ that is,
$(A_\lambda)_*T=T$ for all $\lambda\in\C^*.$
We identify $V$ with the zero section of $\E.$  
Let $\pi_0:\ \overline \E:=\P(\E\oplus\C)\to V$ be its canonical compactification. 

Let $T$ be a positive closed $(p,p)$-current on $X.$

When  $V=\{\text{a single point $x$}\},$ we can assume that $T$ is defined on an open neighborhood of $0$ in $\C^k.$  In this  context $\E=\C^k,$  and by  Harvey's exposition \cite{Harvey},
 when $\lambda$ goes to infinity, the domain of definition
of the current $T_\lambda := (A_\lambda)_* (T)$ converges to $\C^k .$ This family of currents
is relatively compact,  and any limit current $T_\infty$ for $\lambda \to\infty,$ is  
called a {\it tangent current} to $T.$ 
A tangent  current  is defined on the  whole $\C^k,$ and it is 
conic. 

 When $V$ has
positive dimension, 
We  expect as  in  Harvey's exposition \cite{Harvey} that  one can define  $T_\lambda$ and  obtain a   tangent  current  $T_\infty$ living on $\E.$
However,  a basic difficulty arises: in general, no neighbourhood of $V$ in $X$ is biholomorphic
to a neighbourhood of $V$ in $\E.$

To  encounter  this  difficulty,  Dinh and  Sibony propose  a  softer  notion:  {\it the  admissible maps.}  More  precisely, let $\tau$  be a diffeomorphism between a neighbourhood of $V$ in $X$ and a
neighbourhood of $V$ in $\E$ whose restriction to $V$ is identity. We assume that
$\tau$ is admissible in the sense that the endomorphism of $\E$ induced by the
differential of $\tau$  when restricted to $V$ is the identity map from $\E$ to $\E.$

 Fix $0\leq p\leq  k$ and set
\begin{equation}\label{e:m}\upm:= \min(l,k-p)\qquad\text{and}\qquad
  \lowm:=\max(0,l-p).
  \end{equation}
Here is  the main result of  Dinh and  Sibony.
\begin{theorem}\label{T:Dinh-Sibony}{\rm  (Dinh-Sibony \cite[Theorems 1.1, 4.6   and Definition 4.8]{DinhSibony18})}  Let $X,$ $V,$ $\E,$  $\overline \E,$ $A_\lambda$ and $\tau$ be as above. 
Let $T$ be  a positive closed $(p,p)$-current on $X.$
Assume in addition  that $X$ is  K\"ahler and  $\supp(T)\cap V$ is  compact. Then:
\begin{enumerate} \item The family
of currents $T_\lambda:= (A_\lambda)_* \tau_* (T )$ is relatively compact and any limit current, for
$\lambda\to\infty,$ is a positive closed $(p, p)$-current on $\E$ whose trivial extension is a
positive closed $(p, p)$-current on $\overline \E.$   Such a limit current $R$ is  called a {\rm tangent current to $T$ along $V.$} 
\item If $R$ is a tangent current to $T$ along $V$, then it is
$V$-conic, i.e., invariant under $(A_\lambda)_* ,$  and its de Rham cohomology class $\{R\}$ in the cohomology group with compact support
$H^{2p}_\comp (\overline\E, \C)$ does not depend on the choice of $\tau$ and $R.$  We denote $\{R\}$ by $\bfc^\DS(T,V)$ (or simply by $\bfc^\DS(T)$  if there is  no confusion on $V$), and call it  {\rm the  total  tangent class of $T$  along $V$}, or equivalently, {\rm  Dinh-Sibony (total) cohomology class of $T$ along $V.$}
In fact, $\bfc^\DS(T,V)\in H^{p,p}_\comp(\overline\E,\C).$

\item  Let $-h_{\overline \E}$ denote the  tautological class of the bundle $\pi_0:\  \overline \E\to  V.$ Then, by  Leray-Hirsch  theorem (see for example \cite{BottTu}), we have the following   decomposition
of Dinh-Sibony   class $\bfc^\DS(T,V)$:
$$
\bfc^\DS(T,V)=\sum_{j=\lowm}^\upm \pi_0^*(\bfc^\DS_j(T,V))\smile  h_{\overline \E}^{j-l+p},
$$
where  $\bfc^\DS_j(T,V)$  is a class in  $H^{2l-2j}_\comp(V,\C).$ Moreover, this decomposition is  unique.

\end{enumerate}
\end{theorem}

 When $V$  has positive dimension $l,$
 according   to  Dinh  and Sibony,  the    notion of Lelong number of the current $T$ at a single point should be replaced by
 the  family of  cohomology classes $\{ \bfc^\DS_j(T,V):\  \lowm \leq j\leq \upm\}$  given by Theorem \ref{T:Dinh-Sibony} (3) above. 
This is an important and original viewpoint  of  Dinh and Sibony. 
 
To prove  their theorem,  Dinh and  Sibony develop a cohomological calculus on positive closed currents $T$ such that   $\supp(T)\cap V$ is  compact.
It is  worth noting that later on  Vu \cite{Vu21a} weakens the assumption of K\"ahlerian  on $X.$

\begin{definition}\label{D:hor-dim} {\rm (Dinh-Sibony \cite[Definition 3.7]{DinhSibony18}) } \rm 
Let $X,V$ be  as  in Theorem \ref{T:Dinh-Sibony}.
 Let $T$  be a positive closed  $(p,p)$-current on $X$ such that $\supp(T)\cap V$ is compact.
 If $\bfc^\DS(T,V)\not=0,$ then we define {\rm the horizontal dimension} (or  {\rm h-dimension} for short) of $T$  along $V$ to be  the  maximal $j:\ \lowm\leq j\leq \upm$  such that $\bfc^\DS_j(T,V)\not=0,$   otherwise (i.e.  $\bfc^\DS(T,V)=0$) 
 we define  by convention the h-dimension of $T$ to be simply $\lowm.$
 Denote by  $\hbar$   the h-dimension of $T.$
\end{definition}

\begin{theorem}\label{T:Dinh-Sibony-hor}{\rm (Dinh-Sibony \cite[Lemma  3.8]{DinhSibony18}) } Let $X,V$ be  as  in Theorem \ref{T:Dinh-Sibony}.
Let $T$  be a positive closed  $(p,p)$-current on $X$  such that $\supp(T)\cap V$ is compact. Let $\hbar$  be the h-dimension of $T$ along $V.$
Let $R$ be a tangent  current to $T$ along $V.$
 Then the following  assertions hold:
 \begin{itemize}
 \item  If  $R\not=0,$ then  $\hbar$  is also  the maximal  $j$ such that $R\wedge \pi^*\omega^j\not=0,$
 otherwise  (i.e. $R$=0) we have  $\hbar=\lowm.$
 \item 
The class $\bfc^\DS_\hbar(T,V)$ is pseudo-effective  compactly supported in $V$, i.e., it contains a positive closed current  compactly supported in $V.$
\end{itemize}
\end{theorem}

\begin{theorem}\label{T:Dinh-Sibony-Siu}{\rm (Dinh-Sibony  \cite[Theorem   4.11]{DinhSibony18}) }
Let $X,V$ be  as  in Theorem \ref{T:Dinh-Sibony}.
Let $W$ be an  open subset in $X$  such that $W\cap V$ is relatively compact in $V.$
Let $T_n$ and $T$  be a positive closed  $(p,p)$-currents on $X$  with support in $W$  such that $T_n\to T.$ Let $\hbar$  be the h-dimension of $T$ along $V.$
 Then the following  assertions hold:
 \begin{itemize}
 \item  If $j$ is  an integer with  $ j>\hbar,$  then $\bfc_j^\DS(T_n,V)\to 0.$
 \item  If $\bfc_\hbar$  is a limit  class of the  sequence 
 $\bfc^\DS_\hbar(T_n,V),$  then  the class  $\bfc_\hbar$ and $ \bfc^\DS_\hbar(T,V)-\bfc_\hbar$ are pseudo-effective compactly supported in $V.$  
\end{itemize}
\end{theorem}

\begin{theorem}\label{T:Dinh-Sibony-product}{\rm (Dinh-Sibony  \cite[Lemma 3.4 and Definition 5.9]{DinhSibony18})}
 Let $X$ be a  compact K\"ahler manifold. For $1\leq j\leq m,$  let $T_j$ be a positive closed  current of bidegree $(p_j,p_j)$ on $X.$
 Consider $\T:=T_1\otimes\ldots\otimes T_m$ which is  a positive  closed  $(p,p)$-current on  $X^m,$
 where $p:=p_1+\cdots+p_m.$
 Let $\Delta:=\{(x,\ldots,x):\ x\in X\}$  be the diagonal of $X^m.$ 
 Let $\pi:\ \E\to\Delta$ be the normal bundle to $\Delta$ in $X^m.$
  Suppose that 
  
  \begin{enumerate}
   \item There exists a unique tangent current $\T_\infty$ to $\T$ along $\Delta$;
   
   \item The horizontal dimension of $\T$ along $\Delta$ is minimal, i.e. $\hbar= k-p.$
   
  \end{enumerate}
   Then there exists a unique positive  closed  $(p,p)$-current $S$ on $\Delta$  such that $\T_\infty=\pi^* S.$ 
\end{theorem}
  
  \begin{definition}  \label{D:DS-wedge-prod}\rm  Under the  hypothesis and  the  inclusion  of Theorem  \ref{T:Dinh-Sibony-product},
  we  say that {\it $T_1,\ldots,T_m$ are  wedgeable in the sense  of the   theory of tangent  currents.} Identifying $\Delta$ to $X,$
  Dinh-Sibony \cite[Definition 5.9]{DinhSibony18} define  $T_1\curlywedge\ldots  \curlywedge T_m:=S.$
\end{definition}

\subsection{Alessandrini-Bassanelli theory} 

 On the other hand,  Alessandrini and Bassanelli introduce in  \cite{AlessandriniBassanelli96} a remarkable  notion of Lelong number  of   positive plurisubharmonic currents. Recall  that  a  positive  current  $T$ is  said to be  plurisubharmonic if the current  $\ddc T$ is positive. 
 So all positive closed currents are positive plurisubharmonic.
 In the context of    Alessandrini and Bassanelli,     
$(X,V)$  is a  very special pair of manifolds (affine manifolds), however  they  allow a domain $B\Subset V$ and formulate their Lelong number along $B.$
This means that $\supp(T)\cap B$ may be noncompact in $B.$

\begin{theorem}\label{T:AB-1}{\rm  (Alessandrini and Bassanelli \cite[Theorem I and Definition 2.2]{AlessandriniBassanelli96})}
Consider $X=\C^k$ and $V$ is a linear complex subspace of dimension $l\geq 0. $   
  We use the coordinates $(z,w)\in\C^{k-l}\times \C^l$ so  that  $V=\{z=0\}.$ Let $0\leq p<k-l$ and  
  let $T$  be    a positive plurisubharmonic $(p,p)$-current on an open  neighborhood $\Omega$ of $0$ in $\C^k.$ 
Then, for every open ball $B$ in $V,$ $B\Subset \Omega,$ the  following limit  exists and is  finite
$$
 \nu_\AB(T,B):=\lim_{r\to 0+} {1\over r^{2(k-l-p)}}\int_{\Tube(B,r)}  T(z,w)\wedge (\ddc \|z\|^2)^{k-l-p}\wedge (\ddc \|w\|^2)^{l},
 $$
 where the  tube  $\Tube(B,r)$ of radius $r$ over $B$ is  given by 
 \begin{equation} \label{e:Tube-AB}
\Tube(B,r):=\left\lbrace (z,w)\in \C^{k-l}\times\C^{l}:\ \|z\|<r,\ w\in B\right\rbrace .
\end{equation}
 $\nu_\AB(T,B)$ is  called    the {\it  Alessandrini-Bassanelli's  Lelong  number of $T$ along  $B$}.
\end{theorem}

The  important viewpoint of Alessandrini-Bassanelli is  that when $V$ is of positive dimension,   tubular neighborhoods $\Tube(B,r)$ of $B$ and a  mixed form $  (\ddc \|z\|^2)^{k-l-p}\wedge (\ddc \|w\|^2)^{l}$ should replace
the usual balls $\B(x,r)$ around a single point $x$ with  the usual form  $  (\ddc \|z\|^2)^{k-p}.$
When $V$ is  a single point $\{x\}$ and $B=\{x\},$   Alessandrini-Bassanelli's  Lelong  number $\nu_\AB(T,x)$  coincides with the classical Lelong number
$\nu(T,x).$

 Alessandrini-Bassanelli's method relies on some Lelong-Jensen formulas  which can be obtained  from the usual Lelong-Jensen formula (see  \cite{Demailly, Skoda})  by slicing. 
They  also  characterize  this Lelong number geometrically in the sense   of Siu \cite{Siu}.
 
\begin{theorem}\label{T:AB-3} {\rm  (Siu \cite[Section 11]{Siu} for positive closed currents,   Alessandrini-Bassanelli \cite[Theorem II]{AlessandriniBassanelli96} for positive plurisubharmonic  currents)}
 Let $F:\ \Omega\to\Omega'$ be  a  biholomorphic map between open subsets of $\C^k.$
 If $T$ is a  positive  plurisubharmonic $(p,p)$-current on $\Omega$ and $x\in\Omega,$ then
 $$
 \nu(T,x)=\nu(F_*T,F(x)).
 $$
\end{theorem}
Hence, the limit \eqref{e:Lelong-number-point} does not depend on the choice of coordinates even  for   positive  plurisubharmonic currents. 
So,  the Lelong number  of a  positive plurisubhamonic   current at a single point is an intrinsic notion.
 
 Although  the  assumption   on  the  pair of manifolds $(X,V)$ in Theorem  \ref{T:AB-1} is  quite  restrictive and  this theorem   provides only one  Lelong number, Alessandrini-Bassanelli theory may be regarded as the first  effort to elaborate  the notion of  numerical Lelong numbers  when the dimension of $V$ is positive.

 \subsection{Organization of the  article and acknowledgments}
 The paper is  organized  as follows. The main results of this article is  divided into 5 groups.

 In Section  \ref{S:Lelong-numbers-main-results} we first  formulate  the concept of the generalized Lelong number of  a positive plurisubharmonic current  associated to a  closed smooth  $(j,j)$-form on the  base submanifold. This formulation permits us  to       state   the main results shortly afterwards.

  Section \ref{S:Background} lay  out the background  of   the  problems  considered 
 in this  work.     We first recall some basic facts  on currents  and positive  currents.
Next,  we recall  the  notion of strongly  admissible maps  introduced  in \cite{Nguyen21}.
The remainder of the section is  devoted  to recalling   various classes of positive currents introduced  in \cite{Nguyen21}.

Section  \ref{S:Lelong-Jensen}    recalls  several  Lelong-Jensen  type formulas  which  have been obtained  in \cite{Nguyen21}.
These formulas  play the  key role  in  the present  work.

 Section \ref{S:Finiteness}    begins by introducing  standard   settings for further technical  developments. The rest   the  section is  then devoted to
 admissible estimates, mass indicators and their finiteness. The new thing here   is that the $(j,j)$-forms $\omega^{(j)}$ are allowed to be simply closed smooth, whereas  in  our previous work in \cite{Nguyen21} they are  only  some powers of a K\"ahler form.

 In Section \ref{S:Lelong-numbers-for-closed-currents} we  adapt and improve  the technique developed in \cite{Nguyen21}
 in order  to   study  the generalized Lelong numbers for positive closed currents.
 We also  prove      Theorem \ref{T:Nguyen-2}.

 In Section \ref{S:Lelong-numbers-for-plurisubhamonic-currents}
 we  outline the adaptation  and the improvement of  the technique from  \cite{Nguyen21}
 in order  to   study  the generalized Lelong numbers for positive plurisubharmonic currents.
 This is necessary for  the proof of       Theorem \ref{T:Nguyen-1}. This, together with Theorem \ref{T:Nguyen-2},  constitutes  the first main result of the  article.  
 
 Section  \ref{S:Basic-formulas} establishes explicit formulas expressing the generalized Lelong numbers in terms of the tangent currents.
 
 Using these formulas,  Section \ref{S:H-dim-and-Siu} proves the second main result (Theorem  \ref{T:Nguyen-DS} and Theorem \ref{T:Nguyen-Siu}).
 
 Based on  
 Section  \ref{S:Basic-formulas}, Section \ref{S:Dinh-Sibony-vs-Lelong}  discovers the link between 
 the generalized  Lelong numbers and  Dinh-Sibony cohomology classes (Theorem \ref{T:Nguyen-Dinh-Sibony}).
 This is the third  main result of the  article. It is quite accidental and unexpecting that
 in the limiting   case  (that is, when $\lambda\to\infty$ in Theorem \ref{T:Dinh-Sibony}),
 knowing Dinh-Sibony cohomology classes are equivalent to knowing the generalized  Lelong numbers.   Since  the latter characteristic  numbers   are in fact  the limits of some
  generalized Lelong functionals (see, for example,  \eqref{e:Lelong-number-point},   \eqref{e:Lelong-number-point-bisbis(1)}, \eqref{e:Lelong-number-point-bisbis(2)},\eqref{e:Lelong-number-point-bis} in the case of a single point),
the advantage  of these  functionals is   that they help  us to quantify
 the former cohomology classes. This  idea  will be  illustrated  in the last  two sections.
 As an outcome  of Theorem \ref{T:Nguyen-Dinh-Sibony}, we  prove  in Theorem \ref{T:Nguyen-Siu-AB} that in some important situations (including the  context of Dinh-Sibony in Theorem \ref{T:Dinh-Sibony}), the generalized Lelong numbers are totally intrinsic.  It is 
 a natural  generalization of Theorem \ref{T:AB-3}  of Siu and  Alessandrini-Bassanelli.

 Section  \ref{S:Intersection} gives  a  solution  to the  above Problem 1. Namely,
 it provides  an  effective sufficient condition for defining the intersection of $m$ positive  closed  currents in the sense of Dinh-Sibony’s theory of tangent currents on a  compact K\"ahler manifold, see Theorem \ref{T:Nguyen-intersection}.  We  use the  second  main result (Theorem  \ref{T:Nguyen-DS}) and   a recent result  
 of  Nguyen-Truong \cite{NguyenTruong} on the uniqueness of tangent currents.
 
 Section  \ref{S:Continuity} gives  a  solution  to the  above Problem 2. 
  More  specifically,
 it gives us   an  effective sufficient condition for  the continuity of   the intersection of $m$ positive  closed  currents in  the above sense, see Theorem \ref{T:Nguyen-convergence}.  The proof  relies on an explicit formula of  $m$ positive  closed  currents in the sense of Dinh-Sibony’s theory of tangent currents (see Theorem  \ref{T:Nguyen-intersection-formula} below),   which is  inspired by 
 the  work of  Dinh-Nguyen-Vu \cite{DinhNguyenVu} on the theory of super-potentials.
 The  section is  concluded  with a reformulation of Theorem \ref{T:Nguyen-convergence} in terms of the  blow-up along the diagonal   (see Theorem \ref{T:Nguyen-convergence-bis}).

 \smallskip
 
\noindent
{\bf Acknowledgments. }   The    author   acknowledges support from the Labex CEMPI (ANR-11-LABX-0007-01)
and from  the project QuaSiDy (ANR-21-CE40-0016).
The paper was partially prepared 
during the visit of the  author at the Vietnam  Institute for Advanced Study in Mathematics (VIASM). He would like to express his gratitude to this organization 
for hospitality and  for  financial support.

 
 \section{The generalize Lelong numbers and  statement of main results}
 \label{S:Lelong-numbers-main-results}

\subsection{The generalized Lelong numbers}\label{SS:generalized-Lelong-numbers}

Let $X$ be a complex  manifold of dimension $k,$
 $V\subset  X$    a submanifold  of dimension $l,$ and $B\subset V$ a  relatively compact piecewise  $\Cc^2$-smooth open subset.
Let  $V_0$ be  a relatively compact  open subset of $V$ such that  $B\Subset V_0.$ Let $\omega$ be a Hermitian form on $V.$
Let $\tau:\ U\to\tau(U)$ be  an admissible  map along $B$
 from an open neighborhood $U$ of $\overline B$ in $X,$ see Definition \ref{D:admissible-maps} below.

  Denote by $\pi:\ \E\to V$ the canonical projection.
Consider a Hermitian metric  $h=\|\cdot\|$  on the  vector bundle  $\E_{\pi^{-1}(V_0)}$  and    let   $\varphi:\ \E_{\pi^{-1}(V_0)}\to \R^+$ be the function defined by  
\begin{equation}\label{e:varphi-spec}
 \varphi(y):=\|y\|^2\qquad \text{for}\qquad  y\in \pi^{-1}(V_0)\subset \E.
\end{equation}
Consider also the following  closed  $(1,1)$-forms on  $ \pi^{-1}(V_0)\subset \E $
\begin{equation}\label{e:alpha-beta-spec}
 \alpha:=\ddc\log\varphi\quad\text{and}\quad \beta:= \ddc\varphi.
\end{equation} 
So, for every $x\in V_0\subset  X$  the  metric $\| \cdot\|$  on the fiber $\E_x\simeq \C^{k-l}$ is  an Euclidean metric (in a suitable basis). In particular, we have 
\begin{equation}\label{e:varphi_bis-spec}  \varphi(\lambda y)=|\lambda|^2\varphi(y)\qquad\text{for}\qquad  y\in \pi^{-1}(V_0)\subset \E,\qquad\lambda\in\C. 
\end{equation}
For $r>0$ 
consider the  following {\it tube with base $B$ and radius $r$}
\begin{equation}
\label{e:tubular-nbh-0}
\Tube(B,r):=\left\lbrace y\in \E:\    \pi(y)\in B\quad\text{and}\quad  \|y\|<r  \right\rbrace.
\end{equation}
So   this is  a  natural generalization of Euclidean  tubes   considered by  Alessandrini-Bassanelli  in \eqref{e:Tube-AB}.
For for  all $0\leq s<r<\infty,$  define also the  corona tube
\begin{equation}\label{e:tubular-corona-0}\Tube(B,s,r):=\left\lbrace y\in \E:\   \pi(y)\in B\quad\text{and}\quad   s<\|y\|<r \right\rbrace.
\end{equation}
Since  $V_0\Subset V,$ there  is  a constant $c>0$ large enough such that 
$c\pi^*\omega+\beta$
is  positive on $\pi^{-1}(V_0).$  Moreover,   the latter  form defines
a K\"ahler metric there  if  $\omega$ is  K\"ahler on $V_0.$

 Let $\bfr$ be small enough such that $\Tube(B,\bfr)\subset  \tau(U),$  see \eqref{e:tubular-nbh-0}.
 Fix $0\leq p\leq k.$
 Let $T$ be a real  current of degree $2p$ and  of  order $0$    on    $U.$
 For $0\leq j\leq \upm,$  let $\omega^{(j)}$ be  a   closed  smooth real $(j,j)$-form on $V_0$ such that
 $\omega^{(1)}=\omega.$
 For $0\leq j\leq \upm$ and  $0<r\leq \bfr,$   consider 
\begin{equation}\label{e:Lelong-numbers}
 \nu(T,B,\omega^{(j)},r,\tau,h):=  {1\over r^{2(k-p-j)}}\int_{\Tube(B,r)} (\tau_*T)\wedge \pi^*(\omega^{(j)}) \wedge \beta^{k-p-j}.  
\end{equation} 

Let  $0\leq j\leq \upm.$  For  $0<s<r\leq \bfr,$   consider 
\begin{equation}\label{e:Lelong-corona-numbers}
 \kappa(T,B,\omega^{(j)},s,r,\tau,h):=   \int_{\Tube(B,s,r)} (\tau_*T)\wedge \pi^*(\omega^{(j)}) \wedge \alpha^{k-p-j}. 
\end{equation} 
For $0<r\leq\bfr,$   we consider
\begin{equation}\label{e:Lelong-log-numbers-bullet}
 \kappa^\bullet(T,B,\omega^{(j)},r,\tau,h):=\kappa(T,B,\omega^{(j)},0,r,\tau,h)= \int_{\Tube(B,r)\setminus V} (\tau_*T)\wedge \pi^*(\omega^{(j)}) \wedge \alpha^{k-p-j};  
\end{equation}
we  also consider
\begin{equation}\label{e:Lelong-log-numbers}
 \kappa(T,B,\omega^{(j)},r,\tau,h):= \int_{\Tube(B,r)} (\tau_*T)\wedge \pi^*(\omega^{(j)}) \wedge \alpha^{k-p-j},  
\end{equation} 
provided that the RHS side makes sense according  to the  following definitions.

 

\begin{definition}\label{D:Lelong-log-numbers(1)}\rm
We say that  \eqref{e:Lelong-log-numbers} holds in the spirit of \eqref{e:Lelong-number-point-bisbis(1)}
if $T=T^+-T^-$  in an open neighborhood of $\overline B$ in $X$ and $T^\pm\in\SH^{p;m,m'}(B)$ (resp.  $T^\pm\in\PH^{p;m,m'}(B),$  resp.
$T^\pm\in\CL^{p;m,m'}(B)$ for some  suitable integers $0\leq m'\leq m$) with the corresponding sequences of approximating forms   $(T^\pm_n)_{n=1}^\infty,$
and  for any  such  forms  $(T^\pm_n),$  the two  limits  on the  following  RHS exist and  are finite   
\begin{equation}\label{e:Lelong-log-numbers(2)}
\kappa(T,B,\omega^{(j)},r,\tau):= \lim\limits_{n\to\infty} \kappa(T^+_n,B,\omega^{(j)},r,\tau)-\lim\limits_{n\to\infty}\kappa(T^-_n,B,\omega^{(j)} ,r,\tau),
\end{equation}
and the  value on the RHS is  independent of the choice of $(T^\pm_n)_{n=1}^\infty.$
\end{definition}

\begin{definition}\label{D:Lelong-log-numbers(2)}\rm
We say that  \eqref{e:Lelong-log-numbers} holds in the spirit of \eqref{e:Lelong-number-point-bisbis(2)}
if  the    limit  on the  following  RHS exists and  is finite   
\begin{equation}\label{e:Lelong-log-numbers(1)}
\kappa(T,B,\omega^{(j)},r,\tau):= \lim\limits_{\epsilon\to 0+}  \int_{\Tube(B,r)} (\tau_*T)\wedge \pi^*(\omega^{(j)}) \wedge \alpha^{k-p-j}_\epsilon.
\end{equation}
Here, the  smooth  form $\alpha_\epsilon$ is given  by \eqref{e:alpha-beta-eps} below.
\end{definition}

\begin{remark}\rm
In   \cite[Section 3]{Nguyen21} the author introduced  the following particular  case:
\begin{equation}\label{e:standard-omega-j}
 \omega^{(j)}:=\omega^{j}\qquad\text{for}\qquad \lowm\leq j\leq \upm.
\end{equation}
 In this  case,  we  also denote   $\nu(T,B,\omega^j,r,\tau,h)$    by $\nu_j(T,B,\omega,r,\tau,h).$
 Moreover,  when $j=\upm$ we also denote $\nu(T,B,\omega^\upm,r,\tau,h)$    by $\nu_\top(T,B,\omega,r,\tau,h).$
\end{remark}

 \begin{remark}\rm
  When  there is  no confusion about the choice  of  a metric $h$ and an  admissible map  $\tau,$
  we often write $\nu_j(T,B,\omega)$  (resp. $\nu_j(T,B,\omega,r)$) for short  instead of  $\nu_j(T,B,\omega,\tau,h)$   (resp. $\nu_j(T,B,\omega,r,\tau,h)$).
 \end{remark}

 \subsection{First  main result: existence of the generalized Lelong numbers}
 \label{SS:First-Main-Result}

The first theorem  of the first main result of the  article is  the following one which deals  with positive  plurisubharmonic/positive  pluriharmonic  currents.
 
 \begin{theorem}\label{T:Nguyen-1}  
  Let $X,$ $V$ be as  above  and  suppose that $(V,\omega)$ is  K\"ahler,  
  and  that  $B$ is   a piecewise $\Cc^2$-smooth open subset  of $V$ and that  there exists a strongly   admissible map for $B.$
  For $0\leq j\leq \upm,$  let $\omega^{(j)}$ be  a  closed  smooth real $(j,j)$-form on $V_0.$
  Let $T$ be  a   positive plurisubharmonic  $(p,p)$-current  on a neighborhood of $\overline B$ in $X$ such that  $T=T^+-T^-$  for some $T^\pm\in\SH^{p;3,3}( B).$
  Then  the following  assertions  hold:
  \begin{enumerate}
  \item For every $\lowm\leq j\leq \upm,$  the following limit  exists and is  finite
  $$    \nu(T,B,\omega^{(j)},h):=\lim\limits_{r\to 0+}\nu(T,B,\omega^{(j)},r,\tau ,h)                              $$
   for  all strongly  admissible maps  $\tau$     for $B$ and for all  Hermitian  metrics  $h$ on $\E.$   
   
  
  \item  
  The following  equalities hold
  \begin{eqnarray*}\lim\limits_{r\to 0+}\kappa^\bullet(T,B,\omega^{(j)},r,\tau,h)&=&0,\\
  \lim\limits_{r\to 0+}\kappa(T,B,\omega^{(j)},r,\tau,h)&=&\nu(T,B,\omega^{(j)},h),
  \end{eqnarray*}
  for all $\lowm\leq j\leq \upm$ with $j>l-p,$  and  for all strongly  admissible maps $\tau$  for $B$ and for all Hermitian  metrics $h$ on $\E.$

\item The  real numbers  $\nu(T,B,\omega^{(j)},h)$ are intrinsic, that is, they are
  independent of  the choice  of  $\tau.$

  \item  If  $\tau$ is a holomorphic admissible map  and if  $T^\pm$ belong  only  to  the class $\SH^{p;2,2}(\overline B),$
  then the above three assertions (1)--(3) still hold for $j=\upm.$
    
\item If  instead of the above  assumption on $T,$  we assume  that  $T$ is a   positive pluriharmonic  $(p,p)$-current  on a neighborhood of $\overline B$ in $X$ such that  $T=T^+-T^-$  for some $T^\pm\in\PH^{p;2,2}( B),$  then  all the above  assertions still  hold.  
  \end{enumerate}
 \end{theorem}

 \begin{remark}\rm
At the first glance    the condition  $T=T^+-T^-$  looks artificial. However, it is unavoidable in practice, see Theorem \ref{T:strong-admissible-maps} and  \ref{T:approximation} below. 

The  above  results  generalizes  \cite[Tangent Theorem I (Theorem 1.8)]{Nguyen21}  except  the assertions about   the  tangent currents.
\end{remark}
 
 The  second theorem of our first main result deals with positive closed  currents.  
 \begin{theorem}\label{T:Nguyen-2} 
  Let $X,$ $V$ be as  above.   Assume that there is a Hermitian  metric $\omega$ on $V$  for which
  $\ddc \omega^j=0$ for $\lowm\leq j\leq \upm-1.$
  For $\lowm\leq j\leq \upm,$  let $\omega^{(j)}$ be  a    closed smooth real $(j,j)$-form   on $V_0.$ 
 Assume
  also that  $B$ is   a piecewise $\Cc^2$-smooth open subset  of $V$ and that  there exists a strongly   admissible map for $B.$
  Let $T$ be  a   positive closed  $(p,p)$-current  on a neighborhood of $\overline B$ in $X$ such that  $T=T^+-T^-$  for some $T^\pm\in\CL^{p;2,2}( B).$
  Then  the following  assertions  hold:
   \begin{enumerate}

    \item  For  $\lowm\leq j\leq \upm$ and for $0<r_1<r_2\leq\bfr,$
  $$ \nu(T,B,\omega^{(j)},r_2,\tau,h)-\nu(T,B,\omega^{(j)},r_1,\tau,h)=\kappa(T,B,\omega^{(j)}, r_1,r_2,\tau,h)+O(r_2).
  $$
   \item For every $\lowm\leq j\leq \upm,$  the following limit  exists and is  finite
  $$    \nu(T,B,\omega^{(j)},h):=\lim\limits_{r\to 0+}\nu(T,B,\omega^{(j)},r , \tau ,h)                              $$
   for  all strongly  admissible maps  $\tau$     for $B$ and for all  Hermitian  metrics  $h$ on $\E.$   

  \item  The following  equalities hold
  \begin{eqnarray*}
   \lim\limits_{r\to 0+}\kappa^\bullet(T,B,\omega^{(j)},r,\tau,h)&=&0,\\
   \lim\limits_{r\to 0+}\kappa(T,B,\omega^{(j)},r,\tau,h)&=&\nu(T,B,\omega^{(j)},h),
  \end{eqnarray*}
 for all $\lowm\leq j\leq \upm$ with $j>l-p,$  and for all strongly  admissible maps $\tau$  for $B$ and for all Hermitian  metrics $h$ on $\E.$

\item The  real numbers  $\nu(T,B,\omega^{(j)},h)$ are intrinsic, that is, they are
  independent of  the choice  of  $\tau.$

  \item  If  $\tau$ is a holomorphic admissible map  and if  $T^\pm$ belong  only  to  the class $\CL^{p;1,1}(\overline B),$
  then the above four assertions (1)--(4) still hold.

 \item  If  instead of the above  assumption on $\omega$ and  $T,$  we assume  that  the form $\omega$ is  K\"ahler    and $T$ is a   positive closed  $(p,p)$-current  on a neighborhood of $\overline B$ in $X$ such that  $T=T^+-T^-$  for some $T^\pm\in\CL^{p;1,1}( B),$  then all the above  assertions still hold.   
  \end{enumerate}
 \end{theorem}
 
\begin{remark} \rm  The above result generalizes  \cite[Tangent Theorem II (Theorem 1.11)]{Nguyen21}.
 \end{remark}

 \begin{remark}\label{R:optimal-T:Nguyen2}\rm  We keep  the  hypothesis  on $X,V,$ and $\omega$ as  in  Theorem \ref{T:Nguyen-2}. 
 For $\lowm\leq j\leq \upm,$  let $\omega^{(j)}$ be  a    smooth real $(j,j)$-form (not necessarily closed) on $V_0$ such that   $\ddc(\omega^{(j)}\wedge \omega^q)=0$  for $0\leq  q\leq \upm-j-1.$
  Then  assertions  (1)--(5) of Theorem \ref{T:Nguyen-2} still hold.
 \end{remark}

 \begin{definition}\label{D:Lelong-numbers}\rm Let $h$ be a  fixed  Hermitian metric on the vector bundle $\E.$
 The value $\nu(T,B,\omega^{(j)},h)$  is called {\it the  (generalized)  Lelong number of $T$ along $B$ with respect to the form $\omega^{(j)}.$}
 The   set of real numbers  $\{\nu(T,B,\omega^{(j)},h):$\ $\lowm\leq j\leq \upm\}$  are called {\it  the Lelong numbers of $T$ along $B$
  with respect to the set of forms $\{\omega^{(j)}\}.$}  
The   number $\nu(T,B,\omega^{(\upm)},h) $ is  called  {\it  the top Lelong numbers of $T$ along $B$
with respect to the form $\omega^{(\upm)}$},
it is  also denoted by  $\nu_\top(T,B,\omega^{(\upm)},h).$
\end{definition}

 \subsection{Second main results: horizontal dimension and a Siu's upper-semicontinuity type theorem}\label{SS:Second-main-results}

  Let $X,$ $V$ be as  above  and  suppose that $(V,\omega)$ is  K\"ahler,  
  and  that  $B$ is   a piecewise $\Cc^2$-smooth open subset  of $V$ and that  there exists a strongly   admissible map for $B.$
  Let $T$ be  a   positive plurisubharmonic  $(p,p)$-current  on a neighborhood of $\overline B$ in $X$ such that  $T=T^+-T^-$  for some $T^\pm\in\SH^{p;3,3}( B).$
 By Tangent Theorem I (Theorem 3.8) in \cite{Nguyen21}, let $T_\infty$ be a  tangent current to $T$
 along $B,$ that is,  $T_\infty=\lim_{n\to\infty} T_{\lambda_n}$ for some $(\lambda_n)\nearrow\infty,$  where
  $T_\lambda:= (A_\lambda)_* \tau_* (T ).$ 
  We know   by  this  theorem that $T_\infty\wedge \pi^*\omega^\lowm$ is $V$-conic positive pluriharmonic.
  If   $T$ is a   positive pluriharmonic  $(p,p)$-current  on a neighborhood of $\overline B$ in $X$ such that  $T=T^+-T^-$  for some $T^\pm\in\PH^{p;2,2}( B),$
  then we also  know by this theorem that  $T_\infty$ is $V$-conic positive  pluriharmonic.
 
\begin{definition} {\rm The  horizontal dimension $\hbar$ of $T$ along $B$ is the largest integer $j\in[\lowm,\upm]$ such that $T_\infty\wedge \pi^*\omega^j\not=0$  if it exists, otherwise  we set  $\hbar:=\lowm.$}
\end{definition}

\begin{theorem}\label{T:Nguyen-DS}
 We keep  the above  hypothesis and notation.
Then:
\begin{enumerate}

\item $T_\infty\wedge \pi^*\omega^{(j)}$ is $V$-conic  pluriharmonic for every    closed  smooth real  $(j,j)$-form $\omega^{(j)}$ on $V$ with $\lowm \leq j\leq \upm.$

\item    $\nu(T,B,\omega^{(q)},h)=0$ for $\hbar<q\leq \upm$ and for every  closed  smooth real $(q,q)$-form $\omega^{(q)}$ on $V.$

\item  The horizontal dimension of $T$ along $B$ is  also the smallest integer $\hbar\in[\lowm,\upm]$ such that
   $\nu_q(T,B,\omega,h)=0$ for $\hbar<q\leq \upm.$  

\item If  $\omega^{(\hbar)}$ is  strongly  positive, then  $\nu(T,B,\omega^{(\hbar)},h)$ is  nonnegative.
If  $\omega^{(\hbar)}$ is  strictly  positive, then  if $T_\infty\not=0$ then  $\nu(T,B,\omega^{(\hbar)},h)>0,$
  else (that is, $T_\infty=0$)   $\nu(T,B,\omega^{(\hbar)},h)=0.$
 
  \end{enumerate}
\end{theorem}

\begin{theorem}\label{T:Nguyen-Siu}
 We keep  the above  hypothesis and notation. Suppose that $X$ is  K\"ahler. Let $U$ be an open neighborhood of $\overline B$ in $X,$ and $W$ be an open neighborhood of $\partial B$ in $U.$  
Let $T_n, T\in \CL^p(X)$ such that $T_n\to T$ as $n\to\infty$ and $\supp T_n \subset  U $ and  $\supp T_n\cap W=\varnothing.$ Let $\hbar$ be the horizontal dimension  of $T$ along $B.$ 
Then
\begin{enumerate}
 \item If $j>\hbar,$ then  $\lim_{n\to\infty} \nu(T_n,B,\omega^{(j)},h)= 0$ for every  closed  smooth real $(j,j)$-form $\omega^{(j)}$ on $V.$
 \item   For every strongly positive   closed  smooth $(\hbar,\hbar)$-form $\omega^{(\hbar)}$ on $V$ (see Definition \ref{D:strong-pos}), we have   \begin{eqnarray*} \liminf_{n\to\infty} \nu(T_n,B,\omega^{(\hbar)},h)&\geq& 0,\\ 
  \nu(T,B,\omega^{(\hbar)},h) &\geq & \limsup_{n\to\infty} \nu(T_n,B,\omega^{(\hbar)},h).
    \end{eqnarray*}
\end{enumerate}
\end{theorem}

\subsection{Third main result: Dinh-Sibony classes versus  the generalized Lelong numbers
}
   For  a  smooth closed $(p,q)$-form $\gamma$       on  complex manifold $X,$ let $\{\gamma\}$ be  its  class
in $H^{p,q}(X).$   If $\gamma$ is  moreover,  compactly supported,  let $\{\gamma\}$ be  its  class
in $H^{p,q}_\comp(X).$   
The cup-product $\smile$ on $H^{p,q}_\comp (X) \times H^{k-p,k-q} (X)$ is defined by
$$
(\{\gamma\},\{\gamma'\}) \mapsto \{\gamma\} \smile \{\gamma' \} :=\int_X\gamma\wedge\gamma',
$$
where $\gamma$ and $\gamma'$ are  smooth closed forms. The last integral depends only on
the classes of $\gamma$  and $\gamma'.$ The bilinear form $\smile$ is non-degenerate and induces
a canonical  isomorphism  between  $H^{p,q}_\comp (X)$ and its  dual  $H^{k-p,k-q} (X)^*$  (Poincar\'e duality). In
the definition of $\smile$ one can take $\gamma'$ smooth and $\gamma$ a current in the sense of
de Rham.   Assume that $X$ is  compact K\"ahler and $\gamma$ is a 
$\ddc$-closed $(p, q)$-current. 
 Then 
by the $\ddc$-lemma, the integral $\int_X \gamma \wedge \gamma'$ is also independent of the
choice of $\gamma'$ smooth and closed in a fixed cohomology class. So, using the above isomorphism,
one can associate to  such $\gamma$ a class $\{\gamma\}$ in $H^{p,q} (X).$

\begin{theorem}\label{T:Nguyen-Dinh-Sibony}
  Let $X$ be a  K\"ahler manifold,   and $V\subset X$  a    submanifold 
  of dimension $1\leq l<k,$ and $0\leq p\leq k.$ Let $-h_{\overline \E}$ denote the  tautological class of the bundle $\pi_0:\  \overline \E\to  V.$
  For $\lowm\leq j\leq \upm,$  let  $\omega^{(j)}$ be a closed   smooth real $(j,j)$-form on $V.$ Then:
 \begin{enumerate} 
  \item 
For any  current
  $T\in \CL^{p}(X)$ such that $\supp (T)\cap \supp ( \omega^{(\hbar)})\cap V\Subset V,$ we have,  by using the notations in Theorem \ref{T:Dinh-Sibony}:
 \begin{equation}\label{e:Lelong-numbers-vs-DS}
 \nu(T,V,\omega^{(j)},h)=\sum_{i=j}^\upm \pi^*_0 \bfc^\DS_{i}(T,V)\smile \pi^*_0 \{\omega^{(j)}\} \smile h_{\overline\E}^{k-l+i-j}, \qquad\forall  \lowm\leq j\leq \upm.
  \end{equation}
 \item  Suppose that $T\in \PH^p(X)$ and  $X$ is  compact.
 We know  by the  discussion at the  beginning of Subsection \ref{SS:Second-main-results}
   that any tangent current $T_\infty$ to $T$ along $V$  is $V$-conic positive  pluriharmonic in $\E,$
   and  therefore, by the  discussion preceeding this theorem applied to the compact K\"ahler manifold $\overline\E,$  $T_\infty$ 
 defines a  class $\{T_\infty\}$ in $H^{p,p}(\overline\E).$
    We denote $\{T_\infty\}$ by $\bfc^\DS(T,V)$ (or simply by $\bfc^\DS(T)$  if there is  no confusion on $V$), and call it  {\rm the  total  tangent class of $T$  along $V$}, or equivalently, {\rm  Dinh-Sibony (total) cohomology class of $T$ along $V.$}
In fact, $\bfc^\DS(T,V)\in H^{p,p}(\overline\E,\R).$
Therefore, by  Leray-Hirsch  theorem,
  we have the following   unique decomposition
of Dinh-Sibony   class $\bfc^\DS(T,V)$ as in Theorem \ref{T:Dinh-Sibony}:
$$
\bfc^\DS(T,V)=\sum_{j=\lowm}^\upm \pi_0^*(\bfc^\DS_j(T,V))\smile  h_{\overline \E}^{j-l+p},
$$
where  $\bfc^\DS_j(T,V)$  is a class in  $H^{l-j,l-j}(V,\C).$ Moreover,  the following identity holds:
 \begin{equation}\label{e:Lelong-numbers-vs-DS-bis}
 \nu(T,V,\omega^{(j)},h)=\sum_{i=j}^\upm  \pi^*_0\bfc^\DS_i(T,V)\smile  \pi^*_0  \{\omega^{(j)}\} \smile h_{\overline\E}^{k-l+i-j}, \qquad\forall  \lowm\leq j\leq \upm.
  \end{equation}
  In particular, $\bfc^\DS(T,V)$ is  intrinsic, i.e., it does not depend on the choice of  a tangent current $T_\infty.$
  \end{enumerate}
 \end{theorem}

 \begin{remark}
 \rm
  In the context of Dinh-Sibony (that is,  under assertion (1)),   knowing  Dinh-Sibony  cohomology classes of $T$  is {\rm \bf equivalent} to knowing the generalized Lelong numbers of $T.$
 Indeed, we use, for  $\lowm\leq j\leq \upm,$  several forms  $\omega^{(j)}_s$  such that   the classes $\{\omega^{(j)}_s\}$'s  span $H^{j,j}(V)$ (see  the proof of  
 Theorem \ref{T:Nguyen-Dinh-Sibony}).
 \end{remark}
 
 \begin{theorem}\label{T:Nguyen-Siu-AB}
    Under  the assumption  of Theorem \ref{T:Nguyen-Dinh-Sibony}, then    $\nu(T,B,\omega^{(\hbar)},h)$  is  totally intrinsic, i.e.  it is independent of  the choice  of both   $\tau$ and $h.$ So we will denote  $\nu(T,B,\omega^{(\hbar)},h)$ simply by  $\nu(T,B,\omega^{(\hbar)}).$
 \end{theorem}

 \begin{remark}
  \rm In  \cite[Theorem 1.8 (4) (Tangent Theorem I) and Theorem 1.11 (4) (Tangent Theorem II)]{Nguyen21}
  we only proved that the top Lelong number  $\nu_\top(T,B):=\nu_\upm(T,B,\tau,h)$ is  totally intrinsic. But in the latter theorems, there are less assumptions on $T$  and on $X$ than  those of Theorem \ref{T:Nguyen-Siu-AB}.
 \end{remark}

 \subsection{Fourth main results:  intersection theory  and  effective criteria}\label{SS:Fourth-main-results}

 Let $(X,\omega)$  be a compact K\"ahler manifold of dimension $k.$
 Consider $m$ integers $p_1,\ldots, p_m\geq 1$  such that   $p:=p_1+\ldots+p_m\leq k.$
 Let $\Delta:=\{(x,\ldots,x):\ x\in X\}$  be the diagonal of $X^m,$ and $\omega_\Delta$ be a K\"ahler form on $\Delta,$   and $\tau$  be a  strongly admissible map  along $\Delta$ in $X^m.$ 
 Let $\pi:\ \E\to\Delta$ be the normal bundle to $\Delta$ in $X^m.$
 Let $h$ be a Hermitian metric on $\E.$
Let $\dist(\bfx,\Delta)$ be the distance  from a point  $\bfx\in X^m$ to $\Delta.$
We may assume that  $\dist(\cdot,\Delta)\leq 1/2.$  So $-\log\dist(\cdot,\Delta) \cdot\T$
is a positive $(p,p)$-current on $X^m.$

  \begin{theorem}\label{T:Nguyen-intersection}
  Let 
  $T_j\in \CL^{p_j}(X)$ for $1\leq j\leq m$ with $m\geq 2,$ and  
consider $\T:=T_1\otimes\ldots\otimes T_m\in \CL^p(X^m).$
  Suppose that
  \begin{enumerate}
   \item $\kappa^\bullet_j(-\log\dist(\cdot,\Delta)\cdot\T,\Delta,\omega_\Delta,\bfr,\tau,h)<\infty$ for some  $\bfr>0$ and for all  $ k-p< j\leq k-\max_{1\leq i\leq m} p_i$;
   
   \item  $\nu_j(\T,\Delta,\omega_\Delta,\tau,h)=0$ for all  $ k-p<j\leq k-\max_{1\leq i\leq m} p_i.$
   \end{enumerate}
  Then   $T_1\curlywedge\ldots  \curlywedge T_m$  exists in the sense of  the theory of tangent currents (see Definition \ref{D:DS-wedge-prod}). 
  \end{theorem}

\begin{remark}   \label{R:check-condition-independent-of-choices}
  \rm  Assumption (1) can be  checked using  an arbitrary  finite  cover of $\Delta$ by  local holomorphic charts $(B_i)_{i\in I}$ in $\Delta,$ and over  each chart $B_i$ we may use an arbitrary Hermitian  metric $h_i$ for $\E|_{B_i}$     (see Remark \ref{R:four-equivalent-conds} and Lemma \ref{L:comparison-two-tau-two-h}). 
  More  concretely, assumption (1) is equivalent  to  the  following: $\kappa^\bullet_j(-\log\dist(\cdot,\Delta)\cdot\T,B_i,\omega_\Delta,\bfr,\tau_i,h_i)<\infty$ for some  $\bfr>0$ and for all  $ k-p< j\leq k-\max_{1\leq s\leq m} p_s$ and $i\in I.$

  By  \cite{DinhNguyenVu}, if    $m=2$ and if the super-potential of $T_1$ is  continuous,
then   $T_1 \curlywedge T_2$  exists in the sense of the theory of tangent currents (see Definition \ref{D:DS-wedge-prod}) for all $T_2\in \CL(X).$ Theorem  \ref{T:Nguyen-intersection} provides a more general  result  than  the last one.
\end{remark}
  
\begin{remark}
  \rm At  the end of the article, we  will  discuss another  version  of Theorem \ref{T:Nguyen-intersection} using the blow-up  along the diagonal $\Delta$ in $X^m,$  see Theorem \ref{T:Nguyen-intersection-bis} below.
 \end{remark}
 
 \subsection{Fifth main results:  continuity  of  Dinh-Sibony intersection}\label{SS:Fifth-main-results}
  
  \begin{theorem}\label{T:Nguyen-convergence}
  Let $m\geq 2$ be an integer, and let 
  $(T_{j,n})_{n=0}^\infty$ be a sequence  of currents  in $ \CL^{p_j}(X)$ for $1\leq j\leq m.$ Consider, for $n\in\N,$   
 $\T_n:=T_{1,n}\otimes\ldots\otimes T_{m,n}\in \CL^p(X^m).$ 
For $0<r\leq\bfr$  and $j$ with $ k-p<j\leq k-\max_{1\leq i\leq m} p_i,$ set
\begin{equation}\label{e:vartheta_j }\vartheta_j(r):=\sup_{n\in\N}\kappa_j(-\log\dist(\cdot,\Delta)\cdot\T_n,\Delta,\omega_\Delta,r,\tau,h).
\end{equation}
  Suppose that \begin{enumerate}
  \item  $T_{j,n}\to  T_j$ weakly as $n$ tends to infinity for $1\leq j\leq m;$   
   \item $\lim_{r\to 0}\vartheta_j(r)=0$ for $ k-p<j\leq k-\max_{1\leq i\leq m} p_i.$
   \end{enumerate}
  Then    $T_{1,n}\curlywedge\ldots  \curlywedge T_{m,n}$   and  $T_1\curlywedge\ldots  \curlywedge T_m$  exist in the sense of Dinh-Sibony's theory of tangent currents (see Definition \ref{D:DS-wedge-prod}).  Moreover, $T_{1,n}\curlywedge\ldots  \curlywedge T_{m,n}$   converge weakly to $T_1\curlywedge\ldots  \curlywedge T_m$ as $n$ tends to infinity.
  \end{theorem}

\begin{remark}   
  \rm  Assumption (2) can be  checked using  an arbitrary  finite  cover of $\Delta$ by  local holomorphic charts $(B_i)_{i\in I}$ in $\Delta,$  and over  each chart $B_i$ we may use an arbitrary Hermitian  metric $h_i$ for $\E|_{B_i}.$    See Proposition \ref{P:two-equivalent-conds} and Lemma \ref{L:comparison-two-tau-two-h} below,  see also Remark \ref{R:check-condition-independent-of-choices} above.  More  concretely, assumption (2) is equivalent  to  the  following:
  $$\lim_{r\to 0}\vartheta^{(i)}_j(r)=0\quad\text{for}\quad k-p<j\leq k-\max_{1\leq s\leq m} p_s\quad\text{and}\quad i\in I.$$
  Here,  
  $$\vartheta^{(i)}_j(r):= \sup_{n\in\N} \kappa^\bullet_j(-\log\dist(\cdot,\Delta)\cdot\T,B_i,\omega_\Delta,\bfr,\tau_i,h_i) .$$
\end{remark}
  
 \begin{remark}
  \rm At  the end of the article, we  will  discuss another  version  of Theorem \ref{T:Nguyen-convergence} using the blow-up  along the diagonal $\Delta$ in $X^m,$  see Theorem \ref{T:Nguyen-convergence-bis} below.
 \end{remark}

 
 \section{Background}\label{S:Background}

\subsection{Currents  and positive  currents}\label{SS:currents}


 Let $X$ be a  complex manifold of dimension $k$ endowed with  a Hermitian metric $\beta.$
 \begin{definition}\rm \label{D:strong-pos} A  
 real $(p,p)$-current $T$ on $X$ is said to be {\rm strongly positive} if 
for every smooth test $(n-p,n-p)$-form $\phi$ compactly supported in $X$ which is also  a positive current in the  sense of Definition \ref{D:positive-currents},
we have $\langle T,\phi\rangle\geq 0.$
  \end{definition}
  
  Let  $T$ be  a positive  $(p,p)$-current on  $X.$ Then  
$T\wedge \beta^{k-p}$ is a positive measure on
$X$. The mass of  $T\wedge \beta^{k-p}$
on a measurable set $A$ is denoted by $\|T\|_A$ and is called {\it the mass of $T$ on $A$}.
{\it The mass} $\|T\|$ of $T$ is the total mass of  $T\wedge \beta^{k-p}$ on $X.$
A $(p, p)$-current $T$ on $X$  is {\it strictly positive} if  we have locally $T \geq \epsilon \beta^p ,$ i.e., $T -\epsilon \beta^p$
is positive, for some constant $\epsilon > 0.$ The definition does not depend on the
choice of $\beta.$

 Let $R$ be  a current with measure coefficients (or equivalently,  of order $0$) on an open set $\Omega$  in  $X.$ 
 Let $W$ be a relatively compact open subset of $\Omega$ and $\Phi$ a smooth test form  on $\Omega,$ we will  write
 \begin{equation}\label{e;cut-off}
 \int_W  R\wedge \Phi:=  \langle  R,\ind_W\Phi\rangle,
 \end{equation}
 where $\ind_W$ is  the  characteristic  function of $W.$ Let $(R_n)_{n=1}^\infty$ be a  sequence of positive  currents on $\Omega$ such that
 $\lim_{n\to\infty}  R_n=R$  weakly on  $\Omega,$ then  we  see that
 \begin{equation}\label{e:continuity-cut-off}
 \lim_{n\to\infty}\int_W  R_n\wedge \Phi=\int_W  R\wedge \Phi
 \end{equation}
 for every smooth test form $\Phi$  on $\Omega$ and  every relatively compact open subset $W\subset \Omega$ with $\|R\|(\partial W)=0.$
 Here, $\partial W$ is  the topological boundary of $W$ and $\|R\|$ is  the mass-measure of $R.$
 Consequently, if $K$  is a compact subset of $\Omega$ and $(W_i)_{i\in I}$ is a family of open subsets of $\Omega$ such that
 $K\subset  W_i$ for all $i\in I$ and $\partial W_i\cap \partial W_{j}=\varnothing$ for $i\not=j,$ then  we have
 \begin{equation} \label{e:except-countable}
   \lim_{n\to\infty}\int_{W_i}  R_n\wedge \Phi=\int_{W_i}  R\wedge \Phi
 \end{equation}
for every smooth test form $\Phi$  on $\Omega$ and  every $i\in I$  except  for a countable subset of $I.$

 In this  article we are concerned with the following  notion of weak convergence of quasi-positive currents.
 \begin{definition}\label{D:quasi-positivity}\rm
  We say that   a    current $R$ defined on $\Omega$  is {\it  quasi-positive} if, for every $x\in \Omega,$ there are an open neighborhood $\Omega_x$ of $x$  in $\Omega$ and  a $\Cc^1$-diffeomorphism $\tau_x$
 of $\Omega_x$ such that $\tau_x^*R$ is a  positive current.
 
 We say that   a  sequence of  currents $(R_n)_{n=1}^\infty$ {\it converge  in the sense of quasi-positive currents  on $\Omega$ to a current $R$}
  if   for every $x\in \Omega,$ there are an open neighborhood $\Omega_x$ of $x$  in $\Omega$ and  a $\Cc^1$-diffeomorphism $\tau_x$
 of $\Omega_x$  and two sequences of  positive  currents $(T^\pm_n)_{n=1}^\infty$  on $\Omega_x$   such that all currents $\tau_x^*(R_n-R)= T^+_n-T^-_n$   and that both sequences $T_n^\pm$ converge  weakly to  a common positive current  $T$ on $\Omega_x.$  
 \end{definition}
 
 The relevance of this notion is justified by the following simple result.
 \begin{lemma}{\rm (see \cite[Lemma 2.4]{Nguyen21})}\label{L:quasi-positivity}
  If  a  sequence of  currents $(R_n)_{n=1}^\infty$  converge  in the sense of quasi-positive currents  on $\Omega$ to a current $R,$ then 
  both  \eqref{e:continuity-cut-off} and \eqref{e:except-countable}  hold.
 \end{lemma}
\subsection{Normal bundle and admissible maps}\label{SS:admissible-maps}

The  following notion,  introduced  by  Dinh-Sibony \cite{DinhSibony18}, plays a vital role in their tangent theory for positive closed currents.  
\begin{definition}\label{D:admissible-maps} {\rm  (See \cite[Definitions 2.15 and  2.18]{DinhSibony18}) } \rm   
 Let $B$ be a relatively compact nonempty open subset of $V.$
An admissible map   along  $B$  is   a $\Cc^1$-smooth diffeomorphism $\tau$  from an open neighborhood  $U$ of    $\overline{B}$ in $X$
onto an open neighborhood of $B\subset V$ in $\E$ (where $V$ is  identified with the zero section $0_\E$) which is identity on an open neighborhood of  $\overline B\subset V$ such that the  endomorphism   on $\E$ induced  by the restriction
of the differential $d\tau$ to ${\overline B}$ is identity.  More  specifically, a $\Cc^1$-smooth diffeomorphism $\tau: \ U\to\tau(U)$ 
is  admissible  if for every $x\in V\cap U,$ the following two conditions hold:

\begin{itemize}
 \item[(i)] $\tau(x)=x;$ 
\item[(ii)] Let $d\tau_x:\ \Tan_x(X)\to \Tan_x(\E)$ be the  differential $d\tau$ at $x.$ By writing  $\Tan_x(X)=\Tan_x(V)\oplus \E_x$ and 
$\Tan_x(\E)=\Tan_x(V)\oplus \E_x,$  $d\tau_x$ induces in a natural way a $\R$-endomorphism on $\E_x.$  This second  condition  requires that this endomorphism is identity $\id_{\E_x}.$
 \end{itemize}
  
In local coordinates, we can describe an admissible
map $\tau$ as follows: for every  point $x\in V\cap U,$  for every local
chart $y=(z,w)$  on a neighborhood $W$ of   $x$ in $U$  with  $V\cap W=\{z=0\}$, 
 we have
\begin{equation}\label{e:admissible-maps}
\tau (y) = \big( z + O(\|z\|^2), w+ O(\|z\|)\big) ,
\end{equation}
and
\begin{equation}
d\tau  (y) = \big( dz + \widetilde O(\|z\|^2), dw + \widetilde O(\|z\| ) \big),
\end{equation}
as $z \to  0$ where for every positive integer $m,$  $\widetilde O(\|z \|^m )$ denotes the sum of $1$-forms with
$O(\|z\|^m )$-coefficients and a linear combination of $dz,$ $d\bar z$  with $O(\|z\|^{m-1} )$-coefficients.
\end{definition}
It  is worth noting that in  \cite{DinhSibony18}   Dinh-Sibony  use the terminology {\it almost-admissible}  for those
maps satisfying  Definition \ref{D:admissible-maps}.
In general, $\tau$ is not holomorphic. When $U$ is a small enough local chart,
we can choose a holomorphic  admissible map by using suitable holomorphic coordinates
on $U .$   For the  global situation, the following  result gives a positive answer.
\begin{theorem}{\rm (\cite[Lemma 4.2]{DinhSibony18})}
For every  compact  subset $V_0\subset V,$ there always exists  an admissible map $\tau$ defined on a small enough tubular neighborhood  $U$ of $V_0$ in $X.$
\end{theorem}
  
In order  to  develop   a quantitative theory of  tangent  and density currents for  positive plurisubharmonic currents, the following notion, inspired
by Dinh-Sibony \cite[Proposition 3.8]{DinhSibony18b},  was  introduced  in \cite[Definition 2.7]{Nguyen21}.  
\begin{definition}
 \label{D:Strongly-admissible-maps}\rm  Let $B$ be a relatively compact nonempty open subset of $V.$
   A {\it strongly admissible} map  along  $B$ is   a $\Cc^2$-smooth diffeomorphism $\tau$  from an open neighborhood $U$ of $\overline B$  in $X$
onto an open neighborhood of $V \cap U$ in $\E$  such that for every  point $x\in V\cap U,$  for every local
chart $y=(z,w)$  on a neighborhood $W$ of   $x$ in $U$  with  $V\cap W=\{z=0\}$,
 we have   
 \begin{eqnarray*}
 \tau_{j}(z,w)&= &  z_j+\sum_{p,q=1}^{k-l}  a_{pq}(w) z_pz_{q}  +O(\|z\|^3)\quad\text{for}\quad 1\leq j\leq  k-l,\\
 \tau_{j}(z,w)&= &  w_{j-(k-l)} +\sum_{p=1}^{k-l} b_p(w)z_p 
 +O(\|z\|^2)\quad\text{for}\quad k-l< j\leq  k.
 \end{eqnarray*}
 Here,  we  write $\tau(y)=(\tau_{1}(y),\ldots,\tau_{k-l}(y),\tau_{k-l+1}(y),  \ldots \tau_{k}(y))\in\C^k,$  and $a_{pq},$  $b_p,$ 
 are $\Cc^2$-smooth functions depending only on $w.$
 In other words, if we write $\tau(z,w)=(z',w')\in\C^{k-l}\times \C^l,$ then
 \begin{eqnarray*}
 z'&= &  z+  z A z^T  +O(\|z\|^3),\\
 w'&= &  w + Bz 
 +O(\|z\|^2),
 \end{eqnarray*}
 where  $A$ is a $(k-l)\times(k-l)$-matrix and $B$ is  a $l\times(k-l)$-matrix whose entries are $\Cc^2$-smooth functions in $w,$ $z^T$ is the transpose of $z,$
\end{definition}
 
Observe  that  a  strongly  admissible map is necessarily  admissible in the sense of Definition \ref{D:admissible-maps}.
 On the  other hand,  holomorphic  admissible maps  are always strongly   admissible. Roughly speaking,  strongly  admissible maps
 are  those  admissible maps which are in some   sense {\it nearly}  holomorphic.

\begin{theorem}\label{T:strong-admissible-maps}{\rm (see \cite[Theorem 1.19 (i)]{Nguyen21}) }
  Let $X$ be a K\"ahler  manifold of dimension $k.$ Let  $V\subset X$ be  a submanifold
  of dimension $l$ and  $B\subset  V$  a relatively compact piecewisely $\Cc^2$-smooth open subset.
   Then there  exists a  strongly  admissible map  for $B.$
 \end{theorem}

 \subsection{Classes $\CL^p,$   $\PH^{p},$ $\SH^{p}(X)$ and   $\CL^{p;m,m'},$  $\PH^{p;m,m'}, $ $\SH^{p;m,m'}$ ($0\leq m'\leq m$)  of positive  currents}\label{SS:Classes}
 
 Let $X$ be a complex manifold of dimension $k$ and $p$ be an integer with $0\leq p\leq k.$
 Denote by $\CL^p(X)$  (resp.  $\PH^{p}(X)$)   (resp.  $\SH^{p}(X)$)  the class of positive
 closed  $(p,p)$-currents (resp. the class of positive
 pluriharmonic $(p,p)$-currents) (resp.  the class of positive plurisubharmonic $(p,p)$-currents  $(p,p)$-currents) on $X.$

Let $V\subset  X$ be   a  K\"ahler submanifold  of dimension $l,$ and  $\omega$ a K\"ahler form on $V,$ 
and $B\subset V$
a  relatively compact piecewise  $\Cc^2$-smooth open subset.

\begin{definition}\rm \label{D:Class}{\rm  (\cite[Subsection 1.5]{Nguyen21})}
Let $m,m'\in\N$ with $0\leq m'\leq m.$
Let $T$ be a positive   current of bidegree $(p,p)$ on $X.$
We say that {\it $T$ is  approximable along $B$ with $\Cc^m$-smooth positive  closed forms (resp. positive pluriharmonic forms, resp.  positive plurisubharmonic  forms)   with $\Cc^{m'}$-control on boundary,}
if, there are an open neighborhood $U$ of $\overline B$ in $X,$   and an open neighborhood $W$ of $\partial  B$ in $X,$  and  a sequence of positive closed  $\Cc^m$-smooth forms   (resp. positive pluriharmonic $\Cc^m$-smooth forms, resp.  positive plurisubharmonic  $\Cc^m$-smooth forms)           $T_n$ defined on $U,$    such that 
 \begin{enumerate}
\item the masses  $\|T_n\|$ on $U$ are uniformly bounded;

\item $T_n$ converge weakly  to $T$ on $U$ as $n$ tends to infinity;

\item     The $\Cc^{m'}$-norms of $T_n$'s on $W$ are uniformly bounded, that is,  
$$\sup_{n\geq 1} \|T_n\|_{\Cc^{m'}(W)} <\infty.$$
 \end{enumerate}
 Let  $ \CL^{p;m,m'}(B)$ denote the class of all positive  closed currents on $X$
 which are approximable  along $B$ by  $\Cc^m$-smooth positive closed forms.
  
  We denote  by  $\CL^{p;m,m'}( B)$  the class of all positive $(p,p)$-currents which are approximable  along $B$ by  $\Cc^m$-smooth positive closed forms with  $\Cc^{m'}$-control on boundary.

  Analogously, we denote  by  $\PH^{p;m,m'}( B)$  the class of all positive $(p,p)$-currents which are approximable  along $B$ by  $\Cc^m$-smooth positive pluriharmonic forms with  $\Cc^{m'}$-control on boundary.

  We denote  by  $\SH^{p;m,m'}( B)$  the class of all positive $(p,p)$-currents which are approximable  along $B$ by  $\Cc^m$-smooth positive plurisubharmonic forms with  $\Cc^{m'}$-control on boundary.
 \end{definition}

 The relevance of these classes of positive currents is illustrated by the  following result which is  inspired by the  work in \cite{DinhSibony04}.

\begin{theorem}\label{T:approximation}{\rm (see \cite[Theorem 1.19 (ii)]{Nguyen21}) }
  Let $X$ be a  K\"ahler manifold of dimension $k.$ Let  $V\subset X$ be  a submanifold
  of dimension $l$ and  $B\subset  V$  a relatively compact piecewisely $\Cc^2$-smooth open subset.
  Let $m,m'\in\N$ with $m\geq m'.$  
  Let $T$ be a   positive plurisubharmonic (resp.  positive  pluriharmonic, resp.  positive closed) $(p,p)$-current   on $X$ which satisfies the following conditions (i)--(ii):
   \begin{itemize}
   \item[(i)] $T$ is  of class $\Cc^{m'}$ near  $\partial B;$
   \item[(ii)] There is   a relatively compact open subset $\Omega$ of $X$ with  $B\Subset \Omega$
   and $dT$ is of class $\Cc^0$ near $\partial\Omega.$
   \end{itemize}
  Then  $T$ can be written  in an open neighborhood of $\overline B$ in $X$ as $T=T^+-T^-$  for some $T^\pm  \in  \SH^{p;m,m'}(B)$ (resp. $T^\pm  \in  \PH^{p;m,m'}(B),$ $T^\pm  \in  \CL^{p;m,m'}(B)$).
 \end{theorem}
  
 \begin{corollary}\label{C:approximation}
  Let $X$ be a compact  K\"ahler manifold of dimension $k.$   Let  $V\subset X$ be  a submanifold
  of dimension $l$ and  $B\subset  V$  a relatively compact piecewisely $\Cc^2$-smooth open subset.
  Let $m,m'\in\N$ with $m\geq m'.$  
  Let $T$ be a      positive  pluriharmonic (resp.  positive closed) $(p,p)$-current   on $X$  such that  $T$ is  of class $\Cc^{m'}$ near  $\partial B$ in  $X.$
  Then  $T$ can be written  in an open neighborhood of $\overline B$ in $X$ as $T=T^+-T^-$  for some  $T^\pm  \in  \PH^{p;m,m'}(B)$  (resp. $T^\pm  \in  \CL^{p;m,m'}(B)$).
 \end{corollary}

 \section{Lelong-Jensen formulas for holomorphic vector bundles}\label{S:Lelong-Jensen}
 
 We keep the notation introduced in Subsection \ref{SS:generalized-Lelong-numbers}.
 Let $B\subset  V$ be a domain and  $0<r<\infty.$
 Consider the tube $\Tube(B,r)$ defined  by \eqref{e:tubular-nbh-0}.
The boundary   $\partial\Tube(B,r)$ can be decomposed as the disjoint union  of the {\it  vertical boundary}   $\partial_\ver\Tube(B,r)$ and
the  {\it horizontal  boundary}    $\partial_\hor\Tube(B,r)$, where
\begin{eqnarray*}
 \partial_\ver\Tube(B,r)&:=&  \left\lbrace y\in \E:\   \pi(y)\in \partial B\quad\text{and}\quad  \varphi(y)\leq r^2  \right\rbrace ,\\
 \partial_\hor\Tube(B,r)&:=&  \left\lbrace y\in \E:\   \pi(y)\in B\quad\text{and}\quad  \varphi(y)=r^2  \right\rbrace .
\end{eqnarray*}
We see  easily that $\Tube(B,r)$ is a  manifold  with piecewise $\Cc^2$-smooth  boundary
for  every $r\in (0,\bfr]$ except a  countable set of values. Here the constant $\bfr>0$ was defined in Subsection \ref{SS:generalized-Lelong-numbers}. When  $\partial B=\varnothing,$  we have $\partial_\ver\Tube(B,r)=\varnothing.$

\begin{figure}[h]
 \begin{center}

\begin{tikzpicture}
\draw[->] (0,0) -- (10,0) node[right] {$\C^l$};
\node [below left=0.25cm] at (10,0) {base};
\draw[->] (0,0) -- (0,5) node[above] {$\C^{k-l}$};
\node [above right=0.25cm] at (0,4) {fibers};

\draw[->] (2,0) -- (2,4);

\draw[->] (5,0) -- (5,4);

\draw[->] (8,0) -- (8,4);
\draw[black, ultra thick, fill=gray!20] (2,0) rectangle (8,2.5);
\draw[color=red, decorate,decoration={brace,raise=0.1cm}]
(8,0) -- (2,0) node[below=0.25cm,pos=0.5] {$B$};
\draw[color=blue, decorate,decoration={brace,raise=0.7cm}]
(10,0) -- (0,0) node[below=1cm,pos=0.5] {$V$};

\draw[color=red, decorate,decoration={calligraphic brace,raise=0.1cm}]
(2,2.5) -- (8,2.5) node[above=0.25cm,pos=0.5] {$\partial_\hor \Tube(B,r)$};

\node [color=blue, above left] at (2,1.25) {$\partial_\ver \Tube(B,r)$};

\node [color=blue, above right] at (8,1.25) {$\partial_\ver \Tube(B,r)$};
\node  at (2,0) {$[$};
\node [above left] at (2,0) {$\partial B$};
\node  at (8,0) {$]$};
\node [above right] at (8,0) {$\partial B$};
 \draw [<->, > = stealth, red] (5, 0) -- (5, 2.5) node [midway, fill = white] {$r$};
 
 \node[right=-0.15cm] at (12.5,2) {\textbf{Key}};
\draw[black, thick,fill=gray!20] (12,0.5) rectangle (13,1)
node[pos=0.5,right=0.6cm] {$\Tube(B,r)$};
 \draw[color=red, decorate,decoration={calligraphic brace,raise=0.1cm}]
(13,0.5) -- (12,0.5) node[below=0.05cm,pos=0.5] {$B$};
 \draw[color=orange,decorate,decoration={calligraphic brace,raise=0.1cm}]
(12,0.5) -- (12,1) node[above=0.2cm,pos=0.5,sloped] {$r$};

\draw[->, blue, ultra thick] (1,1.3) to 
(2,1); 

\draw[->, blue, ultra thick] (9,1.3) to 
(8,1);
\end{tikzpicture}
\caption{Illustrations of a Tube  $\Tube(B,r)$  with base $B$ and radius $r,$ its horizontal boundary $\partial_\hor \Tube(B,r)$ and  its vertical boundary $\partial_\ver \Tube(B,r).$  The boundary $\partial B$ of $B$ is represented by $[$ and $]$ on the base $\overline B.$ }
\end{center}
\end{figure}
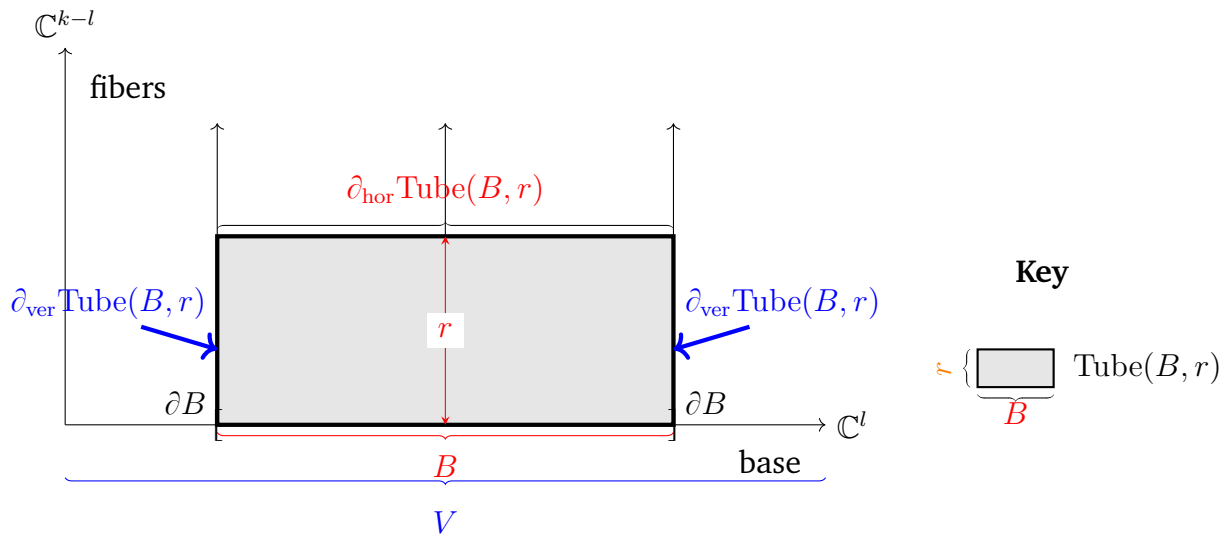

 
 
 
 
 
\begin{notation}
 \label{N:principal}\rm 
  Let $S$ be  a current of bidegree $2p$  defined on $\Tube(B,\bfr)\subset \E.$
  We  denote by $S^\sharp$  or equivalently $(S)^\sharp$ its component of bidegree $(p,p).$
\end{notation}

In this section we  recall  from \cite[Section 4]{Nguyen21}  some  Lelong-Jensen  formulas   for  vector bundles.   These formulas  play
a key role  throughout   this  work.

\begin{theorem}\label{T:Lelong-Jensen} {\rm  (see \cite[Theorem 4.2]{Nguyen21})}
Let $r\in(0,\bfr]$ and   $B\Subset V$  a relatively compact  open set with piecewice $\Cc^2$-smooth boundary.
Let $S$ be a real current of dimension $2q$ on  a neighborhood of $\overline\Tube(B,r)$ such that $S$ and $\ddc S$ are of order  $0$
and  that  $S$ is of class $\Cc^1$ near $\partial_\ver \Tube( B, r).$  
Suppose that there is a sequence of $\Cc^2$-smooth forms of dimension $2q$ $(S_n)_{n=1}^\infty$ defined on  a neighborhood of $\overline\Tube(B,r)$ 
such  that
\begin{itemize}
 \item[(i)] 
$S_n$ converge to $S$  in the sense of quasi-positive currents   on a neighborhood of  $\overline\Tube(B,r)$ as $n$ tends to infinity (see Definition \ref{D:quasi-positivity});
\item [(ii)]  $\ddc S_n$ converge to $\ddc S$  in the sense of quasi-positive currents   on a neighborhood of  $\overline\Tube(B,r)$ as $n$ tends to infinity;
\item[(iii)] there is an  open neighborhood of $\partial_\ver\Tube(B,r)$  on which  $S_n$  converge to $S$  in  $\Cc^1$-norm.
\end{itemize}
Then,   for  all $r_1,r_2\in  (0,r]$ with  $r_1<r_2$  except for  a  countable set of   values,  we have that  
\begin{equation}\label{e:Lelong-Jensen}\begin{split}
 {1\over  r_2^{2q}} \int_{\Tube(B,r_2)} S\wedge \beta^q- {1\over  r_1^{2q}} \int_{\Tube(B,r_1)} S\wedge \beta^q
 = \lim\limits_{n\to\infty}\Vc(S_n,r_1,r_2)+   \int_{\Tube(B,r_1,r_2)} S\wedge \alpha^q\\
 +  \int_{r_1}^{r_2} \big( {1\over t^{2q}}-{1\over r_2^{2q}}  \big)2tdt\int_{\Tube(B,t)} \ddc S\wedge \beta^{q-1} 
 +  \big( {1\over r_1^{2q}}-{1\over r_2^{2q}}  \big) \int_{0}^{r_1}2tdt\int_{\Tube(B,t)} \ddc S\wedge \beta^{q-1}.
 \end{split}
\end{equation} 
Here the vertical boundary term  $\Vc(S,r_1,r_2)$  for a $\Cc^1$-smooth form $S$  is  given by the following formula,  where $S^\sharp$ denotes, 
according to Notation \ref{N:principal}, the component of bidimension $(q,q)$ of the current $S:$
\begin{equation}\label{e:vertical-boundary-term}
\begin{split}
 \Vc(S,r_1,r_2)&:=-\int_{r_1}^{r_2} \big( {1\over t^{2q}}-{1\over r_2^{2q}}  \big)2tdt\int_{\partial_\ver\Tube(B,t)} \dc S^\sharp\wedge \beta^{q-1} 
 \\ &-  \big( {1\over r_1^{2q}}-{1\over r_2^{2q}}  \big) \int_{0}^{r_1}2tdt\int_{\partial_\ver\Tube(B,t)} \dc S^\sharp\wedge \beta^{q-1}
 +{1\over r_2^{2q}} \int_{\partial_\ver \Tube(B,r_2)}\dc\varphi\wedge S^\sharp\wedge \beta^{q-1} \\&-{1\over r_1^{2q}} \int_{\partial_\ver \Tube(B,r_1)}\dc\varphi\wedge S^\sharp\wedge \beta^{q-1}
  -\int_{\partial_\ver \Tube(B,r_1,r_2)}\dc\log\varphi\wedge S^\sharp\wedge \alpha^{q-1}.
  \end{split}
\end{equation}
  \end{theorem}

The next theorem deals with the special case where  the current  is  approximable by smooth {\bf closed}  forms with  control on the boundary.
Here, we gain the  smoothness.

\begin{theorem}\label{T:Lelong-Jensen-closed}{\rm  (see \cite[Theorem 4.5]{Nguyen21})}
Let $r\in (0,\bfr]$ and  let $S$ be a real closed current of dimension $2q$ on  a neighborhood of $\overline\Tube(B,r).$ 
%
%
Suppose that there is a sequence of $\Cc^1$-smooth closed forms of dimension $2q:$  $(S_n)_{n=1}^\infty$ defined on  a neighborhood of $\overline\Tube(B,r)$ 
such  that
$S_n$ converge to $S$  in the sense of quasi-positive currents   on a neighborhood of  $\overline\Tube(B,r)$ as $n$ tends to infinity (see Definition \ref{D:quasi-positivity}).

Then the following  two assertions  hold:
\begin{enumerate} \item   
For  all $r_1,r_2\in  (0,r]$ with  $r_1<r_2$  except for  a  countable set of   values, we have that  
\begin{equation}\label{e:Lelong-Jensen-closed}
 {1\over  r_2^{2q}} \int_{\Tube(B,r_2)} S\wedge \beta^q- {1\over  r_1^{2q}} \int_{\Tube(B,r_1)} S\wedge \beta^q
 =\lim\limits_{n\to\infty} \Vc(S_n,r_1,r_2)+   \int_{\Tube(B,r_1,r_2)} S\wedge \alpha^q.
\end{equation} 
Here the vertical boundary term  $\Vc(S,r_1,r_2)$ for  a  continuous form $S$  is  given by
\begin{equation}\label{e:vertical-boundary-term-closed}
\begin{split}
 \Vc(S,r_1,r_2)&:=
 {1\over r_2^{2q}} \int_{\partial_\ver \Tube(B,r_2)}\dc\varphi\wedge S\wedge \beta^{q-1}-{1\over r_1^{2q}} \int_{\partial_\ver \Tube(B,r_1)}\dc\varphi\wedge S\wedge \beta^{q-1}\\
  &-\int_{\partial_\ver \Tube(B,r_1,r_2)}\dc\log\varphi\wedge S\wedge \alpha^{q-1}.
  \end{split}
\end{equation}
\item  If $S$ is a closed $\Cc^1$-smooth form, then identity \eqref{e:Lelong-Jensen-closed} (with $S_n:= S$ for $n \geq 1$) holds
for  all $r_1,r_2\in  [0,r]$ with  $r_1<r_2.$
\end{enumerate}
\end{theorem}

Next, we  consider the  following variant  of $\varphi$  in the spirit of 
\eqref{e:varphi-spec}:
for  every $\epsilon>0,$  set
\begin{equation}\label{e:varphi_eps}
 \varphi_\epsilon(y):=\|y\|^2+\epsilon^2\qquad \text{for}\qquad  y\in \E.
\end{equation}
In this case  where we have $r_0=\epsilon$ and $\bfr=\infty.$
Following the model \eqref{e:alpha-beta-spec}, consider also the following  closed  $(1,1)$-form  for each $\epsilon >0$ on  $\U:$
\begin{equation}\label{e:alpha-beta-eps}
 \alpha_\epsilon:=\ddc\log\varphi_\epsilon\quad\text{and note that }\quad \beta= \ddc\varphi_\epsilon.
\end{equation}  
The  following   result which  will play a key role  for proving logarithmic   interpretation  version in the spirit of  \eqref{e:Lelong-number-point-bisbis(2)}.

\begin{theorem}\label{T:Lelong-Jensen-eps}{\rm  (see \cite[Theorem 4.10]{Nguyen21})}
Let $\bfr\in\R^+_*$  and   $B\Subset V$  a relatively compact  open set with piecewice $\Cc^2$-smooth boundary.
Let $S$ be a real current of dimension $2q$ on  a neighborhood of $\overline\Tube(B,\bfr).$ 
%
%
Suppose that there is a sequence of $\Cc^2$-smooth forms of dimension $2q$ $(S_n)_{n=1}^\infty$ defined on  a neighborhood of $\overline\Tube(B,\bfr)$ 
such  that
\begin{itemize}
 \item[(i)] 
$S_n$ converge to $S$  in the sense of quasi-positive currents   on a neighborhood of  $\overline\Tube(B,\bfr)$ as $n$ tends to infinity (see Definition \ref{D:quasi-positivity});
\item [(ii)]  $\ddc S_n$ converge to $\ddc S$  in the sense of quasi-positive currents   on a neighborhood of  $\overline\Tube(B,\bfr)$ as $n$ tends to infinity.
\end{itemize}
Then  the following  two assertions hold:
\begin{enumerate}
 \item 
For  all  $r\in (0,\bfr)$ and   $\epsilon\in  (0,r)$      except for a  countable  set of  values, we have that
\begin{multline*}
 {1\over  (r^2+\epsilon^2)^q} \int_{\Tube(B,r)} S\wedge \beta^q
 =\lim\limits_{n\to\infty} \Vc_\epsilon(S_n,r)+   \int_{\Tube(B,r)} S\wedge \alpha^q_\epsilon\\
 +  \int_{0}^{r} \big( {1\over (t^2+\epsilon^2)^q}-{1\over (r^2+\epsilon^2)^q}  \big)2tdt\int_{\Tube(B,t)} \ddc S\wedge \beta^{q-1} .
\end{multline*} 
Here the vertical boundary term  $\Vc_\epsilon(S,r)$ for  a  $\Cc^1$-smooth form $S$ is  given by
\begin{equation}\label{e:vertical-boundary-term-eps}
\begin{split}
 \Vc_\epsilon(S,r)&:=-\int_{0}^{r} \big( {1\over (t^2+\epsilon^2)^q}-{1\over (r^2+\epsilon^2)^q}  \big)2tdt\int_{\partial_\ver\Tube(B,t)} \dc S^\sharp\wedge \beta^{q-1} 
  \\
 &+{1\over  (r^2+\epsilon^2)^q} \int_{\partial_\ver \Tube(B,r)}\dc\varphi\wedge S^\sharp\wedge \beta^{q-1}
  -\int_{\partial_\ver \Tube(B,r)}\dc\log\varphi_\epsilon\wedge S^\sharp\wedge \alpha_\epsilon^{q-1}.
  \end{split}
\end{equation}

\item  If $S$  is a $\Cc^2$-smooth form, then the above identity (with $S_n:=S$ for $n\geq 1$) holds  for  all  $r\in (0,\bfr)$ and 
$\epsilon\in  (0,r).$ 

\end{enumerate}
\end{theorem}

We record a version of Theorem \ref{T:Lelong-Jensen-eps} for {\bf closed} currents.

\begin{theorem}\label{T:Lelong-Jensen-closed-eps}{\rm  (see \cite[Theorem 4.11]{Nguyen21})}
Let $\bfr\in\R^+_*$ and  let $S$ be a real closed current of dimension $2q$ on  a neighborhood of $\overline\Tube(B,\bfr).$ 
%
%
Suppose that there is a sequence of $\Cc^1$-smooth closed forms of dimension $2q:$  $(S_n)_{n=1}^\infty$ defined on  a neighborhood of $\overline\Tube(B,\bfr)$ 
such  that
%
$S_n$ converge to $S$  in the sense of quasi-positive currents   on a neighborhood of  $\overline\Tube(B,\bfr)$ as $n$ tends to infinity (see Definition \ref{D:quasi-positivity}).
%
%
Then the  following two  assertions hold:
\begin{enumerate}
 \item 
For  all $r\in (0,\bfr]$ and   $\epsilon\in  (0,r)$  except for a countable set of values,   we have that 
\begin{equation}\label{e:Lelong-Jensen-closed-eps}
 {1\over  (r^2+\epsilon^2)^{q}} \int_{\Tube(B,r)} S\wedge \beta^q
 = \lim\limits_{n\to\infty} \Vc_\epsilon(S_n,r)+   \int_{\Tube(B,r)} S\wedge \alpha^q_\epsilon.
\end{equation} 
Here the vertical boundary term  $\Vc_\epsilon(S,r)$ for a  continuous   form $S$  is  given by
\begin{equation}\label{e:vertical-boundary-term-closed-eps}
 \Vc_\epsilon(S,r):=
 {1\over (r^2+\epsilon^2)^{q}} \int_{\partial_\ver \Tube(B,r)}\dc\varphi\wedge S\wedge \beta^{q-1}-
  \int_{\partial_\ver \Tube(B,r)}\dc\log\varphi_\epsilon\wedge S\wedge \alpha^{q-1}_\epsilon.
\end{equation}

\item  If $S$  is a closed  $\Cc^1$-smooth form, then the above identity (with $S_n:=S$ for $n\geq 1$) holds  for  all  $r\in (0,\bfr)$ and 
$\epsilon\in  (0,r).$
\end{enumerate}
\end{theorem}

Here is a version of   Theorem  \ref{T:Lelong-Jensen-closed} for smooth closed forms when the  minor radius $r_1$  becomes  infinitesimally small.
\begin{theorem}\label{T:Lelong-Jensen-smooth-closed}{\rm  (see \cite[Theorem 4.16]{Nguyen21})}
Let $\bfr\in\R^+_*$  and let $S$ be a $\Cc^1$-smooth closed form of dimension $2q$ on  a neighborhood of $\overline\Tube(B,\bfr).$
Suppose that $q\leq k-l.$
\begin{enumerate}
\item
Then,  for  all $0<r\leq\bfr,$
\begin{equation}\label{e:Lelong-Jensen-smooth-closed}
 {1\over  r^{2q}} \int_{\Tube(B,r)} S\wedge \beta^q- \lim_{s\to 0+}{1\over  s^{2q}} \int_{\Tube(B,s)} S\wedge \beta^q
 = \Vc(S,r)+   \int_{\Tube(B,r)} S\wedge \alpha^q.
\end{equation} 
Here the vertical boundary term  $\Vc(S,r)$ is  given by
\begin{equation}\label{e:vertical-boundary-term-closed-bis}
\begin{split}
 \Vc(S,r)&:=\big(
 {1\over r^{2q}} \int_{\partial_\ver \Tube(B,r)}\dc\varphi\wedge S\wedge \beta^{q-1}-\lim_{s\to 0+}{1\over s^{2q}} \int_{\partial_\ver \Tube(B,s)}\dc\varphi\wedge S\wedge \beta^{q-1}\big)\\
  &-\int_{\partial_\ver \Tube(B,r)}\dc\log\varphi\wedge S\wedge \alpha^{q-1}.
  \end{split}
\end{equation}
\begin{itemize}
\item If  $q<k-l$ then $\lim_{s\to 0+} {1\over  s^{2q}} \int_{\Tube(B,s)} S\wedge \beta^q=0.$
\item If  $q=k-l$ and  $S(y)$  is a  positive  form for all $y\in \overline B,$   then  the last  limit is nonnegative.
\end{itemize}
\item  Suppose  in addition that  $\supp(S) \cap \partial_\ver \Tube( B, \bfr)=\varnothing.$   
Then,  for  all $0<r<\bfr,$
\begin{equation*}
 {1\over  r^{2q}} \int_{\Tube(B,r)} S\wedge \beta^q- \lim_{s\to 0+}{1\over  s^{2q}} \int_{\Tube(B,s)} S\wedge \beta^q
 =   \int_{\Tube(B,r)} S\wedge \alpha^q.
\end{equation*} 
\end{enumerate}
\end{theorem}

Here is a version of   Theorem  \ref{T:Lelong-Jensen} for smooth  forms when the  minor radius $r_1$  becomes  infinitesimally small.
\begin{theorem}\label{T:Lelong-Jensen-smooth}{\rm  (see \cite[Theorem 4.15]{Nguyen21})}
Let $\bfr\in\R^+_*$  and let $S$ be a $\Cc^2$-smooth  form of dimension $2q$ on  a neighborhood of $\overline\Tube(B,\bfr).$
Suppose that $q\leq k-l.$
\begin{enumerate}
\item
Then,  for  all $0<r\leq\bfr,$
\begin{equation}\label{e:Lelong-Jensen-smooth}\begin{split}
& {1\over  r^{2q}} \int_{\Tube(B,r)} S\wedge \beta^q -\lim_{s\to 0+} {1\over  s^{2q}} \int_{\Tube(B,s)} S\wedge \beta^q
 =   \int_{\Tube(B,r)} S\wedge \alpha^q\\
 &+  \int_{0}^{r} \big( {1\over t^{2q}}-{1\over r^{2q}}  \big)2tdt\int_{\Tube(B,t)} \ddc S\wedge \beta^{q-1} 
 + \Vc(S,r).
 \end{split}
\end{equation} 
Here, the vertical boundary term  $\Vc(S,r)$ is  given by 
\begin{equation}\label{e:vertical-boundary-term-bis}
\begin{split}
 \Vc(S,r)&:=-\int_{0}^{r} \big( {1\over t^{2q}}-{1\over r^{2q}}  \big)2tdt\int_{\partial_\ver\Tube(B,t)} \dc S^\sharp\wedge \beta^{q-1} 
 -\int_{\partial_\ver \Tube(B,r)}\dc\log\varphi\wedge S^\sharp\wedge \alpha^{q-1}\\
 &+\big( {1\over r^{2q}} \int_{\partial_\ver \Tube(B,r)}\dc\varphi\wedge S^\sharp\wedge \beta^{q-1}-\lim_{s\to 0+} {1\over s^{2q}} \int_{\partial_\ver \Tube(B,s)}\dc\varphi\wedge S^\sharp\wedge \beta^{q-1}\big)
  .
  \end{split}
\end{equation}
\begin{itemize}
\item If   $q<k-l,$ then $\lim_{s\to 0+} {1\over  s^{2q}} \int_{\Tube(B,s)} S\wedge \beta^q=0.$
 \item If   $q=k-l$ and  $S(y)$ is a positive form  for all $y\in \overline B,$   then  the last  limit is nonnegative.
\end{itemize}

\item   Suppose  in addition that  $\supp(S) \cap \partial_\ver \Tube( B, \bfr)=\varnothing.$   
Then,  for  all $0<r<\bfr,$
\begin{eqnarray*}
 &&{1\over  r^{2q}} \int_{\Tube(B,r)} S\wedge \beta^q-\lim_{s\to 0+} {1\over  s^{2q}} \int_{\Tube(B,s)} S\wedge \beta^q\\
& =&\int_{\Tube(B,r)} S\wedge \alpha^q
 + \int_{0}^{r} \big( {1\over t^{2q}}-{1\over r^{2q}}  \big)2tdt\int_{\Tube(B,t)} \ddc S\wedge \beta^{q-1} 
 .
\end{eqnarray*}
\end{enumerate}
\end{theorem}

Finally,  we conclude the section with two  asymptotic   Lelong-Jensen formulas.

\begin{theorem}\label{T:vertical-boundary-terms}{\rm  (see \cite[Theorem 4.17]{Nguyen21})}
Let $\bfr\in\R^+_*$ and   $0\leq q\leq k-l.$ Let $S$ be a real current of dimension $2q$ on  a neighborhood of $\overline\Tube(B,\bfr)$ such that $S$ and $\ddc S$ are of order  $0.$
Suppose that there is a sequence of $\Cc^2$-smooth forms of dimension $2q$ $(S_n)_{n=1}^\infty$ defined on  a neighborhood of $\overline\Tube(B,\bfr)$ 
such  that
\begin{itemize}
 \item[(i)] 
$S_n$ converge to $S$  in the sense of quasi-positive currents   on a neighborhood of  $\overline\Tube(B,\bfr)$ as $n$ tends to infinity (see Definition \ref{D:quasi-positivity});
\item [(ii)]  $\ddc S_n$ converge to $\ddc S$  in the sense of quasi-positive currents   on a neighborhood of  $\overline\Tube(B,\bfr)$ as $n$ tends to infinity;
\item[(iii)] there is an  open neighborhood of $\partial_\ver\Tube(B,\bfr)$  on which the $\Cc^1$-norms of  $S_n$  are uniformly bounded.
\end{itemize}
Then,  for  all $s,r\in  (0,\bfr]$ with  $s<r$  except for  a  countable set of   values,   formula \eqref{e:Lelong-Jensen} for $r_1:=s,$ $r_2:=r$
(resp.  formula  \eqref{e:Lelong-Jensen-smooth})
holds with 
\begin{equation*}
 |\Vc(S,s,r)|\leq cr\qquad\big(\text{resp.}\qquad  |\Vc(S,r)|\leq  cr\quad\big),
\end{equation*}
where $c$ is a constant  independent of $s,r.$
\end{theorem}

\begin{theorem}\label{T:vertical-boundary-closed}{\rm  (see \cite[Theorem 4.18]{Nguyen21})}
Let $\bfr\in\R^+_*$ and   $0\leq q\leq k-l.$  Let $S$ be a real closed current of dimension $2q$ on  a neighborhood of $\overline\Tube(B,\bfr).$ 
Suppose that there is a sequence of $\Cc^1$-smooth closed forms of dimension $2q:$  $(S_n)_{n=1}^\infty$ defined on  a neighborhood of $\overline\Tube(B,\bfr)$ 
such  that
\begin{itemize}
 \item[(i)] 
$S_n$ converge to $S$  in the sense of quasi-positive currents   on a neighborhood of  $\overline\Tube(B,\bfr)$ as $n$ tends to infinity (see Definition \ref{D:quasi-positivity});
\item[(ii)] there is an  open neighborhood of $\partial_\ver\Tube(B,\bfr)$  on which the $\Cc^m$-norms of  $S_n$  are uniformly bounded, where
$m=0$ if  $q<k-l$ and $m=1$ if $q=k-l.$ 
\end{itemize}
Then,  for  all $s,r\in  (0,\bfr]$ with  $s<r$  except for  a  countable set of   values,  
formula \eqref{e:Lelong-Jensen-closed}  for $r_1:=s,$ $r_2:=r$  (resp. formula \eqref{e:Lelong-Jensen-smooth-closed}) holds
  with 
\begin{equation*}
 |\Vc(S,s,r)|\leq cr\qquad\big(\text{resp.}\qquad  |\Vc(S,r)|\leq  cr\quad\big),
\end{equation*}
where $c$ is a constant independent of $s,r.$
\end{theorem}

\begin{theorem}\label{T:vertical-boundary-eps}{\rm  (see \cite[Theorem 4.19]{Nguyen21})}
 We keep the hypothesis and the notation of Theorem \ref{T:Lelong-Jensen-eps} (resp.  Theorem \ref{T:Lelong-Jensen-closed-eps}).
 Then there  is  a  constant $c$ depending only on $S$ such that 
 for all  $r\in (0,\bfr]$  and $\epsilon\in(0,r),$ the  following assertions  hold:
 \begin{enumerate}
  \item If $q<k-l,$ then 
 $
 |\Vc_\epsilon(S,r)|\leq  cr.$
 
 \item If $  q=k-l$ and  we  are in the assumption of Theorem \ref{T:Lelong-Jensen-eps}, then
 
 $$ \big|\Vc_\epsilon(S,r)-{
1\over r^{2q}} \int_{\partial_\ver \Tube(B,r)}\dc\varphi\wedge S^\sharp\wedge \beta^{q-1}\big|
\leq cr 
.$$
\item If $  q=k-l$ and  we  are in the assumption of Theorem \ref{T:Lelong-Jensen-closed-eps}, then
 
 $$ \big|\Vc_\epsilon(S,r)-{
1\over r^{2q}} \int_{\partial_\ver \Tube(B,r)}\dc\varphi\wedge S\wedge \beta^{q-1}\big|
\leq cr 
.$$
\end{enumerate}

\end{theorem}

\section{Admissible estimates, mass indicators and their finiteness}\label{S:Finiteness}

In the first six  subsections we  recall  standard   settings for further technical  developments, and preliminary results established in  \cite[Sections 5 and 7]{Nguyen21}. The last subsection is devoted to new results  on  mass indicators and their finiteness.

\subsection{Forms $\alpha_\ver$ and $\beta_\ver$}
 
 Since  the  transition  functions  of the holomorphic vector bundle $\E$ are holomorphic, 
 the vertical operators $\partial_\ver,$ $\dbar_v$ 
which are the  restrictions of the usual operators $\partial$ and $\dbar$ on fibers of $\E$  are well-defined. More precisely, for a smooth
form $\Phi$ on an open set $\Omega$ in $\E,$ we can define
\begin{equation}\label{e:d-and-dbar-ver}
 \partial_\ver\Phi(y):=\partial |_{\E_{\pi(y)}}\Phi (y) \qquad\text{and}\qquad \dbar_\ver\Phi(y):=\dbar|_{\E_{\pi(y)}}\Phi(y)\qquad\text{for}\qquad y\in \Omega.
\end{equation}
So  the vertical operators  $d_\ver$ and  $\ddcv$ are  also well-defined  by the formulas \begin{equation}
\label{e:ddc-ver}
d_\ver \Phi:=\partial_\ver \Phi+\dbar_\ver \Phi\quad\text{and}\quad \ddcv\Phi:={i\over \pi}  \partial_\ver\dbar_\ver\Phi.                           
                          \end{equation}  
Consider for $ y\in \E,$
\begin{equation}\label{e:alpha-beta-ver}
 \alpha_\ver(y):=\ddcv \log\varphi (y)=\ddc|_{\E_{\pi(y)}}\log\varphi (y) \qquad\text{and}\qquad \beta_\ver(y):=\ddcv \varphi(y)=\ddc|_{\E_{\pi(y)}}\varphi(y),
\end{equation}
where  $\ddc|_{\E_{\pi(y)}}$  is restriction  of  the operator $\ddc$ on the fiber $\E_{\pi(y)}.$  
Observe  that  both
 $\alpha_\ver$  and $\beta_\ver$  are positive  $(1,1)$-forms on $\E.$
However, they are not necessarily closed.

\subsection{Positive forms $\hat\alpha$ and $\hat\beta$}
Recall  from  Subsection  \ref{SS:generalized-Lelong-numbers} the constant $\bfr>0.$ 
Recall from \cite[Lemma 5.1]{Nguyen21}  the construction of   positive  currents/forms   $\hat\alpha,$ $\hat\alpha'$  and $\hat\beta.$  
This positivity plays  a  crucial  role    in the sequel.

\begin{lemma}\label{L:hat-alpha-beta} 
\begin{enumerate}
\item  There is  a  constant $c_1>0$ large enough such that 
\begin{equation}\label{e:hat-beta}
  \hat\beta:=c_1\varphi\cdot  \pi^*\omega+\beta
 \end{equation}
is  positive  on $\pi^{-1}(V_0)$ and is strictly positive on $\pi^{-1}(V_0)\setminus V_0,$
and 
\begin{equation}\label{e:hat-alpha'}
\hat\alpha':= c_1 \pi^*\omega+ \alpha
\end{equation}
satisfies 
\begin{equation}\label{e:hat-alpha'-vs-alpha_ver}
c_1\alpha_\ver+c^2_1\pi^*\omega\geq\hat\alpha'\geq  c^{-1}_1\alpha_\ver.
\end{equation}
In particular,   $\hat\alpha'$ is positive  on $\pi^{-1}(V_0).$
\item 
  There are constants $c_2,c_3>0$ such that on $\Tube(V_0,\bfr)\setminus V_0,$  
\begin{equation}\label{e:hat-alpha}
\hat\alpha:=\hat\alpha'+c_2\beta= c_1 \pi^*\omega+ \alpha+c_2\beta
\end{equation}
is   strictly positive, and 
\begin{equation}\label{e:hat-alpha-vs-alpha_ver}
 \hat\alpha\geq  c^{-1}_1\alpha_\ver,
\end{equation}
and
\begin{equation}\label{e:hat-alpha-vs-hat-beta}
\varphi\hat\alpha\leq  c_3\hat\beta.
\end{equation}
\item 
  There are constants $c_3>0$ such that  on $\Tube(V_0,\bfr),$
  \begin{equation}\label{e:hat-beta-vs-beta_ver}
\hat\beta\geq  c^{-1}_1\beta_\ver,
\end{equation}
  and on $\Tube(V_0,\bfr)\setminus V_0,$  
\begin{equation}\label{e:tilde-alpha-vs-hat-beta}
\varphi\alpha_\ver\leq c_3\hat\beta.
\end{equation}
\end{enumerate}

\end{lemma}

\subsection{Analysis  in local coordinates}\label{SS:Anal-local-coord}
  The following local model studied in   \cite{Nguyen21} is   useful. Consider  an open set $V_0\Subset  V$ and  let us  study $\E$ near a given point $y_0\in V_0.$
We use the coordinates $(z,w)\in\C^{k-l}\times \C^l$  around a neighborhood $U$ of $y_0$  such that $y_0=0$  in  these coordinates. 
We may assume that  $U$ has the form $U=U'\times U'',$ where $U'$ (resp. $U'')$ is an  open neighborhood of $0'$ in $\C^{k-l}$ (resp. of  $0''$ in $\C^l$),
and  $V=\{z=0\}\simeq U''.$ Moreover,  we  may assume  that $U''=(2\D)^l.$  
Consider  the trivial  vector bundle $\pi:\ \E \to  U''$ with  $\E\simeq  \C^{k-l}\times U''.$ 
There is a smooth function  $A:\  \D^l\to \GL(\C,k-l)$ such that
\begin{equation}\label{e:varphi-new-exp}
\varphi(z,w)= \| A(w)z\|^2\qquad\text{for}\qquad  z\in\C^{k-l},\ w\in \D^l,
\end{equation}
where $\varphi$ is   defined in \eqref{e:varphi-spec}.
It follows  from \eqref{e:alpha-beta-ver} and  \eqref{e:varphi-new-exp} that  
\begin{equation}\label{e:tilde-alpha-beta-local-exp}
\alpha_\ver(z,w)=  A(w)^* [\ddc \log{\|z\|^2}]\quad\text{and}\quad    \beta_\ver(z,w)=  A(w)^* [\ddc \|z\|^2]\quad\text{for}\quad  z\in\C^{k-l},\ w\in \D^l.
\end{equation}
 
\subsection{Extended Standing Hypothesis}\label{SS:Ex-Stand-Hyp}

Let $B$ be a  relatively compact piecewise $\Cc^2$-smooth open
subset. Let $V_0$ be a relatively compact open subset of $V$ such that $B \Subset V_0 .$
Consider  a  strongly admissible map $\tau:\ \bfU\to\tau(\bfU)$  along $B,$  with $\bfU$ a neighborhood of $\overline B$ in $X.$
By shrinking $\bfU$ if necessary, we may   fix a finite collection $\Uc=(\bfU_\ell,\tau_\ell)_{1\leq \ell\leq \ell_0} ,$
 of holomorphic admissible maps for $\bfU.$  More precisely,  
there is a  finite cover of $\overline \bfU$ by open subsets $\bfU_\ell,$ $1\leq \ell\leq  \ell_0,$ of $X$
such that   there is   a holomorphic coordinate system on $\overline \bfU_\ell$ in $X$ and $\bfU_\ell$  is  biholomorphic  to  $\U_\ell:=\tau_\ell(\bfU_\ell)\subset \E$
by a  holomorphic admissible map $\tau_\ell.$ By  choosing $\bfr>0$ small enough, we may assume  without loss of generality that $\overline\Tube(B,\bfr)\Subset \U:=\bigcup_{\ell=1}^{\ell_0} \U_\ell.$ 
Fix  a partition of unity  $(\theta_\ell)_{1\leq \ell\leq \ell_0}$ subordinate   to the open cover  $(\bfU_\ell\cap V)_{1\leq \ell\leq \ell_0}$   of $\overline{\bfU\cap V}$  in $V$  such that $\sum_{1\leq \ell\leq \ell_0}  \theta_\ell=1$ on an open neighborhood of $\overline {\bfU\cap V} \subset V.$ We  may  assume  without loss of generality  that there are open  subsets
$\widetilde V_\ell\subset V$  for $1\leq \ell\leq \ell_0$ such that
\begin{equation}\label{e:supp} \supp(\theta_l)\subset  \widetilde V_\ell\Subset \bfU_\ell\cap V \quad\text{and}\quad 
\tau( \widetilde V_\ell)\Subset  \U_\ell\quad\text{and}\quad \pi^{-1}(\supp(\theta_\ell))\cap \U\subset \U_\ell.
\end{equation}
 For $1\leq \ell\leq \ell_0$ set 
 \begin{equation}\label{e:tilde-tau_ell}\tilde \tau_\ell:=\tau\circ\tau_\ell^{-1}.
 \end{equation}
  So  $\tilde \tau_\ell$ defines a map  from  $\U_\ell\subset \E$  onto  $\tau(\bfU_\ell)\subset \E.$

  We also assume  that  for every $1\leq\ell\leq\ell_0,$  there is a local  coordinate system $y=(z,w)$  on $\U_\ell$  with  $V\cap \U_\ell=\{z=0\}.$
  
  $\Uc=(\bfU_\ell,\tau_\ell)_{1\leq \ell\leq \ell_0}$ is  said to be  a  {\it covering  family of  holomorphic admissible maps for $B.$}

Recall that $\pi:\  \E\to V$ is  the canonical projection.
\begin{definition}\label{D:T-hash}\rm
Let  $ T$ be a current  defined on $\bfU.$ Consider the current  $T^{\#}$ defined  on $\U$ by the   following formula:
\begin{equation}\label{e:T-hash}
T^{\#}:=\sum_{\ell=1}^{\ell_0}  (\pi^*\theta_\ell)\cdot (\tau_\ell)_*( T|_{\bfU_\ell}).
\end{equation}
By  \eqref{e:supp},  $T^{\#}$ is  well-defined.

Let  $ T$ be a current  defined on $\bfU$  and $0\leq s<r\leq \bfr.$  Consider the currents  $T^{\#}_r$  and $T^{\#}_{s,r}$  defined  on $\U$ as  follows:
\begin{equation}\label{e:T-hash_r}
T^{\#}_r:=\sum_{\ell=1}^{\ell_0}  (\pi^*\theta_\ell)\cdot (\ind_{\Tube(B,r)}\circ \tilde \tau_\ell )\cdot(\tau_\ell)_*( T|_{\bfU_\ell})\quad\text{and}\quad
T^{\#}_{s,r}:=\sum_{\ell=1}^{\ell_0}  (\pi^*\theta_\ell)\cdot (\ind_{\Tube(B,s,r)}\circ \tilde \tau_\ell )\cdot(\tau_\ell)_*( T|_{\bfU_\ell}).
\end{equation}
 \end{definition}
 
\subsection{Admissible estimates}\label{SS:admissible-estimates}

Admissible estimates   are those  estimates which are related to   admissible maps.
This  subsection provides  necessary admissible estimates.  

Let $\U$ be an open neighborhood of $0$ in $\C^k.$
We use the local coordinates $y=(z,w)\in\C^{k-l}\times \C^l$ on $\U.$

The following  notion  will be  needed  in order  to obtain   admissible estimates.
\begin{definition}
 \label{D:precsim}
 \rm  Let $\Gamma$ be  a form of degree $2$ and $S$ a positive $(1,1)$-form  defined on $\U.$ 
 For  $(p,q)\in \{(0,2),(1,1),(2,0)\},$ $\Gamma^{p,q}$ denotes the component of bidegree $(p,q)$ of $\Gamma.$
 So $\Gamma^{1,1}=\Gamma^\sharp$  according to  Notation \ref{N:principal}.
 
 We  write $\Gamma\trianglelefteq  S$  
 if there is a constant $c>0$ such that   the  following   two inequalities  hold for $y\in\U:$
 $$
  \Gamma^{0,2}(y)\wedge \overline{\Gamma^{0,2}}(y)\leq c S^2(y)\quad\text{and}\quad  \Gamma^{2,0}(y)\wedge \overline{\Gamma^{2,0}}(y)\leq cS^2(y).
 $$
 \end{definition}
\begin{notation}\label{N:pm}
 \rm  Let $\Gamma$ and  $S$ be  two  real $(1,1)$-forms   defined on $\U.$
 
 We  write  $\Gamma\lesssim S$ if there is a constant $c>0$ such that  $\Gamma\leq c S.$
 We write  $\pm \Gamma\lesssim   S$  if  we  have both $ \Gamma\lesssim  S$ and  $- \Gamma\lesssim  S.$
 
 We  write $\Gamma\approx  S$ if we  have both $ \Gamma\lesssim  S$ and  $S\lesssim  \Gamma.$
\end{notation}

\begin{definition}\label{D:Hc}
 \rm Let $\Hc=\Hc(\U)$ be the class of all   real $(1,1)$-forms $H$ on $\U$  which  can be written as
 $$
 H=\sum f_{pq'} dz_p\wedge d\bar w_{q'}+ \sum g_{p'q} d\bar z_{p'}\wedge d w_q,
 $$
 where  $f_{pq'}$ and  $g_{p'q}$ are  bounded functions. 
\end{definition}

Now  we place ourselves under the Extended Standing  Hypothesis introduced in Subsection \ref{SS:Ex-Stand-Hyp}.
Since $\tau$ is  strongly  admissible, we infer  from  Definition \ref{D:Strongly-admissible-maps} that the  following    estimates of $1$-forms  for 
the components of $\tau=(s_1,\ldots, s_k)$ in the  local coordinates $y=(z,w).$
Note  that $s_j= \tau^*z_j$  for $1\leq j\leq k-l$  and  $s_j= \tau^*w_{j-k+l}$  for  $k-l<j\leq k.$
 \begin{equation}\label{e:diff-dz}
 d(\tau^*z_j)-dz_j=\sum_{p=1}^{k-l}  O(\|z\|) dz_p +  O(\|z\|^2)\quad\text{and}\quad 
 d(\tau^*\bar z_j)-d\bar z_j=\sum_{p=1}^{k-l}  O(\|z\|) d\bar z_p +  O(\|z\|^2).
 \end{equation}
 \begin{equation}\label{e:diff-dw}
 d(\tau^*w_m)-dw_m=\sum_{p=1}^{k-l}  O(1) dz_p + O(\|z\|)\quad\text{and}\quad
 d(\tau^*\bar w_m)-d\bar w_m=\sum_{p=1}^{k-l}  O(1) d\bar z_p +  O(\|z\|).
 \end{equation}
Using  this we  infer   the  following  estimates for the    change  under $\tau$ of  a $\Cc^1$-smooth function, of  a $\Cc^1$-smooth $(1,1)$-form, and of the basic $(1,1)$-forms  $\pi^*\omega,$ $\beta,$ $\hat\beta.$ 
\begin{proposition}\label{P:basic-admissible-estimates-I}
Let $\gamma$ be a $\Cc^1$-smooth $(1,1)$-form on $V_0.$ Then
there  are constants $c_3,c_4>0$ such that   $ c_3 \pi^*\omega +c_4\beta\geq 0$  on  $\pi^{-1}(V_0)\subset \E$  and that
for every $1\leq \ell \leq \ell_0,$
 the following inequalities hold  on $\U_\ell \cap \Tube(B,\bfr):$
\begin{enumerate}

\item   $|\tilde \tau_\ell^*(\varphi)-\varphi|\leq c_3 \varphi^{3\over 2},$ and   $|\tilde \tau_\ell^*(f)-f|\leq c_3 \varphi^{1\over 2}$  for every  $\Cc^1$-smooth  function $f$ on  $\Tube(B,\bfr);$
 \item $\pm\big(   \tilde \tau_\ell^*(\pi^*\gamma)  -\pi^*\gamma-H\big)^\sharp \lesssim  c_3 \varphi^{1\over 2}  \pi^*\omega +c_4\varphi^{1\over 2} \beta,   $
  and 
 $\tilde\tau_\ell^*(\pi^*\omega)  -\pi^*\omega\trianglelefteq   c_3 \varphi^{1\over 2}  \pi^*\omega +c_4\varphi^{1\over 2} \beta;$
 \item $\pm\big(   \tilde \tau_\ell^*(\beta) -\beta  \big)^\sharp \lesssim c_3 \phi^{3\over 2}\cdot  \pi^*\omega +c_4\phi^{1\over 2}\cdot \beta,$ and 
 $\pm\big(   \tilde \tau_\ell^*(\beta) -\beta  \big) \trianglelefteq c_3 \phi^{3\over 2}\cdot  \pi^*\omega +c_4\phi^{1\over 2}\cdot \beta;$
   \item $\pm\big(   \tilde \tau_\ell^*(\hat\beta) -\hat\beta  \big)^\sharp \lesssim c_3 \phi^{3\over 2}\cdot  \pi^*\omega +c_4\phi^{1\over 2}\cdot \hat\beta,$ and 
 $\pm\big(   \tilde \tau_\ell^*(\hat \beta) -\hat\beta  \big) \trianglelefteq c_3 \phi^{3\over 2} \cdot  \pi^*\omega +c_4\phi^{1\over 2}\cdot \hat\beta.$
\end{enumerate}
Here, in the first  inequalities of  (2)-(3)-(4), $H$ is  some  form in the class $\Hc$ given in Definition  \ref{D:Hc}.
\end{proposition}
 \proof
 The case  where $\gamma:=\omega$ has been proved in Proposition 7.8 and  Proposition 7.9 in \cite{Nguyen21}. However,  the proof therein is also  valid  for  the  general case of  a $\Cc^1$-smooth $(1,1)$-form $\gamma$ on $V_0.$
 \endproof
 The  following technical lemma is very often  needed.

\begin{lemma}\label{L:basic-positive-difference-bis}
 Let $T$ be a  positive current of bidgree $(p,p)$ on $\bfU.$
 Let  $R_1,\ldots,R_{k-p}$  be  real $(1,1)$-currents on $\Tube(B,\bfr)\subset \E,$ let   
 $S_1,\ldots,S_{k-p}$ and $S'_1,\ldots,S'_{k-p}$  be  positive $(1,1)$-currents on $\Tube(B,\bfr)\subset \E,$
 and for each $1\leq\ell\leq \ell_0,$ 
 let $H_{\ell,1},\ldots,H_{\ell,k-p}$ be    real $(1,1)$-forms  in the class $\Hc$ on $\U_\ell$    such that
 \begin{eqnarray*}
 \varphi^{1\over 2}R_j&\lesssim& S_j\quad  \text{and}\quad \varphi^{1\over 2}R_j\lesssim S'_j  \quad\mbox{on $\Tube(B,\bfr)$  for  $ 1\leq j\leq k-p;$}\\
 \pm[(\tilde \tau_\ell)^* R_j-R_j-H_{\ell,j}]^\sharp&\lesssim & S_j \quad\mbox{on $\U_\ell$  for  $1\leq \ell\leq \ell_0$ and $ 1\leq j\leq k-p;$}\\
 \big( (\tilde \tau_\ell)^* R_1-R_1,\ldots,    (\tilde \tau_\ell)^* R_{k-p}-R_{k-p}) &\trianglelefteq &\big(S'_1,\ldots, S'_{k-p}\big)\quad \mbox{on $\U_\ell$  for  $1\leq \ell\leq \ell_0.$}
 \end{eqnarray*}
Let $0<s<r\leq \bfr$ and set  $R:=R_1\wedge \ldots \wedge R_{k-p}.$ 
Suppose in addition that there are constants $0<c_5<1$ and $c_6>1$ and    positive $(1,1)$-forms $R'_1,\ldots,R'_{k-p}$  such that 
\begin{itemize}
 \item $R'_j\geq R_j\geq -R'_j$  for $1\leq j\leq k-p;$ 
\item if  $y\in  \U_\ell$  with  $0<\theta_\ell(y)< c_5,$  then  we may find  $1\leq \ell'\leq \ell_0$ 
and  an open neighborhood $\U_y$  of $y$ in  $\U$  such that for $x\in\U_y,$ we have that $\theta_{\ell'}(x)>  c_5$  and that $$ -c_6(\tilde\tau_{\ell'}\circ \tilde\tau^{-1}_\ell)^* R'_j(x)\leq R_j(x)\leq  c_6(\tilde\tau_{\ell'}\circ \tilde\tau^{-1}_\ell)^* R'_j(x)$$ and  that $S_j(x)\leq  c_6(\tilde\tau_{\ell'}\circ \tilde\tau^{-1}_\ell)^* S_j(x)$  and that
 $S'_j(x)\leq  c_6(\tilde\tau_{\ell'}\circ \tilde\tau^{-1}_\ell)^* S'_j(x).$   
\end{itemize}
Then there is a constant 
$c$ that depends on $c_5,c_6$  and $\ell_0$ such that
 \begin{multline*}
\left   |                \langle \tau_*T ,\ind_{\Tube(B,r)}R\rangle  -\langle T^\hash_r ,R\rangle \right |^2\\
\leq  c\cdot\, \sum_{\ell=1}^{\ell_0}  \sum\limits_{I,J,K}\sum_{j=0}^{|I|} \big( \int
 (\ind_{\Tube(B,r+c_0r^2)}\circ \tilde\tau_\ell) (\pi^*\theta_\ell)(\tau_\ell)_* T 
\wedge R'_K\wedge  \pi^*\omega^j\wedge \hat\beta^{|I|-j} \wedge S_J\wedge S'_{(I\cup J\cup K)^\bfc}\big)\\
   \cdot    \big(  \int  (\ind_{\Tube(B,r+c_0r^2)}\circ \tilde\tau_\ell) (\pi^*\theta_\ell)(\tau_\ell)_* T        \wedge R'_K \wedge \pi^*\omega^{|I|-j}\wedge \hat\beta^{j} \wedge S_J\wedge S'_{(I\cup J\cup K)^\bfc}\big).
   \end{multline*} 
\begin{multline*}
\left   |                \langle \tau_*T ,\ind_{\Tube(B,s,r)}R\rangle  -\langle T^\hash_{s,r} ,R\rangle \right |^2\\
\leq 
  c\cdot\, \sum_{\ell=1}^{\ell_0}  \sum\limits_{I,J,K}\sum_{j=0}^{|I|} \big( \int
 (\ind_{\Tube(B,s-c_0s^2,r+c_0r^2)}\circ \tilde\tau_\ell) (\pi^*\theta_\ell)(\tau_\ell)_* T 
\wedge R'_K\wedge  \pi^*\omega^j\wedge \hat\beta^{|I|-j} \wedge S_J\wedge S'_{(I\cup J\cup K)^\bfc}\big)\\
   \cdot    \big( \int   (\ind_{\Tube(B,s-c_0s^2,r+c_0r^2)}\circ \tilde\tau_\ell) (\pi^*\theta_\ell)(\tau_\ell)_* T        \wedge R'_K \wedge \pi^*\omega^{|I|-j}\wedge \hat\beta^{j} \wedge S_J\wedge S'_{(I\cup J\cup K)^\bfc}\big).
\end{multline*}
Here, the  sum $\sum_{I,J,K}$ is  taken over all $I,J,K\subset \{1,\ldots,k-p\}$  such that    $H_j\not\equiv 0$ for $j\in I,$ and that $I,J,K$ are mutually disjoint,
   and $|(I\cup J\cup K)^\bfc|$ is  even, and  $K\not=\{1,\ldots,k-p\}.$
Moreover, 
$$
S_J:=\bigwedge_{j\in J} S_j,\qquad R'_K:=\bigwedge_{j\in K} R'_j,\qquad S'_J:=\bigwedge_{j\in J} S'_j.
$$
   \end{lemma}

\proof
The  case  where     $R_1,\ldots,R_{k-p}$   are  positive $(1,1)$-currents on $\Tube(B,\bfr)\subset \E$  has been established in \cite [Lemma  7.22]{Nguyen21}. The proof therein is  still valid  for the   general case 
using   the     additional lower estimates on $R_j$  in the above two $\bullet.$
\endproof
The  following   convergence test  established in \cite[ Lemma 5.2 ]{Nguyen21}  will be  needed.
\begin{lemma}\label{L:convergence-test}
Let $0<r_1<r_2\leq\bfr.$   Consider  two functions
 $f:\ (0,\bfr]\to \R$  and  $\epsilon:\   [\bfr^{-1},\infty)\to (0,\infty),$ $ \lambda\mapsto\epsilon_\lambda$
  such that
 \begin{itemize}\item[(i)]   there are  two constants $c>0$ and $N\in\N$
 such that  if $2^n\leq  \lambda <2^{n+1}$ and $2^{n-N}>\bfr^{-1},$   then 
 $\epsilon_\lambda\leq c\sum_{j=- N}^N \epsilon_{2^{n+j}};$ 
 \item[(ii)]  $\sum_{n\in\N: 2^n\geq \bfr^{-1}}\epsilon_{2^n}<\infty;$
 \item[(iii)] For $r\in(r_1,r_2), $  we have $
 f({r\over \lambda})-f({r_1\over\lambda})\geq  -\epsilon_\lambda.$  
 \end{itemize}
 \begin{enumerate}
  \item Then we have 
$\lim_{r\to 0} f(r)=\liminf_{r\to 0} f(r)\in\R\cup\{-\infty\}.$

\item If instead of condition (iii) we  have the following stronger  condition (iii'): 
$$ |f({r_2\over \lambda})-f({r_1\over\lambda})|\leq  \epsilon_\lambda,$$
 then $\lim_{r\to 0} f(r)=\liminf_{r\to 0} f(r)\in\R,$ that is, the last limit is  finite.

\end{enumerate}
\end{lemma}

\subsection{Local and global mass indicators for positive currents}\label{SS:Local-global-mass-indicators}

We use the notation  introduced  in    Subsection \ref{SS:Ex-Stand-Hyp}. 
 Following the  model  formula \eqref{e:T-hash},
 we introduce  the following  mass indicators   for a  positive  current $T$  of bidegree $(p,p)$ defined on $X.$
 For $0\leq j\leq k$ and $0\leq q\leq k-l$ and  $1\leq\ell\leq \ell_0,$ and  for $0<s<r\leq \bfr,$ 
 \begin{equation}\label{e:local-mass-indicators-bis}
 \begin{split}
 \Mc_j(T,r,\tau_\ell)&:={1\over  r^{2(k-p-j)}}\int (\ind_{\Tube(B,r)}\circ \tilde\tau_\ell)  (\pi^*\theta_\ell)\cdot (\tau_\ell)_*( T|_{\bfU_\ell})
 \wedge\pi^*\omega^j\wedge  \hat\beta^{k-p-j},\\
 \Kc_{j,q}(T,r,\tau_\ell)&:=\int(\ind_{\Tube(B,r)}\circ \tilde\tau_\ell)  (\pi^*\theta_\ell)\cdot (\tau_\ell)_*( T|_{\bfU_\ell})
 \wedge\pi^*\omega^j\wedge  \hat\beta^{k-p-q-j}\wedge \hat \alpha^q,\\
 \Kc_{j,q}(T,s,r,\tau_\ell)&:=\int (\ind_{\Tube(B,s,r)}\circ \tilde\tau_\ell) (\pi^*\theta_\ell)\cdot (\tau_\ell)_*( T|_{\bfU_\ell})
 \wedge\pi^*\omega^j\wedge  \hat\beta^{k-p-q-j}\wedge \hat \alpha^q.
 \end{split}
 \end{equation}
\begin{remark}\rm A consideration on  bidegree  using  the  formula for $\upm$ given in \eqref{e:m}
shows 
that $(\tau_\ell)_*( T|_{\bfU_\ell})
 \wedge\pi^*\omega^\upm$ is of full bidegree $(l,l)$ in $\{dw,d\bar w\},$ see also 
Corollary 4.8 in \cite{Nguyen21}.
Consequently,
by the bidegree reason, we  deduce  that  $\Mc_j(T,r,\tau_\ell),$ $\Kc_{j,q}(T,r,\tau)$ and $\Kc_{j,q}(T,s,r,\tau)$
are equal to $0$ provided that $j>\upm.$ 
\end{remark}
 
  We define   the following  global mass indicators.
 \begin{equation}\label{e:global-mass-indicators}
 \begin{split}
 \Mc_j(T,r)=\Mc_j(T,r,\Uc)&:=\sum_{\ell=1}^{\ell_0}\Mc_j(T,r,\tau_\ell) ,\\
 \Mc^\tot(T,r)=\Mc^\tot(T,r,\Uc)&:=\sum_{j=0}^\upm \Mc_j(T,r),\\
 \Kc_{j,q}(T,r)=    \Kc_{j,q}(T,r,\Uc) &:=\sum_{\ell=1}^{\ell_0}\Kc_{j,q}(T,r,\tau_\ell),\\
 \Kc_{j,q}(T,s,r)=    \Kc_{j,q}(T,s,r,\Uc) &:=\sum_{\ell=1}^{\ell_0}\Kc_{j,q}(T,s,r,\tau_\ell).
 \end{split}
 \end{equation}

 \begin{lemma}
 \label{L:Mc_j}
  \begin{equation*}
 \begin{split}
 \Mc_j(T,r)&= {1\over r^{2(k-p-j)}}\int  T^\hash_r\wedge\pi^*\omega^j\wedge  \hat\beta^{k-p-j} ,\\
 \Kc_{j,q}(T,r) &=     \int  T^\hash_r\wedge
 \pi^*\omega^j\wedge  \hat\beta^{k-p-q-j}\wedge \hat\alpha^q ,\\
 \Kc_{j,q}(T,s,r) &=\int  T^\hash_{s,r}
 \wedge\pi^*\omega^j\wedge  \hat\beta^{k-p-q-j}\wedge \hat\alpha^q.
 \end{split}
 \end{equation*}
  \end{lemma}

 \subsection{Inequalities of  mass indicators}\label{SS:Inequa-mass-indicators}

Let $\omega$ be  a Hermitian  metric  on $V.$

Let $\bfj=(j_1,j_2,j_3,j_4)$ with $j_1,j_3,j_4\in\N $ and $j_2\in{1\over 4}\N,$ and 
$k-p-j_1-j_3\geq 0.$ Let $0\leq s<r\leq\bfr.$ Consider a positive current $T$  on $\bfU.$  Define
\begin{equation}\label{e:I_bfj}\begin{split}  I_{\bfj}(s,r)&:=\int_{\Tube(B,s,r)}\tau_*T\wedge \varphi^{j_2}(c_1-c_2\varphi)^{j_4}\hat\beta^{k-p-j_1-j_3}\wedge
(\pi^*\omega)^{j_3}\wedge \hat\alpha^{j_1},\\ 
I^\hash_{\bfj}(s,r)&:=\int_{\Tube(B,s,r)}T^\hash_{s,r}\wedge \varphi^{j_2}(c_1-c_2\varphi)^{j_4}\hat\beta^{k-p-j_1-j_3}\wedge
(\pi^*\omega)^{j_3}\wedge \hat\alpha^{j_1}.
\end{split}
\end{equation}
We define  $I_\bfj(r)$ and $I^\hash_\bfj(r)$ similarly  replacing the current $T^\hash_{s,r}$ (resp.  the domain of integration $\Tube(B,s,r)$) by
$T^\hash_r$  (resp. $\Tube(B,r)$).
\begin{remark}\label{R:I_bfj}
 \rm  Observe that $\Kc_{j,q}(T,r)=I^\hash_{(q,0,j,0)}(r)$ and $\Kc_{j,q}(T,s,r)=I^\hash_{(q,0,j,0)}(s,r).$
\end{remark}
 Let $\omega^{(j_3)}$ be a  $\Cc^1$-smooth  real  $(1,1)$-form  on $V_0.$ Define also
\begin{equation}\label{e:wide_I_bfj}\begin{split}  \widetilde I_{\bfj}(s,r)&:=\int_{\Tube(B,s,r)}\tau_*T\wedge \varphi^{j_2}(c_1-c_2\varphi)^{j_4}\hat\beta^{k-p-j_1-j_3}\wedge
(\pi^*\omega^{(j_3)})\wedge \hat\alpha^{j_1},\\ 
\widetilde I^\hash_{\bfj}(s,r)&:=\int_{\Tube(B,s,r)}T^\hash_{s,r}\wedge \varphi^{j_2}(c_1-c_2\varphi)^{j_4}\hat\beta^{k-p-j_1-j_3}\wedge
(\pi^*\omega^{(j_3)})\wedge \hat\alpha^{j_1}.
\end{split}
\end{equation}
We define  $ \widetilde I_\bfj(r)$ and $ \widetilde I^\hash_\bfj(r)$ similarly  replacing the current $T^\hash_{s,r}$ (resp.  the domain of integration $\Tube(B,s,r)$) by
$T^\hash_r$  (resp. $\Tube(B,r)$).

 Here is  the first  main technical   result of this  section.
\begin{lemma}\label{L:spec-wedge}  There are   constants $c,c_0$ independent of $T$    such that the   following inequalities holds for  $0\leq s<r<\bfr:$
 \begin{equation*}\begin{split}
 |\widetilde I_\bfj(r)- \widetilde I^\hash_\bfj(r)|^2 &\leq c\big(\sum_{\bfj'} I^\hash_{\bfj'}(r+c_0r^2)\big)\big ( \sum_{\bfj''} I^\hash_{\bfj''}(r+c_0r^2) \big),\\
|\widetilde I_\bfj(s,r)- \widetilde I^\hash_\bfj(s,r)|^2 &\leq c\big(\sum_{\bfj'} I^\hash_{\bfj'}(s-c_0s^2,r+c_0r^2)\big)\big ( \sum_{\bfj''} I^\hash_{\bfj''}(s-c_0s^2,r+c_0r^2) \big).  
\end{split}
\end{equation*}
Here, on the RHS:
\begin{itemize} \item[$\bullet$] the first sum  is taken over a finite number of multi-indices    $\bfj'=(j'_1,j'_2,j'_3,j'_4)$ as above  such that  $j'_1\leq  j_1$  and $j'_2\geq j_2;$ and either ($j'_3\leq j_3$) or ($j'_3>j_3$ and $j'_2\geq j_2+{1\over 2}$).
\item  the second sum   is taken over  a finite number of multi-indices $\bfj''=(j''_1,j''_2,j''_3,j''_4)$ as above   such that   either  ($j''_1< j_1$)
or ($j''_1=j_1$ and $j''_2\geq {1\over 4}+j_2$) or ($j''_1=j_1$ and $j''_3<j_3$).
\end{itemize}
\end{lemma}

\proof  We only give the proof of the  second inequality since the proof of the first one is  similar. 
  Using the compactness of $\overline B$ in $V,$ we can  write 
  $$\omega^{(j)}= \sum_{I\in\Ic}  \gamma_{I1}\wedge  \ldots\wedge \gamma_{Ij},$$ 
  where
 $\gamma_{I1},\ldots, \gamma_{Ij}$  are $\Cc^1$-smooth forms of bidegree $(1,1)$ compactly supported in  $V_0,$ and $\Ic$ is a nonempty  finite  index set.
 Observe that
 $$
 \widetilde I_\bfj(s,r)= \sum_{I\in\Ic} \widetilde I_{I,\bfj}(s,r)\quad\text{and}\quad \widetilde I^\hash_\bfj(s,r)=\sum_{I\in\Ic} \widetilde I^\hash_{I,\bfj}(s,r),
 $$
 where $ \widetilde I_{I,\bfj}(s,r)$ (resp. $\widetilde I^\hash_{I,\bfj}(s,r)$) are obtained from
  $\widetilde I_\bfj(s,r)$  (resp. $\widetilde I^\hash_\bfj(s,r)$) by replacing the form $\omega^{(j)}$ by the form  $\gamma_{I1}\wedge  \ldots\wedge \gamma_{Ij} .$
 Since 
 $$
 |\widetilde I_\bfj(s,r)- \widetilde I^\hash_\bfj(s,r)|^2\leq (\# \Ic) \big(\sum_{I\in\Ic} |\widetilde I_{I,\bfj}(s,r)- \widetilde I^\hash_{I,\bfj}(s,r)|^2 \big),    
 $$
  where  $\#\Ic$ is the cardinal of  the set$\Ic,$  
  we are reduced  to proving  the lemma  for the case where $$\omega^{(j)}= \gamma_{1}\wedge  \ldots\wedge \gamma_{j},$$ where
 $\gamma_{1},\ldots, \gamma_{j}$  are $\Cc^1$-smooth forms of bidegree $(1,1)$ compactly supported in  $V_0.$
 We argue as in the  proof of \cite[ Lemma 8.6]{Nguyen21} using  Proposition \ref{P:basic-admissible-estimates-I} (resp.  Lemma  \ref{L:basic-positive-difference-bis}) instead of   Proposition 7.8 (resp.  Lemma  7.22) therein.
\endproof

 We conclude  the section with an inequality  
 which will be useful later on.
 This is  the last  main technical   result of this  section.
 
\begin{proposition}\label{P:Lc-finite-bis}
  For  $0<r_1<r_2\leq \bfr,$ there is  a constant $c_8>0$ such that for every $q\leq  \min(k-p,k-l)$ and 
    every positive    current $T,$     we have the following estimate for $  \lambda>1:$
 $$|\kappa_{k-p-q}(T,{r_1\over \lambda},{r_2\over \lambda},\tau)|<c_8\sum\limits_{0\leq q'\leq q,\ 0\leq j'\leq \min(\upm,k-p-q')} \Kc_{j',q'}(T,{r_1\over \lambda}-c_0\big({r_1\over \lambda}\big)^2,{r_2\over \lambda}
 +c_0\big({r_2\over \lambda}\big)^2).$$
 \end{proposition}
\proof  
Fix  $0\leq q_0\leq \min(k-p,k-l)$ and  set $j_0:=k-p-q_0.$
Set  $s:={r_1\over \lambda}$ and  $r:={r_2\over \lambda}.$  In the remainder of the proof,  we  use  formula \eqref{e:wide_I_bfj} for $ I_{\bfj}(s,r)$ and 
$\widetilde I^\hash_{\bfj}(s,r)$, where  $\omega^{(j_3)}$  with $j_3\geq j_0$ is defined by
$$
 \omega^{(j_3)}:=  \omega^{(j_0)}\wedge \omega^{j_3-j_0}. 
 $$
Recall from \eqref{e:hat-alpha} and \eqref{e:hat-beta} that
\begin{equation*} \alpha=\hat\alpha- c_1 \pi^*\omega-c_2\beta
=\hat\alpha -c_2\hat \beta +(c_2\varphi-c_1)\pi^*\omega
\quad\text{and}\quad  \beta=\hat\beta-c_1\varphi\cdot  \pi^*\omega.
 \end{equation*}
So we get  that  
\begin{eqnarray*} 
 \alpha^{q_0}
&=&(\hat\alpha -c_2\hat \beta +(c_2\varphi-c_1)\pi^*\omega)^{q_0} \\
&=&
\hat \alpha^{q_0}
+\sum_{j_1,j'_1}^q {q_0 \choose j_1}{q_0-j_1 \choose j'_1}(-c_2)^{j_1'}\cdot\,\hat\beta^{j'_1}\wedge 
((c_2\varphi-c_1)\pi^*\omega)^{q_0-j_1-j'_1}\wedge \hat\alpha^{j_1},
\end{eqnarray*}
where the last sum is  taken over all $(j_1,j'_1)$   such that   $0\leq j_1,j'_1\leq q_0$  and
$j_1+j'_1\leq  q_0$ and 
$j_1\not=q_0.$  Using this 
and the first equality of  \eqref{e:I_bfj}, we  have 
\begin{equation}\label{e:P-Lc-bullet_finite(5)}
 \begin{split}
\kappa(T,B,\omega^{(j_0)},s,r)&=\int_{\Tube(B,s,r)}\tau_*T\wedge \pi^*(\omega^{(j_0)})\wedge\alpha^{q_0}
  = I_{q_0,0,j_0,0}(T,s,r)\\&+\sum_{j_1,j'_1} {q_0 \choose j_1}
  {q_0-j_1 \choose j'_1}(-1)^{q_0-j_1-j'_1}(-c_2)^{j'_1}
 I_{j_1, 0, q_0+j_0-j_1-j'_1,q_0-j_1-j'_1}(T,s,r).
\end{split}
\end{equation}
 We rewrite   \eqref{e:P-Lc-bullet_finite(5)} as 
 \begin{equation}\label{e:P-Lc-bullet_finite(6)}
 \Ic_{j_0}(T,s,r)
 =\Ic_1+\Ic_2+\Ic_3,
\end{equation}
where
\begin{equation*}
 \begin{split}
 \Ic_1&:=I^\hash_{q_0,0,j_0,0}(T,s,r)+\sum_{j'_1,j_1}{q_0 \choose j_1}
  {q_0-j_1 \choose j'_1}(-1)^{q_0-j_1-j'_1}(-c_2)^{j'_1}\\
 &\cdot I^\hash_{j_1, 0, q_0+j_0-j_1-j'_1,q_0-j_1-j'_1}(T,s,r), \\
 \Ic_2&:= I_{q_0,0,j_0,0}(T,r)-I^\hash_{q_0,0,j_0,0}(T,s,s), \\
 \Ic_3&:=  \sum_{j'_1,j_1} {q_0 \choose j_1}
  {q_0-j_1 \choose j'_1}(-1)^{q_0-j_1-j'_1}(-c_2)^{j'_1}\\
  &\cdot\big(I_{j_1, 0, q_0+j_0-j_1-j'_1,q_0-j_1-j'_1}(T,s,r)-I^\hash_{j_1, 0, q_0+j_0-j_1-j'_1,q_0-j_1-j'_1}(T,s,r)\big)  .
\end{split}
\end{equation*}
Here the  sums $\Ic_1$ and $\Ic_3$  are  taken over all $(j_1,j'_1)$   such that   $0\leq j_1,$ $0<j'_1$ and
$j_1+j'_1\leq  q_0.$ 

Observe that  $\Ic_1$ is bounded from above  by a constant times $$\sum_{0\leq q'\leq q_0,\ 0\leq j'\leq \min(\upm,k-p-q')} \Kc_{j',q'}(T,{r_1\over \lambda},{r_2\over \lambda} ).$$
 
 Applying Lemma \ref{L:spec-wedge}  to  each  difference term  in $\Ic_2$ and $\Ic_3$ 
yields that
\begin{equation} \label{e:P-Lc-bullet_finite(7)}
|I_\bfi(r)- I^\hash_\bfi(s,r)|^2 \leq c\big(\sum_{\bfi'} I^\hash_{\bfi'}(s-c_0s^2,r+c_0r^2)\big)\big ( \sum_{\bfi''} I^\hash_{\bfi''}(s-c_0s^2,r+c_0r^2) \big).  
\end{equation}
Here, on the LHS  $\bfi=(i_1,i_2,i_3,i_4)$  is  either $(q_0,0,j_0,0)$ or $(j_1,0, q_0+j_0-j_1-j'_1,q_0-j_1-j'_1)$  with $j_1,j'_1$ as above, 
and on the RHS:
\begin{itemize} \item[$\bullet$] the first sum  is taken over a finite number of multi-indices    $\bfi'=(i'_1,i'_2,i'_3,i'_4)$ as above  such that  $i'_1\leq  i_1$  and $i'_2\geq i_2;$ 
\item  the second sum   is taken over  a finite number of multi-indices $\bfi''=(i''_1,i''_2,i''_3,i''_4)$ as above   such that   either  ($i''_1< i_1$)
or ($i''_1=i_1$ and $i''_2\geq {1\over 4}+i_2$) or ($i''_1=i_1$ and $i''_3<i_3$).
\end{itemize}
 Observe that when $\bfr$ is  small  enough,  $c_1-c_2\varphi\approx  1$ and $\varphi \lesssim \bfr^2\ll 1$ on $\Tube(B,\bfr).$
 Therefore, $I_{i_1,i_2,i_3,i_4}(T,s,r)\leq c I_{i_1,0,i_3,0}(T,s,r)$ for a constant $c>0$ independent of $T$ and $0<r\leq \bfr.$
 Consequently, each  sum on the RHS of \eqref{e:P-Lc-bullet_finite(7)} is  bounded from above by a constant times  
 $$\sum_{0\leq q'\leq q_0,\ 0\leq j'\leq \min(\upm,k-p-q')} \Kc_{j',q'}(T,{r_1\over \lambda},{r_2\over \lambda} ).$$
 This, combined with \eqref{e:P-Lc-bullet_finite(5)}--\eqref{e:P-Lc-bullet_finite(6)}
 gives the result.
\endproof

 \section{The generalized Lelong numbers for positive closed currents}
 \label{S:Lelong-numbers-for-closed-currents}

Using   the results developed in Section \ref{S:Finiteness} and the arguments in the proof of \cite[Theorem 3.7]{Nguyen21},  we will prove   Theorem \ref{T:Nguyen-2} in this  section.

\subsection{Finiteness of the mass indicator  $\Kc_{j,q}$}

Fix  an open neighborhood $\bfW$ of $\partial B$ in $X$ with $\bfW\subset \bfU.$

\begin{definition}\label{D:sup}\rm
Fix an open neighborhood $\bfU$ of $\overline B$ and an open neighborhood $\bfW$ of $\partial B$ in $X$ with $\bfW\subset \bfU.$  Let $m,m'\in\N$ with $m\geq m'.$
Let $\widetilde\CL^{p;m,m'}(\bfU,\bfW)$ be the  set of all $T\in \CL^{p;m,m'}(\bfU,\bfW)$  whose  a sequence of approximating  forms $(T_n)_{n=1}^\infty$ (see Definition  \ref{D:Class})
satisfies the following   condition:
 \begin{equation}\label{e:unit-CL-1,1} \|T_n\|_{\bfU}\leq  1\qquad\text{and}\qquad  \| T_n\|_{\Cc^1(\bfW)}\leq 1.\end{equation}
 \end{definition} 
  Recall  that $\omega$ is a  Hermitian  metric on $V.$
 \begin{theorem}\label{T:Lc-finite}
  Assume that one of the following two  conditions   is fulfilled:
  \begin{enumerate}
   \item $(m,m')=(1,1)$ and $\omega$ is  K\"ahler.
  \item  $(m,m')=(2,2)$ and  $\ddc \omega^j=0$ for $\lowm\leq j\leq \upm-1.$
 \end{enumerate}
   Then there is a constants $c_7>0$ such that for every positive closed  current $T$ on $\bfU$ belonging to the class $\widetilde\CL^{p;m,m'}(\bfU,\bfW),$  we have  
   $\Kc_{j,q}(T,\bfr)<c_7.$
 \end{theorem}
 \proof
 Under  condition (1) the result follows from \cite[Theorem 8.7]{Nguyen21}. Under  condition (2) the result follows from \cite[Theorem 14.10]{Nguyen21}.
 \endproof
 
\subsection{Existence of Lelong numbers}

This  subsection  is  devoted to the proof of assertions (1)--(3) of  Theorem \ref{T:Nguyen-2}.

\proof[Proof of assertion (1) of Theorem \ref{T:Nguyen-2}]
First assume  that the current $T$ is a closed $\Cc^1$-smooth form. 
 Since $\omega^{(j)}$ is  closed,
 we have for $1\leq j\leq   \upm$ that
\begin{equation*}
d [(\tau_*T)\wedge \pi^*\omega^{(j)}]=d(\tau_*T)\wedge \pi^*\omega^{(j)}
=(\tau_*dT)\wedge \pi^*\omega^{(j)}=0 .
\end{equation*}
Applying Theorem  \ref{T:Lelong-Jensen-closed}  to  $\tau_*T\wedge \pi^*(\omega^{(j)})$  and using   the  above  equality,  we get  that
\begin{equation}\label{e:Lelong-closed-all-degrees-bis-asser(1)-Kaehler}
  \nu(T,B,\omega^{(j)}, r_2,\tau)- \nu(T,B,\omega^{(j)},r_1,\tau)= \int_{\Tube(B,r_1,r_2)}\tau_*T\wedge \pi^*(\omega^{(j)})\wedge\alpha^{k-p-j}+\Vc(\tau_*T\wedge \pi^*(\omega^{(j)}),r_1,r_2).
\end{equation}
On  the other hand,  since $j\geq \lowm$ we  get that  $k-p-j\leq  k-l.$  Therefore, we can apply
Theorem \ref{T:vertical-boundary-closed} to the current $\tau_*T\wedge \pi^*(\omega^{(j)}),$ which gives  that  $\Vc(\tau_*T\wedge \pi^*(\omega^{(j)}),r_1,r_2)=O(r_2).$ This proves  assertion (1) in  the special case where $T$ is $\Cc^1$-smooth.

Now  we consider the general case where $T$ is a general positive closed  $(p,p)$-current such that $T=T^+-T^-,$ where $T^\pm$ are approximable
along $B\subset V$ by positive closed $\Cc^1$-smooth  $(p,p)$-forms $(T^\pm_n)$  with $\Cc^1$-control on boundary. 
So $T^+_n\to T^+$ and $T^-_n\to T^-$ as $n$ tends to infinity.  
By the  previous case applied  to $T^\pm_n,$ we get that
\begin{equation*}
\nu( T^\pm_n,B,\omega^{(j)},r_2,\tau)- \nu(T^\pm_n,B,\omega^{(j)},r_1,\tau)=\kappa(T^\pm_n,B,\omega^{(j)},r_1,r_2,\tau) +O(r_2).
\end{equation*}
Letting  $n$ tend to infinity, we infer that
\begin{equation*}
\nu( T^\pm,B,\omega^{(j)},r_2,\tau)- \nu(T^\pm,B,\omega^{(j)},r_1,\tau)=\kappa(T^\pm,B,\omega^{(j)},r_1,r_2,\tau)+O(r_2).
\end{equation*}
This implies assertion (1)  since $T=T^+-T^-.$
\endproof

\proof[Proof of assertion (2) of Theorem \ref{T:Nguyen-2}]  Let $q:=k-p-j.$ Fix $r_1,\ r_2\in  (0,\bfr]$ with $r_1<r_2/2.$
Applying   Proposition  \ref{P:Lc-finite-bis}  yields for $ \lambda>1$ that
 \begin{equation}\label{e:T:Lelong-closed(1)}|\kappa_{j}(T^\pm,B, \omega^{(j)}, {r_1\over \lambda} ,{r_2\over \lambda},\tau)|<c_8\sum_{0\leq q'\leq q,\ 0\leq j'\leq \min(\upm,k-p-q')} \Kc_{j',q'}\big(T,{r_1\over \lambda}-c_0\big({r_1\over \lambda}\big)^2,{r_2\over \lambda}+c_0\big({r_2\over \lambda}\big)^2 \big).
 \end{equation}
 On the other  hand,  since there is an  $M\in\N$ such that
 $$  1\leq \#\left\lbrace n\in\N:\  y\in  \Tube\big (B, {r_1\over 2^n}-c_0 {r^2_1\over 4^n},{r_2\over 2^n}+c_0 {r^2_2\over 4^n}\big)\right\rbrace \leq M\quad\text{for}\quad y\in\Tube(B,\bfr),$$
 it follows that
 $$
 \sum_{n=1}^\infty \sum_{0\leq q'\leq q,\ 0\leq j'\leq \min(\upm,k-p-q')} \Kc_{j',q'}
 \big (T, {r_1\over 2^n}-c_0 {r^2_1\over 4^n},{r_2\over 2^n}+c_0 {r^2_2\over 4^n}\big)
 \leq M
  \sum_{0\leq q'\leq q,\ 0\leq j'\leq \upm} \Kc_{j',q'}(T,\bfr ).
 $$
 By Theorem \ref{T:Lc-finite} the  RHS is  finite.  Therefore, we infer from \eqref{e:T:Lelong-closed(1)} that
\begin{eqnarray*}
\sum_{n=1}^\infty\big|\kappa_{j}\big(T,B, \omega^{(j)},{r_1\over 2^n}-c_0 {r^2_1\over 4^n},{r_2\over 2^n}+c_0 {r^2_2\over 4^n},\tau\big)\big|&\leq&  \sum_{n=0}^\infty\big|\kappa_{j}\big(T^+,B, \omega^{(j)},{r_1\over 2^n}-c_0 {r^2_1\over 4^n},{r_2\over 2^n}+c_0 {r^2_2\over 4^n},\tau\big)\big|\\
&+&\sum_{n=0}^\infty\big|\kappa_{j}\big(T^-
,B, \omega^{(j)},
{r_1\over 2^n}-c_0 {r^2_1\over 4^n},{r_2\over 2^n}+c_0 {r^2_2\over 4^n},\tau\big)\big|\\
&\leq &   M c_8\sum_{0\leq q'\leq q,\ 0\leq j'\leq \min(\upm,k-p-q')} \Kc_{j',q'}(T,\bfr )
<\infty.
\end{eqnarray*}
Now  we apply  Lemma  \ref{L:convergence-test} (2)  to functions $f^\pm$  defined by
$
f^\pm(r):=\nu(T^\pm,B, \omega^{(j)},r,\tau)$  and
$$
\epsilon^\pm_\lambda:=2c_8    
 \lambda^{-1}+c_8\sum_{0\leq q'\leq q,\ 0\leq j'\leq \min(\upm,k-p-q')} \Kc_{j',q'}\big (T, {r_1\over 2^n}-c_0 {r^2_1\over 4^n},{r_2\over 2^n}+c_0 {r^2_2\over 4^n}\big).
$$
By  assertion  (1) and inequality \eqref{e:T:Lelong-closed(1)},  we have  by increasing the constant $c_8$  if necessary:
$$
|f^\pm({r_2\over \lambda})-f^\pm({r_1\over\lambda})|=|\kappa_j(T^\pm,B, \omega^{(j)},{r_1\over \lambda},{r_2\over \lambda},\tau )+O(\lambda^{-1})|\leq \epsilon_\lambda.
$$
Hence, assertion (2) follows.
\endproof

\proof[Proof of assertion (3) of Theorem \ref{T:Nguyen-2}]
We start  with the proof of the first part of assertion (3). 
By  \eqref{e:Lelong-log-numbers-bullet} and   assertion (1), we have 
\begin{eqnarray*}\
  \kappa^\bullet(T,B,\omega^{(j)},r,\tau)= \limsup\limits_{s\to0+}   \kappa(T,B,\omega^{(j)},s,r,\tau)&=& \nu(T,B,\omega^{(j)},r,\tau)-\liminf_{s\to0+} \nu(T,B,\omega^{(j)},s,\tau)\\
 &=&  \nu(T,B,\omega^{(j)},r,\tau)- \nu(T,B,\omega^{(j)},\tau),
\end{eqnarray*} 
where the last equality holds by assertion (2).   Consequently, we  infer  from  assertion  (2) again that
\begin{equation*}
 \lim\limits_{r\to0+}\kappa^\bullet(T,B,\omega^{(j)},r,\tau)= \lim\limits_{r\to0+}\nu(T,B,\omega^{(j)},r,\tau)- \nu(T,B,\omega^{(j)},\tau)=0.
\end{equation*}

We turn  to  the proof of the  second  part  of assertion (3). 
  
First, 
we  will prove  the interpretation of assertion (3)  in the spirit  of \eqref{e:Lelong-number-point-bisbis(1)}.  Since  $q:=k-p-j_0<k-l,$
we infer from   Theorems  \ref{T:Lelong-Jensen-smooth-closed}  and \ref{T:vertical-boundary-closed}  that
$$
\kappa(T^\pm_n,B,\omega^{(j_0)},r,\tau)=\nu( T^\pm_n,B,\omega^{(j_0)},r,\tau)+O(r).
$$
Consequently,
\begin{eqnarray*}
\kappa(T,B,\omega^{(j_0)},r,\tau)&:=& \lim\limits_{n\to\infty} \kappa(T^+_n-T^-_n,B,\omega^{(j_0)},r,\tau)\\
&=&\lim\limits_{n\to\infty}\nu(T^+_n,B,\omega^{(j_0)},r,\tau)-\lim\limits_{n\to\infty}\nu(T^-_n,B,\omega^{(j_0)},r,\tau)+O(r)\\
&=&\nu(T^+,B,\omega^{(j_0)},r,\tau)-\nu(T^-,B,\omega^{(j_0)},r,\tau)+O(r) =\nu(T,B,\omega^{(j_0)},r,\tau)+O(r).
\end{eqnarray*}
This  implies  the  desired interpretation  according to  Definition \ref{D:Lelong-log-numbers(1)}.

Second,  we  will prove  the interpretation of assertion (3)  in the spirit  of \eqref{e:Lelong-number-point-bisbis(2)}.
To start  with,  we  fix $0<r<\bfr$ and let  $0<\epsilon<r.$  Theorem \ref{T:Lelong-Jensen-closed-eps} applied  to  $\tau_*T\wedge \pi^*(\omega^{(j_0)})$  and using $\ddc \lbrack \tau_*T\wedge \pi^*(\omega^{(j_0)})\rbrack=0$  gives
 \begin{eqnarray*}
 {1\over  (r^2+\epsilon^2)^{k-p-j_0}} \int_{\Tube(B,r)} \tau_*T\wedge \pi^*(\omega^{(j_0)})\wedge \beta_\epsilon^{k-p-j_0}
 &= &\lim\limits_{n\to\infty} \Vc_\epsilon( \tau_*T_n\wedge \pi^*(\omega^{(j_0)}),r)\\
 &+&   \int_{\Tube(B,r)} \tau_*T\wedge \pi^*(\omega^{(j_0)})\wedge \alpha^{k-p-j_0}_\epsilon.
\end{eqnarray*}
Now  we let $\epsilon$ tend to $0.$  Then the LHS  tends to $\nu(T,B,\omega^{(j_0)},r,\tau).$  On the other hand,  we deduce  from \eqref{e:vertical-boundary-term-eps} and  Theorem \ref{T:vertical-boundary-eps}  that $\Vc_\epsilon( \tau_*T_n\wedge \pi^*(\omega^{(j_0)}),r)=O(r).$  Consequently,    the second  term on the RHS  tends to $\nu(T,B,\omega^{(j_0)},r,\tau)+O(r).$
This  proves the  desired  interpretation  according to  Definition \ref{D:Lelong-log-numbers(2)}.
\endproof

\subsection{Independence of admissible maps}
\label{SS:Independence}

In this  subsection we will prove assertions (5) and (6) of Theorem \ref{T:Nguyen-2}.

Fix $j_0$  with  $\lowm \leq j_0\leq \upm$ and  fix a closed   smooth real $(j_0,j_0)$-form $\omega^{(j_0)}$ on $V_0.$  Writing $\omega^{(j_0)}=(\omega^{(j_0)}+c\omega^{j_0})- c\omega^{j_0}$ for a large constant $c>0$ so that $\omega^{(j_0)}+c\omega^{j_0}$ is a strictly positive  form,
we are reduced to working on two  strictly positive  forms 
  $\omega^{(j_0)}+c\omega^{j_0}$ and $ c\omega^{j_0}.$  Consequently, we may assume without loss of generality that $\omega^{(j_0)}$  is strictly positive  form on $V_0.$ 
 In the remainder of the subsection,  we  use  formula \eqref{e:wide_I_bfj} for $ I_{\bfj}(s,r)$ and 
$\widetilde I^\hash_{\bfj}(s,r)$, where  $\omega^{(j_3)}$  with $j_3=j\geq j_0$ is defined by
$$
 \omega^{(j)}:=  \omega^{(j_0)}\wedge \omega^{j-j_0}. 
 $$
 Using the compactness of $\overline B$ in $V,$ we can  write 
  \begin{equation}\label{e:decomposition-omega-j_0} \omega^{(j_0)}= \sum_{I_0\in\Ic_0}  \gamma_{I_01}\wedge  \ldots\wedge \gamma_{I_0j_0},
  \end{equation}
  where
 $\gamma_{I1},\ldots, \gamma_{Ij}$  are $\Cc^1$-smooth forms of bidegree $(1,1)$ compactly supported in  $V_0,$ and $\Ic_0$ is a nonempty  finite  index set.
 
Consider the following mass indicators, for $j$ with $j_0 \leq j\leq \upm:$
 \begin{equation}\label{e:other-mass-indicator} \hat\nu_j(T,r):={1\over  r^{2(k-p-j)}}\int\limits_{\Tube(B,r)}
 \tau_* T\wedge (\beta+c_1r^2\pi^*\omega)^{k-p-j}\wedge  \pi^*\omega^{(j)}.
 \end{equation}
 We also consider the following new mass indicators, where  $T^\hash$ and $T^\hash_r$ are given in \eqref{e:T-hash} and \eqref{e:T-hash_r}:
\begin{equation}\label{e:other-mass-indicators-bis}
 \begin{split}
 \widehat\Mc^\hash_j(T,r)&:={1\over   r^{2(k-p-j)}}\int\limits_{\Tube(B,r)}
 T^\hash \wedge (\beta +c_1r^2\pi^*\omega)^{k-p-j}\wedge  \pi^*\omega^{(j)},\\
 \Mc^\hash_j(T,r)&:={1\over r^{2(k-p-j)}}\int
 T^\hash_r \wedge (\beta+c_1r^2\pi^*\omega)^{k-p-j}\wedge  \pi^*\omega^{(j)}
 .
 \end{split}
 \end{equation}

\begin{lemma} \label{L:Mc-hash-vs-nu} There is a  constant $c>0$ such that for every $\lowm\leq j\leq \upm$ and $0<r\leq\bfr:$
$$|\Mc^\hash_j(T,r)- \hat\nu_j(T,r)| 
\leq    cr\sum^{\upm}_{q=\lowm}  \Mc^\hash_q(T,r).$$
\end{lemma}

\proof
By Propositions \ref{P:basic-admissible-estimates-I}, for every  smooth  real $(1,1)$-form $\gamma$ on $V_0,$
there  are constants $c_3,c_4>0$ such that    $ c_3 r^2\pi^*\omega +c_4\beta\geq 0$  on  $\Tube(B,r)$  for $0<r\leq \bfr,$ and that
for every $1\leq \ell \leq \ell_0,$
 the following inequalities hold  on $\U_\ell \cap \Tube(B,r)$  for $0<r\leq \bfr:$
\begin{equation}\label{e:admissible-estimates-bis}
\begin{split}
  \pm\big(   \tilde \tau_\ell^*(\pi^*\gamma)  -\pi^*\gamma-H\big)^\sharp &\lesssim  c_3 r  \pi^*\omega +c_4r( \beta+c_1r^2\pi^*\omega),\\   
   \pm\big(   \tilde \tau_\ell^*(\beta+c_1r^2\pi^*\omega) - (\beta+c_1r^2\pi^*\omega) \big)^\sharp &\lesssim c_3 r^3  \pi^*\omega +c_4r (\beta+c_1r^2\pi^*\omega).
                \end{split}
\end{equation}
Here,  on the LHS of the  first line,  $H$ is some  form 
in the class $\Hc$ given in Definition  \ref{D:Hc}.  We also have  that
  \begin{equation}\label{e:admissible-estimates-bis(2)}
\begin{split}
  &\left\lbrace \big(\tilde\tau_\ell^*(\pi^*\gamma)  -\pi^*\gamma\big),   \big(    \tilde \tau_\ell^*(\beta+c_1r^2\pi^*\omega) -(\beta+c_1r^2\pi^*\omega)  \big) \right\rbrace\\
  &\trianglelefteq    \left\lbrace  \big( c_3 r  \pi^*\omega +c_4r (\beta+c_1r^2\pi^*\omega)\big),\big(c_3 r^3  \pi^*\omega +c_4r (\beta+c_1r^2\pi^*\omega)\big)\right\rbrace.
  \end{split}
\end{equation}
   Arguing  as  in the proof  of   Lemma  \ref{L:spec-wedge},
   we are in the  position to  apply Lemma \ref{L:basic-positive-difference-bis}. We come back the statement of this lemma. Consider  the $(1,1)$-forms $\gamma_{I1},\ldots,\gamma_{Ij_0}$ appearing in \eqref{e:decomposition-omega-j_0}.
 For  every $I_0\in\Ic_0,$ let $R_1,\ldots,R_{k-p}$ be the $k-p$ forms among $\{\pi^*\omega,\beta+c_1r^2\pi^*\omega, \gamma_{I_01},\ldots,\gamma_{I_0j_0}\}$  which appear in the integral of $\hat\nu_j(T,r)$ in 
 \eqref{e:other-mass-indicator} where $\omega^{(j_0)}$ is replaced by$  \gamma_{I_01}\wedge  \ldots\wedge \gamma_{I_0j_0}.$  So  setting  $R^{I}:=R_1\wedge \ldots \wedge R_{k-p}.$ we get $$\hat\nu_j(T,r)=\sum_{I_0\in \Ic_0} {1\over  r^{2(k-p-j)}}\int_{\Tube(B,r)}\tau_*T\wedge R^{(I_0)}.$$
 Now  we define  $R'_1,\ldots,R'_{k-p}$ as  follows. If  $R_j\in \{\pi^*\omega,\gamma_{I_01},\ldots,\gamma_{I_0j_0}\},$ set $R'_j:= c_1\pi^*\omega+\beta,$
  otherwise   $R_j =\beta+c_1r^2\pi^*\omega$ and set $R'_j:=R_j.$
 Let   $S_1,\ldots,S_{k-p}$ be  the  corresponding positive $(1,1)$-form associated to $R_1,\ldots,R_{k-p}$ respectively on the RHS  of \eqref{e:admissible-estimates-bis}.
 Let $S'_1,\ldots,S'_{k-p}$  be  the  corresponding positive $(1,1)$-form associated to $R_1,\ldots,R_{k-p}$ respectively on the RHS   of \eqref{e:admissible-estimates-bis(2)}.
 Let $H_1,\ldots,H_{k-p}$ be    the  corresponding  real $(1,1)$-forms  associated to $R_1,\ldots,R_{k-p}$ respectively on the  LHS of  each  inequality of \eqref{e:admissible-estimates-bis}.
  
  The rest of the proof follows along the  same lines as those given in the proof of \cite[Lemma 8.10]{Nguyen21} making the obviously necessary adaptation.
 \endproof

 \begin{proposition}\label{P:comparison-Mc}
 For $j_0 \leq j\leq \upm,$  we have that $$ \lim_{r\to 0+} \Mc^\hash_j(T,r)=\lim_{r\to 0+}  \hat\nu(T,r)=\sum\limits_{q=0}^{k-p-j}  {k-p-j\choose q} c^q_1\nu(T,B,\omega^{(j+q)},\tau).$$
 \end{proposition}
\proof Using formula \eqref{e:other-mass-indicator}
 and arguing as in the proof of Lemma 6.1 in \cite{Nguyen21}, we obtain  the following   identity:
 \begin{equation*}
  \hat\nu_j(T,r)=\sum\limits_{q=0}^{k-p-j}  {k-p-j\choose q} c^q_1\nu(T,B,\omega^{(j+q)},r,\tau).  
 \end{equation*}
 Next,  letting  $r$ tend to $0$
  in this identity, we infer
  from  Theorem  \ref{T:Nguyen-2} (2) that   
 \begin{equation}\label{e:Mc_j-vs-nu-bis}
  \lim_{r\to 0+} \hat\nu_j(T,r)=\sum\limits_{q=0}^{k-p-j}  {k-p-j\choose q} c^q_1\nu(T,B,\omega^{(j+q)},\tau).  
 \end{equation}
 This proves the second identity of the  proposition.
 
 It remains  to show  the first  identity. 
 Applying  Lemma   \ref{L:Mc-hash-vs-nu}  yields that 
   there is a constant $c>0$ such that for $0<r\leq\bfr,$
 \begin{equation}\label{e:Mc_j-induction}
 \big|\sum_{j=\lowm}^\upm  \Mc^\hash_j(T,r)-\sum_{j=\lowm}^\upm\hat\nu_j(T,r)\big|\leq  c r \sum_{j=\lowm}^{\upm}\Mc^\hash_j(T,r).
 \end{equation}
This, combined  with   \eqref{e:Mc_j-vs-nu-bis}, implies that there is a constant  $c>0$ such  that
$$
 \sum_{j=\lowm}^{\upm}\Mc^\hash_j(T,r)\leq c\qquad\text{for}\qquad 0<r\leq\bfr.  
$$
Therefore, we infer from Lemma   \ref{L:Mc-hash-vs-nu} that
$
|\Mc^\hash_j(T,r)- \hat\nu_j(T,r)|\leq  cr
$ for $\lowm\leq j\leq \upm.$  Letting $r$ tend to $0,$  
the first  identity of the proposition  
 follows.
 \endproof

 \proof[Proof of assertion (5) of Theorem \ref{T:Nguyen-2}]
 As  in the beginning of the  subsection, fix  $j_0$ with $\lowm\leq j_0\leq \upm.$ 
 Let $\tau$ and $\tau'$ be two   strongly admissible  maps.
 For $1\leq \ell\leq \ell_0$  we define  $\tilde\tau'_\ell:=\tau'\circ \tau^{-1}_\ell$  according to formula \eqref{e:tilde-tau_ell}.
 So   $\tilde\tau'_\ell$ is defined in the same was as  $\tilde\tau_\ell$
 using $\tau'$ instead of $\tau.$
 Similarly, we define $T^{'\hash}$ and $T^{'\hash}_r$ according to formulas \eqref{e:T-hash} and \eqref{e:T-hash_r}  by using $\tilde\tau'_\ell$ instead of $\tilde\tau_\ell.$ 
 Similarly, we  define
 $
 \widehat\Mc'^\hash_j(T,r)$ and 
 $\Mc'^\hash_j(T,r)$  according to formula \eqref{e:other-mass-indicators-bis} by using  $T^{'\hash}$ and $T^{'\hash}_r$  instead of $T^{\hash}$ and $T^{\hash}_r.$

 We need to show  that 
 \begin{equation}\label{e:tau-tau'}\nu(T,B,\omega^{(j)},\tau)=\nu(T,B,\omega^{(j)},\tau')\qquad\text{for}\qquad j_0 \leq j\leq \upm.
 \end{equation}
 By \eqref{e:admissible-estimates-bis} 
there  are constants $c_3,c_4>0$ such that    $ c_3 r^2\pi^*\omega +c_4\beta\geq 0$  on  $\Tube(B,r)$  for $0<r\leq \bfr,$ and that
for every $1\leq \ell \leq \ell_0,$
 the following inequalities hold  on $\U_\ell \cap \Tube(B,r)$  for $0<r\leq \bfr:$
\begin{equation}\label{e:admissible-estimates-bis-bis}
\begin{split}
  \pm\big(   \tilde \tau_\ell^*(\pi^*\gamma)  -(\tilde \tau'_\ell)^*(\pi^*\gamma)-H\big) &\lesssim  c_3 r  \pi^*\omega +c_4r( \beta+c_1r^2\pi^*\omega),\\   
\pm\big(   \tilde \tau_\ell^*(\beta+c_1r^2\pi^*\omega) - (\tilde \tau'_\ell)^*(\beta+c_1r^2\pi^*\omega)  \big) &\lesssim c_3 r^3  \pi^*\omega +c_4r (\beta+c_1r^2\pi^*\omega).
                \end{split}
\end{equation}
Here,    $H$ is some  form 
in the class $\Hc$ given in Definition  \ref{D:Hc}.

 By \eqref{e:admissible-estimates-bis(2)} for every $1\leq \ell \leq \ell_0,$
 the following inequality holds  on $\U_\ell \cap \Tube(B,r)$  for $0<r\leq \bfr:$
\begin{equation}\label{e:admissible-estimates-bis-bis(2)}
\begin{split}
& \left\lbrace \big(\tilde\tau_\ell^*(\pi^*\gamma)  -(\tilde \tau'_\ell)^*(\pi^*\gamma)\big),\big(    \tilde \tau_\ell^*(\beta+c_1r^2\pi^*\omega) -(\tilde \tau'_\ell)^*(\beta+c_1r^2\pi^*\omega)  \big)  \right\rbrace\\
&\trianglelefteq \left\lbrace \big(  c_3 r  \pi^*\omega +c_4r (\beta+c_1r^2\pi^*\omega)\big),\big(  c_3 r^3  \pi^*\omega +c_4r (\beta+c_1r^2\pi^*\omega)\big)\right\rbrace.
 \end{split}
 \end{equation}
Using \eqref{e:admissible-estimates-bis-bis}--\eqref{e:admissible-estimates-bis-bis(2)} and arguing as in the  proof of Lemma  \ref{L:Mc-hash-vs-nu} we can show that 
there is a  constant $c>0$ such that for every $j_0\leq j\leq \upm$ and $0<r\leq\bfr:$
$$|\Mc^\hash_j(T,r)- {\Mc'}^\hash_j(T,r)|\leq   cr\sum_{j=j_0}^{\upm}  \Mc^\hash_j(T,r).$$
Thus by  Proposition  \ref{P:comparison-Mc},  $|\Mc^\hash_j(T,r)- {\Mc'}^\hash_j(T,r)|\leq  cr.$ 
 So by this  proposition  again, we get that
\begin{equation*}
 \lim_{r\to 0}{\Mc'}^\hash_{j}(T,r)  =\lim_{r\to 0+}  \hat\nu_{j_0}(T,r)=\sum\limits_{q=0}^{k-p-j_0}  {k-p-j_0\choose q} \nu(T,B,\omega^{(j_0+q)},\tau).  
 \end{equation*}
 Hence, for $j_0\leq j\leq \upm,$ we have 
 $$
 \sum\limits_{q=0}^{k-p-j}  {k-p-j\choose q} \nu(T,B,\omega^{(j_0+q)},\tau)=\sum\limits_{q=0}^{k-p-j}  {k-p-j\choose q} \nu(T,B,\omega^{(j_0+q)},\tau').
 $$
 These  equalities imply \eqref{e:tau-tau'}.
 The proof is  thereby completed.
 \endproof

 \section{The generalized Lelong numbers for positive plurisubhamonic currents}
 \label{S:Lelong-numbers-for-plurisubhamonic-currents}
 In Subsection 
\ref{SS:m-negligible-form-and-abstract-estimates} we adapt to the present more general context the results obtained in \cite[Sections 11, 12]{Nguyen21}.
Using this    Subsection \ref{SS:basic-estimates-and-finitness-Kc-Lc}
first presents  suitable adaptation  of \cite[Section 13]{Nguyen21}. The remainder of the  subsection is devoted to  the proof of Theorem \ref{T:Nguyen-1}.
 
\subsection{$m$-negligible  test forms and  abstract  estimates }
\label{SS:m-negligible-form-and-abstract-estimates}

  Recall  from Subsection \ref{SS:Ex-Stand-Hyp} that  for every $1\leq\ell\leq\ell_0,$  there is a local  coordinate system $y=(z,w)$  on $\U_\ell$  with  $V\cap \U_\ell=\{z=0\}.$
  
 \begin{definition}
  \label{D:negligible}
  \rm  Let $S$ be  a differential  form (resp.  a  current) defined on $\Tube(B,r)\subset \E$  for some $0<r\leq\bfr.$
   So we can write  in a  local representation of $S$ in coordinates $y=(z,w)\in\C^{k-l}\times\C^l: $
   \begin{equation}\label{e:S_IJKL} S=\sum_{M=(I,J;K,L)} S_Mdz_I\wedge d\bar z_J\wedge dw_K\wedge d\bar w_L,
   \end{equation}
where the $S_M=S_{I,J;K,L}(z,w)$  are  the component  functions  (resp.  component  distributions), and the sum is taken over   $M=( I,J;K,L)$ with $I,J\subset\{1,\ldots,k-l\}$ and $K,L\subset \{1,\ldots,l\}.$

For  $M=( I,J;K,L)$ as above, we also write $dy_M$ instead of $dz_I\wedge d\bar z_J\wedge dw_K\wedge d\bar w_L.$

   Given $0\leq  m\leq 2l,$ we say that a bounded differential form $S$ is  {\it $m$-negligible}  if  in the  above  representation,
  for every $I,J,K,L,$ 
 it holds that  $S_{I,J;K,L}$ is  smooth  outside $V$ and   $S_{I,J;K,L}(z,w)=O(\|z\|^{\bfe(m,K,L)}),$
 where
 \begin{equation*}
  \bfe(m,K,L):=\max\big(0,|K|+|L|-m  \big)\in\N.
 \end{equation*}
 \end{definition}

Let $T$ be  positive plurisubharmonic  current $T$ of bidegree $(p,p)$ on $\bfU.$ 
Consider the integers 
\begin{equation}\label{e:m+}
 \lowm^+:=\max(0,l-p-1)\qquad\text{and}\qquad  \upm^+:=\min(l,k-p-1).
\end{equation}
In other words,    $\lowm^+,$ $\upm^+$ are    associated  to the  $(p+1,p+1)$-current
$\ddc T$ in the same  way  as    $\lowm,$ $\upm$ are    associated  to the  $(p,p)$-current
$T$ in formula \eqref{e:m}. 

Following  the model of \eqref{e:global-mass-indicators}, consider the following   mass indicators, for $0<r\leq\bfr,$
\begin{equation}\label{e:global-mass-indicators-bis}
 \begin{split}
 \Mc^\tot( T,r)&:=\sum_{j=0}^{\upm} \Mc_j( T,r),\quad    \Mc^\tot( \ddc T,r):=\sum_{j=0}^{\upm^+} \Mc_j( \ddc T,r) \\
 \Nc(T,r)&:= \Mc^\tot(T,r)+\Mc^\tot(\ddc T,r)=\sum_{j=0}^{\upm} \Mc_j( T,r)+\sum_{j=0}^{\upm^+} \Mc_j(\ddc T,r).
 \end{split}
 \end{equation}
where  
the $\Mc_j$'s   are defined  in \eqref{e:global-mass-indicators}.

In this section 
following   Definition \ref{D:sup}, we  introduce the following class of currents.


\begin{definition}\label{D:sup-bis}\rm
Fix an open neighborhood $\bfU$ of $\overline B$ and an open neighborhood $\bfW$ of $\partial B$ in $X$ with $\bfW\subset \bfU.$
Let $\widetilde\SH^{p;3,3}_p(\bfU,\bfW)$ be the  set of all $T\in \SH^{p;3,3}(\bfU,\bfW)$  whose  a sequence of approximating  forms $(T_n)_{n=1}^\infty$
satisfies the following   condition:
 \begin{equation}\label{e:unit-SH-3,3} \|T_n\|_{\bfU}\leq  1 \quad\text{and}\quad \| \ddc T_n\|_{\bfU}\leq 1  \quad\text{and}\quad  \| T_n\|_{\Cc^3(\bfW)}\leq 1.\end{equation}

Given   a  class of currents  $\Fc$ and  a mass indicator $\Mc(T)$   for all currents $T\in\Fc,$  
 We denote  by $\sup_{T\in\Fc}\Mc(T)$  the supremum of $\Mc(T)$  when  $T$ is taken over $\Fc.$
\end{definition}

Recall  the notation from  the  Standing  Hypothesis in Subsection \ref{SS:Ex-Stand-Hyp}.
 Fix an integer $j$ with $\lowm\leq j\leq \upm$ and  a smooth $(j,j)$-form  $\omega^{(j)}$ on $V_0.$ 
Consider  the  forms on $\bfU$:
\begin{equation}\label{e:partition-can-forms-j,l}
\Phi:= \pi^*(\omega^{(j)})\wedge \beta^{k-p-j-1} \quad\text{and}\quad \Phi^{(\ell)}:= (\pi^*\theta_\ell)\cdot\pi^*(\omega^{(j)})\wedge \beta^{k-p-j-1}\quad\text{for}\quad
1\leq\ell\leq \ell_0.
\end{equation}
So we have
\begin{equation}\label{e:sum-forms-j,l}
\Phi=\sum_{\ell=1}^{\ell_0} \Phi^{(\ell)}\qquad\text{on}\qquad \bfU.
\end{equation}
For $\ell$  with $1\leq \ell\leq \ell_0$ and   set $\tilde\tau:=\tilde\tau_\ell.$
For $r\in(0,\bfr],$   set $\H_r:=\Tube (\widetilde V_\ell,r)\subset \E.$

Let $T$ be a positive  plurisubharmonic   current on $\bfU$ in the class  $\widetilde\SH^{p;3,3}(\bfU,\bfW).$ 
Consider the  current
\begin{equation}\label{e:S_ell}
 S^{(\ell)}:=(\tau_\ell)_*(T|_{\bfU_\ell})   .
\end{equation}
By \eqref{e:T-hash} we get that
\begin{equation}
 T^\hash =\sum_{\ell=1}^{\ell_0} (\pi^*\theta_\ell)\cdot S^{(\ell)}.
\end{equation}
Note that the current $S^{(\ell)}$ is  positive plurisubharmonic on $\H_\bfr.$ Moreover, arguing  as in the  proof of  
\cite[Lemma 11.6]{Nguyen21},  we can check that  $\Phi^{(\ell)}$  is a $2j$-negligible smooth form. 
By \cite[Proposition  11.41]{Nguyen21},
there are 
 \begin{itemize}
 \item [$\bullet$]  two  functions $\Ic^{(\ell)}_1,\ \Ic^{(\ell)}_2 :\ (0,\bfr]\to\R;$ 
  \item [$\bullet$]   three differential operators $D^{(\ell)}_{10},$ $D^{(\ell)}_{11},$  $D^{(\ell)}_{12}$ in  the class $\widehat\Dc^0_\ell;$
   and three differential operators $D^{(\ell)}_{20},$ $D^{(\ell)}_{21},$  $D^{(\ell)}_{22}$   in  the class $\Dc^0_\ell;$
  \item[$\bullet$]  three smooth $2q$-forms $\Phi^{(\ell)}_{10}$   which is  $(2j-1)$-negligible,
   $\Phi^{(\ell)}_{11}$  which  is  $2j$-negligible,  $\Phi^{(\ell)}_{12}$  which  is  $(2j-1)$-negligible;
   and three smooth $2q$-forms $\Phi^{(\ell)}_{20}$   which is  $2j$-negligible,
   $\Phi^{(\ell)}_{21}$  which  is  $(2j+1)$-negligible,  $\Phi^{(\ell)}_{22}$  which  is  $2j$-negligible;
 \end{itemize}
  such that 
 every $0<r_1<r_2\leq\bfr$ and  every smooth  function $\chi$ on $(0,\bfr],$  we have for $\nu\in\{1,2\},$
 \begin{equation}\label{e:Stokes-ddc-difference-bis}
 \begin{split}
  \int_{r_1}^{r_2}\chi(t) \Ic^{(\ell)}_\nu(t)dt&= 
\int_{\Tube(B,r_1,r_2)}\chi(\|y\|) (D^{(\ell)}_{\nu 1}S^{(\ell)}\wedge \Phi^{(\ell)}_{\nu 1})(y)+\int_{\Tube(B,r_1,r_2)}\chi'(\|y\|) (D^{(\ell)}_{ \nu 2}S^{(\ell)}\wedge \Phi^{(\ell)}_{\nu 2})(y)\\
&+ \int_{\partial_\hor\Tube(B,r_2) }\chi(r_2)( D^{(\ell)}_{\nu 0}S^{(\ell)}\wedge \Phi^{(\ell)}_{\nu 0})(y)- \int_{\partial_\hor\Tube(B,r_1)} \chi(r_1)(D^{(\ell)}_{\nu 0}S^{(\ell)}\wedge \Phi^{(\ell)}_{\nu 0})(y)  ,
 \end{split}
 \end{equation}
 and  that the  following inequality holds 
 for all $  0<t\leq \bfr:$
  \begin{equation}\label{e:inequal-ddc-difference-bis}
 {1\over  r^{2(k-p-j)} } \int_{r\over 2}^r\big| \langle  \ddc [(\tilde\tau_\ell)_* S^{(\ell)}] -(\tilde\tau_\ell)_*(\ddc  S^{(\ell)}),\Phi^{(\ell)}  \rangle_{\tilde \tau (\H_t)}
 -  \Ic^{(\ell)}_1(t) -\Ic^{(\ell)}_2(t)\big|dt
  \leq 
\sum_{m=\lowm}^\upm \nu_m( S^{(\ell)},B,r,\id).
 \end{equation} 
 Recall  from  \cite[Lemma  12.2]{Nguyen21}  the following  result.
 \begin{lemma}\label{L:partition}
  The  following   equalities hold:
  \begin{eqnarray*}
     (\tilde\tau_\ell)_* S^{(\ell)}&=&\tau_*T \qquad \text{and}\quad (\tilde\tau_\ell)_*(\ddc  S^{(\ell)})=\tau_*(\ddc T)\quad\text{on}\quad \bfU_\ell,\\
     \sum_{\ell=1}^{\ell_0}  \ddc[ (\tilde\tau_\ell)_* S^{(\ell)}]\wedge \Phi^{(\ell)}  &=&\ddc (\tau_*T)\wedge \Phi\quad\text{and}\quad
      \sum_{\ell=1}^{\ell_0} (\tilde\tau_\ell)_*(  \ddc S^{(\ell)})\wedge \Phi^{(\ell)}  = \tau_*(\ddc T)\wedge \Phi
      \quad\text{on}\quad \bfU.
  \end{eqnarray*}
 \end{lemma}

Recall  from  \cite[Lemma  12.3]{Nguyen21}  the following  result.
 \begin{lemma}\label{L:Ic-0}
  Under  the above  hypotheses and  notations,  there is a constant $c$ independent of $T$  such that for $\nu\in \{1,2\}$ and for all $1\leq \ell\leq \ell_0$ and for all $0<r\leq \bfr:$
  \begin{equation*}
   {1\over  r^{2(k-p-j)} } \int_{r\over 2}^r\big|  \int_{\partial_\hor\Tube(B,t) }( D^{(\ell)}_{\nu0}S^{(\ell)}\wedge \Phi^{(\ell)}_{\nu0}) \big|dt\leq  cr^2\Mc^\tot(T,r).
  \end{equation*}
 \end{lemma}
 
 Consider  two functions $\chi_1,\chi_2: (0,r]\to \R^+$ defined by    
  \begin{equation}\label{e:chi_1-chi_2}\chi_1(t):= {t\over  r^{2(k-p-j)} }\quad\text{and}\quad
  \chi_2(t):={1\over  t^{2(k-p-j)-1} } \quad\text{for}\quad t\in (0,r].
  \end{equation}
  
Recall  from  \cite[Lemma  12.4]{Nguyen21}  the following  result.
  \begin{lemma}\label{L:Ic-1-and-2}
Under  the above  hypotheses and  notations, let  $0<r\leq \bfr.$
  Then  there is a constant $c$ independent of $T$ and $r$ such that for $\nu\in \{1,2\}$ and for all $1\leq \ell\leq \ell_0$ and for
  all $0<s<r:$
  \begin{eqnarray*}
  \big |\int_{\Tube(B,s,r)}\chi(\|y\|) (D^{(\ell)}_{\nu 1}S^{(\ell)}\wedge \Phi^{(\ell)}_{\nu1})(y)\big|  &\leq & c\sum_{n=0}^\infty {r\over 2^n}\Mc^\tot(T,{r\over 2^n}),\\
  \big|  \int_{\Tube(B,s,r)}\chi'(\|y\|) (D^{(\ell)}_{\nu2}S\wedge \Phi^{(\ell)}_{\nu2})(y) \big|&\leq & c\sum_{n=0}^\infty {r\over 2^n}\Mc^\tot(T,{r\over 2^n}).
  \end{eqnarray*}
  Here $\chi$ is  either the function $\chi_1$ or the function $\chi_2$ given in \eqref{e:chi_1-chi_2}.
 \end{lemma}

\subsection{Basic  estimates  and finiteness of the mass indicators $\Kc_{j,q}$ and $\Lc_{j,q}$}
\label{SS:basic-estimates-and-finitness-Kc-Lc}
 
 Let $\lowm\leq i\leq \upm$ and  let $\omega^{(i)}$ be a strictly positive  closed smooth $(i,i)$-form on $V_0.$
 
 Using 
Subsection 
\ref{SS:m-negligible-form-and-abstract-estimates}, 
 we first adapt to  the  present  more general context the results of \cite[Section 13]{Nguyen21}.
 We only state  these results.
\begin{lemma}\label{L:ddc-difference} {\rm  (See  \cite[Lemma 13.1]{Nguyen21}). } There is a constant $c>0$  such that  for every
$j$ with $\lowm\leq j\leq \upm,$  and every $m$ with $i\leq j-m,$ and every positive plurisubharmonic  current $T$   in the class  $\widetilde\SH^{p;3,3}(\bfU,\bfW),$
 there exists  a function  $(0,\bfr]\ni r\to \tilde r$ (depending on $T$)  with ${r\over 2}\leq \tilde r\leq r$
such  that the  following  two inequalities hold for  $0<s<r\leq \bfr:$
 \begin{multline*} \big|\int_{\tilde s}^{\tilde r} \big({1\over t^{2(k-p-j)}} -{1\over r^{2(k-p-j)}}\big)2tdt  \int_{\Tube(B,t)}  \big( \ddc (\tau_*T)- \tau_*(\ddc T)\big)\wedge \pi^*(\omega^{(i)}\wedge\omega^{j-m-i})\wedge \beta^{k-p-j+m-1}\big|\\
 \leq cr^{2m+1}\Mc^\tot(T,r),
 \end{multline*}

\end{lemma}

\begin{lemma}\label{L:ddc-difference-m-1}  {\rm  (See  \cite[Lemma 13.2]{Nguyen21}). }
There is a constant $c>0$  such that  for every
$j$ with $\lowm\leq j\leq \upm,$ and every $m$ with $1\leq m$ and $i\leq j-m,$  and every positive plurisubharmonic  current $T$   in the class  $\widetilde\SH^{p;3,3}(\bfU,\bfW),$
 the  following properties  hold for every $0<s<r\leq \bfr:$
 \begin{multline*} \Big|\int_{\tilde s}^{\tilde r} \big({1\over t^{2(k-p-j)}} -{1\over r^{2(k-p-j)}}\big)2tdt  \int_{\Tube(B,t)}   \tau_*(\ddc T)\wedge \pi^*(\omega^{(i)}\wedge\omega^{j-m-i})\wedge \beta^{k-p-j+m-1} \Big|\\
\leq cr^{2m}\Mc^\tot(  \ddc T,r ).
 \end{multline*}

\end{lemma}

As  an immediate consequence of  Lemmas \ref{L:ddc-difference} and \ref{L:ddc-difference-m-1}, we obtain  the following result.

\begin{corollary}\label{C:ddc-estimate-m-1}  {\rm  (See  \cite[Corollary 13.3]{Nguyen21}). }There is a constant $c>0$  such that  for every
$j$ with $\lowm\leq j\leq \upm,$ and every $m$ with $1\leq m$ and $i\leq j-m,$ and every positive plurisubharmonic  current $T$   in the class  $\widetilde\SH^{p;3,3}(\bfU,\bfW),$
 the  following inequality  holds:
 \begin{multline*} \Big |\int_{\tilde s}^{\tilde r} \big({1\over t^{2(k-p-j)}} -{1\over r^{2(k-p-j)}}\big)2tdt  \int_{\Tube(B,t)}   \ddc (\tau_*T)\wedge \pi^*(\omega^{(i)}\wedge\omega^{j-m-i})\wedge \beta^{k-p-j+m-1}\Big|\\
 \leq cr^{2m}\Nc(  T,r ).
 \end{multline*}

\end{corollary}


Fix  an open neighborhood $\bfW$ of $\partial B$ in $X$ with $\bfW\subset \bfU.$
Recall  the class $\widetilde\SH^{p;3,3}(\bfU,\bfW)$  given in Definition
 \ref{D:sup-bis}.
 
   For $0<r\leq \bfr$ and $0 \leq q\leq k-l$ and $0\leq i\leq  j\leq k-p-q,$ consider  following  global mass indicator
  \begin{equation}\label{e:mass-indicators-log}
 \Lc_{i,j,q}(T,r) :=     \int_0^r {2dt\over t^{2q-1}}    \big( \int_{\Tube(B,t)}(\ddc T)^\hash\wedge
 \pi^*(\omega^{(i)}\wedge \omega^{j-i})\wedge  (\beta+c_1t^2\pi^*\omega)^{k-p-j-1}\big).
 \end{equation}
 Since  $\beta+c_1t^2\pi^*\omega$ is a positive form on $\Tube(B,t),$ it follows that $\Lc_{j,q}(T,r)\geq 0.$
  \begin{lemma}\label{L:Lc-vs-analog} {\rm  (See  \cite[Lemma 13.4]{Nguyen21}). }
  There is a constant  $c>0$ such that for every positive plurisubharmonic   current $T$ on $\bfU$ belonging to the class $\widetilde\SH^{p;3,3}(\bfU,\bfW),$  and  every $0<r\leq \bfr$ and $i,j,q\geq 0$ as  above, we have
  \begin{equation*}
   \big|  \int_0^{\tilde r} {2dt\over t^{2q-1}}   \big( \int_{\Tube(B,t)}\tau_*(\ddc T)\wedge
 \pi^*(\omega^{(i)}\wedge \omega^{j-i})\wedge  (\beta+c_1t^2\pi^*\omega)^{k-p-j-1}\big)-\Lc_{j,q}(T,r) \big|\leq cr\Mc^\tot(\ddc T,r).
   \end{equation*}
  
  \end{lemma}

  \begin{lemma}\label{L:Lc-vs-analog-bis} {\rm  (See  \cite[Lemma 13.5]{Nguyen21}). }
  There is a constant  $c>0$ such that for every positive plurisubharmonic   current $T$ on $\bfU$ belonging to the class $\widetilde\SH^{p;3,3}(\bfU,\bfW),$  and  every $0<r\leq \bfr$ and $i,j,q\geq 0$ as  above, we have
  \begin{multline*}
   \big|  \int_0^{\tilde r} {2dt\over t^{2q-1}}   \big( \int_{\Tube(B,t)}\tau_*(\ddc T)\wedge
 \pi^*(\omega^{(i)}\wedge \omega^{j-i})\wedge  (\beta+c_1t^2\pi^*\omega)^{k-p-j-1}\big)\\
 - \int_0^{\tilde r} {2dt\over t^{2q-1}}   \big( \int_{\Tube(B,t)}\ddc (\tau_*T)\wedge
 \pi^*(\omega^{(i)}\wedge \omega^{j-i})\wedge  (\beta+c_1t^2\pi^*\omega)^{k-p-j-1}\big)
 \big|\leq cr\Nc( T,r).
   \end{multline*}
  
  \end{lemma}

  \begin{lemma}\label{L:Lc-vs-analog-bis'} {\rm  (See  \cite[Lemma 13.6]{Nguyen21}). }
  There is a constant  $c>0$ such that for every positive plurisubharmonic   current $T$ on $\bfU$ belonging to the class $\widetilde\SH^{p;3,3}(\bfU,\bfW),$  and  every $0<r\leq \bfr$ and $i,j,q\geq 0$ as  above, we have
  \begin{multline*}
   \big|  \int_0^{\tilde r} {2dt\over t^{2q-1}}   \big( \int_{\Tube(B,t)}\tau_*(\ddc T)\wedge
 \pi^*(\omega^{(i)}\wedge \omega^{j-i})\wedge  \beta^{k-p-j-1}\big)\\
 - \int_0^{\tilde r} {2dt\over t^{2q-1}}   \big( \int_{\Tube(B,t)}(\ddc T)^\hash\wedge
 \pi^*(\omega^{(i)}\wedge \omega^{j-i})\wedge  \beta^{k-p-j-1}\big)
 \big|\leq cr\Mc^\tot( \ddc T,r).
   \end{multline*}
  
  \end{lemma}
  
 \begin{lemma}\label{L:Lc-vs-analog-bisbis}  {\rm  (See  \cite[Lemma 13.7]{Nguyen21}). } There is a constant  $c>0$ such that for every $j$ with $\lowm^+\leq j\leq \upm^+,$ and for every positive plurisubharmonic   current $T$ in the class $\widetilde\SH^{p;3,3}(\bfU,\bfW),$  
  and  for every $0<r\leq \bfr,$ we have
  \begin{multline*}
   \big|\Lc_{j,q}(T,\tilde r)-  \int_0^{ \tilde r} {2dt\over t^{2q-1}} \big(  \int_{\Tube(B,t)}  
   (\ddc T)^\hash\wedge \pi^*(\omega^{j})\wedge \beta^{k-p-j-1}\big)  \big|\\
   \leq 
   c\sum_{j'=1}^{\upm-j}  \Lc_{i,j+j',q-j'}(T,r)
   +
   c\Nc(T,r).
  \end{multline*}
 \end{lemma}

  \begin{lemma}\label{L:Lc-vs-analog-bis-new} {\rm  (See  \cite[Lemma 13.8]{Nguyen21}). }
  There is a constant  $c>0$ such that for every positive plurisubharmonic   current $T$ on $\bfU$ belonging to the class $\widetilde\SH^{p;3,3}(\bfU,\bfW),$  and  every $0<r\leq \bfr$ and $i,j,q\geq 0$ as  above, we have
  \begin{multline*}
   \big|  \int_0^{\tilde r} \big({1\over t^{2q}} -  {1\over \tilde r^{2q}}\big)2tdt   \big( \int_{\Tube(B,t)}\tau_*(\ddc T)\wedge
 \pi^*(\omega^{(i)}\wedge \omega^{j-i})\wedge  \beta^{k-p-j-1}\big)\\
 - \int_0^{\tilde r}  \big({1\over t^{2q}} -  {1\over \tilde r^{2q}}\big)2tdt \big( \int_{\Tube(B,t)}\ddc (\tau_*T)\wedge
 \pi^*(\omega^{(i)}\wedge \omega^{j-i})\wedge  \beta^{k-p-j-1}\big)
 \big|\leq cr\Nc( T,r).
   \end{multline*}
  
  \end{lemma}

  \begin{lemma}\label{L:Lc-vs-analog-bis'-new} {\rm  (See  \cite[Lemma 13.9]{Nguyen21}). }
  There is a constant  $c>0$ such that for every positive plurisubharmonic   current $T$ on $\bfU$ belonging to the class $\widetilde\SH^{p;3,3}(\bfU,\bfW),$  and  every $0<r\leq \bfr$ and $i,j,q\geq 0$ as  above, we have
  \begin{multline*}
   \big|  \int_0^{\tilde r}   \big({1\over t^{2q}} -  {1\over \tilde r^{2q}}\big)2tdt  \big( \int_{\Tube(B,t)}\tau_*(\ddc T)\wedge
 \pi^*(\omega^{(i)}\wedge \omega^{j-i})\wedge  \beta^{k-p-j-1}\big)\\
 - \int_0^{\tilde r}  \big({1\over t^{2q}} -  {1\over \tilde r^{2q}}\big)2tdt  \big( \int_{\Tube(B,t)}(\ddc T)^\hash\wedge
 \pi^*(\omega^{(i)}\wedge \omega^{j-i})\wedge  \beta^{k-p-j-1}\big)
 \big|\leq cr\Mc^\tot( \ddc T,r).
   \end{multline*}
  
  \end{lemma} 
 \begin{lemma}\label{L:Lc-vs-analog-bisbis-new}  {\rm  (See  \cite[Lemma 13.10]{Nguyen21}). } There is a constant  $c>0$ such that for every $j$ with $\lowm^+\leq j\leq \upm^+,$ and for every positive plurisubharmonic   current $T$ in the class $\widetilde\SH^{p;3,3}(\bfU,\bfW),$  
  and  for every $0<r\leq \bfr,$ we have
  \begin{multline*}
   \big|\Lc_{i,j,q}(T,\tilde r)-  \int_0^{ \tilde r} \big({1\over t^{2q}} -  {1\over \tilde  r^{2q}}\big)2tdt  \big(  \int_{\Tube(B,t)}  
   (\ddc T)^\hash\wedge \pi^*(\omega^{(i)}\wedge \omega^{j-i})\wedge \beta^{k-p-j-1}\big)  \big|\\
   \leq 
   c\sum_{j'=1}^{\upm-j}  \Lc_{i,j+j',q-j'}(T,\tilde r)
   +
   c\Nc(T,\tilde r).
  \end{multline*}
 \end{lemma}

  \begin{lemma}\label{L:Lc-vs-analog-bis-new-second-term}
  There is a constant  $c>0$ such that for every positive plurisubharmonic   current $T$ on $\bfU$ belonging to the class $\widetilde\SH^{p;3,3}(\bfU,\bfW),$  and  every $0<s\leq \bfr$  and $0<\tilde s< r\leq \bfr$ and    $i,j,q\geq 0$ as  above, we have
  \begin{multline*}
   \big|  \big({1\over \tilde s^{2q}} -  {1\over r^{2q}}\big) \int_0^{\tilde s}2tdt   \big( \int_{\Tube(B,t)}\tau_*(\ddc T)\wedge
 \pi^*(\omega^{(i)}\wedge \omega^{j-i})\wedge  \beta^{k-p-j-1}\big)\\
 -   \big({1\over \tilde s^{2q}} -  {1\over r^{2q}}\big)\int_0^{\tilde s}2tdt \big( \int_{\Tube(B,t)}\ddc (\tau_*T)\wedge
 \pi^*(\omega^{(i)}\wedge \omega^{j-i})\wedge  \beta^{k-p-j-1}\big)
 \big|\leq cr\Nc( T,r).
   \end{multline*}
  
  \end{lemma}

  \begin{lemma}\label{L:Lc-vs-analog-bis'-new-second-term}
  There is a constant  $c>0$ such that for every positive plurisubharmonic   current $T$ on $\bfU$ belonging to the class $\widetilde\SH^{p;3,3}(\bfU,\bfW),$  and  every $0<s\leq \bfr$ and $0<\tilde s< r\leq \bfr$ and  $i,j,q\geq 0$ as  above, we have
  \begin{multline*}
   \big|    \big({1\over \tilde s^{2q}} -  {1\over r^{2q}}\big) \int_0^{\tilde s}2tdt  \big( \int_{\Tube(B,t)}\tau_*(\ddc T)\wedge
 \pi^*(\omega^{(i)}\wedge \omega^{j-i})\wedge  \beta^{k-p-j-1}\big)\\
 -   \big({1\over \tilde s^{2q}} -  {1\over r^{2q}}\big)\int_0^{\tilde s}2tdt  \big( \int_{\Tube(B,t)}(\ddc T)^\hash\wedge
 \pi^*(\omega^{(i)}\wedge \omega^{j-i})\wedge  \beta^{k-p-j-1}\big)
 \big|\leq cr\Mc^\tot( \ddc T,r).
   \end{multline*}
  
  \end{lemma}

\begin{lemma}\label{L:Lelong-smooth-forms-psh} {\rm  (See  \cite[Lemma 13.11]{Nguyen21}). }
 Let $T$ be a positive plurisubharmonic   $\Cc^2$-smooth $(p,p)$-form   on $\bfU.$
 Then  for every $\lowm \leq j\leq \upm,$  we have  $\nu_j(T,B,\tau)=0$ if $j\not=l-p$ and $\nu_j(T,B,\tau)\geq 0$ if $j=l-p.$ 
\end{lemma}

Recall  from \cite[Theorem 13.12]{Nguyen21} the finiteness of the mass indicators $\Kc_{j,q} $ and   $\Lc_{j,q} .$  
 \begin{theorem}\label{T:Lc-finite-psh}
   There is  a constant $c_{10}>0$ such that for every positive plurisubharmonic   current $T$ on $\bfU$ belonging to the class $\widetilde\SH^{p;3,3}(\bfU,\bfW),$ 
   and 
   for $0\leq q\leq  k-l$ and   $0\leq   j\leq k-p-q,$ we have 
   In particular,
   $$\Kc_{j,q}(T,\bfr)<c_{10}\quad\text{and}\quad  \Lc_{j,q}(T,\bfr)<c_{10} .
   $$
   \end{theorem}
   
  Recall also from \cite[Corollary 13.13]{Nguyen21} the following finiteness of $\ddc T.$  

\begin{corollary}\label{C:ddc-positive-finite}
There is a constant $c_{11}>0$  such that for every  positive plurisubharmonic   current $T$ in the class $\widetilde\SH^{p;3,3}(\bfU,\bfW),$    and every $q,j$ with $0\leq q\leq  \min(k-l,k-p-1)$ and  $0\leq j\leq k-p-q-1,$ we have  
$$
 \int_0^\bfr{dt\over t^{2q-1}}  \big( \int_{\Tube(B,t)}(\ddc T)^\hash\wedge
 \pi^*\omega^j\wedge  \hat\beta^{k-p-j-1}\big)
<c_{11}.
$$ 
\end{corollary}

\begin{theorem}\label{T:vanishing-Lelong}{\rm  (See  \cite[Theorem 13.14]{Nguyen21}). }
 For every positive plurisubharmonic   current $T$ such that $T=T^+-T^-$ on an open neighborhood of $\overline B$ in $X$  with $T^\pm$  in the class $\SH^{p;3,3}(B),$    and every $\lowm^+\leq j\leq \upm^+,$
 we have  $(\ddc T,B,\omega^{(j)},\tau)=0.$
\end{theorem}
\proof 
 We use
 Lemmas  \ref{L:Lc-vs-analog-bis'} and  \ref{L:Lc-vs-analog-bisbis} instead of \cite[Lemmas 13.6 and 13.7]{Nguyen21}.
\endproof

 \begin{proposition}\label{P:Lc-finite-psh-bis}{\rm  (See  \cite[Proposition 13.15]{Nguyen21}). }
  For  $0<r_1<r_2\leq \bfr,$ there is  a constant $c_{11}>0$ such that for every $q\leq  \min(k-p,k-l)$ and 
    every positive plurisubharmonic   current $T$ in the class $\widetilde\SH^{p;3,3}(\bfU,\bfW),$   we have the following estimate:
 $$|\kappa(T,B,\omega^{(k-p-q)},{r_1\over \lambda},{r_2\over \lambda},\tau)|<c_{11}\sum\limits_{0\leq q'\leq q,\ 0\leq j'\leq \min(\upm,k-p-q')} \Kc_{j',q'}(T,{r_1\over \lambda},{r_2\over \lambda} )\qquad\text{for}\quad  \lambda>1.$$
 \end{proposition}

\subsection{End of the proof of Theorem \ref{T:Nguyen-1}}
We follows  along the same lines as those given in \cite[Subsection 13.3]{Nguyen21} making the obviously necessary changes. There are however  some small  remarks below.
\proof[Proof of assertion (1)]
   We use Lemma \ref{L:Lc-vs-analog-bisbis-new} instead of \cite[Lemma 13.10]{Nguyen21}.
\endproof

 \proof[Proof of assertion (3)]
   We use   Lemmas \ref{L:Lc-vs-analog-bis'-new} and \ref{L:Lc-vs-analog-bisbis-new} instead of \cite[Lemmas 13.9 and  13.10]{Nguyen21}.
\endproof

 
 \section{Basic formulas for the  generalized Lelong numbers}\label{S:Basic-formulas}
 
 In this  section we obtain explicit formulas representing   the generalized  Lelong numbers 
 in terms of  the tangent currents. 
\subsection{A basic formula}
 \begin{theorem}\label{T:Lelong-numbers-vs-T-infty}
 Let $X,$ $V$ be as  above  and  suppose that $(V,\omega)$ is  K\"ahler,  
  and  that  $B$ is   a piecewise $\Cc^2$-smooth open subset  of $V$ and that  there exists a strongly   admissible map  for $B.$
  \begin{enumerate} \item If $k-p\leq l,$ we have 
  $$
  \nu(  T,B,\omega^{(k-p)})= \int_{B}  \ind_B T \wedge \omega^{(k-p)}.$$
  \item For $j\in\N$ such that $0\leq j\leq \upm$ and that $j<k-p,$ let $\omega^{(j)}$ be  a  closed  smooth real $(j,j)$-form on $V_0.$
  Let $T$ be  a   positive plurisubharmonic  $(p,p)$-current  on a neighborhood of $\overline B$ in $X$ such that  $T=T^+-T^-$  for some $T^\pm\in\SH^{p;3,3}( B).$
  By \cite[Tangent Theorem I (Theorem 1.8)]{Nguyen21}, let $T_\infty$ be  a  tangent  current to $T$ along $B.$
 Then   the following   identity holds:
 $$
  \nu(  T,B,\omega^{(j)})={1\over r^2} \int_{\Tube(B,0,r)}  T_\infty\wedge \omega^{(j)}\wedge \beta \wedge \alpha^{k-p-j-1}\quad\text{for}\quad 0<r\leq \bfr .$$
 \end{enumerate}
 \end{theorem}

 \proof 
 The first assertion is straightforward. We come to the proof of the second one.  
 Let $\tau$ be a strongly   admissible map  for $B.$
 Consider  
 $$
 S:=\tau_*T\wedge \omega^{(j)}\wedge \beta\qquad \text{and}\qquad q:=k-p-j-1.
 $$
 Applying  Theorem \ref{T:Lelong-Jensen-eps} to such   $S$ and $q,$   we obtain  
 for  all  $r\in (0,\bfr)$ and   $\epsilon\in  (0,r)$      except for a  countable  set of  values that
\begin{multline*}
 {1\over r^2}{1\over  (r^2+\epsilon^2)^q} \int_{\Tube(B,r)} \tau_* T\wedge \omega^{(j)}\wedge \beta^{q+1}
 ={\Vc_\epsilon(S,r)\over r^2}+  {1\over r^2} \int_{\Tube(B,r)} \tau_* T\wedge \omega^{(j)}\wedge \beta \wedge \alpha^q_\epsilon\\
 +   {1\over r^2}\int_{0}^{r} \big( {1\over (t^2+\epsilon^2)^q}-{1\over (r^2+\epsilon^2)^q}  \big)2tdt\int_{\Tube(B,t)} \ddc (\tau_*) T\wedge \omega^{(j)}\wedge \beta^{q+1} .
\end{multline*} 
\begin{lemma}\label{L:T:Lelong-numbers-vs-T-infty} The following estimates hold:
 ${\Vc_\epsilon(S,r)\over r^2}=O(r)$  and  for $0<\epsilon\ll r,$ $$ {1\over r^2}\int_{0}^{r} \big( {1\over (t^2+\epsilon^2)^q}-{1\over (r^2+\epsilon^2)^q}  \big)2tdt\int_{\Tube(B,t)} \ddc (\tau_*) T\wedge \omega^{(j)}\wedge \beta^{q+1} =O(r).$$
\end{lemma}
\proof[Proof of  Lemma \ref{L:T:Lelong-numbers-vs-T-infty}]
The first estimates follows from an application of  Theorem   \ref{T:vertical-boundary-eps}.
The  second one  can be proved  arguing as in the proof of Lemma \ref{L:Lc-vs-analog-bis'-new}  
and Corollary \ref{C:ddc-positive-finite}.  
\endproof
Applying Lemma \ref{L:T:Lelong-numbers-vs-T-infty} to the last equality,     we   obtain that
 \begin{equation}\label{e:Lelong-alpha_eps-r}
 {r^{2q}\over  (r^2+\epsilon^2)^q}\,\cdot \nu(  T,B,\omega^{(j)},r)
 =   {1\over r^2} \int_{\Tube(B,r)} \tau_* T\wedge \omega^{(j)}\wedge \beta \wedge \alpha^q_\epsilon
 +   O(r)\qquad\text{for}\qquad  0<\epsilon<r.
\end{equation}
Fix an  arbitrarily small $0<\delta \ll 1.$
Applying the above  inequality to  $\delta r$ in place of $r,$ it follows that
\begin{equation}\label{e:Lelong-alpha_eps-r-bis}
 {(\delta r)^{2q}\over  ((\delta r)^2+\epsilon^2)^q}\,\cdot \nu(  T,B,\omega^{(j)},\delta r)
 =    {1\over \delta^2 r^2}\int_{\Tube(B,\delta r)} \tau_* T\wedge \omega^{(j)}\wedge \beta \wedge \alpha^q_\epsilon
 +   O(\delta r)\,\,\text{for}\,\,  0<\epsilon<\delta r.
\end{equation}
Since   $\lim_{r\to 0}\nu(  T,B,\omega^{(j)},r)=\lim_{r\to 0}\nu(  T,B,\omega^{(j)},\delta r)=\nu(  T,B,\omega^{(j)})$  and   ${r^{2q}\over  (r^2+\epsilon^2)^q}=1+O(\delta)$ for  $0<\epsilon<\delta^2r,$ it follows from \eqref{e:Lelong-alpha_eps-r}--\eqref{e:Lelong-alpha_eps-r-bis} that
\begin{equation}\label{e:Lelong-alpha_eps-r,delta}
  \big| \nu(  T,B,\omega^{(j)})- {1\over r^2}\int_{\Tube(B,\delta r,r)}  \tau_*T\wedge \omega^{(j)}\wedge \beta \wedge \alpha^q_\epsilon \big|<\delta,\quad\text{for}\quad 0<r<\bfr,\   0<\epsilon<\delta^2r.
\end{equation}
Fix   $r\in[0,\bfr).$ Then  there is a smooth  positive form  $\gamma=\gamma_{r,\delta}$  such that for 
all $0<\epsilon<\delta^2r,$ 
$$
-\gamma\leq  \alpha_\epsilon^q\leq  \gamma \qquad\text{on}\qquad  \Tube(B, \delta r/ 2,\bfr). 
$$
Since   
$ \lim_{\epsilon\to 0+}  \alpha_\epsilon^q=\alpha^q$ on  $\Tube(B,\delta r,r),$  it follows from Lebesgue's dominated  convergence   that 
$$
\int_{\Tube(B,\delta r,r)}  \tau_*T\wedge \omega^{(j)}\wedge \beta \wedge \alpha^q_\epsilon\to \int_{\Tube(B,\delta r,r)}  \tau_*T\wedge \omega^{(j)}\wedge \beta \wedge \alpha^q\quad\text{as}\quad \epsilon\to 0+.
$$
Hence,  \eqref{e:Lelong-alpha_eps-r,delta} becomes
\begin{equation}\label{e:Lelong-alpha_eps-r,delta-bis}
  \big| \nu(  T,B,\omega^{(j)})- {1\over r^2}\int_{\Tube(B,\delta r,r)}  \tau_*T\wedge \omega^{(j)}\wedge \beta \wedge \alpha^q \big|<\delta,\quad\text{for}\quad 0<r<\bfr.
\end{equation}
For  $0<r<\bfr,$ by letting $\lambda:={\bfr \over r}$ and by acting  the homothetic map $A_\lambda$
and  using  that $A^*_\lambda(\alpha)=\alpha,$  $ A^*_\lambda(\beta)=|\lambda|^2\beta ,$ we obtain that
\begin{equation*}{1\over r^2}\int_{\Tube(B,\delta r,r)}   \tau_*T\wedge \omega^{(j)}\wedge \beta \wedge \alpha^q 
={1\over \bfr^2}\int_{\Tube(B,\delta \bfr,\bfr)}  (A_\lambda)_* (\tau_*T)\wedge \omega^{(j)}\wedge \beta \wedge \alpha^q. 
\end{equation*}
Since   $\pi^* \omega^{(j)},$ $ \beta$ and $\alpha$ are all smooth forms on $\Tube(B,\delta \bfr,\bfr)$
and  $(A_\lambda)_* (\tau_*T)\to T_\infty$  as $\lambda\to\infty.$
the RHS tends to ${1\over \bfr^2}\int_{\Tube(B,\delta \bfr,\bfr)}  T_\infty\wedge \omega^{(j)}\wedge \beta \wedge \alpha^q$ as $r\to 0+.$  Putting this together with \eqref{e:Lelong-alpha_eps-r,delta-bis}, we get that
\begin{equation*} 
  \big| \nu(  T,B,\omega^{(j)})-{1\over \bfr^2} \int_{\Tube(B,\delta \bfr,\bfr)}  T_\infty\wedge \omega^{(j)}\wedge \beta \wedge \alpha^q  \big|<\delta,\quad\text{for}\quad 0<r<\bfr.
\end{equation*}
Letting $\delta\to 0$ in the last line, the theorem follows for $r:=\bfr.$   Clearly,  the above proof also holds for $0<r\leq \bfr.$
\endproof
We conclude this subsection with the following variant of Theorem \ref{T:Lelong-numbers-vs-T-infty-bis}.
\begin{theorem}\label{T:Lelong-numbers-vs-T-infty-bis}\begin{enumerate}
 \item If  instead of the above  assumption on $T,$  we assume  that  $T$ is a   positive pluriharmonic  $(p,p)$-current  on a neighborhood of $\overline B$ in $X$ such that  $T=T^+-T^-$  for some $T^\pm\in\PH^{p;2,2}( B),$  then  the conclusion of  Theorem \ref{T:Lelong-numbers-vs-T-infty} still  holds.   
  \item If  instead of the above  assumption on $T,$  we assume  that  $T$ is a   positive closed  $(p,p)$-current  on a neighborhood of $\overline B$ in $X$ such that  $T=T^+-T^-$  for some $T^\pm\in\CL^{p;1,1}( B),$  then  the conclusion of  Theorem \ref{T:Lelong-numbers-vs-T-infty} still  holds.  
  \item
  Let $X,$ $V$ be as  above.   Assume that there is a Hermitian  metric $\omega$ on $V$  for which
  $\ddc \omega^j=0$ for $\lowm\leq j\leq \upm-1.$   If  instead of the above  assumption on $T,$  we assume  that  $T$ is a   positive  closed  $(p,p)$-current  on a neighborhood of $\overline B$ in $X$ such that  $T=T^+-T^-$  for some $T^\pm\in\CL^{p;2,2}( B),$  then  the conclusion of  Theorem \ref{T:Lelong-numbers-vs-T-infty} still  holds.   
\end{enumerate}
\end{theorem}
\proof
Assertion (1) can be proved similarly as in the proof of Theorem \ref{T:Lelong-numbers-vs-T-infty-bis}
using  Theorem \ref{T:Nguyen-1} (5).

The proof of assertion (2) is  basically similar to  and even  simpler than  the proof of Theorem \ref{T:Lelong-numbers-vs-T-infty-bis}. We make use of   Theorem \ref{T:Nguyen-2} (6). We    utilize Theorem \ref{T:Lelong-Jensen-closed-eps} instead of Theorem  \ref{T:Lelong-Jensen-eps}. We do not need  Lemma \ref{L:T:Lelong-numbers-vs-T-infty}.

The proof of assertion (3) is  basically similar to  the proof of assertion (2). We make use of   Theorem \ref{T:Nguyen-2}.
\endproof

 
\subsection{Tangent currents versus  the generalized Lelong numbers}

Suppose that $V$ is a K\"ahler manifold of dimension $l,$ not necessarily compact, and
let $\omega$ be a Hermitian  metric on $V .$ Let $\E$ be the normal  bundle to $V$ in $X$
 and denote by $\P(\E)$ its projectivization. The complex manifold $\P(\E)$
is of dimension $k - 1.$ Denote by $\pi_{\P(\E)} : \P(\E) \to V$ the canonical projection.
The map $\pi_{\P(\E)}$ deﬁnes a regular ﬁbration over $V$ with $\P^{k-l-1}$ ﬁbers.

Consider a Hermitian metric $h$ on $\E$ and denote by $\omega_{\P(\E)}$ the closed
$(1, 1)$-form on $\P(\E)$ induced by $\ddc \log\|y\|^2_h$  with $y\in \E.$ The restriction of
$\omega_{\P(\E)}$ to each ﬁber of $\P(\E)$ is the Fubini-Study form on this ﬁber. So $\omega_{\P(\E)}$
is strictly positive in the ﬁber direction. It follows that given an open set
$V_0 \Subset V $ there is a constant $c > 0$ large enough such that $c\pi_{\P(E)}^*  (\omega) + \omega_{\P(\E)}$ is
positive on $\pi^{-1} (V_0 ).$  If  $\omega$ is  a K\"ahler form on $V,$ then   the last sum deﬁnes a K\"ahler metric there.

Let $\pi_\FS:\ \E\setminus V\to\P(\E)$ be the canonical projection.

\begin{theorem}\label{T:tangent-currents-vs-Lelong-numbers} Let $X,$ $V$ be as  above  and  suppose that $(V,\omega)$ is  K\"ahler,  
  and  that  $B$ is   a piecewise $\Cc^2$-smooth open subset  of $V$ and that  there exists a strongly   admissible map  for $B.$
  For $0\leq j\leq \upm,$  let $\omega^{(j)}$ be  a  closed  smooth real $(j,j)$-form on $V_0.$
  Let $T$ be  a   positive plurisubharmonic  $(p,p)$-current  on a neighborhood of $\overline B$ in $X$ such that  $T=T^+-T^-$  for some $T^\pm\in\SH^{p;3,3}( B).$
  By \cite[Tangent Theorem I(Theorem 1.8)]{Nguyen21}, let $T_\infty$ be  a  tangent  current to $T$ along $B.$ We will know  from Theorem \ref{T:Nguyen-DS} that $T_\infty\wedge \pi^*\omega^{(j)}$ is   positive pluriharmonic in $\E|_B$ and $V$-conic along $B.$
Hence, $T_\infty\wedge \pi^*\omega^{(j)}=\pi^*_\FS(\T^{(j)}_\infty),$  where $\T^{(j)}_\infty$ is a  positive pluriharmonic  current living on $\P(\E).$
   Then,   the following   identity holds:
 $$
  \nu(  T,B,\omega^{(j)})=\int_{\P(\E|_B)}\T^{(j)}_\infty\wedge \omega_{\P(\E)}^{k-p-j-1}.
 $$
 \end{theorem}
 \proof
Since $V_0\Subset  V,$ we only need to  prove  a local result near a given point $y_0\in V_0.$  We make an analysis  in local coordinates.
We use the coordinates $(z,w)\in\C^{k-l}\times \C^l$  around a neighborhood $U$ of $y_0$  such that $y_0=0$  in  these coordinates. So $y=(z,w).$ 
We may assume that  $U$ has the form $U=U'\times U'',$ where $U'$ (resp. $U'')$ are open neighborhood of $0'$ in $\C^{k-l}$ of  ($0''$ in $\C^l$)
and  $V=\{z=0\}\simeq U''.$ Moreover,  we  may assume  that $U''=(2\D)^l.$  
Consider  the trivial  vector bundle $\pi:\ \E \to  U''$ with  $\E\simeq  \C^{k-l}\times U''.$ 
So the canonical projection  $\pi_\FS|_U:\ (\C^{k-l}\setminus \{0\})\times \D^l\to  \P^{k-l-1}\times \D^l$ reads as  $(z,w)\mapsto  \pi_\FS(z,w):=([z],w),$  where  $ \C^{k-l}\setminus \{0\}\to  \P^{k-l-1}$ is  the canonical  projection mapping $z$ to $[z]:=[z_1:\ldots:z_{k-l}].$
Observe that
\begin{equation}\label{e:FS}
\pi^*_\FS  (\omega_\FS([z],w))=\ddc (\log{\|(z,w)\|^2_h})=\alpha(z,w)\qquad\text{for}\qquad  (z,w)\in(\C^{k-l}\setminus \{0\})\times\D^l.
\end{equation}
We place ourselves on  an open set of $\C^{k-l}$ defined by $z_{k-l}\not=0.$
We   may assume without loss of generality that
\begin{equation}\label{e:max-coordinate} 2|z_{k-l}| > \max\limits_{1\leq j\leq k-l}|z_j|.
\end{equation}
and use the projective coordinates
\begin{equation}\label{e:homogeneous-coordinates}
\zeta_1:={z_1\over z_{k-l}},\ldots, \zeta_{k-l-1}:={z_{k-l-1}\over z_{k-l}},\quad \zeta_{k-l}=z_{k-l}.
\end{equation}
In the coordinates  $\zeta=(\zeta_1,\ldots,\zeta_{k-l})=(\zeta',\zeta_{k-l}),$ the form $\omega_\FS([z],w)$  can be  rewritten as  
\begin{equation}\label{e:FS-zeta} \omega_\FS([z],w)= \ddc \log{ \|(\zeta_1,\ldots,\zeta_{k-l-1},1,w)\|^2_h}= \ddc \log{ \|(\zeta',1,w)\|^2_h}.
\end{equation}
Let $T_\infty$ be  a  tangent current to $T$ along $B.$

 Fix $0<r\leq \bfr.$ By Theorem \ref{T:Lelong-numbers-vs-T-infty}, 
  we  obtain that
 $$
  \nu(  T,B,\omega^{(j)})={1\over r^2} \int_{\Tube(B,0,r)}  T_\infty\wedge \omega^{(j)}\wedge \beta \wedge \alpha^{k-p-j-1}.$$
 Using  a partition of unity on $V_0,$ the RHS can be rewritten as  
$$ {1\over r^2} \int_{  \left\lbrace(\zeta',\zeta_{k-l},w)\in\D^{k-l-1}\times\D\times \D^l:\ \|(\zeta',1,w)\|_h|\zeta_{k-l}|<r \right\rbrace}  T_\infty\wedge \pi^*\omega^{(j)}\wedge \beta \wedge \omega_{\P(\E)}^{k-p-j-1}.$$
 Since $\T_\infty\wedge \pi^*\omega^{(j)}\wedge \omega_{\P(\E)}^{k-p-j-1}$ is a current of full degree $(k-1,k-1)$ in $(d\zeta',d\bar\zeta',dw,d\bar w),$     the  form $\beta=\ddc \varphi$ only contributes
 $\ddc |\zeta_{k-l}|^2$ to the wedge-product $T_\infty\wedge \omega^{(j)}\wedge \beta \wedge \alpha^{k-p-j-1}.$
  Therefore,by Fubini's theorem,  we can write the above RHS as
 $$  \int_{(\zeta',w) }  \T_\infty^{(j)}\wedge \omega_{\P(\E)}^{k-p-j-1}\big({1\over r^2}\int_{\left\lbrace (\zeta',\zeta_{k-l},w)\in\D^{k-l-1}\times\D\times \D^l:\                \|(\zeta',w)\|_h |\zeta_{k-l}|<r\right\rbrace}\|(\zeta',1,w)\|^2_h \ddc |\zeta_{k-l}|^2 \big) .$$
 Since by  the change of variable  $t:=\|(\zeta',1,w)\|\zeta_{k-l},$   the inner integral is  equal to
 $$
 {1\over r^2}\int_{|t|<r}  \ddc |t|^2  
 =1
 ,$$   the  conclusion of the theorem  follows.
 \endproof

We conclude this subsection with the following variant of Theorem \ref{T:tangent-currents-vs-Lelong-numbers}.
\begin{theorem}\label{T:tangent-currents-vs-Lelong-numbers-bis}\begin{enumerate}
 \item Instead of the above  assumption on $T,$  we assume  that  $T$ is a   positive pluriharmonic  $(p,p)$-current  on a neighborhood of $\overline B$ in $X$ such that  $T=T^+-T^-$  for some $T^\pm\in\PH^{p;2,2}( B).$ By \cite[Tangent Theorem I (Theorem 1.8)  (7)]{Nguyen21}, let $T_\infty$ be  a  tangent  current to $T$ along $B.$ We know  from  this theorem that $T_\infty$ is positive pluriharmonic in $\E|_B$ and  $V$-conic  along $B.$
Hence, $T_\infty=\pi^*_\FS(\T_\infty),$  where $\T_\infty$ is a  positive pluriharmonic  current living on $\P(\E).$
 Then  the conclusion of  Theorem \ref{T:tangent-currents-vs-Lelong-numbers} still  holds, more concretely, 
 \begin{equation}\label{e:T:tangent-currents-vs-Lelong-numbers-bis} 
  \nu(  T,B,\omega^{(j)})=\int_{\P(\E|_B)}\T_\infty\wedge \pi^*\omega^{(j)}\wedge \omega_{\P(\E)}^{k-p-j-1}.
 \end{equation}
  \item Instead of the above  assumption on $T,$  we assume  that  $T$ is a   positive closed  $(p,p)$-current  on a neighborhood of $\overline B$ in $X$ such that  $T=T^+-T^-$  for some $T^\pm\in\CL^{p;1,1}( B).$ By \cite[Tangent Theorem II (Theorem 1.11)]{Nguyen21}, let $T_\infty$ be  a  tangent  current to $T$ along $B.$ We know  from  this theorem that $T_\infty$ is  positive closed in $\E|_B$ and  $V$-conic along $B.$
Hence, $T_\infty=\pi^*_\FS(\T_\infty),$  where $\T_\infty$ is a  positive closed  current living on $\P(\E).$
  Then  identity \eqref{e:T:tangent-currents-vs-Lelong-numbers-bis}  still  holds.  
  \item
  Let $X,$ $V$ be as  above.   Assume that there is a Hermitian  metric $\omega$ on $V$  for which
  $\ddc \omega^j=0$ for $\lowm\leq j\leq \upm-1.$   If  instead of the above  assumption on $T,$  we assume  that  $T$ is a   positive  closed  $(p,p)$-current  on a neighborhood of $\overline B$ in $X$ such that  $T=T^+-T^-$  for some $T^\pm\in\CL^{p;2,2}( B).$  Then  the conclusion of  assertion (2) still  holds.   
\end{enumerate}
\end{theorem}
 We leave the proof to the interested reader.
 
\section{Horizontal dimension and  a Siu's upper-semicontinuity type theorem}
\label{S:H-dim-and-Siu}
 
 In the first three  subsections we prove Theorem \ref{T:Nguyen-DS}. 
 The last  subsection is devoted to the proof of Theorem \ref{T:Nguyen-Siu}. 
\subsection{Pluriharmonicity and $V$-conicity}
We prove  assertion (1) of Theorem \ref{T:Nguyen-DS}. Writing $\omega^{(j)}=(\omega^{(j)}+c\omega^j)-c\omega^j$ for a large enough constant $c>0,$  we may assume  without loss of generality that $\omega^{(j)}$ is  strictly positive.
 By  \cite[Theorem   15.5]{Nguyen21},
$T_\infty$ is   the tangent current  to $T$ along $B$ associated   to a sequence
 $(\lambda_n ) \subset \C^*$ converging to $\infty.$   
Fix $r_1,\ r_2\in  (0,\bfr)$ with $r_1<r_2.$  Let $\lambda\in\R$ with $\lambda \geq 1.$  

For every $j$ with  $\lowm\leq j\leq \upm,$  applying Theorem  \ref{T:Lelong-Jensen}  to  $(A_{\lambda_n})_*(\tau_*T)\wedge \pi^*(\omega^{(j)})$ yields that
\begin{multline*}
  \nu(T,B,\omega^{(j)},{r_2\over |\lambda_n|},\tau)- \nu(T,B,\omega^{(j)},{r_1\over |\lambda_n|},\tau)=  \Vc\big((A_{\lambda_n})_*(\tau_*T)\wedge \pi^*(\omega^{(j)}),r_1,r_2\big)\\ +\int_{\Tube(B,r_1,r_2)}(A_{\lambda_n})_*(\tau_*T)\wedge \pi^*(\omega^{(j)})\wedge\alpha^{k-p-j}\\ 
   +  \int_{r_1}^{r_2} \big( {1\over t^{2(k-p-j)}}-{1\over r_2^{2(k-p-j)}}  \big)2tdt\int_{\Tube(B,t)}  \ddc (A_{\lambda_n})_*(\tau_*T)\wedge \pi^*(\omega^{(j)})\wedge \beta^{k-p-j-1} \\
  +  \big( {1\over r_1^{2(k-p-j)}}-{1\over r_2^{2(k-p-j)}}  \big) \int_{0}^{r_1}2tdt\int_{z\in \Tube(B,t)} \ddc (A_{\lambda_n})_*(\tau_*T)\wedge \pi^*(\omega^{(j)})\wedge \beta^{k-p-j-1}.
\end{multline*}
When  $n$ tends to infinity, the LHS  tends to $0$   since by
Theorem \ref{T:Nguyen-1} (1), $\lim_{n\to\infty}\nu(T,B,\omega^{(j)},{r\over |\lambda_n|},\tau)=\nu(T,B,\omega^{(j)},\tau)\in\R$ for $0<r\leq\bfr.$
By  Theorem \ref{T:vertical-boundary-terms},  $\Vc\big((A_{\lambda_n})_*(\tau_*T)\wedge \pi^*(\omega^{(j)}),r_1,r_2\big)\to 0$ as $n\to\infty.$
Therefore, we obtain that 
\begin{equation}\label{e:sum_integrals_equal_zero}\begin{split}
  0&=   \int_{\Tube(B,r_1,r_2)}T_\infty\wedge \pi^*(\omega^{(j)})\wedge\alpha^{k-p-j}\\
  & +  \int_{r_1}^{r_2} \big( {1\over t^{2(k-p-j)}}-{1\over r_2^{2(k-p-j)}}  \big)2tdt\int_{\Tube(B,t)}  \ddc T_\infty\wedge \pi^*(\omega^{(j)})\wedge \beta^{k-p-j-1} \\
  &+  \big( {1\over r_1^{2(k-p-j)}}-{1\over r_2^{2(k-p-j)}}  \big) \int_{0}^{r_1}2tdt\int_{z\in \Tube(B,t)} \ddc T_\infty\wedge \pi^*(\omega^{(j)})\wedge \beta^{k-p-j-1}.
  \end{split}
\end{equation}
Next,   
consider  a small neighborhood $V(x_0)$ of  an arbitrary  point $x_0\in  \Tube(B, r_0),$  where in a local chart $V(x_0)\simeq \D^l$ and  $\E|_{V(x_0)}\simeq \C^{k-l}\times \D^l.$
For $x\in \E|_{V(x_0)},$ write $x=(z,w).$  Since   $\upm=\min(l,k-p)$ and $T_\infty$ is  of bidegree $(p,p)$   we  see that  $T_\infty\wedge \pi^*\omega^{(\upm)}$ is  of full bidegree $(l,l)$  in $dw,$ $d\bar w.$
Consequently, we infer from \eqref{e:hat-alpha'} that
\begin{eqnarray*}
T_\infty\wedge \pi^*(\omega^{(\upm)})\wedge\alpha^{k-p-\upm}
&=&T_\infty\wedge \pi^*(\omega^{(\upm)})\wedge(\hat\alpha')^{k-p-\upm},\\
\ddc T_\infty\wedge \pi^*\omega^{(\upm)}\wedge\beta^{k-p-\upm}&=&\ddc T_\infty\wedge \pi^*\omega^{(\upm)}\wedge\hat\beta^{k-p-\upm}.
\end{eqnarray*}
This,  combined  with \eqref{e:sum_integrals_equal_zero} for $j:=\upm,$  implies  that 
\begin{eqnarray*}
  0&&= \int_{\Tube(B,r_1,r_2)}T_\infty\wedge \pi^*(\omega^{(\upm)})\wedge(\hat\alpha')^{k-p-\upm}\\
  &+&  \int_{r_1}^{r_2} \big( {1\over t^{2(k-p-\upm)}}-{1\over r_2^{2(k-p-\upm)}}  \big)2tdt\int_{\Tube(B,t)}  \ddc T_\infty\wedge (\pi^*\omega^{(\upm)})\wedge \hat\beta^{(k-p-\upm)-1} \\
 &+&  \big( {1\over r_1^{2(k-p-\upm)}}-{1\over r_2^{2(k-p-\upm)}}  \big) \int_{0}^{r_1}2tdt\int_{z\in \Tube(B,t)} 
 \ddc T_\infty\wedge (\pi^*\omega^{(\upm)})\wedge \hat\beta^{(k-p-\upm)-1}. 
\end{eqnarray*}
Since we know  by \cite[Tangent Theorem I (Theorem 1.8)]{Nguyen21} that $T_\infty$ is  positive  plurisubharmonic,  both $T_\infty$ and $\ddc T_\infty$  are   positive  currents.
Moreover,  $\omega,$ $\hat\alpha',$  $\hat\beta$ are positive forms.  Consequently,  all integrals of  the RHS of the last line  are $\geq 0.$
On the ther hand, their sum  is  equal to $0.$  So  all integrals  are  $0,$ that is,
\begin{equation*}
  \int_{z\in \Tube(B,r_2)} \ddc T_\infty\wedge \pi^*(\omega^{(\upm)})\wedge \hat\beta^{k-p-\upm-1}=0\quad\text{and}\quad  \int_{\Tube(B,r_1,r_2)}T_\infty\wedge \pi^*(\omega^{(\upm)})\wedge\alpha^{k-p-\upm}=0.
\end{equation*}
Note that  $\hat\beta$  and $\pi^*\omega$ are  smooth  strictly positive $(1,1)$ forms on $\Tube(B,\bfr),$ and that for every   smooth   positive $(1,1)$ form $H$ on $\Tube(B,\bfr),$ we can find a constant $c>0$ such that $H\leq c (\hat\beta+\pi^*\omega)$ on $\Tube(B,\bfr).$
Since $0<r_1<r_2\leq\bfr$ are arbitrarily chosen, 
 we infer that  the following  equality holds for all $j$ with $\upm\leq j\leq k$:
 \begin{equation}\label{e:fact-partial-plurihar}
  \ddc T_\infty\wedge\pi^*(\omega^{(j)})=0\quad\text{on}\quad\Tube(B,\bfr)\quad\text{and}\quad T_\infty\wedge \pi^*(\omega^{(j)})\wedge\alpha^{k-p-j}=0 \quad\text{on}\quad\Tube(B,\bfr)\setminus B. 
 \end{equation}
Suppose  that \eqref{e:fact-partial-plurihar} holds for all $j$  with $j_0<j\leq\upm,$ where $j_0$ is a given integer
with $\lowm\leq j_0<\upm.$  We need to  prove \eqref{e:fact-partial-plurihar}  for  $j=j_0.$ 

 Using  \eqref{e:fact-partial-plurihar} for all $j$  with $j_0<j\leq k,$ we  infer from \eqref{e:hat-alpha'}  that
\begin{eqnarray*}
T_\infty\wedge \pi^*(\omega^{(j_0)})\wedge\alpha^{k-p-j_0}
&=&T_\infty\wedge \pi^*(\omega^{(j_0)})\wedge(\hat\alpha')^{k-p-j_0},\\
\ddc T_\infty\wedge \pi^*\omega^{(j_0)}\wedge\beta^{k-p-j_0}&=&\ddc T_\infty\wedge \pi^*\omega^{(j_0)}\wedge\hat\beta^{k-p-j_0}.
\end{eqnarray*}
This,  combined  with \eqref{e:sum_integrals_equal_zero} for $j:=j_0,$  implies  that 
\begin{eqnarray*}
  0&&= \int_{\Tube(B,r_1,r_2)}T_\infty\wedge \pi^*(\omega^{(j_0)})\wedge(\hat\alpha')^{k-p-{j_0}}\\
  &+&  \int_{r_1}^{r_2} \big( {1\over t^{2(k-p-j_0)}}-{1\over r_2^{2(k-p-j_0)}}  \big)2tdt\int_{\Tube(B,t)}  \ddc T_\infty\wedge (\pi^*\omega^{(j_0)})\wedge \hat\beta^{(k-p-j_0)-1} \\
 &+&  \big( {1\over r_1^{2(k-p-j_0)}}-{1\over r_2^{2(k-p-j_0)}}  \big) \int_{0}^{r_1}2tdt\int_{z\in \Tube(B,t)} 
 \ddc T_\infty\wedge (\pi^*\omega^{(j_0)})\wedge \hat\beta^{(k-p-j_0)-1}. 
\end{eqnarray*}
We repeat the above argument   using that   both $T_\infty$ and $\ddc T_\infty$  are   positive  currents and that $\pi^*\omega,$  $\hat\alpha',$  $\hat\beta$ are positive forms.  Consequently,    all integrals on the RHS  are  $0.$
Therefore, \eqref{e:fact-partial-plurihar} holds for  $j=j_0.$
Hence,  the proof of \eqref{e:fact-partial-plurihar} is  completed.
In particular,  $\ddc T_\infty\wedge\pi^*(\omega^{(\lowm)})=0$ on $\Tube(B,\bfr).$
Since we will prove shortly below that  $T_\infty\wedge \pi^*(\omega^{(\lowm)})$ is $V$-conic,
it follows that   $\ddc T_\infty\wedge\pi^*(\omega^{(j)})=0$ on $\pi^{-1}(B)\subset \E.$

  Recall from  \eqref{e:hat-alpha'} that  $\hat\alpha'=\alpha+c_1\pi^*\omega,$ and from
  \eqref{e:hat-alpha'-vs-alpha_ver} that  $
\hat\alpha'\geq  c^{-1}_1\alpha_\ver\geq 0.$  Moreover, $T_\infty$ is   a positive current.
 Therefore, we infer from the  second identity of \eqref{e:fact-partial-plurihar} that
  \begin{equation}\label{e:T_infty-proj_iden-bis} T_\infty\wedge (\pi^*(\omega^{(j)})\wedge \alpha_\ver^{k-p-j}=0\quad\text{on}\quad \Tube(B,\bfr)\setminus B\quad\text{for}\quad \lowm\leq j\leq k.
\end{equation}
Consider  the  positive  pluriharmonic  current $\Theta:=T_\infty\wedge\pi^*(\omega^{(j)}).$
We need  to show  that $\Theta$ is  $V$-conic. 
   We argue as  in the proof of assertion (3) of \cite[Theorem 16.3]{Nguyen21}.
\endproof
\subsection{Horizontal dimension}
We prove  assertions (2), (3) and (4)  of Theorem \ref{T:Nguyen-DS}. 

\proof[Proof of  assertion (2)  of Theorem \ref{T:Nguyen-DS}]
Let $q>\hbar$ and let $\omega^{(q)}$ be a    closed  smooth $(q,q)$-for on $V.$ 
  By Theorem \ref{T:tangent-currents-vs-Lelong-numbers}, the  pluriharmonic current $\T^{(q)}_\infty,$
  defined  by $\pi_\FS^*( \T^{(q)}_\infty)=T_\infty  \wedge \pi^*\omega^{(q)},$ satisfies
 $$
  \nu(  T,B,\omega^{(q)})=\int_{\P(\E|_B)}\T^{(q)}_\infty\wedge \omega_{\P(\E)}^{k-p-q-1}.
 $$
 Fix a large constant $c>0$ such that $c\omega^q+ \omega^{(q)}$ is strictly positive.
 We infer from the assumption  $q>\hbar$  that $T_\infty \wedge \pi^*(c\omega^q+ \omega^{(q)})=0$
 and  $T_\infty \wedge \pi^*(c\omega^q )=0.$ Therefore,    $T_\infty  \wedge \pi^*\omega^{(q)}=0,$
 hence $\T^{(q)}_\infty=0,$ which in turn implies that $
  \nu(  T,B,\omega^{(q)})=0.$
  \endproof

\proof[Proof of  assertion (3)  of Theorem \ref{T:Nguyen-DS}]
 By  assertion (2), we  obtain that $
  \nu_q(  T,B,\omega)=0$ for all $\hbar <q\leq\upm.$ Therefore, 
 we
 only need to show that  if $j\in[\lowm,\upm]$ is an integer  satisfying   $\nu_q(  T,B,\omega)=0$ for all $j <q\leq\upm,$  
 and  if $j>\lowm$ then $\nu_j(T,B,\omega)\not=0,$ then $j=\hbar.$

 There is a constant $c>0$ such that  $c\pi^*\omega+\omega_{\P(\E)}$ is a  K\"ahler form on  $\P(\E).$
 By Theorem \ref{T:tangent-currents-vs-Lelong-numbers}, 
 the  pluriharmonic current $\T^{(j)}_\infty$ on $\E|_B \setminus B$
  defined  by $\pi_\FS^*( \T^{(j)}_\infty)=T_\infty  \wedge \pi^*\omega^{(j)},$ satisfies
 $$
  \nu_j(  T,B,\omega)=\int_{\P(\E|_B)}\T^{(j)}_\infty\wedge \omega_{\P(\E)}^{k-p-j-1}.
 $$ 
 Using that  $\nu_q(  T,B,\omega)=0$ for all $j <q\leq\upm,$ we infer that 
 $$  \nu_j(  T,B,\omega)=\int_{\P(\E|_B)}\T^{(j)}_\infty\wedge (\omega_{\P(\E)}+c\pi^*\omega)^{k-p-j-1}.
 $$
 Since $T_\infty$ is a positive current, and $  \pi^*\omega^{(j)}$ and $\omega_{\P(\E)}+c\pi^*\omega$ are   positive smooth forms, we infer that   $\nu_j(  T,B,\omega)\geq 0.$ By hypothesis,  $\nu_j(  T,B,\omega)\not= 0.$
 So $\nu_j(  T,B,\omega)>0 $ and $T_\infty\wedge \pi^*\omega^{j}\not=0.$ Hence, $j=\hbar.$
\endproof

\proof[Proof of  assertion (4)  of Theorem \ref{T:Nguyen-DS}]
Let  $\omega^{(\hbar)}$ be a   strongly  positive   closed smooth $(\hbar,\hbar)$-form on $V.$ 
By Theorem \ref{T:tangent-currents-vs-Lelong-numbers}, the  pluriharmonic current $\T^{(\hbar)}_\infty,$  on $\E|_B \setminus B$
  defined  by $\pi_\FS^*( \T^{(\hbar)}_\infty)=T_\infty  \wedge \pi^*\omega^{(\hbar)},$ satisfies
 $$
  \nu(  T,B,\omega^{\hbar)})=\int_{\P(\E|_B)}\T^{(\hbar)}_\infty\wedge \omega_{\P(\E)}^{k-p-\hbar-1}.
 $$ 
 Using that  $ T_\infty\wedge \pi^*(\omega^{(\hbar)}\wedge \omega^{q-\hbar})=0$ for all $\hbar<q\leq\upm,$ we infer that $\T^{(\hbar)}_\infty\wedge \omega^{q-h}=0,$ and hence,
 $$  \nu(  T,B,\omega^{(\hbar)})=\int_{\P(\E|_B)}\T^{(\hbar)}_\infty\wedge (\omega_{\P(\E)}+c\pi^*\omega)^{k-p-\hbar-1}.
 $$
 Since $T_\infty$ is a positive current, and $  \pi^*\omega^{(\hbar)}$ and $\omega_{\P(\E)}+c\pi^*\omega$ are   positive smooth forms, we infer that   $\nu(  T,B,\omega^{(\hbar)})\geq 0.$

To complete the proof, let   $\omega^{(\hbar)}$ be a   strictly  positive and $T_\infty\not=0.$ 
there is a constant $c'>0$ such that $\omega^{(\hbar)}\geq c'\omega^\hbar,$ and by  the above  argument, we infer that   the current $\T^{\langle \hbar\rangle}_\infty$  on $\E|_B \setminus B$ defined by   $\pi_\FS^*( \T^{\langle \hbar\rangle}_\infty)=T_\infty  \wedge \pi^*\omega^{\hbar},$ satisfies 
$$\nu_\hbar(T,B,\omega)=\int_{\P(\E|_B)}\T^{\langle \hbar\rangle}_\infty\wedge (\omega_{\P(\E)}+c\pi^*\omega)^{k-p-\hbar-1}.
 $$
 Since $T_\infty\wedge \pi^*\omega^{\hbar}$ is a nonzero positive  current and $\omega_{\P(\E)}+c\pi^*\omega$ is a  K\"ahler form, we infer that $\T^{\langle \hbar\rangle}_\infty$ is a nonzero positive  current, and  hence   $\nu_\hbar(  T,B,\omega)> 0.$ This implies that $\nu(  T,B,\omega^{(\hbar)})\geq  c'\nu_\hbar(  T,B,\omega)> 0.$
\endproof

\subsection{A Siu's upper-semicontinuity type theorem}
\label{SS:Siu}
 
We adapt the method  of the proof of Theorem \ref{T:Dinh-Sibony-Siu} in Dinh-Sibony's article  \cite{DinhSibony18}   in order to prove   Theorem \ref{T:Nguyen-Siu}. 

For a positive closed current $T,$ let $T_\infty$ be a tangent current to $T$ along $B,$ and 
$\T$ be  the associated current  on $\E|_B \setminus B$ defined by $T_\infty=\pi^*_\FS(\T).$
 
 If $T$ has positive mass on $V ,$  then the horizontal dimension of $T$
along $V$ is maximal, i.e., equal to $k- p.$ The theorem is trivial. Suppose  now
that $T$ gives no mass to $V .$ We deduce that the mass of $T_n$ on $V$ tends to 0.
So removing from $T_n$ its restriction to $V$ allows us to suppose  that $T_n$ has no
mass on $V$ for every $n.$
Denote by $\widehat T$ and $\widehat T_n$ the strict transforms of $T$ and $T_n$ with respect to the
the  blow-up $\sigma :\ \widehat X\to X $ along $V .$  We identify $\widehat V$ with $\P(\E). $ By \cite[Lemma 4.7]{DinhSibony18} and the last
assertion of \cite[Lemma 4.10]{DinhSibony18}, we have
\begin{equation*}
\{\T\} = \{\widehat\T\}= \{\widehat T\}|_{\P(E)}\quad\text{and}\quad 
\{\T_n\} = \{\widehat\T_n\}= \{\widehat T_n\}|_{\widehat \P(E)}.
\end{equation*}
Passing to  a subsequence if necessary, we may suppose  that $\widehat T_n$ converges to a current $\widehat T' .$ 
Write $\widehat T' = \widehat T + \widehat R,$
 where $\widehat R$
is the restriction of
$\widehat T'$ to $V .$ 
So there is a   positive closed $(p-1,p-1)$-current $R$ on $\P(E)$ such that  $\widehat R=\pi^*R\wedge [\widehat V].$
Let $\{R\}$ be the class of $R$ in $H^{p-1,p-1}(\P(\E)).$
Let  ${\widehat R}'$ denote the class of $\widehat R$ in $H^{2p}_{\widehat V} (\widehat X, \C).$
 Therefore, we 
obtain
\begin{equation}\label{e:class_T_n-T}\lim\limits_{n\to\infty} \{\T_n\} - \{\T\} = {\widehat R}'|_{\P(\E)}=-\{R\} \smile h_{\P(\E)},
\end{equation}
where the last equality follows  from \cite[Lemma 3.17]{DinhSibony18}.

Recall from the proof of  \cite[Theorem 4.11]{DinhSibony18}
  that the horizontal dimension of $\widehat R$ is at most equal to $\hbar.$
Consequently, the horizontal  of the limit   $\lim\limits_{n\to\infty} \{\T_n\}$ does not exceed $\hbar.$

 To prove  assertion (1),  fix $\lowm\leq j\leq\upm$ with $j>\hbar$  and  let   $\omega^{(j)}$
 be a  closed  smooth $(j,j)$-form on $V.$
 Fix  an arbitrary $n\in\N.$  By the hypothesis on $T_n$ and on $U,W,$ and multipling $T$ by a  small psoitive constant if necessary,  we apply Theorem \ref{T:approximation} in order to find  two sequences $(T^\pm_{n,N})_{N\in\N}\in \widetilde{\CL}^{p;1,1}(U,W)$
 such that $T^\pm_{n,N}\to T^\pm_n$ as $N\to\infty$ and  $T_n=T^+_n-T^-_n.$
 
By Theorem \ref{T:tangent-currents-vs-Lelong-numbers}, we have that
$$
  \nu(  T^\pm_n,B,\omega^{(j)})=\int_{\P(\E|_B)}\T^\pm_n\wedge \pi^*\omega^{(j)}\wedge h_{\P(\E)}^{k-p-j-1}.
 $$
Since  $T_n=T^+_n-T^-_n,$ it follows  that $\T_n=\T^+_n-\T^-_n,$ and 
$$
  \nu(  T_n,B,\omega^{(j)})=\int_{\P(\E|_B)}\T_n\wedge \pi^*\omega^{(j)}\wedge h_{\P(\E)}^{k-p-j-1}.
 $$
As $\supp \T_n\cap  \E_{\partial B}=\varnothing,$ we rewrite   this as  follows:
 $$
  \nu(  T_n,B,\omega^{(j)})=\int_{\P(\E|_B)}\{\T_n\}\smile\{ \pi^*\omega^{(j)}\}\wedge h_{\P(\E)}^{k-p-j-1}.
 $$
Since  the horizontal  of the limit   $\lim\limits_{n\to\infty} \{\T_n\}$ does not exceed $\hbar,$
 we infer that $\lim_{n\to\infty}\{\T_n\}\smile\{ \pi^*\omega^{(j)}\}=0,$ and hence
 $\lim_{n\to\infty} \nu(T_n,B,\omega^{(j)})= 0.$

  To prove assertion (2),  let $\omega^{(\hbar)}$ be a  strongly positive   closed  smooth $(\hbar,\hbar)$-for on $V.$ 
  There is a constant $c>0$ such that  $c\pi^*\omega+\omega_{\P(\E)}$ is a  K\"ahler form on  $\P(\E).$
  By Theorem \ref{T:tangent-currents-vs-Lelong-numbers}, we have that
 $$
  \nu(  T_n,B,\omega^{(j)})=\int_{\P(\E|_B)}\{\T_n\}\smile\{ \pi^*\omega^{(j)}\}\smile h_{\P(\E)}^{k-p-j-1}.
 $$
 This, combined with  assertion (1), implies that
 \begin{eqnarray*}
 \liminf_{n\to\infty}  \nu(  T_n,B,\omega^{(j)})&=&\liminf_{n\to\infty} \int_{\P(\E|_B)}\{\T_n\}\smile\{ \pi^*\omega^{(j)}\}\wedge \{\omega_{\P(\E)}\}^{k-p-j-1}\\
 &=&\liminf_{n\to\infty} \int_{\P(\E|_B)}\{\T_n\}\smile\{ \pi^*\omega^{(j)}\}\smile (h_{\P(\E)}+c\{\pi^*\omega\})^{k-p-j-1}.
 \end{eqnarray*}
 Since $\T_n$ is a positive current, and $  \pi^*\omega^{(j)}$ and $\omega_{\P(\E)}+c\pi^*\omega$ are   positive smooth forms, we infer that  
  $ \liminf_{n\to\infty} \nu(T_n,B,\omega^{(\hbar)})\geq 0.$

 By Theorem \ref{T:tangent-currents-vs-Lelong-numbers}, we have that
  \begin{eqnarray*}
  \nu(T,B,\omega^{(\hbar)}) - \limsup_{n\to\infty} \nu(T_n,B,\omega^{(\hbar)})
  &=& \liminf\limits_{n\to\infty}\int_{\P(\E|_B)}(\{T\}-\{\T_n\})\smile\{ \pi^*\omega^{(\hbar)}\}\smile h_{\P(\E)}^{k-p-\hbar-1}\\
  &=&
 \int_{\P(\E|_B)} \{R\} \smile h_{\P(\E)}\smile\{ \pi^*\omega^{(\hbar)}\}\smile h_{\P(\E)}^{k-p-\hbar-1},
\end{eqnarray*}  
where the last equality holds by \eqref{e:class_T_n-T}.

Since the horizontal dimension  of $\{R\}$ does not exceed $\hbar,$ we see that the last integral is  equal to
$$
\int_{\P(\E|_B)} \{R\} \smile\{ \pi^*\omega^{(\hbar)}\}\smile(c\{\pi^*\omega\}+ h_{\P(\E)})^{k-p-\hbar}.
$$
The last expression is $\geq 0$ because $R$ is a positive current, and $  \pi^*\omega^{(\hbar)}$ and $\omega_{\P(\E)}+c\pi^*\omega$ are   positive smooth forms.
This  completes the proof that $\nu(T,B,\omega^{(\hbar)}) \geq \limsup_{n\to\infty} \nu(T_n,B,\omega^{(\hbar)}).$              \hfill $\square$

 
 \section{Dinh-Sibony classes versus generalized Lelong numbers}\label{S:Dinh-Sibony-vs-Lelong}
  
This section is devoted to the proof of Theorems   \ref{T:Nguyen-Dinh-Sibony} and \ref{T:Nguyen-Siu-AB}.
Consider a Hermitian metric $\bar h$ on $\overline \E:=\P(\E\oplus \C),$ and denote by $\omega_{\overline \E}$ the closed
$(1, 1)$-form on $\overline \E$ induced by $\ddc \log\|y\|^2_h$  with $y\in \overline\E.$ The restriction of
$\omega_{\overline\E}$ to each ﬁber of $\overline\E$ is the Fubini-Study form on this ﬁber. So $\omega_{\overline\E}$
is strictly positive in the ﬁber direction. 

Recall $\pi_\FS:\ \E\to\P(\E)$  the canonical projection. Consider $\iota:\ \E\hookrightarrow \P(\E)$ the canonical injection.
  \begin{lemma}
   \label{L:PE-to-bar-E}
   For every  smooth volume form $\Theta$ on  $\P(E),$  we have
   $$
   \int_{\P(\E)}\Theta=\int_{\overline\E} \pi_\FS^*\Theta\wedge \omega_{\overline \E}. 
   $$
  \end{lemma}
\proof
It  suffices to prove  the  theorem for the case where $V$ is  a  single point, that is, $\E=\C^k.$ 
In this  case, the lemma follows by a straightforward computation.
\endproof
  \proof[End of the proof of Theorem  \ref{T:Nguyen-Dinh-Sibony}]
  We only give  the proof of assertion (2) for the sake of clarity. Indeed although the proof of assertion (1) is  quite similar to that of assertion (2), it is  somehow technically  a bit more complicated.
  
  Let $T_\infty\in \PH^p(\overline\E)$ be a tangent  current to $T$ along $V.$
  By \cite[Proposition 3.10]{DinhSibony18} (the proof therein goes through in  the  present context of positive pluriharmonic  currents), we obtain
  $$
  T_\infty=\pi^*_\FS(\T_\infty)+\ind_V T.
  $$
  Let $T_0\in \PH^{l-k+p}(V)$ be such that $ \ind_V T=(\iota_V)_* T_0,$  where $\iota_V:\ V\hookrightarrow \overline\E$ is  the canonical  injection. When $l<k-p,$ we have $T_0=0$ and $ \ind_V T=0.$
By Theorem \ref{T:Dinh-Sibony},  we have  $\{T_\infty\}=\bfc^\DS(T,V)$  is  the  total  tangent class of $T$  along $V$, or equivalently,   Dinh-Sibony (total) cohomology class of $T$ along $V.$
Let $-h_{\overline \E}$ denote the  tautological class of the bundle $\pi_0:\  \overline \E\to  V.$ Then  we have the following   decomposition
$$
\{ \T_\infty\}=\sum_{i=\lowm}^\upm \pi_0^*(\bfc^\DS_i(T,V))\smile  h_{\overline \E}^{i-l+p},
$$
where  $\bfc^\DS_j(T,V)$  is a class in  $H^{2l-2j}(V,\C).$ Moreover, for $\lowm\leq j\leq \upm,$
$$\bfc^\DS_{k-p}(T,V)=\{T_0\},\qquad\text{and for $j<k-l,$}\quad  \bfc^\DS_{j}(T,V)=\bfc^\DS_{j}(\T_\infty,V). $$
This decomposition is, in fact, unique.

On the other hand,  since $\omega_{\P(\E)}$ (resp. $\omega_{\overline\E}$) is  the  curvature form of the  line bundle $\Oc_{\P(\E)}$ (resp. the line bundle  $\Oc_{\overline \E}$), we see  that 
$$ \iota_*( \pi^*_\FS \{\omega_{\P(\E)}\})=h_{\overline \E}=\{\omega_{\overline\E}\}.$$

By Theorem \ref{T:tangent-currents-vs-Lelong-numbers}, we   obtain for $\lowm\leq j\leq\upm$  that  
\begin{equation*}
  \nu(  T,V,\omega^{(j)})=
  \begin{cases}
   \int_V T_0\wedge \omega^{(j)}, &\text{if}\ j=k-p;\\ 
   \int_{\P(\E)}\T_\infty\wedge \pi^*\omega^{(j)}\wedge \omega_{\P(\E)}^{k-p-j-1}, &\text{otherwise.}
  \end{cases}
  \end{equation*}
  Putting together the  above   equalities and  applying  Lemma   \ref{L:PE-to-bar-E}, it follows that 
   \begin{equation*}
\nu(  T,V,\omega^{(j)}) = \big(\sum_{i=\lowm}^\upm \pi_0^*(\bfc^\DS_i(T,V))\smile  h_{\overline \E}^{i-l+p} \big) \smile\pi^*\{\omega^{(j)}\}   \smile  h_{\overline \E}^{k-p-j-1}  \smile  h_{\overline \E}.
 \end{equation*}
By  a consideration of bidegree, we  see that on $V$:
$$
\bfc^\DS_i(T,V)\smile\{\omega^{(j)}\} =0 \qquad\text{for}\qquad  i<j.
$$
 Therefore,  we infer that  
 \begin{equation*}
  \nu(  T,B,\omega^{(j)}) 
 = \big(\sum_{i=j}^\upm \pi_0^*(\bfc^\DS_i(T,V))\smile  h_{\overline \E}^{i-l+p} \big) \smile\pi^*\{\omega^{(j)}\}   \smile  h_{\overline \E}^{k-p-j}.
 \end{equation*}
 This proves   formula  \eqref{e:Lelong-numbers-vs-DS-bis} of assertion (2).
 
To complete the proof of assertion (2), we need to show that the total cohomology class $
\bfc^\DS_i(T,V)$ is  independent of the choice of $\tau$ and  $T_\infty.$ We fix a  Hermitian metric on $\E.$
  By the Leray-Hirsch theorem \cite{BottTu},  the cohomology $H^\ast(\P(\overline\E))$ is
a free module over $H^\ast(V)$ with basis $\{1, x, \ldots, x^{k-l}\},$  where $x:= h_{\overline \E}.$ 
So  $x\in H^{1,1}(\P(\overline\E)).$ Moreover,  $x^{k-l+1}$ can be written
uniquely as a linear combination of $1, x, \ldots, x^{k-l}$ with coefficients in
$H^\ast(V);$ these coefficients are by definition the Chern classes of the complex
vector bundle $\overline\E:$
\begin{equation}\label{e:min-poly}
x^{k-l+1}+c_1(\overline \E)x^{k-l}+\ldots+ c_{k-l+1}(\overline \E)=0,\qquad c_j(\overline\E)\in H^{j,j}(V).
\end{equation}
In this  equation we identify $c_j(\overline\E)$ to $\pi_0^*(c_j(\overline\E)).$
 With this definition of the Chern classes, we see that
the ring structure of the cohomology of $\P(\overline\E)$ is given by
 \begin{equation}\label{e:ring-structure}
 H^\ast(\P(\E)) = H^*(V)[x]/\big(x^{k-l+1}+c_1(\overline \E)x^{k-l}+\ldots+ c_{k-l+1}(\overline \E)\big).
 \end{equation}
  Let $T_\infty$ and $T'_\infty$ be two tangent current to $T$ along $V.$ By Tangent Theorem I (Theorem 3.8) in \cite{Nguyen21},  there are  strongly admissible maps $\tau,$ $\tau'$ along $V$ in $X$  and   two sequences $(\lambda_n)_{n\in\N},$  $(\lambda'_n)_{n\in\N}\subset\C$  such that $(\lambda_n)\nearrow\infty,$ 
  $(\lambda'_n)\nearrow\infty,$ and that 
  $$T_\infty=\lim_{n\to\infty} T_{\lambda_n}\qquad\text{and}\qquad T'_\infty=\lim_{n\to\infty} T'_{\lambda'_n},$$ 
  where
  $T_\lambda:= (A_\lambda)_* \tau_* (T )$ and   $T'_\lambda:= (A_\lambda)_* \tau'_* (T ).$
  We know that $T_\infty$ and $T'_\infty$  are both positive pluriharmonic  $(p,p)$-currents on $\overline\E.$
  Consider the  pluriharmonic  current $S:= T_\infty- T'_\infty$ on $\overline\E.$ To complete the proof of assertion (2),  we only need to show that
 $\{T_\infty\}=\{ T'_\infty\},$ i.e. $\{S\}=0.$
   By  Leray-Hirsch  theorem,
  we have the following   unique decomposition:
\begin{equation}\label{e:decomposition-S}
\{S\}=\sum_{j=\lowm}^\upm \pi_0^*(\{S_j\})\smile  h_{\overline \E}^{j-l+p},
\end{equation}
where   $\{S_j\}$  is a class in  $H^{l-j,l-j}(V,\C).$

 On the other  hand,  by   \eqref{e:Lelong-numbers-vs-DS-bis}, 
 we have for $\lowm\leq j\leq \upm$ and for every  closed   smooth real $(j,j)$-form   $\omega^{(j)}$ on $V$ that 
 \begin{eqnarray*}
  \sum_{i=j}^\upm \pi^*_0 \{(T_\infty)_i\}\smile \pi^*_0 \{\omega^{(j)}\} \smile h_{\overline\E}^{k-l+i-j}&=&\nu(T,V,\omega^{(j)},\tau,h),\\
  \sum_{i=j}^\upm \pi^*_0 \{(T'_\infty)_i\}\smile \pi^*_0 \{\omega^{(j)}\} \smile h_{\overline\E}^{k-l+i-j}&=&\nu(T,V,\omega^{(j)},\tau',h).
 \end{eqnarray*}
Since we have by Theorem \ref{T:Nguyen-1} (3) that  $\nu(T,V,\omega^{(j)},\tau,h)=\nu(T,V,\omega^{(j)},\tau',h),$
 substracting the  second line  from the first one and  using    $x:= h_{\overline \E}$ yield that
 the class $ \{S\}\in H^{p,p}(\overline\E)$ satisfies
  the  equation
 \begin{equation}\label{e:coho-classes-relations}
 \sum_{i=j}^\upm \pi^*_0 \{S_i\}\smile \pi^*_0 \{\omega^{(j)}\} \smile x^{k-l+i-j}=0, 
  \end{equation}
  for $\lowm\leq j\leq \upm$ and for every  closed   smooth real $(j,j)$-form   $\omega^{(j)}$ on $V.$
  It follows from  \eqref{e:min-poly} and \eqref{e:ring-structure} that the map $\theta:\ H^*(V)\to H^*(\P(\overline\E)),$ given by
  $\theta(\{\gamma\}):=\{\gamma\}\smile x^{k-l}$ for $\{\gamma\}\in H^*(V), $ is  one-to-one.
  
  Applying  \eqref{e:coho-classes-relations}  for $j=\upm$ yields that  $\pi^*_0 \{S_\upm\}\smile \pi^*_0 \{\omega^{(\upm)}\} \smile x^{k-l}=0.$
  The injectivity of $\theta$ implies that $\{S_\upm\}\smile \{\omega^{(\upm)}\}=0$ in $H^\ast(V).$  By choosing $\{\omega^{(\upm)}\}$ among
  a basis of $H^{\upm,\upm}(V),$ it follows from the Poincar\'e duality that   $\{S_\upm\}=0.$
  
  If $\upm=\lowm$ then  we stop. Otherwise, using  $\{S_\upm\}=0$ and  applying  \eqref{e:coho-classes-relations}  for $j=\upm-1$ and arguing as the  above paragraph  yields that $\{S_{\upm-1}\}=0.$
  We continue this process  until we show that  $\{S_\lowm\}=0,\ldots, \{S_\upm\}=0.$ Hence,  by \eqref{e:decomposition-S} $\{S\}=0.$
This completes the proof of the  theorem.
\endproof

\proof[End of the proof of Theorem  \ref{T:Nguyen-Siu-AB}]
By Theorem \ref{T:Dinh-Sibony} and Theorem \ref{T:Nguyen-Dinh-Sibony},  the classes $\bfc^\DS_i(T,V)$'s are independent of the  choice of an admissible map $\tau.$ Therefore, the RHS's  in the  formulas  of  Theorem \ref{T:Nguyen-Dinh-Sibony} are independent of  the  choice of an admissible map $\tau$ and a Hermitian metric $h$ on $\overline\E.$   By these formulas, so are the LHS's. 
\endproof
 
 \section{Applications in 
 intersection theory}\label{S:Intersection}
 
 \subsection{A criterion for  the uniqueness of the tangent currents}
Our first main result  provides a relevant ``local'' sufficient  condition for    the uniqueness of the tangent currents.
It is very convenient in practice.

For a $(p,p)$-current $T$ on $X$ consider the 
current $\widehat T$ on $X$ defined by
\begin{equation}\label{e:widehat-T} \widehat T(x)=-\log\dist(x,V)\cdot T(x),\qquad x\in X,
\end{equation}
where $\dist(x,V)$ is  the distance from $x$ to $V$ with respect to a fixed smooth Hermitian metric on $X.$

\begin{theorem}\label{T:unique-tangent-current} Let $X$ be a complex manifold of dimension $k.$
Let $V\subset  X$ be   a  K\"ahler submanifold  of dimension $l,$ and  $\omega$ a K\"ahler form on $V,$ 
and $B\subset V$
a  relatively compact piecewise  $\Cc^2$-smooth open subset  admitting a  strongly admissible map.
Suppose in addition   that   $\overline B$ can be covered  by a finite collection of    $\Cc^2$-smooth  domains $(B_i)_{i\in I}$ with $B_i\Subset V$  such that $B_i\cap B$ is  a  $\Cc^2$-smooth open subset in $V.$  Let $\tau_i$ be a  strongly admissible  map for $B_i,$ and $h_i$   a Hermitian metric  on $\E|_{\overline B_i}$ for $i\in I.$
Let $T$ be a positive closed currents of bidegree $(p,p)$  
 such that     
     $T =T^+-T^-$ on a  neighborhood of $\overline B$ in $X,$
 where $T^\pm\in \CL^{p;1,1} (B),$ and that 
    for all $i\in I,$    
    \begin{equation}\label{e:Nguyen}\kappa^\bullet_j( \widehat T,B_i\cap B,\omega, \bfr,\tau_i,h_i) <\infty \qquad\text{ for all}\qquad \lowm\leq j\leq \upm    .
    \end{equation}
Then $T$ admits a  unique  tangent  current along $B$.  
\end{theorem}

Prior  to the  proof of Theorem  \ref{T:unique-tangent-current}, we need  some preparatory  results.
Recall from \cite[Proposition 7.9]{Nguyen21} the  following: 
\begin{proposition}\label{P:basic-admissible-estimates-II}
There  are constants $c_3,c_4>0$  such that   the conclusion of Proposition \ref{P:basic-admissible-estimates-I}  holds and that
for every $1\leq \ell \leq \ell_0,$
 the following inequalities hold  on $\U_\ell \cap \Tube(B,\bfr):$
\begin{enumerate}

\item  
$\pm\big(   \tilde \tau_\ell^*(\alpha) -\alpha  \big)^\sharp\lesssim   c_3 \pi^*\omega +c_4\beta+c_3\varphi^{1/2}\alpha$ and 
$  \big(   \tilde \tau_\ell^*(\alpha) -\alpha  \big)\trianglelefteq   c_3 \pi^*\omega +c_4\beta+c_3\varphi^{1/4}\alpha ;$
\item 
$\pm\big(  \tilde \tau_\ell^*(\hat\alpha) -\hat\alpha -H\big)^\sharp \lesssim  c_3 \pi^*\omega +c_4\hat\beta +c_3\varphi^{1/2}\hat\alpha$ and 
$\big(  \tilde \tau_\ell^*(\hat\alpha) -\hat\alpha  \big) \trianglelefteq   c_3 \pi^*\omega +c_4\hat\beta +c_3\varphi^{1/4}\hat\alpha.$   Here,  $H$ is  some  form in the class $\Hc$ given in Definition  \ref{D:Hc}.           
\end{enumerate} 
\end{proposition}

The following result is needed.
\begin{lemma}\label{L:comparison-two-tau-two-h} Let $\tau'_i$ be a  strongly admissible  map for $B_i,$ and $h'_i$   a Hermitian metric  on $\E|_{\overline B_i}$ for $i\in I.$ Let $T$ be a positive  $(p,p)$-current on $X.$
Then the following assertions hold:
\begin{enumerate}
 \item There is a constant $c>0$ independent of $T $ such that such for every $i\in I,$
$$c^{-1}\Kc_{j,k-p-j}(T,B,0,c^{-1} r,\tau'_i,h')\leq \Kc_{j,k-p-j}(T,B,0,r,\tau_i,h)\leq c\Kc_{j,k-p-j}(T,B,0,cr,\tau'_i,h').$$
\item  The  following integrals  are  either simultaneously finite or  simultaneously infinite:  $$ \sum_{j=\lowm}^\upm \int_0^{\bfr} {\Kc_{j,k-p-j}(T,B_i\cap B ,0,r,\tau_i,h)\over  r}dr\quad\text{and}\quad  \sum_{j=\lowm}^\upm \int_0^{\bfr} {\Kc_{j,k-p-j}(T,B_i\cap B,0,r,\tau'_i,h')\over  r}dr.$$
\end{enumerate}

\end{lemma}
\proof   Assertion (2) is an  immediate consequence of assertion (1). For the sake of simplicity we  omit  the sub-index $i$  and  suppose that $B_i=B.$
We only need to prove  assertion (1)  in the following  two cases:

\noindent {\bf  Case I:} $h'=h.$

By \eqref{e:local-mass-indicators-bis} and \eqref{e:global-mass-indicators}, we see easily that there  is a constant  $c>0$ such that
$$c^{-1}\Kc_{j,k-p-j}(T,B,0,r/2,\tau',h)\leq \Kc_{j,k-p-j}(T,B,0,r,\tau,h)\leq c\Kc_{j,k-p-j}(T,B,0,2r,\tau',h).$$
This  inequality implies  the result.

\noindent {\bf  Case II:} $\tau'=\tau.$

It follows from Lemma \ref{L:hat-alpha-beta} that
\begin{equation}\label{e:hat-alpha-beta-with-ver-exp}\hat\beta\approx c_1\varphi\cdot \pi^*\omega+c_2\beta_\ver\quad\text{and}\quad
\hat\alpha\approx c_1\varphi\cdot \pi^*\omega+c_2\beta_\ver +c_1\alpha_\ver.
\end{equation}
Denote by $\alpha^{h},$ $\beta^{h}$ (resp. $\alpha^{h'},$ $\beta^{h'}$) the  forms $\alpha,$ $\beta$ associated to the metric $h$ (resp. $h'$). 
On the  other hand,  using  \eqref{e:alpha-beta-ver} and \eqref{e:tilde-alpha-beta-local-exp},
we  see that 
$$\alpha^h_\ver\approx \alpha^{h'}_\ver\quad\text{and}\quad \beta^h_\ver\approx \beta^{h'}_\ver.$$
Putting these estimates together implies that 
 $$\alpha^h\approx \alpha^{h'}\quad\text{and}\quad \beta^h\approx \beta^{h'}.$$
Using this and combining 
 \eqref{e:local-mass-indicators-bis} and \eqref{e:global-mass-indicators}, the last inequality implies  the result.
 \endproof

We keep the  hypothesis and  the notation in 
 Subsection \ref{SS:Inequa-mass-indicators}.

\begin{lemma}\label{L:normal-ves-hash-widehat-T}  Let $\lowm \leq j\leq \upm$ and set $\bfj:=(k-j,0,j,0)$  and  write, according \eqref{e:I_bfj},
$$I_j(\widehat T,r):= I_\bfj(\widehat T,r),\ I_j(\widehat T,s,r):= I_\bfj(\widehat T,s,r)\quad \text{and}\quad I^\hash_j(\widehat T,r) := I^\hash_\bfj(\widehat T,r),\ I^\hash_j(\widehat T,s,r) := I^\hash_\bfj(\widehat T,s,r).$$
Let $0\leq s<r\leq \bfr.$
Then  there  is  a constant $c$ independent of $T$, $s,r$ such that
\begin{equation*}
\begin{split}
 |I_{j} (\widehat T, r)- I^\hash_{j}(\widehat T,r)|&\leq cr^{1\over 4}(\log r)^{1/2}\big(\sum_{i:\ i\geq j} I^\hash_i(\widehat T,r+c_0r^2)\big),\\
 |I_{j} (\widehat T, s,r)- I^\hash_{j}(\widehat T,s,r)|&\leq cr^{1\over 4}(\log r)^{1/2}\big(\sum_{i:\ i\geq j} I^\hash_i(\widehat T,s-c_0s^2,r+c_0r^2)\big).  
\end{split}\end{equation*}
\end{lemma}
\proof  We only give the proof of the second inequality since the proof of the first one is
similar.
Applying  Lemma 
\ref{L:spec-wedge} yields  a constant $c$ independent of $T$ and $s,r$   such that the   following inequality holds
 \begin{equation*}
|I_\bfj(\widehat T,s,r)- I^\hash_\bfj(\widehat T,s,r)|^2 \leq c\big(\sum_{\bfj'} I^\hash_{\bfj'}(\widehat T,s-c_0s^2,r+c_0r^2)\big)\big ( \sum_{\bfj''} I^\hash_{\bfj''}(\widehat T,s-c_0s^2,r+c_0r^2) \big).  
\end{equation*}
Here, on the RHS:
\begin{itemize} \item[$\bullet$] the first sum  is taken over a finite number of multi-indices    $\bfj'=(j'_1,j'_2,j'_3,j'_4)$ as above  such that  $j'_1\leq  j_1$  and $j'_2\geq j_2;$ and either ($j'_3\leq j_3$) or ($j'_3>j_3$ and $j'_2\geq j_2+{1\over 2}$).
\item  the second sum   is taken over  a finite number of multi-indices $\bfj''=(j''_1,j''_2,j''_3,j''_4)$ as above   such that   either  ($j''_1< j_1$)
or ($j''_1=j_1$ and $j''_2\geq {1\over 4}+j_2$) or ($j''_1=j_1$ and $j''_3<j_3$).
\end{itemize} 
By   the  first item $\bullet$, the  first sum $\sum_{\bfj'} I^\hash_{\bfj'}(\widehat T, s-c_0s^2,r+c_0r^2)$ is bounded by a constant times 
$\sum_{i:\ i\geq j} I^\hash_j(\widehat T,s-c_0s^2,r+c_0r^2).$

 By  Lemma
  \ref{L:Kc_j,q-log} below,  the  second item $\bullet$, the  second  sum  is bounded by a constant times 
 $$\sum_i \Kc_{i,q}(\varphi^{1/4}\widehat T, s-c_0s^2,r+c_0r^2)+  \sum_{i: i<k-p-q} \Kc_{i,q}
 (\widehat T, s-c_0s^2,r+c_0r^2).$$ 
We can  prove   that  the second sum in the last line  is bounded by a constant times $r^2\log r.$ 
In all, the   the second sum  $\sum_{\bfj''} I^\hash_{\bfj''}(s-c_0s^2,r+c_0r^2)$ is bounded by a constant times $cr^{1\over 2}(\log r)\big(\sum_{i:\ i\geq j} I^\hash_j(\widehat T,s-c_0s^2,r+c_0r^2)\big).$
 This, combined with  the previous estimate on   $\sum_{\bfj'} I^\hash_{\bfj'}(s-c_0s^2,r+c_0r^2),$  implies the result.
\endproof

Recall from  \cite[Theorem 1.6]{NguyenTruong}  the following   result.
 
\begin{theorem}\label{T:NguyenTruong} Let $X,\, V,\,  B$ and $I,$  $(B_i)_{i\in I},$ $(\tau_i)_{i\in I},$ $(h_i)_{i\in I}$ be as
in Theorem \ref{T:unique-tangent-current}. Let $T$ be a positive closed currents of bidegree $(p,p)$  
 such that     
     $T =T^+-T^-$ on a  neighborhood of $\overline B$ in $X,$
 where $T^\pm\in \CL^{p;1,1} (B),$ and that 
    for all $i\in I,$
 \begin{equation} \label{e:Nguyen-Truong}  \int_0^{\bfr} {\kappa^\bullet_j(T,B_i\cap B,\omega, r,\tau_i,h_i)\over  r}dr <\infty \qquad\text{ for all}\qquad \lowm\leq j\leq \upm    .
 \end{equation}
Then $T$ admits a  unique  tangent  current along $B$.  
\end{theorem}
\proof
Since assumption \eqref{e:Nguyen-Truong}   is  condition (b-i) in    \cite[Theorem 1.6]{NguyenTruong}, the result follows.
\endproof

\proof[Proof of Theorem  \ref{T:unique-tangent-current}] 

Applying  Lemma \ref{L:normal-ves-hash-widehat-T}, assumption \eqref{e:Nguyen}
implies that  $I^\hash_{j}(\widehat T,B_i\cap B,0,r,\tau_i,h_i)<\infty,$ that is, 
\begin{equation}\label{e:Kc-widehat-inequ}\Kc_{j,k-p-j}(\widehat T,B_i\cap B,0,r,\tau_i,h_i)<\infty\qquad\text{for all}\qquad\lowm\leq j\leq \upm.
\end{equation}
Fix  $j$ with $\lowm\leq j\leq \upm .$  Write   
 \begin{equation}\label{e:int_kappa-i_j-over-r}    \int_0^{\bfr} {\Kc_{j,k-p-j}(T,B_i\cap B,\omega,r,\tau_i,h_i)\over  r}dr =  \sum_{n=0}^\infty \int_{\bfr\over 2^{n+1}}^{\bfr\over 2^n} {\Kc_{j,k-p-j}(T,B_i\cap B,\omega,r,\tau_i,h_i)\over  r}dr.
 \end{equation}
 Since  we have,  for  $r\in[{\bfr\over 2^{n+1}},{\bfr\over 2^n} ], $
\begin{eqnarray*}
\Kc_{j,k-p-j}\big(T,B_i\cap B,\omega,{\bfr\over 2^{n+1}} ,\tau_i,h_i\big)&\leq& \Kc_{j,k-p-j}\big(T,B_i\cap B,\omega,r,\tau_i,h_i\big)\\
&\leq& \Kc_{j,k-p-j}\big(T,B_i\cap B,\omega,{\bfr\over 2^n} ,\tau_i,h_i\big),
\end{eqnarray*}
 the RHS of \eqref{e:int_kappa-i_j-over-r} is  dominated by  
\begin{equation*}    \sum_{n=0}^\infty \int_{\bfr\over 2^{n+1}}^{\bfr\over 2^n} {\Kc_{j,k-p-j}\big(T,B_i\cap B,\omega,{\bfr\over 2^n},\tau_i,h_i\big)\over  r}dr= \ln{2}\cdot  \sum_{n=0}^\infty  \Kc_{j,k-p-j}\big(T,B_i\cap B,\omega,{\bfr\over 2^n},\tau_i,h_i\big) .
 \end{equation*}
 Rewrite the sum on the RHS as
 \begin{equation*}      \sum_{n=0}^\infty \sum_{q=n}^\infty \Kc_{j,k-p-j}\big(T,B_i\cap B,\omega,{\bfr\over 2^{q+1}},{\bfr\over 2^{q}},\tau_i,h_i\big)=  \sum_{q=0}^\infty (q+1) \Kc_{j,k-p-j}\big(T,B_i\cap B,\omega,{\bfr\over 2^{q+1}},{\bfr\over 2^{q}},\tau_i,h_i\big).
 \end{equation*}
 Since for each $q\in\N,$  $|\log\dist(y,V)-q|\leq  c$ for  $y\in\Tube(B_i\cap B, {\bfr\over 2^{q+1}},{\bfr\over 2^{q}}),$
 it follows that
  \begin{equation*}     (q+1)  \Kc_{j,k-p-j}\big(T,B_i\cap B,\omega,{\bfr\over 2^{q+1}},{\bfr\over 2^{q}},\tau_i,h_i\big)\lesssim
    \Kc_{j,k-p-j}\big(\widehat T,B_i\cap B,\omega,{\bfr\over 2^{q+1}},{\bfr\over 2^{q}},\tau_i,h_i\big).
 \end{equation*}
 Hence,  the last sum is  dominated by  $\Kc_{j,k-p-j}\big(\widehat T,B_i\cap B,\omega,0,\bfr,\tau_i,h_i\big),$ which is  in turn finite by   \eqref{e:Kc-widehat-inequ}. By  \cite{NguyenTruong},
 condition \eqref{e:Nguyen-Truong} is  equivalent to 
 the following appearing in \eqref{e:int_kappa-i_j-over-r}:
  \begin{equation*}     \int_0^{\bfr} {\Kc_{j,k-p-j}(T,B_i\cap B,\omega,r,\tau_i,h_i)\over  r}dr <\infty, \qquad\text{ for all}\qquad \lowm\leq j\leq \upm   .
 \end{equation*}
 Therefore, condition \eqref{e:Nguyen}
 is  equivalent to condition \eqref{e:Nguyen-Truong}.
\endproof
\begin{remark}\rm
 \label{R:four-equivalent-conds}
 In fact,  the  method of the proof of the above theorem also shows that for every $i\in I$ and for every $j_0\in\N$ with 
 $\lowm\leq j_0\leq \upm,$ the  following four conditions  are  equivalent:
 \begin{enumerate}
  \item  \begin{equation*}
  \int_0^{\bfr} {\kappa^\bullet_j(T,B_i\cap B,\omega, r,\tau_i,h_i)\over  r}dr <\infty \qquad\text{ for all}\qquad j_0\leq j\leq \upm    .
 \end{equation*}
 \item \begin{equation*}
\kappa^\bullet_j( \widehat T,B_i\cap B,\omega, \bfr,\tau_i,h_i) <\infty \qquad\text{ for all}\qquad j_0\leq j\leq \upm    .
    \end{equation*}
 \item $$ \sum_{j=j_0}^\upm \int_0^{\bfr} {\Kc_{j,k-p-j}(T,B_i\cap B,\omega,0,r,\tau,h)\over  r}dr<\infty.$$
 \item $$\sum_{j=j_0}^\upm  \Kc_{j,k-p-j}(\widehat T,B_i\cap B,\omega,0,r,\tau,h)\leq \infty.$$
 \end{enumerate}
\end{remark}
We give here  the  following consequence   which captures the  essential points of the above  Theorem \ref{T:unique-tangent-current} in the  special but important context where the ambient manifold $X$ is K\"ahler. We expect that  this  explicit statement  will be useful in practice.
\begin{corollary}\label{C:unique-tangent-current}
 Let $X,V$ be as above   and  suppose that $X$ is K\"ahler and that $(V,\omega)$ is K\"ahler  and  that $B\subset V$ is 
a  relatively compact piecewise  $\Cc^2$-smooth open subset.  Let $T$ be a  positive closed current of bidegree $(p,p)$ on $X$ such that there is a neighborhood $W$ of $\partial B$ in $X$ such that the restriction of $T$ on $W$ is    a $\Cc^1$-smooth form and that 
    for all $i\in I,$    $$ \kappa^\bullet_j( \widehat   T,B_i\cap B,\bfr,\tau_i,h_i) <\infty \qquad\text{ for all}\qquad \lowm\leq j\leq \upm    .$$
Then $T$  admits a unique  tangent current  along $B.$  
\end{corollary}

 The following  auxiliary result will be  useful in the next subsection.
\begin{lemma}
  \label{L:Kc_j,q-log} 
Let $T\in \CL^{p;1,1} (B).$
Let  $\lowm \leq j\leq \upm$ and  $0\leq q\leq k-p-j-1.$
Then   $\Kc_{j,q}(\widehat T, r,\tau,h)\leq  cr^{2(k-p-q -j)}|\log{r}|.$
\end{lemma}
\proof
We argue as in the proof of \cite[Proposition 4.3]{NguyenTruong} (which in turn relies on  Theorem \ref{T:Lc-finite}) making  the obviously necessary changes and taking into account expression \eqref{e:widehat-T}.
Consequently, we can show  that $\Kc_{j,q}(\widehat S,0, r,\tau,h)\leq  cr^{2(k-p-q -j)}|\log{r}|$
for all smooth   positive  closed $(p,p)$-currents $S\in \widetilde \CL^{p;1,1}(\bfU,\bfW),$ where $c>0$ is a constan independent of $S.$ 
Since $S$ is smooth, we  infer that
$$\Kc_{j,q}(\widehat S, r,\tau,h) = \Kc_{j,q}(\widehat S,0, r,\tau,h)  \leq  cr^{2(k-p-q -j)}|\log{r}|,$$
Using an approximation of $T$ by smooth  positive  closed $(p,p)$-currents $S$'s, the result follows.
\endproof

We conclude  this  subsection with the  following result which will be needed in Section  \ref{S:Continuity}.
\begin{proposition}\label{P:two-equivalent-conds}
 For every $i\in I$ and for every $j_0\in\N$ with 
 $\lowm\leq j_0\leq \upm,$ the  following two conditions  are  equivalent:
 \begin{enumerate}
 \item $
\kappa_j( \widehat T,B_i\cap B,\omega, \bfr,\tau_i,h_i) <\infty$   for all $ j_0\leq j\leq \upm    .$
 \item $$\sum_{j=j_0}^\upm  \Kc_{j,k-p-j}(\widehat T,B_i\cap B,\omega,r,\tau,h)\leq \infty.$$
 \end{enumerate}
 In particular, the  condition that $
\kappa_j( \widehat T,B_i\cap B,\omega, \bfr,\tau_i,h_i) <\infty$   for all $i\in I$ $ j_0\leq j\leq \upm$ is independent of the choice of $\tau_i, $ $h_i.$
\end{proposition}
\proof  
Applying  the first  inequality of Lemma \ref{L:normal-ves-hash-widehat-T} and applying Theorem  \ref{T:Lc-finite} and  Lemma 
  \ref{L:Kc_j,q-log}, the first  assertion  follows.
  
  The second assertion holds by combining  the first one and  Lemma \ref{L:comparison-two-tau-two-h}.
\endproof
 \subsection{A criterion for the wedgeablity in the sense of Dinh-Sibony}




We keep the  hypothesis and  the notationn  in 
 Subsection \ref{SS:Fourth-main-results}. 
 Let $(X,\omega)$  be a compact K\"ahler manifold of dimension $k.$
 Consider $m \geq 2$ integers $p_1,\ldots, p_m\geq 1$  such that   $p:=p_1+\ldots+p_m\leq k.$
 Let $\Delta:=\{(x,\ldots,x):\ x\in X\}$  be the diagonal of $X^m,$ and $\omega_\Delta$ be a K\"ahler form on $\Delta,$   and $\tau$  be a  strongly admissible map  along $\Delta$ in $X^m.$ 
 Let $\pi:\ \E\to\Delta$ be the normal bundle to $\Delta$ in $X^m.$
 Let $h$ be a Hermitian metric on $\E.$ 
 Using formula \eqref{e:alpha-beta-spec} we define  $\alpha$ and $\beta$ on $\E.$
 Using formula  \eqref{e:alpha-beta-ver} we define  $\alpha_\ver$ on $\E.$
Let $\dist(\bfx,\Delta)$ be the distance  from a point  $\bfx\in X^m$ to $\Delta.$
We may assume that  $\dist(\cdot,\Delta)\leq 1/2.$  So 
$$\widehat\T:=-\log\dist(\cdot,\Delta) \cdot\T$$
is a positive $(p,p)$-current on $X^m.$
 
 Let  $T_j\in \CL^{p_j}(X)$ for $1\leq j\leq m$  with  $p:=p_1+\ldots+p_m\leq k=\dim(X).$ Consider $\T:=T_1\otimes\ldots\otimes T_m\in \CL^p(X^m).$
 In this context, $\lowm=k-p$ and  $\upm=\min(k,mk-p).$

\begin{lemma}\label{L:big-power-alpha_ver-equals-zero} It holds that  $\alpha_\ver^{(m-1)k}=0.$
\end{lemma}
\proof
We use  the  local model recalled in Subsection \ref{SS:Anal-local-coord}.  
 It follows  from \eqref{e:tilde-alpha-beta-local-exp}  that  
\begin{equation}
\alpha_\ver(z,w)=  A(w)^* [\ddc \log{\|z\|^2}] \quad\text{for}\quad  z\in\C^{(m-1)k},\ w\in \D^k,
\end{equation}
where   $A:\  \D^l\to \GL(\C,k-l)$ is a smooth function. 
Since  $(\ddc \log{\|z\|^2})^{(m-1)k}=0,$ the result follows.
\endproof 

 The  last lemma of the  section shows  that the finiteness condition of the minimal dimension is superflous. 
 \begin{lemma}\label{L:min-dim-finite-cond-is-superflous}
  Suppose that assumption (1) of  Theorem
  \ref{T:Nguyen-intersection} is  fulfilled, that is,
  $$\kappa^\bullet_j(-\log\dist(\cdot,\Delta)\cdot\T,\Delta,\omega_\Delta,\bfr,h)<\infty\quad\text{for some  $\bfr>0$ and for all}\quad   k-p< j\leq k-\max_{1\leq i\leq m} p_i.$$ 
  Then the above  inequality also holds for $j=k-p,$ that is,
  $$\kappa^\bullet_{k-p}(-\log\dist(\cdot,\Delta)\cdot\T,\Delta,\omega_\Delta,\bfr,h)<\infty.$$ 
 \end{lemma}
 \proof
 First  consider an index  $j$ with $k-\max_{1\leq i\leq m}p_i<j\leq \upm=\min(k,mk-p).$ 
   Using  a local computation  and  \eqref{e:local-mass-indicators-bis} and \eqref{e:global-mass-indicators}, we see easily that
   $\Kc_{j,k-p-j}\big(T,B_i\cap B,\omega,0,r,\tau_i,h_i\big)=0.$ 
   Therefore,  applying   Lemma \ref{L:Kc_j,q-log} and Lemma \ref{L:normal-ves-hash-widehat-T} yields that
   $$
   \lim_{r\to 0+} \kappa^\bullet_j(-\log\dist(\cdot,\Delta)\cdot \T,\Delta,\omega_\Delta, r,\tau,h)=0.
   $$
 Hence, we are only concerned with  indices  $j$ with $\lowm=k-p\leq j\leq k-\max_{1\leq i\leq m}p_i.$ 
 By Lemma  \ref{L:big-power-alpha_ver-equals-zero} $\alpha_\ver^{(m-1)k}=0.$
 On the other hand,
 by \eqref{e:hat-alpha-beta-with-ver-exp}, we  have that
 $
\hat\alpha\lesssim  c_1\varphi\cdot \pi^*\omega+c_2\hat\beta +c_1\alpha_\ver.
$
 Using  these two estimates,
 we can show that
$$
\hat\alpha^{(m-1)k}\leq c\sum_{(i,j)} \hat\alpha^i \wedge \pi^*\omega^{j} \wedge \hat\beta^{(m-1)k-i-j},
$$ 
where  the sum is  taken over all $(i,j)\in\N^2$ with $i\leq (m-1)k-1$ and $i+j\leq (m-1)k.$
Using  this and  \eqref{e:local-mass-indicators-bis} and \eqref{e:global-mass-indicators}, and applying Lemma    \ref{L:Kc_j,q-log} and Lemma \ref{L:normal-ves-hash-widehat-T}  we see easily that there  is a constant  $c>0$ such that
$$\Kc_{\lowm,0}(T,B,0,r,\tau,h)\leq c\sum_{j=\lowm+1}^{\upm}\Kc_{j,k-p-j}(T,B,0,r,\tau,h) +cr^2|\log r|.$$ By  Remark \ref{R:four-equivalent-conds} for $j_0=\lowm+1,$  the result follows.
 \endproof

 \proof[End of the  proof of Theorem \ref{T:Nguyen-intersection}]
 By Lemma \ref{L:min-dim-finite-cond-is-superflous}, condition (1) of the theorem implies 
 condition \eqref{e:Nguyen}. Therefore, by Theorem \ref{T:unique-tangent-current}, there exists a unique  tangent current to $\T$ along  $\Delta.$
 
 On the  other hand, by Theorem  \ref{T:Nguyen-DS} (3),  condition (2) of the theorem implies that the horizontal dimension of $\T$ along  $\Delta$ is  minimal, i.e. $\hbar=\lowm=k-p.$ Then  by Theorem  \ref{T:Dinh-Sibony-product} and  by Definition \ref{D:DS-wedge-prod}, $T_1\curlywedge\ldots\curlywedge T_m$ exists in the sense of  the theory of tangent current.
This  completes the proof.
 \endproof
\section{Continuity of  Dinh-Sibony intersection}
\label{S:Continuity}
 
 In this  section we combine  the technique of the generalized Lelong numbers developed in the previous  sections  with that of   Dinh-Nguyen-Vu \cite{DinhNguyenVu} on the super-potential theory in order to  study the continuity of  Dinh-Sibony intersection.
\subsection{Super-potentials}

Super-potentials are functions which play the role of quasi-potentials for positive closed
currents of arbitrary bi-degree.  For simplicity, we will not introduce this notion in full generality but limit ourselves in the necessary setting. 
Let $X$ be a compact K\"ahler manifold of dimension $k.$ We also fix  a K\"ahler form $\omega$ on $X$. If $T$ is a positive or negative $(p,p)$-current on $X$, its mass is given by $\|T\|:=\langle T,\omega^{k-p}\rangle$ or  $\|T\|:=-\langle T,\omega^{k-p}\rangle$ respectively. Let $\mD^q(X)$ (or $\mD^q$ for short) denote the real vector space spanned by positive closed $(q,q)$-currents on $X$. 
Define the {\it $\ast$-norm} on this space by $\|R\|_\ast:=\min ( \|R^+\|+\|R^-\|)$, where $R^\pm$ are positive closed $(q,q)$-currents satisfying $R=R^+-R^-$ and $\|\ \|$ denotes the mass of a current. We consider this space of currents with the following {\it topology} : a sequence $(R_n)_{n\geq 0}$ in $\mD^q(X)$ converges in this space to $R$ if $R_n\to R$ weakly and if $\|R_n\|_\ast$ is bounded independently of $n$. On any $\ast$-bounded set of $\mD^q(X)$, this topology coincides with the classical weak topology for currents. 
It was shown in \cite{DinhSibony04} that the subspace $\widetilde\mD^q(X)$ of real closed smooth $(q,q)$-forms is dense in $\mD^q(X)$ for the considered topology.

Let $\mD^q_0(X)$ and $\widetilde\mD^q_0(X)$ denote the linear subspaces in $\mD^q(X)$ and $\widetilde\mD^q(X)$ respectively of currents whose cohomology classes in $H^{q,q}(X,\R)$ vanish. Their co-dimensions are equal to the dimension of $H^{q,q}(X,\R)$ which is finite. Fix a real smooth and closed $(p,p)$-form $\gamma_T$ in the same cohomology class with $T$ in $H^{p,p}(X,\R)$.
We will consider in this paper the super-potential of $T$ which is the real function $\mU_T$ on $\widetilde\mD_0^{k-p+1}$ defined by 
\begin{equation} \label{e:mU} \mU_T(R):=\langle T-\gamma_T, U_R\rangle \qquad \text{for} \qquad R\in \widetilde\mD_0^{k-p+1},
\end{equation}
where $U_R$ is any smooth form of bi-degree $(k-p,k-p)$ such that $\ddc U_R=R$. This form always exists because the cohomology class of $R$ vanishes. Note that since the cohomology class of $T-\gamma_T$ vanishes, we can write $T-\gamma_T=\ddc U_T$ for some current $U_T$. By Stokes theorem, we have  
$$\mU_T(R)=\langle \ddc U_T,U_R\rangle =\langle U_T,\ddc U_R\rangle =\langle U_T, R\rangle.$$
We deduce from this identities that $\mU_T(R)$ does not depend on the choice of $U_R$ and $U_T$. However,  $\mU_T$  depends on the  reference form $\gamma.$ Note also that
if $T$ is smooth, it is not necessary to take $R$ and $U_R$ smooth.

We will not consider other super-potentials of $T$. They are some affine extensions of $\mU_T$ to $\mD^{k-p+1}$ or its extensions to some subspaces. The following notions do not depend on the choice of super-potential nor on the reference form $\alpha$. We say that $T$ has a {\it bounded super-potential} if $\mU_T$ is bounded on each $\ast$-bounded subset of $ \widetilde\mD_0^{k-p+1}$. 
We say that $T$ has a {\it continuous super-potential} if $\mU_T$ can be extended to a continuous function on $\mD_0^{k-p+1}$. 
Recall that $\mD^p(X)$ is a metric space. 

 
\begin{definition} \label{D:DS-wedge-prod-super} \rm {(Dinh-Sibony \cite{DinhSibony09,DinhSibony10})}
Consider now two positive closed currents $T$ and $S$ on $X$ of bi-degree $(p,p)$ and $(q,q)$ respectively. Assume that $p+q\leq k$ and that $T$ has a continuous super-potential. So $\mU_T$ is defined on whole $\mD_0^{k-p+1}$. We can define the wedge-product $T\wedge S$ by 
\begin{equation} \label{e:Formula-wedge-product-special}\langle T\wedge S,\phi\rangle := \langle \gamma_T\wedge S,\phi\rangle +\mU_T(S\wedge\ddc\phi)
\end{equation}
for every smooth real test form $\phi$ of bi-degree $(k-p-q,k-p-q)$.
We say that {\it $T\wedge S$ is defined in the sense of Dinh-Sibony's super-potential theory,} or equivalently,
{\it $T$ and $S$  are  wedgeable in the sense of Dinh-Sibony's super-potential theory}.
\end{definition}
Note that $S\wedge\ddc\phi$ belongs to $\mD_0^{k-p+1}$ because it is equal to $\ddc(S\wedge \phi)$. It is not difficult to check that $T\wedge S$ is equal to the usual wedge-product of $T$ and $S$ when one of them is smooth. The current $T\wedge S$ is positive and closed, see \cite{DinhNguyenVu, DinhSibony09, DinhSibony10,Vu16}.

Recall that the projections $\Pi_j:\widehat{X\times X}\to X$ are submersions, see e.g. \cite{DinhSibony04}.  
We have seen   that the definition of super-potential involves the solutions of the equation $\ddc U_R=R$  for $R\in \widetilde \mD_0^{k-p+1}.$ We will recall here the construction of kernel solving this
equation and refer to \cite{DinhSibony10} for details. 

By Blanchard's theorem \cite{Blanchard}, $\widehat{X\times X}$ is a K\"ahler manifold. So we fix
a K\"ahler form $\widehat\omega$ on $\widehat{X\times X}$ and we
can apply Hodge theory to this manifold. 
By K\"unneth's formula, the cohomology class $\{\Delta\}$ of $[\Delta]$ in $H^{k,k}(X \times X,\R)$ can be represented by a real smooth closed $(k,k)$-form $\gamma_{\Delta}$  which is a finite sum of forms of type $\pi_1^*(\phi_1) \wedge \pi_2^*(\phi_2)$ where $\phi_1$ and $\phi_2$ are closed smooth forms of suitable bi-degrees on $X$. By  \cite[1.2.1]{BGS} and \cite[Ex. 2.3.1]{DinhSibony10}, there is a real smooth closed $(k-1,k-1)$-form $\hat \eta$ on $\widehat{X\times X}$ such that $\hat\eta \wedge [\widehat\Delta]$ is cohomologous to $\Pi^*(\gamma_\Delta)$ and 
\begin{align} \label{e:Delta}
\Pi_{*}(\hat\eta \wedge [\widehat\Delta])=[\Delta].
\end{align}

Choose a real smooth closed $(1,1)$-form $\hat\gamma$ on $\widehat{X\times X}$ which is cohomologous to $[\widehat\Delta]$. So we can write $[\widehat\Delta]-\hat\gamma=\ddc \hat u$, where $\hat u$ is a quasi-p.s.h. function on $\widehat{X\times X}$. This equation implies that $ \hat u$ is smooth outside $\widehat\Delta$ and $\hat u-\log\dist(\cdot,\widehat\Delta)$ is a smooth function near $\widehat\Delta$. Subtracting from $\hat u$ a constant allows us to assume that $\hat u$ is negative.
Observe that since $\hat\gamma\wedge\hat\eta$ is cohomologous to $\Pi^*(\alpha_\Delta)$, there is a real smooth $(k-1,k-1)$-form $\hat\gamma'$ on $\widehat{X\times X}$ such that 
\begin{align} \label{e:hat-gamma-0}
\ddc \hat\gamma'=\hat\gamma\wedge\hat\eta -\Pi^*(\gamma_\Delta).
\end{align}
Adding to $\hat\gamma'$ a constant times $\widehat\omega^{k-1}$ allows us to assume that $\hat\gamma'$ is positive. 
For $\epsilon>0$, denote by $\widehat\Delta_\epsilon$ the set of points in $\widehat{X\times X}$ with distance less than $\epsilon$ to $\widehat\Delta$. 

\begin{proposition} \label{P:Kernel} 
\begin{itemize}
\item[(1)] {\rm (See \cite[Proposition 2.3]{DinhNguyenVu})}  If $R$ is in  $\mD^0_q(X)$ and 
$U_R:=(\Pi_1)_*\big((\hat u\hat\eta+\hat\gamma')\wedge\Pi_2^*(R)\big)$, then
$\ddc U_R=R$.
\item[(2)] {\rm (See \cite[Lemma 2.4]{DinhNguyenVu})}
Let $T,\gamma_T$ and $\mU_T$ be as above. Then for every $R$ in $\widetilde\mD^0_{k-p+1}(X)$, we have 
$$\mU_T(R)=  \int_{\widehat{X\times X}}   (\hat u\hat\eta+ \hat\gamma') \wedge \Pi_1^*(T-\gamma_T)\wedge \Pi_2^*(R).$$
\end{itemize}
\end{proposition}

The following proposition gives us a characterization of currents with bounded super-potentials.

\begin{proposition}\label{P:sp-bounded}  {\rm (See \cite[Proposition 2.5]{DinhNguyenVu})}
Let $T$ be a positive closed $(p,p)$-current on $X$. Then $T$ has a bounded super-potential if and only if there is a constant $c>0$ such that for every smooth positive closed $(k-p+1,k-p+1)$-form $R$ on $X$ with $\|R\|\leq 1$ we have
$$-\int_{\widehat{X\times X}} \hat u\,\widehat\omega^{k-1} \wedge \Pi_1^*(T)\wedge \Pi_2^*(R)\leq c.$$
\end{proposition}

We give now a criterion to check if a current has a continuous super-potential. Let $T$ be a positive closed $(p,p)$-current as above. Recall that
for $\epsilon>0$, $\widehat\Delta_\epsilon$ denotes the set of points in $\widehat{X\times X}$ with distance less than $\epsilon$ to $\widehat\Delta$. 
Consider the following quantity
\begin{equation}\label{e:vartheta-sup}
\vartheta_T(\epsilon):=\sup_R \int_{\widehat\Delta_\epsilon} -\hat u\widehat\omega^{k-1}\wedge \Pi_1^*(T)\wedge \Pi_2^*(R),
\end{equation}
where  the supremum is taken over all smooth positive closed forms $R$ on $X$, of bi-degree $(k-p+1,k-p+1)$, such that $\|R\|\leq 1$.

\begin{proposition} \label{P:sp-continuous}{\rm (See \cite[Proposition 2.7]{DinhNguyenVu})}
Let $T$ be a positive closed $(p,p)$-current on $X$. Then $T$ has a continuous super-potential if and only if $\vartheta_T(\epsilon)$ tends to $0$ as $\epsilon$ tends to $0$. 
\end{proposition}

\subsection{Generalization to $m$-fold  products}\label{SS:m-fold-products}
Fix an integer $m\geq 2.$ Let $\X:=X^m$ and let $\Delta:=\{ (x,\ldots,x):\ x\in X\}$ be the diagonal of $\X.$
 Let $\Pi:\widehat\X \to \X$ be the blow-up of $\X$ 
along the diagonal $\Delta$ and let $\widehat\Delta:=\Pi^{-1}(\Delta)$ be the exceptional hypersurface. Denote by $\pi_j$ the projections from $\X$ onto its $j$-th factors for $1\leq j\leq m,$ and define $\Pi_j:=\pi_j\circ\Pi$.
Recall that the projections $\Pi_j:\widehat{\X}\to X$ are submersions, see e.g. \cite{DinhSibony04}.  
By Blanchard's theorem \cite{Blanchard} again, $\widehat{\X}$ is a K\"ahler manifold. So we fix
a K\"ahler form $\widehat\omega$ on $\widehat{\X}$ and we
can apply Hodge theory to this manifold. 
By K\"unneth's formula, the cohomology class $\{\Delta\}$ of $[\Delta]$ in $H^{k,k}(\X,\R)$ can be represented by a real smooth closed $((m-1)k,(m-1)k)$-form $\gamma_{\Delta}$  which is a finite sum of forms of type $\pi_1^*(\phi_1) \wedge\ldots\wedge  \pi_m^*(\phi_1)$ where $\phi_1,\ldots,\phi_m$  are closed smooth forms of suitable bi-degrees on $X$. By  \cite[1.2.1]{BGS} and \cite[Ex. 2.3.1]{DinhSibony10}, there is a real smooth closed $((m-1)k-1,(m-1)k-1)$-form $\hat \eta$ on $\widehat{\X}$ such that $\hat\eta \wedge [\widehat\Delta]$ is cohomologous to $\Pi^*(\gamma_\Delta)$ and 
\begin{align} \label{e:Delta-bis}
\Pi_{*}(\hat\eta \wedge [\widehat\Delta])=[\Delta].
\end{align}

Choose a real smooth closed $(1,1)$-form $\hat\gamma$ on $\widehat{\X}$ which is cohomologous to $[\widehat\Delta]$. So we can write
\begin{equation} \label{e:hat-u}[\widehat\Delta]-\hat\gamma=\ddc \hat u,
\end{equation} where $\hat u$ is a quasi-p.s.h. function on $\widehat{\X}$. This equation implies that $ \hat u$ is smooth outside $\widehat\Delta$ and $\hat u-\log\dist(\cdot,\widehat\Delta)$ is a smooth function near $\widehat\Delta$. Subtracting from $\hat u$ a constant allows us to assume that $\hat u$ is negative.
Observe that since $\hat\gamma\wedge\hat\eta$ is cohomologous to $\Pi^*(\gamma_\Delta)$, there is a real smooth $((m-1)k-1,(m-1)k-1)$-form $\hat\gamma'$ on $\widehat{\X}$ such that 
\begin{align} \label{e:gamma}
\ddc \hat\gamma'=\hat\gamma\wedge\hat\eta -\Pi^*(\gamma_\Delta).
\end{align}
Adding to $\hat\gamma'$ a constant times $\widehat\omega^{(m-1)k-1}$ allows us to assume that $\hat\gamma'$ is positive. 
For $\epsilon>0$, denote by $\widehat\Delta_\epsilon$ the set of points in $\widehat{\X}$ with distance less than $\epsilon$ to $\widehat\Delta$. 
Let $\iota:\ X\to \Delta$ be the natural map $x\mapsto (x,\ldots,x)$ sending $X$ to $\Delta$. 
We also  denote  by $\iota:\ X\to X^m=\X$  the natural map $x\mapsto (x,\ldots,x)$ sending $X$ to $X^m=\X$.
 Let $\pi':\ \X=X^m\to X^{m-1}$ be the  canonical  projection onto the first $m-1$ factors.
 Set $\Pi':=\pi'\circ \Pi:\ \widehat \Pi\to X^{m-1}.$
 Denote  also by $\pi_1$ the  canonical projection $\pi_1:\ X^{m-1}\to X$ onto the  first factor.
 Let $\Delta'$ be  the  diagonal  of $X^{m-1},$ i.e. $\Delta'=\{(x,\ldots,x)\in X^{m-1}:\ x\in X\}.$
 For a   smooth form $\gamma$ on $X^{m-1}$ let $\gamma|_{\Delta'}$  be the restriction of $\gamma$ to $\Delta'.$
\begin{proposition} \label{P:kernel}
 Let  $R_2\in\mD^{p_2},\ldots, R_m\in\mD^{p_m}$ be  continuous  and  assume that  $\gamma_R$ is   a closed  smooth $(p',p')$-form  in  the  class of $R_2\otimes\ldots \otimes R_m$  in  $H^{p',p'}(X^{m-1}).$ 
 \begin{enumerate}
  \item 
Set  
$U_R:=(\Pi_1)_*\big((\hat u\hat\eta+\hat\gamma')\wedge\big(\Pi_2^*(R_2)\wedge\ldots\wedge \Pi_m^*(R_m)-(\Pi')^*\gamma_R\big)\big)$ on $X.$ Then
$\ddc U_R=R-\gamma_R|_{\Delta'}$.
\item  Set  
$\widehat U_{R}:=(\hat u\hat\eta+\hat\gamma')\wedge\big(\Pi_2^*(R_2)\wedge\ldots\wedge \Pi_m^*(R_m)\big)$ on $\widehat \X.$ Then for every smooth form $\phi$  of bidegree $(k-p',k-p')$   on $\X,$ we have that
$$\langle  \widehat U_R, \Pi^*(\ddc\phi)\rangle_{\widehat \X} =\langle R,\iota^*\phi\rangle_X -\langle   (\pi_1)_*(\gamma_\Delta)\wedge   (\pi_1)_*(\pi')^*(\gamma_R),   \iota^*\phi\rangle_{X}.$$
\end{enumerate}
\end{proposition}
 \proof  It is enough to consider the case where $R_2,\ldots,R_m$ are smooth.

 We  start with the proof of assertion (1). A direct computation using \eqref{e:gamma} and the definition of $\hat u$ in \eqref{e:hat-u} gives
\begin{eqnarray*}\ddc U_R&=&(\Pi_1)_*\big([\widehat\Delta]\wedge\hat\eta\wedge\Pi_2^*(R_2)\wedge\ldots\wedge \Pi_m^*(R_m)-(\Pi')^*\gamma_R\big)\big)\\
 &-& (\Pi_1)_*\big(\Pi^*(\gamma_\Delta)\wedge\Pi_2^*(R_2)\wedge\ldots\wedge \Pi_m^*(R_m)-(\Pi')^*\gamma_R\big)\big).
 \end{eqnarray*}
Observe  that the restriction of $\Pi_j^*(R_j)$ to $\widehat\Delta$ is equal to that of $\Pi_1^*(R_j)$ for $j=2,\ldots,m,$ and that the restriction of $(\Pi')^*(\gamma_R)$ to $\widehat\Delta$ is equal to that of $\Pi_1^*(\gamma_R|_{\Delta'}).$  Therefore, using \eqref{e:Delta} and the identity $\Pi_1=\pi_1\circ\Pi$, we see that the first term in the RHS of 
the last equation is equal to
$$(\Pi_1)_*([\widehat\Delta]\wedge\hat\eta)\wedge (R- \gamma_R|_{\Delta'})= (\pi_1)_*[\Delta]\wedge (R- (\pi_1)_*(\gamma_R|_{\Delta'}))=R-(\pi_1)_*(\gamma_R|_{\Delta'}).$$

It remains to check that the second term vanishes. Using $\Pi_j=\pi_j\circ\Pi$, we see that 
\begin{eqnarray*} &&(\Pi_1)_*\big(\Pi^*(\gamma_\Delta)\wedge \Pi_2^*(R_2)\wedge\ldots\wedge \Pi_m^*(R_m)-(\Pi')^*\gamma_R\big)\big)\\
&=&(\pi_1)_*\Pi_*\big(\Pi^*(\gamma_\Delta)\wedge\Pi_2^*(R_2)\wedge\ldots\wedge \Pi_m^*(R_m)-(\Pi')^*\gamma_R\big)\big)\\&=&(\pi_1)_*\big(\gamma_\Delta\wedge\pi_2^*(R_2)\wedge\ldots\wedge \pi_m^*(R_m)-(\pi')^*\gamma_R\big)\big).\end{eqnarray*}
Now, if $\Phi$ is a smooth test form of the right bi-degree on $X$ and if $(x^{(1)},\ldots, x^{(m)})$ denotes the coordinates of points in $\X$, we have
\begin{multline*}\big\langle (\pi_1)_*\big(\gamma_\Delta\wedge \big( \pi_2^*(R_2)\wedge\ldots\wedge \pi_m^*(R_m)
-(\pi')^*\gamma_R\big)\big),\Phi\big\rangle \\
= \int_{\X} \Phi(x^{(1)})\wedge \gamma_\Delta(x^{(1)},\ldots, x^{(m)})\wedge\big( R_2(x^{(2)})\wedge\ldots\wedge  R_m(x^{(m)}) -\gamma_R(x^{(2)},\ldots,x^{(m)}  ) \big) \big).
\end{multline*}
Since the cohomology class of $R-\gamma_R$ vanishes in $X^{m-1},$ it is an exact form. Recall from the choice of  $\gamma_\Delta$ that it has  a nice  property of variable separation. Therefore, using Stokes and Fubini's theorems, we see that the last integral vanishes when we first integrate in variables $x^{(2)},\ldots, x^{(m)}.$ So
\begin{equation}\label{e:vanish-gamma-Delta-vs-gamma-R}\big\langle (\pi_1)_*\big(\gamma_\Delta\wedge \big( \pi_2^*(R_2)\wedge\ldots\wedge \pi_m^*(R_m)
-(\pi')^*\gamma_R\big)\big),\Phi\big\rangle.
\end{equation}

This completes the proof of  assertion (1).

 To prove assertion (2), a direct computation using \eqref{e:gamma} and the definition of $\hat u$ in \eqref{e:hat-u} gives
\begin{eqnarray*}\ddc \widehat U_R&=&[\widehat\Delta]\wedge\hat\eta\wedge\Pi_2^*(R_2)\wedge\ldots\wedge \Pi_m^*(R_m)\\
 &-& \Pi^*(\gamma_\Delta)\wedge\Pi_2^*(R_2)\wedge\ldots\wedge \Pi_m^*(R_m).
 \end{eqnarray*}
Observe  that the restriction of $\Pi_j^*(R_j)$ to $\widehat\Delta$ is equal to that of $\Pi_1^*(R_j)$ for $j=2,\ldots,m,$ and that the restriction of $(\Pi)^*(\phi)$ to $\widehat\Delta$ is equal to that of $\Pi_1^*(\iota^*(\phi)).$  Therefore, by Stokes' formula,  we have that 
\begin{eqnarray*}\langle  \widehat U_{R,\phi}, \Pi^*(\ddc\phi)\rangle_{\widehat \X}&=&\langle \ddc  \widehat U_{R,\phi}, \Pi^*(\phi)\rangle_{\widehat \X}\\
&=& \langle (\Pi_1)_*( [\widehat\Delta]\wedge\hat\eta) \wedge (\Pi_1)_*(\Pi_1^*(R_2)\wedge\ldots\wedge \Pi_1^*(R_m)) , 
\Pi_1^*(\iota^*\phi)\rangle_{\widehat \X}\\
&-&  \langle (  (\Pi_1)_*([\widehat\Delta]\wedge\hat\eta)\wedge  (\Pi_1)_*(\Pi^*(\gamma_\Delta))\wedge   (\Pi_1)_*(\Pi_1^*(R_2)\wedge\ldots\wedge \Pi_1^*(R_m)),   \Pi_1^*(\iota^*\phi)\rangle_{\widehat \X}.
\end{eqnarray*}
 Using \eqref{e:Delta}  we see that the last expression    is  equal to 
 $$
  \langle R_2\wedge\ldots\wedge R_m , \iota^*\phi\rangle_{X}
-  \langle   (\pi_1)_*(\gamma_\Delta)\wedge   R_2\wedge\ldots\wedge R_m,   \iota^*\phi\rangle_{X}.
 $$
 Using \eqref{e:vanish-gamma-Delta-vs-gamma-R} for the second term, assertion (2) follows.
\endproof

The following theorem  allows us to compute the value of a wedge-product  in terms of  super-potentials.

\begin{theorem} \label{T:Formula-wedge-product}
\begin{enumerate}
\item Let $T_1\in\mD^{p_1}$  and let $\gamma_1$ be a  closed smooth $(p_1,p_1)$-form
such that $\gamma_1\in \{T_1\}.$
  Let $\mU_T$ be as above. Then for every $R_2\in \mD^{p_2},\ldots,R_m\in \mD^{p_m}$ and every closed
  smooth  $(p',p')$-form $\gamma_R$ on $X^{m-1}$ such that  $\gamma_R\in \{R_2\otimes\ldots\otimes R_m\},$ we have 
$$\mU_T(R_2\wedge \ldots \wedge R_m-\gamma_R)=  \int_{\widehat{\X}}   (\hat u\hat\eta+ \hat\gamma') \wedge \Pi_1^*(T_1-\gamma_1)  \wedge\big( \Pi_2^*(R_2)\wedge\ldots\wedge \Pi_m^*(R_m)-(\Pi')^*(\gamma_R)\big).$$
\item  For every $R_1\in \mD^{p_1},\ldots,R_m\in \mD^{p_m}$ and every smooth test form $\phi$ on $\X,$  we have
\begin{equation} \label{e:Formula-wedge-product}\begin{split} \langle T_1\curlywedge\ldots\curlywedge T_m,\iota^*\phi  \rangle&=  \int_{\widehat{\X}}   (\hat u\hat\eta+ \hat\gamma') \wedge \Pi_1^*(T_1)  \wedge\Pi_2^*(T_2)\wedge\ldots\wedge \Pi_m^*(T_m)\wedge \Pi^*( \ddc\phi)\\
&-   \int_{\widehat{\X}}   (\hat u\hat\eta+ \hat\gamma') \wedge \Pi_1^*(\gamma_1)  \wedge(\Pi')^*(\gamma_R)\wedge \Pi^*( \ddc\phi)+                     \langle \gamma_1\wedge (\pi_1)_*(\gamma_R|_{\Delta'}),\iota^*\phi  \rangle.
\end{split}
\end{equation}
\end{enumerate}
\end{theorem}
\proof  
Let $U_R$ be as in Proposition \ref{P:kernel}. Observe that it is smooth since $R$ is smooth. By definition of super-potential in the introduction, we deduce from the definition of $U_R$ that
$$\mU_T(R-\gamma_R|_{\Delta'})=\langle T_1-\gamma_1, U_{R-\gamma_R|_{\Delta'}}\rangle =\big\langle \Pi_1^*(T_1-\gamma_1), (\hat u\hat\eta+\hat\gamma')\wedge \big((\Pi_2^*(R)\wedge\ldots\wedge \Pi_m^*(R_m)    -(\Pi')^*(\gamma_R)\big)\big\rangle.$$
Assertion (1) follows.

By Proposition \ref{P:kernel} (2) applied to   $R_2:=T_2,\ldots, R_m:=T_m$ and   to $\pi^{*}_1(T_1-\gamma_1)\wedge \phi$ in place of $ \phi,$  we have  that 
\begin{equation*} \langle T_1\wedge T_2\wedge \ldots\wedge T_m,\iota^*\phi\rangle = \langle \gamma_1\wedge T_2\wedge\ldots\wedge T_m,\phi\rangle +\langle \widehat U_{R},\Pi^*( ( \pi^{*}_1(T_1-\gamma_1) \wedge\ddc\phi))\rangle_{\widehat \X}.
\end{equation*}
Write 
$
 \langle \widehat U_{R},\Pi^*( ( \pi^{*}_1(T_1-\gamma_1) \wedge\ddc\phi))\rangle_{\widehat \X}= I-II,
$
where 
\begin{eqnarray*}
I&:=&\int_{\widehat{\X}}   (\hat u\hat\eta+ \hat\gamma') \wedge \Pi_1^*(T_1)  \wedge\big( \Pi_2^*(T_2)\wedge\ldots\wedge \Pi_m^*(T_m)\wedge \Pi^*(\ddc\phi) \big),\\
II&:=&\int_{\widehat{\X}}   (\hat u\hat\eta+ \hat\gamma') \wedge \Pi_1^*(\gamma_1)  \wedge\big( \Pi_2^*(T_2)\wedge\ldots\wedge \Pi_m^*(T_m)\wedge \Pi^*(\ddc\phi) \big).
\end{eqnarray*}
Set  
$U_R:=(\Pi_1)_*\big((\hat u\hat\eta+\hat\gamma')\wedge\big(\Pi_2^*(T_2)\wedge\ldots\wedge \Pi_m^*(T_m)-(\Pi')^*\gamma_R\big)\big).$ Then
by Proposition \ref{P:kernel} (1),
$\ddc U_R= T_2\wedge \ldots \wedge T_m-\gamma_R|_{\Delta'}$. 
Therefore, by Stokes' theorem,
$$
\langle\gamma_1\wedge (T_2\wedge \ldots \wedge T_m-\gamma_R|_{\Delta'}),\phi\rangle=\langle\gamma_1 \wedge \ddc  U_R,\phi\rangle=\langle \gamma_1\wedge U_R,\ddc\phi  \rangle
$$
On the other hand, 
\begin{equation*}
 II=\langle \gamma_1\wedge U_R,\ddc\phi  \rangle +  \int_{\widehat{\X}}   (\hat u\hat\eta+ \hat\gamma') \wedge \Pi_1^*(\gamma_1)  \wedge(\Pi')^*(\gamma_R)\wedge \Pi^*( \ddc\phi).
\end{equation*}
Hence,   we deduce that
\begin{equation*}
 II=\langle\gamma_1\wedge (T_2\wedge \ldots \wedge T_m-\gamma_R|_{\Delta'}),\phi\rangle +  \int_{\widehat{\X}}   (\hat u\hat\eta+ \hat\gamma') \wedge \Pi_1^*(\gamma_1)  \wedge(\Pi')^*(\gamma_R)\wedge \Pi^*( \ddc\phi).
\end{equation*}
Putting  this  expression for $(II)$ into  the  above  expression for $
 \langle \widehat U_{R},\Pi^*( ( \pi^{*}_1(T_1-\gamma_1) \wedge\ddc\phi))\rangle_{\widehat \X},
$
the desired expression for $\langle T_1\wedge T_2\wedge \ldots\wedge T_m,\phi\rangle $ follows.
\endproof

Finally, we introduce  a  quantity  which has  some similarity with that  given in 
 \eqref{e:vartheta-sup}. Let $T_i\in \CL^{p_i}(X)$ for $1\leq i\leq m$ and  set $\T=T_1\otimes \ldots\otimes T_m\in \CL^p(\X).$ 
 For  $0<\epsilon\leq \epsilon_0,$ 
 consider  the  quantity \begin{equation}\label{e:vartheta-single}
 \begin{split}\vartheta(\T,\epsilon):&= \int_{\widehat\Delta_\epsilon} -\hat u\widehat\omega^{mk-p}\wedge \Pi_1^*(T_1)\wedge\ldots\wedge  \Pi_m^*(T_m)\in[0,\infty],\\
 \vartheta^{\bullet}(\T,\epsilon):&= \int_{\widehat\Delta_\epsilon\setminus \widehat\Delta} -\hat u\widehat\omega^{mk-p}\wedge \Pi_1^*(T_1)\wedge\ldots\wedge  \Pi_m^*(T_m)\in[0,\infty]
 .
  \end{split}
\end{equation}
It is worth noting that for $m=2,$ the following relation holds between \eqref{e:vartheta-sup} and \eqref{e:vartheta-single}:
$$\vartheta_{T_1}(\epsilon)=\sup  \vartheta(T_1\otimes T_2,\epsilon),$$
where  the supremum is taken over all smooth positive closed forms $T_2$ on $X$, of bi-degree $(k-p_1+1,k-p_1+1)$, such that $\|T_2\|\leq 1$.
\subsection{Blow-up model  versus normal vector bundle}\label{SS:blow-up-vs-normal-bundle}

Let $\pi:\E\to\Delta$ denote the normal vector bundle to $\Delta$ in $\X$. We identify $\Delta$ with the zero section of $\E$. For $\lambda\in\C^*$, let $A_\lambda:\E\to \E$ be the multiplication by $\lambda$ on the fibers of $\pi$. The diagonal $\Delta$ is invariant under the action of $A_\lambda$.
Consider a strongly  admissible map  $\tau$ along $\Delta$ from an open  neighborhood of $\Delta$ in  $\X$ to a neighborhood of $\Delta$ in $\E.$   
 Let $h$ be a Hermitian metric on $\E.$ 
 Using formula \eqref{e:alpha-beta-spec} we define  $\alpha$ and $\beta$ on $\E.$
 Using formula  \eqref{e:hat-alpha} we define  $\hat\alpha$ on $\E.$
The following  result relates the blow-up along the diagonal  versus  the normal bundle along the diagonal.  We use the  holomorphic admissible maps $\tau_\ell:\ \bfU_\ell\to\U_\ell$ introduced in Subsection \ref{SS:Ex-Stand-Hyp}  for the present context of the submanifold $\Delta$
  in the ambient manifold $\X.$ 

\begin{proposition}\label{P:blowup-vs-normal-bundle}
There  is a constant $c_3>0$ such that  
 the following inequalities hold 
\begin{equation}\label{e:P:blowup-vs-normal-bundle}\begin{split}
\Pi_*(\hat\omega)&\approx \omega+\tau_\ell^*\alpha_\ver+\tau_\ell^*\beta_\ver\approx \tau_\ell^*\hat\alpha\quad\text{on}\quad\bfU_\ell\setminus V,\\
\pm\big(  (\tau^{-1}\circ \tau_\ell)^*(\Pi_*\hat\omega) -\Pi_*\hat\omega -H\big)^\sharp &\lesssim  c_3 \omega +c_4\tau_\ell^*(\hat\beta) +c_3\varphi^{1/2}\Pi_*\hat\omega\quad\text{on}\quad(\bfU_\ell\cap \tau^{-1}(\U_\ell))\setminus V,\\ 
\big(  \tau^{-1}\circ \tau_\ell)^*(\Pi_*\hat\omega) - \Pi_*\hat\omega \big) &\trianglelefteq   c_3 \omega +c_4\tau_\ell^*(\hat\beta) +c_3\varphi^{1/4}\Pi_*\hat\omega\quad\text{on}\quad(\bfU_\ell\cap \tau^{-1}(\U_\ell))\setminus V.
\end{split}
\end{equation}
Here,  $H$ is  some  form in the class $\Hc(\X)$ given in Definition  \ref{D:Hc}.   
\end{proposition}
Note that   $ \Tube(B,\bfr)\subset\bigcup_{1\leq\ell\leq\ell_0} (\bfU_\ell\cap \tau^{-1}(\U_\ell)). $
 \proof
 Consider  the following  local model of the blow-up  $\Pi:\widehat{\X} \to \X.$  More specifically, 
consider a $mk$-dimensional polydisc $\D^{mk}$ in $\X$  with holomorphic coordinates  $$(z,w)=\big(z_1,\ldots,z_{(m-1)k},w_1,\ldots,w_k\big)$$    around  an arbitrary  point of $\Delta$ in $\X.$ We can choose these local coordinates so that 
 $\Delta$ is equal to the linear subspace $\{z_1=\cdots=z_{(m-1)k}=0\}$.
Let $[u_1 : \cdots : u_{(m-1)k}]$  be the homogeneous coordinates on $\P^{(m-1)k-1}.$ Then,  $\widehat{\X}\cap \Pi^{-1}(\D^{mk})$ may be
identified  with  the   complex manifold
\begin{multline*}
\widehat{\D^{mk}}:=\left\lbrace \big (z_1,\ldots,z_{(m-1)k},w_1,\ldots,w_k,[u_1:\cdots:u_{(m-1)k}]\big )\in\D^{mk}\times \P^{(m-1)k-1}:\right.\\
\left. z_iu_j=z_ju_i   
\quad\text{ for }\quad  1\leq i,j\leq (m-1)k\right\rbrace.
\end{multline*}
Observe that   $\Pi$ is induced by the canonical projection from  $\D^{mk}\times\P^{(m-1)k-1}$ onto the factor
$\D^{mk}$.   Let $\Pi'$ be the canonical projection from  $\widehat{\D^{mk}}$ onto the factor $\D^{mk}.$ 
Let $\omega_{mk}$ be the canonical K\"ahler form  on $\C^{mk},$ that is, 
$$\omega_{mk}(z,w):=\ddc \|( z_1,\ldots, z_{(m-1)k})  \|^2 +\ddc \|( w_1,\ldots, w_k)  \|^2.$$
  Let    $\omega_\FS$  be the Fubini-Study form on $\P^{(m-1)k-1}.$ 
Recall that $\omega_\FS$ is induced by the $(1,1)$-form $\ddc\log\|(u_1,\ldots,u_{(m-1)k})\|$ on $\C^{(m-1)k}\setminus\{0\}$.
Since  $\hat\omega$ is a  K\"ahler form on  $\widehat{\D^{mk}},$  
we obtain 
 $$\hat \omega \approx {\Pi'}^*(\omega_{mk})+ {\Pi''}^*(\omega_\FS)       = \ddc \|( z_1,\ldots, z_{(m-1)k})  \|^2  + \ddc \|( w_{1},\ldots, w_{k})  \|^2+ \ddc \log\|( u_1,\ldots, u_{(m-1)k})  \|.$$
Therefore, we infer that on  $\D^{2k}\setminus \Delta,$
  $$\Pi_*(\hat \omega) \approx \Pi_*({\Pi'}^*(\omega_{mk}))+ \Pi_*({\Pi''}^*(\omega_\FS))       =
  \omega_{mk}+ \Pi_*({\Pi''}^*(\omega_\FS)).$$ 
 Since  $[u_1:\cdots:u_{(m-1)k}]=[z_1:\cdots: z_{(m-1)k}]$ outside the exceptional hypersurface $\Pi^{-1}(\Delta)$ of $\widehat {\D^{mk}}$, it  follows that
 $$ \Pi_*(\hat \omega) \approx\ddc \|( z_1,\ldots, z_{(m-1)k})  \|^2  + \ddc \|( w_{1},\ldots, w_{k})  \|^2+ \ddc \log\|( z_1,\ldots, z_{(m-1)k})  \|.$$
 This inequality holds on $\D^{2k}\setminus\Delta$ and hence on $\D^{2k}$ since positive closed $(1,1)$-currents have no mass on subvarieties of codimension $\geq 2$. Now  we come back the K\"ahler manifold  $(X,\omega).$
Since
$$
\ddc \|( z_1,\ldots, z_{(m-1)k})  \|^2\approx \beta_\ver,\quad  \ddc \|( w_{1},\ldots, w_{k})  \|^2\approx \omega (w),\quad  \ddc \log\|( z_1,\ldots, z_{(m-1)k})  \|\approx  \alpha_\ver,
$$
we infer that $\Pi^*\hat\omega\approx \omega+\pi_*(\beta_\ver)+\pi_*(\alpha_\ver).$
By \eqref{e:hat-alpha-beta-with-ver-exp},  we deduce that  $\Pi^*\hat\omega\approx  \pi_*(\hat\alpha).$
Consequently, the result follows  from an application of  Proposition \ref{P:basic-admissible-estimates-II}.
 \endproof

\subsection{Effective criteria for the continuity of  Dinh-Sibony intersection}\label{SS:continuity-DS-intersection}

  As in Theorem \ref{T:Nguyen-convergence}
  let $m\geq 2$ be an integer, and let 
  $(T_{j,n})_{n=0}^\infty$ be a sequence  of currents of uniformly  bounded mass in $ \CL^{p_j}(X)$ for $1\leq j\leq m.$ Consider, for $n\in\N,$   
 $\T_n:=T_{1,n}\otimes\ldots\otimes T_{m,n}\in \CL^p(X^m).$ 
Suppose that $T_{j,n}\to  T_j$ weakly as $n$ tends to infinity for $1\leq j\leq m.$
  
  For $0<r\leq\bfr$   set
\begin{equation}\label{e:Xi_infty} \Xi_\infty(r):=\sum_{j=\lowm}^\upm \Kc_{j,k-p-j}\big(-\log\dist(\cdot,\Delta)\cdot\T,\Delta,\omega_\Delta,r,\tau,h\big).
\end{equation}
We also set
\begin{equation}\label{e:Xi_n-Xi}\begin{split}\Xi_n(r)&:=\sum_{j=\lowm}^\upm \Kc_{j,k-p-j}\big(-\log\dist(\cdot,\Delta)\cdot\T_n,\Delta,\omega_\Delta,r,\tau,h\big)\quad \text{for}\quad n\in\N,\\
\Xi(r)&:=\sup\limits_{n\in\N}\Xi_n(r).
\end{split}
\end{equation}
  By Remark \ref{R:four-equivalent-conds},
  assumption \eqref{e:vartheta-sup} implies  that
  \begin{equation}
   \label{e:lim-Xi(r)=0}
   \lim\limits_{r\to 0+} \Xi_\infty(r)=0\qquad\text{and}\qquad  \lim\limits_{r\to 0+} \Xi(r)=0.
  \end{equation}
  Consider the following family of closed currents of degree $2p$ indexed by $\lambda\in\C^*$
$$\mR_\lambda:=(A_\lambda)_*\tau_*(T_1\otimes \ldots\otimes T_m).$$
Since $\tau$ may not be holomorphic, the current $\mR_\lambda$  may not be of bi-degree $(p,p)$ and we cannot talk about its positivity. However, for any sequence $(\lambda_n)_{n\geq1}$ in $\C^*$ converging to infinity, there is a subsequence $(\lambda_{n_j})_{j\geq 1}$ such that $\mR_{\lambda_{n_j}}$ converges to some positive closed current $\mR$ of bi-degree $(p,p)$ in $\E$, as $j$ tends to infinity. We say that $\mR$ is {\it a tangent current} of $T_1\otimes \ldots\otimes T_m$ along $\Delta$. It may depend on the sequence $\lambda_{n_j}$ but it is independent of the choice of $\tau$. Recall from  Theorem \ref{T:Dinh-Sibony} that $\mR$ is invariant under the action of $A_\lambda$ for every $\lambda\in\C^*$.

Tangent currents can be seen using local holomorphic coordinates near $\Delta$. We will introduce here coordinates which are suitable for the proof of Theorem \ref{T:Nguyen-intersection}. Let $x^{(1)}=(x_1^{(1)},\ldots,x^{(1)}_k)$ denote a local holomorphic coordinate system on a local chart  $U$ of $X$. We consider the natural coordinate system $(x^{(1)},\ldots,x^{(m)})$ on $U^m$ with $ x^{(j)}=(x_1^{(j)},\ldots,x^{(j)}_k)  $, a copy of $x^{(1)}$ for $2\leq j\leq m,$ such that $\Delta_U:=\Delta\cap (U^m)$ is given by the equation $x^{(j)}=x^{(1)}$ for $2\leq j\leq m.$ We will use the coordinates $(w,z)$ for a small neighborhood $W$ of $\Delta_U$ with $z^{(j)}:=x^{(j)}-x^{(1)}$ for $2\leq j\leq m$ and $w:=x^{(1)}$. So $\Delta_U$ is given by the equation $z=0$. 
The restriction of $\E$ to $\Delta_U$ can be identified with $\Delta_U\times \C^k$. In this setting, the projection $\pi$ is just the map $(w,z)\mapsto (w,0)$. The dilation $A_\lambda$ is the map $(w,z)\mapsto (w,\lambda z)$. 

For simplicity, we also identify $W$ with an open subset of $\Delta_U\times (\C^k)^{m-1}$. With all these notations, the currents $\mR,$
$\mR_n$ above satisfy
\begin{eqnarray*} \mR&=&\lim_{j\to\infty} (A_{\lambda_{n_j}})_*(  T_1\otimes \ldots\otimes T_m                ) \quad \text{on} \quad \Delta_U\times (\C^k)^{m-1},\\
 \mR_n&=&\lim_{j\to\infty} (A_{\lambda_{n_j}})_*(  T_{1,n}\otimes \ldots\otimes T_{m,n}                ) \quad \text{on} \quad \Delta_U\times (\C^k)^{m-1},
\end{eqnarray*}
or equivalently, for any real smooth form $\Phi$ of bi-degree $(mk-p,mk-p)$ with compact support in $\Delta_U\times (\C^k)^{m-1}$, we have 
\begin{equation} \label{e:tangent}
\begin{split}
\langle \mR, \Phi\rangle &=\lim_{\lambda \to\infty} \big\langle (A_{\lambda})_*(  T_1\otimes \ldots\otimes T_m), \Phi\big\rangle,\\
\langle \mR_n, \Phi\rangle &=\lim_{\lambda\to\infty} \big\langle (A_{\lambda})_*(  T_{1,n}\otimes \ldots\otimes T_{m,n}), \Phi\big\rangle.
\end{split} 
\end{equation}
In what follows, we only need to consider $\lambda$ such that $|\lambda|\geq 1$. 

So we need to study the pairing in the RHS of the identity \eqref{e:tangent}. 
We only need to consider the case where $$\Phi(w,z)=\Phi_1(w)\otimes \Phi_2(z^{(2)})\otimes \ldots\otimes
\Phi_m(z^{(m)}),$$
where  $\Phi_1$ is  a $(k-p,k-p)$-smooth form compactly supported  in  $\C^k$, and  for $2\leq j\leq m,$ 
  $\Phi_j$ is a $(k,k)$-smooth  form compactly supported  in  $\C^k$ so that  $\Phi(w,z)$ has compact support in $W.$
We have  
\begin{equation} \label{e:Psi-lambda}
 \big\langle (A_\lambda)_*(T_1\otimes \ldots\otimes T_m), \Phi\big\rangle= \big\langle T_1\otimes \ldots\otimes T_m, A_\lambda^*(\Phi)\big\rangle.
\end{equation}

  \begin{theorem}\label{T:Nguyen-intersection-formula} 
   Under the hypothesis and the notation of Theorem  \ref{T:Nguyen-intersection},
  then   the intersection    $T_1\curlywedge\ldots  \curlywedge T_m$  provided  by this theorem satisfies   the  following formula   for  every smooth test form  $\phi$ on $X:$
 \begin{equation} \label{e:Formula-wedge-product-bis}\begin{split} \langle T_1\curlywedge\ldots\curlywedge T_m,\phi  \rangle&=  \int_{\widehat{\X}}   (\hat u\hat\eta+ \hat\gamma') \wedge \Pi_1^*(T_1)  \wedge\Pi_2^*(T_2)\wedge\ldots\wedge \Pi_m^*(T_m)\wedge \Pi^*_1( \ddc\phi)\\
&-   \int_{\widehat{\X}}   (\hat u\hat\eta+ \hat\gamma') \wedge \Pi_1^*(\gamma_1)  \wedge(\Pi')^*(\gamma_R)\wedge \Pi^*_1( \ddc\phi)+                     \langle \gamma_1\wedge \gamma_R|_{\Delta'},\phi  \rangle.
\end{split}
  \end{equation}
  \end{theorem}
  \begin{remark}
   \rm In the case $m=2,$ the result was  proved in \cite[Lemma 2.4]{DinhNguyenVu} using  formula \eqref{e:Formula-wedge-product-special}.
  \end{remark}

  Prior to the  proof, we  develop some preparatory results.
For every $\lambda\in\C$  with $|\lambda|\geq 1$ and $n\in\N^*,$ consider  
\begin{equation} \label{e:Ic-Jc-lambda}
\begin{split}  \Ic^\lambda&:=  \int_{\widehat{\X}}   (\hat u\hat\eta+ \hat\gamma') \wedge \Pi_1^*(T_1)  \wedge\Pi_2^*(T_2)\wedge\ldots\wedge \Pi_m^*(T_m)\wedge \Pi^*( A^*_\lambda (\ddc\Phi))\\
\Jc^\lambda&:=-   \int_{\widehat{\X}}   (\hat u\hat\eta+ \hat\gamma') \wedge \Pi_1^*(\gamma_1)  \wedge(\Pi')^*(\gamma_{T})\wedge \Pi^*(  A^*_\lambda (\ddc\Phi))+                     \langle \gamma_{1}\wedge \gamma_{T}|_{\Delta'},\iota^* (A^*_\lambda (\Phi))  \rangle.
\end{split}
\end{equation}
We will show that the following limits  exist:
\begin{equation}\label{e:Ic-Jc-lim}
 \Ic=\lim_{\lambda\to\infty}\Ic^\lambda,\quad  \Jc=\lim_{\lambda\to\infty}\Jc^\lambda,
\end{equation}
where for $\phi:=(\pi_1)_*(\Phi),$ 
\begin{equation} \label{e:Ic-Jc}\begin{split}  \Ic&:=  \int_{\widehat{\X}}   (\hat u\hat\eta+ \hat\gamma') \wedge \Pi_1^*(T_1)  \wedge\Pi_2^*(T_2)\wedge\ldots\wedge \Pi_m^*(T_m)\wedge \Pi^*_1( \ddc\phi)\\
\Jc&:=-   \int_{\widehat{\X}}   (\hat u\hat\eta+ \hat\gamma') \wedge \Pi_1^*(\gamma_1)  \wedge(\Pi')^*(\gamma_{T})\wedge \Pi^*_1(  \ddc\phi)+                     \langle \gamma_{1}\wedge \gamma_{T}|_{\Delta'},\iota^* (\phi)  \rangle.
\end{split}
\end{equation}

 \begin{lemma}\label{L:lim_Jc-lambda} It holds that  $\Jc=\lim_{\lambda\to\infty}\Jc^\lambda.$
 \end{lemma}
\proof
Since the forms $\gamma_1,$ $\gamma_T,$ $\phi,$ $\hat\eta,$ $\hat\gamma'$  are smooth,
and $\hat u=\log\dist(\cdot,\widehat \Delta)+O(1),$ a straightforward computation gives the result.
\endproof
For $2\leq j\leq m,$ we consider  the $X_1\times X_j$  and  two projections  $\pi_{1,j}$  on the first factor and  $\pi_{2,j}$ on the second factor:
\begin{equation} \label{e:Psi_j,lambda} \Psi_{j,\lambda}:= (\pi_{1,j})_*\big[A_\lambda^*(\Phi_j)\wedge \pi_{2,j}^*(T_j)\big].
\end{equation}
Let $\Delta_{j}$ be the  diagonal of  $X_1\times X_j.$
Let $\iota_{j}: X_j\to\Delta_{j}$ be the canonical injection.
Let $\tilde\pi_{j}:\ \E_{j}\to \Delta_{j} $ be the canonical projection,  where $\E_{j}$
is the  normal vector bundle to $\Delta_{j}$ in  $X_1\times X_j.$
\begin{lemma}\label{L:lim_Psi_j,lambda} For $2\leq j\leq m,$ the following assertions hold:
\begin{enumerate}
\item  $\Psi_{j,\lambda}$ is  a closed  smooth $(p_j,p_j)$-form. 
\item $\Psi_{j,\lambda}$ converges weakly  to $\iota^*_{j}((\tilde\pi_{j})_*(\Phi_j))\wedge T_j$ as $\lambda$ tends to infinity.
\end{enumerate}
\end{lemma}
\proof   
We have \begin{equation} \label{e:Psi_j,lambda-bis} \Psi_{j,\lambda}(w)=\int_{z^{(j)}\in\C^k} \Phi_j(\lambda z^{(j)})\wedge T_j(w+z^{(j)})=   \int_{z^{(j)}\in\C^k} \Phi_j(\lambda (z^{(j)}-w))\wedge T_j(z^{(j)}).
\end{equation}
Since   $\Phi_j$ is  a smooth $(k,k)$-form and $T_j$ is  closed, $\Psi_{j,\lambda}$ is  a  closed smooth $(p_j,p_j)$-form,  which  proves assertion (1).

Assertion (2) follows from  \cite[Lemma 3.2]{DinhNguyenVu}.
\endproof
 By \eqref{e:Psi-lambda} we have that 
\begin{equation}
  \big\langle  A_\lambda^*(T_1\otimes \ldots\otimes T_m),\Phi\big\rangle_\X =  \langle  T_1 \wedge \Psi_{2,\lambda}\wedge \ldots \wedge  \Psi_{m,\lambda},\Phi_1\rangle_X.
 \end{equation}
 By Theorem \ref{T:Formula-wedge-product} (2) (formula  \eqref{e:Formula-wedge-product}), the RHS of the above  equation  is  equal to
 $\Ic^\lambda+\Jc^\lambda,$  where 
\begin{eqnarray*}\Ic^\lambda&:= & \int_{\widehat{\X}}   (\hat u\hat\eta+ \hat\gamma') \wedge \Pi_1^*(T_1)  \wedge\Pi_2^*(\Psi_{2,\lambda})\wedge\ldots\wedge \Pi_m^*(\Psi_{m,\lambda})\wedge \Pi^*_1( \ddc\Phi_1)\\
\Jc^\lambda&:=&-   \int_{\widehat{\X}}   (\hat u\hat\eta+ \hat\gamma') \wedge \Pi_1^*(\gamma_1)  \wedge(\Pi')^*(\gamma_R)\wedge \Pi^*_1( \ddc\Psi_1)+                     \langle \gamma_1\wedge \gamma_R|_{\Delta'},\iota^*\Psi  \rangle.
\end{eqnarray*}
Fix an arbitrary positive number $\delta_0,$ and then   using  \eqref{e:lim-Xi(r)=0}, choose a number $0<\epsilon_0<\delta_0$ small enough such that  
\begin{equation}\label{e:choice-epsilon_0}\Xi_\infty(4\epsilon_0)\leq \delta_0.
\end{equation}
 Let $0\leq \rho\leq 1$ be a smooth function with compact support in $\widehat\Delta_{2\epsilon_0}$ which is equal to 1 on $\widehat\Delta_{\epsilon_0}$.
 Let $\lambda\in\C$ with $|\lambda|\geq 1,$
By  \eqref{e:Ic_n-Jc_n}, we write 
\begin{equation}\label{e:Ic_-epsilon_0^lambda}
 \Ic^\lambda:= \Ic^{',\lambda}_{\epsilon_0}+\Ic^{'',\lambda}_{\epsilon_0} ,
\end{equation}
where
\begin{equation*}
 \begin{split}  \Ic^{',\lambda}_{\epsilon_0} &:=  \int_{\widehat{\X}}   \rho  (\hat u\hat\eta+ \hat\gamma') \wedge \Pi_1^*(T_1)  \wedge\Pi_2^*(\Psi_{2,\lambda})\wedge\ldots\wedge \Pi_m^*(\Psi_{m,\lambda})\wedge \Pi^*_1( \ddc\Phi_1),\\
   \Ic^{'',\lambda}_{\epsilon_0} &:=  \int_{\widehat{\X}}  (1- \rho)  (\hat u\hat\eta+ \hat\gamma') \wedge \Pi_1^*(T_1)  \wedge\Pi_2^*(\Psi_{2,\lambda})\wedge\ldots\wedge \Pi_m^*(\Psi_{m,\lambda})\wedge \Pi^*_1(  \ddc\Phi_1).
   \end{split}
   \end{equation*}
  By \eqref{e:Ic-Jc}, we have
  \begin{equation}\label{e:Ic_-epsilon_0}
 \Ic:= \Ic'_{\epsilon_0}+\Ic''_{\epsilon_0}, 
\end{equation}
where
\begin{equation*}
 \begin{split}
\Ic'_{\epsilon_0} &:=  \int_{\widehat{\X}}   \rho  (\hat u\hat\eta+ \hat\gamma') \wedge \Pi_1^*(T_1)  \wedge\Pi_2^*(T_2)\wedge\ldots\wedge \Pi_m^*(T_m)\wedge \Pi^*_1( \ddc\phi),\\
   \Ic''_{\epsilon_0} &:=  \int_{\widehat{\X}}  (1- \rho)  (\hat u\hat\eta+ \hat\gamma') \wedge \Pi_1^*(T_1)  \wedge\Pi_2^*(T_2)\wedge\ldots\wedge \Pi_m^*(T_m)\wedge \Pi^*_1(  \ddc\phi).
\end{split}
   \end{equation*}
 \begin{lemma}\label{L:Ic'_epsilon_0} There is a constant $c>0$  independent  of $\delta_0,\epsilon_0$  and  the  currents $T_j$'s  such that for every $|\lambda|\geq 1,$
 $$
 | \Ic^{',\lambda}_{\epsilon_0}|\leq c \Xi_\infty(4\epsilon_0)\quad\text{and}\quad | \Ic'_{\epsilon_0}|\leq c \Xi_\infty(4\epsilon_0).
 $$
 \end{lemma}
\proof We only prove  the first  inequality since the proof of the second one is  similar, even easier and simpler than  the first one. 

We may assume  without loss of generality that the $\Phi_j$'s  are  positive volume forms for $2\leq j\leq m.$
Since  $\pm\ddc\Phi_1\lesssim \omega^{k-p}$ and $\hat u=-\log\dist(\cdot,\widehat\Delta)+O(1),$
the  above positivity shows that
\begin{eqnarray*}
 | \Ic^{',\lambda}_{\epsilon_0}|&\lesssim& \int_{\widehat{\X}}   \rho (-\log\dist(\cdot,\widehat\Delta)+O(1)) \wedge \hat\omega^{(m-1)k} \wedge \Pi_1^*(T_1)  \wedge\Pi_2^*(\Psi_{2,\lambda})\wedge\ldots\wedge \Pi_m^*(\Psi_{m,\lambda})\wedge \Pi^*( \omega^{k-p})\\
 &\lesssim& \int_{\X}   (\rho\circ \Pi) (-\log\dist(\cdot,\Delta)+O(1)) \wedge \Pi_*(\hat\omega)^{(m-1)k} \wedge T_1  \wedge\Psi_{2,\lambda}\wedge\ldots\wedge \Psi_{m,\lambda}\wedge  \omega^{k-p}.
\end{eqnarray*}
 By Fubini's theorem and using \eqref{e:Psi_j,lambda-bis}, the last line and  hence $| \Ic^{',\lambda}_{\epsilon_0}|$ satisfies
 \begin{multline*}
 | \Ic^{',\lambda}_{\epsilon_0}|\leq \int_{x\in \Dc_{\epsilon_0} } \Big ( \int_{z^{(2)}\in \C^k} \Phi_2(\lambda z^{(2)}) \ldots \int_{z^{(m)}\in \C^k} 
 \Phi_m(\lambda z^{(m)})
  \log{\big(\sum_{j=2}^m      \| x^{(j)}- x^{(1)}- z^{(j)}\|  \big)  }\Big) \\
  (T_1\otimes T_2\otimes \ldots\otimes T_m)(x)\wedge \hat\omega^{mk-p},
 \end{multline*}
 where  
 $$\Dc_{\epsilon_0}:=\left\lbrace x=(x^{(1)},\ldots, x^{(m)})\in\C^k\times\ldots\times\C^k:\  \|x^{(2)}- x^{(1)}\|+\ldots+\|x^{(m)}- x^{(1)}\|<2\epsilon_0\right\rbrace.$$
 
 The  following elementary  result is  needed.
\begin{lemma}\label{L:elementary-log}  Suppose that the $\Phi_j$'s are positive  smooth $(k,k)$-forms compactly supported in $\D^k.$ Then, for  $x=(x^{(1)},\ldots, x^{(m)})\in\D^k\times\ldots\times\D^k$ and for $\lambda\in\C$ with $|\lambda|\geq 1,$
 \begin{multline*}
\int_{z^{(2)}\in \D^k} \Phi_2(\lambda z^{(2)}) \ldots \int_{z^{(m)}\in \D^k} 
 \Phi_m(\lambda z^{(m)})
  \Big|\log{\big(\sum_{j=2}^m      \| x^{(j)}- x^{(1)}- z^{(j)}\|  \big)  }\Big|\\
  \leq c+ c\min\Big(\log {|\lambda|}, \big| \log{\big(\sum_{j=2}^m      \| x^{(j)}- x^{(1)}\|  \big)  }\big|\Big),
  \end{multline*}
  where $c>0$ is a constant independent of $x.$
\end{lemma}
\proof[Proof of Lemma \ref{L:elementary-log}]
We may assume  without loss of generality that $\Phi_j( z^{(j)})=\Vol_{\D^k}( z^{(j)}),$
where $\Vol_{\D^k}$ is the canonical volume form on $\D^k.$
Write
$$ z':=(z^{(2)},\ldots, z^{(m)})\in\C^{(m-1)k} \quad\text{and}\quad    x':=\big(x^{(2)}- x^{(1)},\ldots, x^{(m)}- x^{(1)}\big)\in\C^{(m-1)k} .$$
Let  $\Vol(z')$ denote  the canonical volume form on  $\C^{(m-1)k}.$ The lemma boils down to 
the  following   inequality 
$$
\lambda^{2(m-1)k}\int\limits_{\|z'\|<\lambda^{-1}}\log{\|y'-z'\|} d\Vol(z')\leq c+c\min(\log\lambda,\log{\|y'\|}\qquad\text{for}\qquad \lambda\geq 1.
$$
Using a suitable unitary matrix of order $(m-1)k,$ we may assume without loss of generality that
$y':=(r,0,\ldots,0),$ where $r\in[0, (m-1)k].$ Let $t$ be the real part of the  first  coordinate
in $z'.$  Since $ \|y'-z'\|= \|r-z'\|\geq  |r-t|,$ an application of Fubini's theorem shows that
it suffices  to prove  the following 
$$
\lambda\int_{|t|<\lambda^{-1}}|\log{ |r-t|}|dt\leq c+ c\min(\log\lambda,|\log r|)\quad\text{for}\quad \lambda\geq 1\quad\text{and}\quad  r \in [0, (m-1)k].
$$
To prove  the  last estimate it is  sufficient to observe that  its LHS $\leq c+c\log r$  if $r\geq 2\lambda^{-1},$ and  that its LHS $\leq c+c\log\lambda$ otherwise.
\endproof
  Using  the definition of  
$\Xi_\infty(4\epsilon)$ given in \eqref{e:Xi_infty}
and  Proposition \ref{P:blowup-vs-normal-bundle}, we see that  $| \Ic^{',\lambda}_{\epsilon_0}|$
 is smaller than $\Xi_\infty(4\epsilon_0) $ times a constant independent of $\delta_0, T_j.$
 This  completes the proof.
\endproof

\begin{lemma}\label{L:lim_Ic^'',lambda_epsilon_0} It holds that
  $\lim\limits_{\lambda\to\infty}\Ic^{'',\lambda}_{\epsilon_0}=\Ic''_{\epsilon_0}.$
\end{lemma}
\proof
Since $(1-\rho) (\hat u\hat\eta+ \hat\gamma')$ is smooth and vanishes near $\widehat\Delta$, it is equal to $\Pi^*(\Theta)$ for some smooth form $\Theta$ on $\X$ which vanishes near $\Delta$. 
This, combined with Lemma \ref{L:lim_Psi_j,lambda}
implies that    
 the considered term is equal to 
 $$ \Pi_1^*(T_1)  \wedge\Pi_2^*(\Psi_{2,\lambda})\wedge\ldots\wedge \Pi_m^*(\Psi_{m,\lambda})\wedge \Pi^*_1(  \ddc\Phi_1\wedge \Theta).
$$ Recall that the tensor product of $m$ currents depends continuously on these currents. Since by Lemma \ref{L:lim_Psi_j,lambda},     the last pairing converges to   
 \begin{eqnarray*}&&\big\langle T_{1}\otimes \iota^*_{2}((\tilde\pi_{2})_*(\Phi_2))\wedge T_2\otimes  \ldots \otimes \iota^*_{m}((\tilde\pi_{m})_*(\Phi_m))\wedge T_m,\Theta\wedge \ddc\Phi_1 \big\rangle\\
 &=&
  \big\langle T_{1}\otimes  \ldots\otimes T_m,\Theta\wedge \ddc\Phi \big\rangle
 \end{eqnarray*}
 as $\lambda$ tends to infinity. By \eqref{e:Ic_-epsilon_0^lambda}, the expression in the last line is equal to $\Ic''_{\epsilon_0}.$
This  proves the lemma.
\endproof

We arrive at the
\proof[End of the proof of Theorem \ref{T:Nguyen-intersection-formula}]
Let $\delta_0>0$ be  an arbitrary number as above. Choose $\epsilon_0 \in(0,\delta_0)$ to  satisfy \eqref{e:choice-epsilon_0}.
By Lemma \ref{L:lim_Ic^'',lambda_epsilon_0}, 
there exists  $\lambda_0 \in\N$  such that for  all $\lambda\in\C$ with $|\lambda|\geq \lambda'_0,$
\begin{equation*}
  |\Ic^{'',\lambda}_{\epsilon_0}-\Ic^{''}_{\epsilon_0}|<\delta_0.
  \end{equation*}
On the  other hand, by Lemma \ref{L:Ic'_epsilon_0}, there is a constant $c>0$  independent  of $\delta_0,\epsilon_0$  and  the  currents $T_j$'s  such that for every $|\lambda|\geq 1,$
 $$
 | \Ic^{',\lambda}_{\epsilon_0}|\leq c \Xi_\infty(4\epsilon_0)\quad\text{and}\quad | \Ic'_{\epsilon_0}|\leq c \Xi_\infty(4\epsilon_0).
 $$
  Putting these estimates together  with \eqref{e:choice-epsilon_0} gives that 
for all $\lambda\in\C$ with $|\lambda|\geq \lambda_0,$ 
$$
\big| (\Ic^{',\lambda}_{\epsilon_0}+\Ic^{'',\lambda}_{\epsilon_0})-(\Ic'_{\epsilon_0}+\Ic''_{\epsilon_0})\big|<(2c+1)\delta_0. 
$$
This,  coupled  with \eqref{e:Ic_-epsilon_0} and \eqref{e:Ic_-epsilon_0^lambda}, implies  the  first convergence in 
\eqref{e:Ic-Jc-lim}. On the other hand, Lemma \ref{L:lim_Jc-lambda}  gives the  second convergence in 
\eqref{e:Ic-Jc-lim}. This ends the proof of the theorem.
\endproof

The remainder is devoted to the proof of Theorem \ref{T:Nguyen-convergence}, which is inspired by
the proof of  \cite[Proposition 2.7]{DinhNguyenVu}.
Set  for $n\in \N^*$ 
\begin{equation} \label{e:Ic_n-Jc_n}\begin{split}  \Ic_n&:=  \int_{\widehat{\X}}   (\hat u\hat\eta+ \hat\gamma') \wedge \Pi_1^*(T_{1,n})  \wedge\Pi_2^*(T_{2,n})\wedge\ldots\wedge \Pi_m^*(T_{m,n})\wedge \Pi^*(  \ddc\phi)\\
\Jc_n&:=-   \int_{\widehat{\X}}   (\hat u\hat\eta+ \hat\gamma') \wedge \Pi_1^*(\gamma_{1,n})  \wedge(\Pi')^*(\gamma_{T,n})\wedge \Pi^*(  \ddc\phi)+                     \langle \gamma_{1,n}\wedge \gamma_{T,n}|_{\Delta'},\iota^* (\phi)  \rangle.
\end{split}
\end{equation}
By  \eqref{e:Formula-wedge-product}, we have that
\begin{equation} \begin{split} \langle T_1\curlywedge\ldots\curlywedge T_m,\iota^*\phi  \rangle&=   \Ic+\Jc,\\
 \langle T_{1,n}\curlywedge\ldots\curlywedge T_{m,n},\iota^*\phi  \rangle&=   \Ic_n+\Jc_n.
\end{split}
\end{equation}
To prove  Theorem \ref{T:Nguyen-convergence},
  we need to show that 
\begin{equation} \label{e:Ic-Jc_converge}
 \lim\limits_{n\to\infty} \Ic_n=\Ic\qquad\text{and}\qquad    \lim\limits_{n\to\infty} \Jc_n=\Jc . 
\end{equation}
Using a partition of unity and the notations as above, we can assume that $\Phi$ has compact support in $U\times \C^k$. Moreover, since  the tangent  currents  are invariant under the action of $A_\lambda$, it is enough to consider the case where the support of $\Phi$ is contained in $W$ as in the situation of the above lemmas.

Fix an arbitrary positive number $\delta_0$ and then choose a number $0<\epsilon_0<\delta_0$ small enough which  verifies \eqref{e:choice-epsilon_0}. Let $0\leq \rho\leq 1$ be a smooth function with compact support in $\widehat\Delta_{2\epsilon_0}$ which is equal to 1 on $\widehat\Delta_{\epsilon_0}$.   
By  \eqref{e:Ic_n-Jc_n}, we write 
\begin{equation}\label{e:Ic_n-epsilon_0}
    \Ic_n:= \Ic'_{n,\epsilon_0}+\Ic''_{n,\epsilon_0},
\end{equation}
where
\begin{equation*}
 \begin{split}  
   \Ic'_{n,\epsilon_0} &:=  \int_{\widehat{\X}}   \rho  (\hat u\hat\eta+ \hat\gamma') \wedge \Pi_1^*(T_{1,n})  \wedge\Pi_2^*(T_{2,n})\wedge\ldots\wedge \Pi_m^*(T_{m,n})\wedge \Pi^*_1( \ddc\phi),\\
   \Ic''_{n,\epsilon_0} &:=  \int_{\widehat{\X}}  (1- \rho)  (\hat u\hat\eta+ \hat\gamma') \wedge \Pi_1^*(T_{1,n})  \wedge\Pi_2^*(T_{2,n})\wedge\ldots\wedge \Pi_m^*(T_{m,n})\wedge \Pi^*_1(  \ddc\phi).
 \end{split}
\end{equation*}

 \begin{lemma}\label{L:Ic'_n_epsilon_0} There is a constant $c>0$  independent  of $\delta,\epsilon_0$  and  the  currents $T_{j,n}$'s  such that for every $n\geq 1,$
 $
 | \Ic'_{n,\epsilon_0}|\leq c \Xi_n(4\epsilon_0).
 $
 \end{lemma}
\proof The  proof of Lemma \ref{L:Ic'_epsilon_0} also works in this context.
\endproof

\begin{lemma}\label{L:lim_Ic^''_n_epsilon_0} It holds that
  $\lim\limits_{n\to\infty}\Ic''_{n,\epsilon_0}=\Ic''_{\epsilon_0}.$
\end{lemma}
\proof
Since $(1-\rho) (\hat u\hat\eta+ \hat\gamma')$ is smooth and vanishes near $\widehat\Delta$, it is equal to $\Pi^*(\Theta)$ for some smooth form $\Theta$ on $\X$ which vanishes near $\Delta$. It follows that the considered term is equal to 
$\big\langle T_{1,n}\otimes\ldots \otimes T_{m,n},\Theta\wedge \pi^*_1(\ddc\phi) \big\rangle$. Recall that the tensor product of $m$ currents depends continuously on these currents. Since $T_{j,n}$ converges weakly to $T_j$, the last pairing converges to   $\big\langle T_1\otimes\ldots \otimes T_n,\Theta\wedge \pi_1^*(\ddc\phi)\big\rangle .$  
This  proves  the result.
\endproof
 \begin{lemma} \label{L:lim_Jc-n}
  It holds that
    $\lim\limits_{n\to\infty} \Jc_n=\Jc . $
 \end{lemma}
 \proof Observe that in 
  \eqref{e:Ic_n-Jc_n} $\gamma_{1,n}$ converges uniformly  to $\gamma_1$  and  $\gamma_{T,n}$ converge uniformly to $\gamma_T$ as $n\to\infty.$
Consequently,
a straightforward computation using   the  expressions of $\Jc_n$ in \eqref{e:Ic_n-Jc_n}
and of $\Jc$ in \eqref{e:Ic-Jc}  gives the result.
 \endproof
 \proof[End of the proof of Theorem \ref{T:Nguyen-convergence}]
By Lemma  \ref{L:lim_Ic^''_n_epsilon_0},  there exists  $N_0 \in\N$  such that for  all $n\geq N_0,$
\begin{equation*}
  |\Ic''_{n,\epsilon_0}-\Ic''_{\epsilon_0}|<\delta_0.
  \end{equation*}
  By Lemma \ref{L:Ic'_n_epsilon_0}, there is a constant $c>0$  independent  of $\delta,\epsilon_0$  and  the  currents $T_{j,n}$'s  such that for every $n\geq 1,$
 $
 | \Ic'_{n,\epsilon_0}|\leq c \Xi_n(4\epsilon_0).$
  Since we know by \eqref{e:Xi_n-Xi} that  $\Xi_n(4\epsilon_0)\leq \Xi(4\epsilon_0),$
  inequality \eqref{e:choice-epsilon_0} implies that $
 | \Ic'_{n,\epsilon_0}|\leq c \delta_0.$
  On the other hand, by Lemma  \ref{L:Ic'_epsilon_0}
  we have  that $| \Ic'_{\epsilon_0}|\leq c \delta_0.$ 
   Combining these  estimates  we  obtain that
for all $n\geq N_0,$ 
$$
\big| (\Ic'_{n,\epsilon_0}+\Ic''_{n,\epsilon_0})-(\Ic'_{\epsilon_0}+\Ic''_{\epsilon_0})\big|<(2c+1)\delta_0. 
$$
This,  coupled  with \eqref{e:Ic_n-epsilon_0} and   \eqref{e:Ic_-epsilon_0},  implies  the  first convergence in 
\eqref{e:Ic-Jc_converge}. On the other hand, Lemma \ref{L:lim_Jc-n}  gives the  second convergence in 
\eqref{e:Ic-Jc_converge}.
This ends the proof of the theorem.
\endproof

 \subsection{Intersection theory via  the blow-up along the  diagonal} \label{SS:intersectionvia-blow-up}

We use the notation  introduced in Subsection \ref{SS:m-fold-products}.
Fix an integer $m\geq 2$ and
 let $(X,\omega)$  be a compact K\"ahler manifold of dimension $k.$ Let $\X:=X^m$ and let $\Delta:=\{ (x,\ldots,x):\ x\in X\}$ be the diagonal of $\X.$
 Consider the  K\"ahler form $\omega_\Delta:=\iota_*\omega$ on $\Delta,$ where $\iota:\ X\to \Delta$ given by $x\ni X \mapsto (x,\ldots, x)\in\Delta$ is the canonical biholomorphism. 
 Let $\Pi:\widehat\X \to \X$ be the blow-up of $\X$ 
along the diagonal $\Delta$ and let $\widehat\Delta:=\Pi^{-1}(\Delta)$ be the exceptional hypersurface. Let $\hat \omega$ be a  K\"ahler form on  $\widehat\X.$ Let $\hat u$ be the  function introduced in  
 \eqref{e:hat-u}.This  is a quasi-p.s.h. function on $\widehat{\X}$.
However we can replace it by an arbitrary function of the form  $\log\dist(\cdot,\widehat\Delta)+O(1).$
The following   definition of a pull-back  of  positive  closed currents by a blow-up is needed.  
\begin{definition}\label{D:pull-back-blow-up} \rm 
Let $T$ be a  positive  closed  $(p,p)$-current on $\X.$ By \cite[Theorem 1.1]{DinhSibony04}  (see also Theorem \ref{T:approximation}), there exist two sequences of positive  closed smooth  $(p,p)$-forms $(T_n^\pm)_{n\in\N}$ in $\X$  with uniformly bounded masses such that  $T^\pm_n\to  T^\pm$ as $n\to\infty$ and $T=T^+-T^-.$ 
By  passing to a subsequence if necessary, we may  assume that $\Pi^*(T^\pm_n)\to \widehat T^\pm$ as $n\to\infty.$  The positive closed $(p,p)$-current 
$$S:= \widehat T^+-\widehat T^-\qquad\text{on}\qquad \widehat \X,$$
is called {\it a member of  the  pull-back $\Pi^*T.$}
It is  worth noting that  $S$  is  not necessarily unique. However, $S$ is uniquely determined  
on the Zariski open set $\widehat\X\setminus \widehat \Delta,$ where $\Pi$ is locally biholomorphic. 
\end{definition}
\begin{remark}\label{R::pull-back-blow-up}\rm We are under  the context of Definition \ref{D:pull-back-blow-up}.
 The restriction of a member $S$ of the  pull-back $\Pi^*T$ to $\widehat \Delta$ defines
a unique current   $S_{\widehat \Delta}\in \CL^{p-1}(\widehat\Delta),$ that is,  $\ind_{\widehat \Delta}S=\iota_* S_{\widehat \Delta},$ where $\iota:\ \widehat \Delta\to \widehat\X$ is the canonical  injection (see Skoda \cite{Skoda}). 
So
$S$ gives mass to $\widehat \Delta$ if and only if $\{S_{\widehat \Delta}\}\not=0,$ where $\{S_{\widehat \Delta}\}$ is the cohomology class of $S_{\widehat \Delta}$ in $H^\ast(\widehat \Delta,\C).$
\end{remark}
 Consider $m$ integers $p_1,\ldots, p_m\geq 1$  such that   $p:=p_1+\ldots+p_m\leq k.$
 Recall from \eqref{e:vartheta-single}  the   functions $\vartheta^\bullet( \cdot,\cdot)$ and  $\vartheta( \cdot,\cdot).$
The following    first  main result of this   section  provides an alternative  to  Theorem \ref{T:Nguyen-intersection}.

  \begin{theorem}\label{T:Nguyen-intersection-bis}
  Let 
  $T_j\in \CL^{p_j}(X)$ for $1\leq j\leq m$ and  
consider $\T:=T_1\otimes\ldots\otimes T_m\in \CL^p(X^m).$
  Suppose that
  $\lim\limits_{r\to 0}\vartheta^{\bullet}(\T,r)=0$ and   a member of the pull-back $\Pi^*\T$  gives no mass to $ \widehat\Delta.$   
  Then   $T_1\curlywedge\ldots  \curlywedge T_m$  exists in the sense of Dinh-Sibony's theory of tangent currents. 
  \end{theorem}
  
  The  proof of Theorem \ref{T:Nguyen-intersection-bis} relies on the following two
  auxiliary results. Prior to their formulation, 
  we  keep the notation  introduced in Subsection \ref{SS:blow-up-vs-normal-bundle} and 
  in Subsection \ref{SS:Ex-Stand-Hyp} and
  present some  new ones.
  Fix a strongly admissible map  $\tau$ along  $\Delta$ in $\X,$ this is in particular  a  diffeomorphism from an open neighborhood of  $\Delta$ in $\X$  onto an  open neighborhood of $\Delta$ in $\E,$ which is  the identity on $\Delta.$ 
  For every $0<r\ll 1,$ consider the  following open neighborhood of $\widehat\Delta$ in $\widehat\X:$
  \begin{equation}\label{e:Dc_r}
   \Dc_r:=\Pi^{-1}(\tau^{-1}(\Tube(\Delta,r)))\subset \widehat \X.
  \end{equation}
  Observe that there is a constant $c>1$    such that 
\begin{equation}\label{e:widehat-Delta-vs-Dc}
  \widehat\Delta_{c^{-1}r} \subset  \Dc_r\subset \widehat\Delta_{cr}\quad \text{for}\quad 0<r\ll 1,\quad\text{and}\quad  \Dc_r\searrow \widehat\Delta\quad \text{as}\quad r\searrow 0+.
  \end{equation}
  We use the  holomorphic admissible maps $\tau_\ell:\ \bfU_\ell\to\U_\ell$ introduced in Subsection \ref{SS:Ex-Stand-Hyp}  for the present context of the submanifold $\Delta$
  in the ambient manifold $\X.$
  By  Proposition \ref{P:blowup-vs-normal-bundle}, we infer that for every $1\leq\ell\leq\ell_0,$ 
\begin{equation}\label{e:push-down-P:blowup-vs-normal-bundle}\begin{split}
(\tau_\ell)_*(\Pi_*\hat\omega)&\approx\pi^*\omega+\alpha_\ver+\beta_\ver\approx\hat\alpha\quad\text{on}\quad\U_\ell\setminus \Delta,\\
\pm\big(  \tau_*(\Pi_*\hat\omega) -(\tau_\ell)_*(\Pi_*\hat\omega) -H\big)^\sharp &\leq c_3 (\tau_\ell)_*\omega +c_4\hat\beta +c_3\varphi^{1/2}(\tau_\ell)_*(\Pi_*\hat\omega)\quad\text{on}\quad(\U_\ell\cap \tau(\bfU_\ell))\setminus \Delta,\\ 
\big(   \tau_*(\Pi_*\hat\omega) -(\tau_\ell)_*( \Pi_*\hat\omega) \big) &\trianglelefteq   c_3 (\tau_\ell)_*\omega +c_4\hat\beta +c_3\varphi^{1/4}(\tau_\ell)_*(\Pi_*\hat\omega)\quad\text{on}\quad(\U_\ell\cap \tau(\bfU_\ell))\setminus \Delta.
\end{split}
\end{equation}
Here,  $H$ is  some  form in the class $\Hc(\E)$ given in Definition  \ref{D:Hc}.   
 Note that $\Tube(\Delta,\bfr)\subset \bigcup_{\ell=1}^{\ell_0}(\U_\ell\cap \tau(\bfU_\ell)).$

  The first auxiliary result  reformulates condition (1)   of Theorem  \ref{T:Nguyen-intersection} in terms of the  function $\vartheta^{\bullet}(\T,\cdot).$
  
    \begin{proposition}\label{P:Nguyen-intersection-bis}
   We  keep the notation introduced in  Theorem  \ref{T:Nguyen-intersection} and Theorem  \ref{T:Nguyen-intersection-bis}. 
   Then the following   two conditions  are  equivalent:
   \begin{enumerate}
    \item $\kappa^\bullet_j(-\log\dist(\cdot,\Delta)\cdot\T,\Delta,\omega_\Delta,\bfr,\tau,h)<\infty$ for some  $\bfr>0$ and for all  $ k-p< j\leq k-\max_{1\leq i\leq m} p_i$;
    \item   $\lim\limits_{r\to 0}  \vartheta^{\bullet}(\T,r)=0.$ 
   \end{enumerate}
 In other words,
  $\T$ satisfies condition (1)   of Theorem  \ref{T:Nguyen-intersection} if and only if  
  $\lim\limits_{r\to 0}  \vartheta^{\bullet}(\T,r)=0.$
  \end{proposition}
  \proof First, we prove  that condition (1) implies condition (2). We deduce  from  condition (1) and Remark \ref{R:four-equivalent-conds} that
  \begin{equation}\label{e:finite-Kc_j,q-P:Nguyen-intersection-bis}
  \Kc_{j,q}(-\log\dist(\cdot,\Delta)\cdot\T,\Delta,\omega_\Delta,0,\bfr,\tau,h)<\infty\quad\text{for}\quad
  j\in\N\quad\text{with}\quad  k-p< j\leq k-\max_{1\leq i\leq m} p_i.\end{equation}   
On the other hand, arguing as in the proof of Lemma \ref{L:min-dim-finite-cond-is-superflous}, the last condition is  also equivalent to the   following one:
\begin{equation*}
  \Kc_{j,q}(-\log\dist(\cdot,\Delta)\cdot\T,\Delta,\omega_\Delta,0,\bfr/2,\tau,h)<\infty\quad\text{for}\quad k-p\leq j\leq k,\,\, 0\leq q\leq mk-p-j.
\end{equation*}
  This, combined with \eqref{e:Dc_r} and the first   estimate of \eqref{e:push-down-P:blowup-vs-normal-bundle}, implies that
  \begin{equation}\label{e:finite-on-Dc_r-P:Nguyen-intersection-bis} \int_{\Dc_\bfr\setminus \widehat\Delta} -\log\dist(\cdot,\Delta)\cdot\widehat\omega^{mk-p}\wedge \Pi_1^*(T_1)\wedge\ldots\wedge  \Pi_m^*(T_m)<\infty.
  \end{equation}
  An alternative way to   prove \eqref{e:finite-on-Dc_r-P:Nguyen-intersection-bis} is  as follows. By Lemma  \ref{L:big-power-alpha_ver-equals-zero} $\alpha_\ver^{(m-1)k}=0.$ This,  combined  with the  first estimate of \eqref{e:push-down-P:blowup-vs-normal-bundle}  and  \eqref{e:finite-Kc_j,q-P:Nguyen-intersection-bis},
  implies \eqref{e:finite-on-Dc_r-P:Nguyen-intersection-bis}.
  
 Putting \eqref{e:finite-on-Dc_r-P:Nguyen-intersection-bis} together with  the second conclusion of  \eqref{e:widehat-Delta-vs-Dc}    yields  condition (2).
 
 To prove   that condition (2) implies condition (1), we argue in  a reverse way.
  \endproof
  
  The  second  auxiliary result is a simple test on the minimality of  horizontal dimension.
  
  \begin{proposition}\label{P:minimal-hor-dim}
  Let 
  $T_j\in \CL^{p_j}(X)$ for $1\leq j\leq m$ and  
consider $\T:=T_1\otimes\ldots\otimes T_m\in \CL^p(X^m).$
  Then   $\T$ has  minimal  horizontal dimension  along $\Delta$   if and only if
 a membre of the pull-back   $\Pi^*\T$  gives no mass to $ \widehat\Delta.$ In particular,  whether a member of $\Pi^*T$ gives mass to $\widehat \Delta$ (or not) is independent of the choice of members.
  \end{proposition}
  \proof Let $\widehat \T$ be a member of the  pull-back $\Pi^*\T.$
   By  Theorem \ref{T:Nguyen-DS}  $\T$ has  minimal  horizontal dimension  along $\Delta$   if and only if $\T$ satisfies condition (2)   of Theorem  \ref{T:Nguyen-intersection}. By \cite[Tangent Theorem II (Theorem 1.11)]{Nguyen21} (3), this is  equivalent to
   \begin{equation}\label{e:check-C:T:Nguyen-intersection-bis}
   \lim\limits_{r\to 0+}\kappa_j(\T,\Delta,\omega_\Delta,r,\tau,h)=0,
  \end {equation}
 for all  $j\in\N$ with  $ k-p< j\leq k-\max_{1\leq i\leq m} p_i.$  
 
On the other hand,  applying Lemma \ref{L:spec-wedge}, the last condition is  also equivalent to the   following one:
\begin{equation*}
  \lim_{r\to 0+}\Kc_{j,q}(\T,\Delta,\omega_\Delta,r,\tau,h)=0\quad\text{for}\quad k-p< j\leq k- \max_{1\leq i\leq m} p_i .
\end{equation*}
By Lemma  \ref{L:big-power-alpha_ver-equals-zero} $\alpha_\ver^{(m-1)k}=0.$ This,  combined  with the  first estimate of \eqref{e:push-down-P:blowup-vs-normal-bundle}  and
 \eqref{e:Dc_r},
the last condition 
  is  equivalent to 
  \begin{equation*} \lim_{r\to 0+}\int_{\Dc_\bfr} \widehat\omega^{mk-p}\wedge \Pi_1^*(T_1)\wedge\ldots\wedge  \Pi_m^*(T_m)=0.
  \end{equation*}
 By     \eqref{e:widehat-Delta-vs-Dc},  the last limit is equal to the mass of the restriction of  $\widehat \T$ to $\widehat\Delta.$ 
This completes the proof.
  \endproof

  \proof[End of the proof of Theorem \ref{T:Nguyen-intersection-bis}]
  It follows  from  Propositions \ref{P:Nguyen-intersection-bis}, \ref{P:minimal-hor-dim} and Theorem \ref{T:Nguyen-intersection}.
    \endproof
    
    Here is an immediate consequence of Theorem \ref{T:Nguyen-intersection-bis}
  \begin{corollary}\label{C:T:Nguyen-intersection-bis}
    We  keep the notation introduced in   Theorem  \ref{T:Nguyen-intersection-bis}.
    Suppose that
  $\lim\limits_{r\to 0}\vartheta(\T,r)=0.$
  Then   $T_1\curlywedge\ldots  \curlywedge T_m$  exists in the sense of Dinh-Sibony's theory of tangent currents.
  \end{corollary}
  \proof
  Since $0\leq \vartheta^{\bullet}(\T,r)\leq \vartheta(\T,r),$ by  Theorem \ref{T:Nguyen-intersection-bis} we only need  to check that  a member of the pull-back $\Pi^*\T$  gives no mass to $ \widehat\Delta.$ 
  By Proposition \ref{P:minimal-hor-dim}
   it is  sufficient to check that   $\T$ has  minimal  horizontal dimension  along $\Delta.$
    So  by the argument  in the  proof of Proposition \ref{P:minimal-hor-dim}
   it remains  to check \eqref{e:check-C:T:Nguyen-intersection-bis}.
     
  Since $\lim\limits_{r\to 0}\vartheta(\T,r)=0,$ it follows  
  from \eqref{e:vartheta-single} that
 \begin{equation*}\lim\limits_{r\to 0+} \int_{\widehat\Delta_r} \widehat\omega^{mk-p}\wedge \Pi_1^*(T_1)\wedge\ldots\wedge  \Pi_m^*(T_m)=0.
  \end{equation*}
  By the proof of Proposition \ref{P:minimal-hor-dim}
   this is equivalent to  \eqref{e:check-C:T:Nguyen-intersection-bis}.
  \endproof
    We discuss the 
   continuity  of  Dinh-Sibony intersection in terms of the blow-up along the  diagonal.
  The following   second main result of this  section      provides an equivalent alternative  to  Theorem \ref{T:Nguyen-convergence}.
  \begin{theorem}\label{T:Nguyen-convergence-bis}
  Let $m\geq 2$ be an integer, and let 
  $T_{j,n}$ be a sequence  of currents   in $ \CL^{p_j}(X)$ for $1\leq j\leq m.$ Consider  
 $\T_n:=T_{1,n}\otimes\ldots\otimes T_{m,n}\in \CL^p(X^m).$ For $0\leq r\leq \bfr,$ set
 $$
 \vartheta(r):=\sup_{n\in\N} \vartheta(\T_n,r).
 $$
  Suppose that
  \begin{enumerate}
   \item $T_{j,n}\to T_j\in\CL^{p_j}(X)$ weakly as $n$ tends to infinity for $1\leq j\leq m;$
  \item $\lim\limits_{r\to 0}\vartheta(r)=0.$
  \end{enumerate}
  Then   $T_{1,n}\curlywedge\ldots  \curlywedge T_{m,n}$   and  $T_1\curlywedge\ldots  \curlywedge T_m$  exist in the sense of Dinh-Sibony's theory of tangent currents.  Moreover, $T_{1,n}\curlywedge\ldots  \curlywedge T_{m,n}$   converge weakly to $T_1\curlywedge\ldots  \curlywedge T_m$ as $n$ tends to infinity.
  \end{theorem}
  \proof
  
Using Lemma \ref{L:normal-ves-hash-widehat-T}  and \eqref{e:push-down-P:blowup-vs-normal-bundle}, one can show that condition (2) of Theorem \ref{T:Nguyen-convergence} is equivalent to the  following: for $ k-p< j\leq k$ and $ 0\leq q\leq \min(mk-p-j,(m-1)k-1),$
\begin{equation*} 
  \lim_{r\to 0+}\Kc_{j,q}(-\log\dist(\cdot,\Delta)\cdot\T,\Delta,\omega_\Delta,r,\tau,h)=0.
\end{equation*}
On the other hand, arguing as in the proof of Lemma \ref{L:min-dim-finite-cond-is-superflous}, the last condition is  also equivalent to the   following one:
\begin{equation*}
  \lim_{r\to 0+}\Kc_{j,q}(-\log\dist(\cdot,\Delta)\cdot\T,\Delta,\omega_\Delta,r,\tau,h)=0\quad\text{for}\quad k-p\leq j\leq k,\,\, 0\leq q\leq mk-p-j.
\end{equation*}
By \eqref{e:Dc_r} and the first   estimate of \eqref{e:push-down-P:blowup-vs-normal-bundle},  the last condition 
  is  equivalent to 
  \begin{equation*} \lim_{r\to 0+}\int_{\Dc_\bfr}-\log\dist(\cdot,\Delta)\cdot \widehat\omega^{mk-p}\wedge \Pi_1^*(T_1)\wedge\ldots\wedge  \Pi_m^*(T_m)=0.
  \end{equation*}
 By     \eqref{e:widehat-Delta-vs-Dc}, this is  equivalent to $\lim\limits_{r\to 0}\vartheta(r)=0.$   
Consequently, the result follows from   Theorem \ref{T:Nguyen-convergence}.
  \endproof

\end{document}